\documentclass[letterpaper]{article}

\usepackage{amsmath}
\usepackage{amssymb}
\usepackage{bm}
\usepackage[mathscr]{eucal}
\usepackage{amsfonts}
\usepackage{latexsym}
\usepackage{amsxtra}
\usepackage{amsbsy}
\usepackage{amsthm}
\usepackage{amscd}
\usepackage{amsopn}
\usepackage{amstext}
\newcounter{z0}

\newtheorem{bbbb}{Proposition}[subsection]
\newtheorem{cccc}{Lemma}[subsection]
\newtheorem{dddd}{Theorem}[subsection]
\newtheorem{ffff}{Corollary}[subsection]
\newtheorem{aaaaa}{Definition}[section]
\newtheorem{bbbbb}{Proposition}[section]

\newtheorem{ddddd}{Theorem}[section]
\newtheorem{fffff}{Corollary}[section]

\theoremstyle{definition}

\theoremstyle{definition}
\newtheorem{eeee}{Remark}[subsection]

\theoremstyle{definition}
\newtheorem{eeeee}{Remark}[section]

\theoremstyle{definition}

\newcommand{\me}{\mathrm{e}}
\newcommand{\mi}{\mathrm{i}}
\newcommand{\md}{\mathrm{d}}
\renewcommand{\Im}{\mathrm{Im}}
\renewcommand{\Re}{\mathrm{Re}}
\allowdisplaybreaks[4]

\DeclareMathOperator*{\ad}{ad}
\DeclareMathOperator{\tr}{tr}
\DeclareMathOperator{\diag}{diag}
\numberwithin{equation}{section}
\begin{document}
\baselineskip=12pt
\frenchspacing
\title{Connection Formulae for Asymptotics of Solutions of the Degenerate
Third Painlev\'{e} Equation: II}
\author{A.~V.~Kitaev\thanks{\texttt{E-mail: kitaev@pdmi.ras.ru}} \\
Steklov Mathematical Institute \\
Fontanka 27 \\
St. Petersburg 191023 \\
Russia \and
A.~Vartanian \\
Department of Mathematics \\
College of Charleston \\
Charleston, South Carolina 29424 \\
U.~S.~A.}
\date{\today}
\maketitle
\begin{abstract}
\noindent
The degenerate third Painlev\'{e} equation, $u^{\prime \prime}(\tau) \! = \!
\frac{(u^{\prime}(\tau))^{2}}{u(\tau)} \! - \! \frac{u^{\prime}(\tau)}{\tau}
\! + \! \frac{1}{\tau}(-8 \epsilon u^{2}(\tau) \! + \! 2ab) \! + \!
\frac{b^{2}}{u(\tau)}$, where $\epsilon \! = \! \pm 1$, $b \! \in \!
\mathbb{R} \setminus \lbrace 0 \rbrace$, and $a \! \in \! \mathbb{C}$, is
studied via the Isomonodromy Deformation Method. Asymptotics of general
regular and singular solutions as $\tau \! \to \! \pm \infty$ and $\tau \!
\to \! \pm \mi \infty$ are derived and parametrized in terms of the monodromy
data of the associated $2 \times 2$ linear auxiliary problem introduced in
\cite{a1}. Using these results, three-real-parameter families of solutions
that have infinite sequences of zeroes and poles that are asymptotically
located along the real and imaginary axes are distinguished: asymptotics of
these zeroes and poles are also obtained.

\vspace{0.40cm}


\vspace{0.23cm}
\textbf{Abbreviated Title.} Degenerate Third Painlev\'{e} Equation: II

\vspace{0.23cm}

method, Stokes phenomena
\end{abstract}
\clearpage
\section{Introduction and the Manifold of Monodromy Data} \label{sec1}
In this paper we continue our study \cite{a1} of the degenerate third
Painlev\'{e} equation,
\begin{equation}
u^{\prime \prime}(\tau) \! = \! \dfrac{(u^{\prime}(\tau))^{2}}{u(\tau)} \! -
\! \dfrac{u^{\prime}(\tau)}{\tau} \! + \! \dfrac{1}{\tau}(-8 \epsilon u^{2}
(\tau) \! + \! 2 ab) \! + \! \dfrac{b^{2}}{u(\tau)}, \qquad \epsilon \! = \!
\pm 1, \label{eq1.1}
\end{equation}
where the prime denotes differentiation with respect to $\tau$, and $b \! \in
\! \mathbb{R} \setminus \lbrace 0 \rbrace$ and $a \! \in \! \mathbb{C}$ are
parameters. For other Painlev\'{e} equations, there exist families of
solutions that have zeroes and poles accumulating at the point at infinity
(see, for example, \cite{a2}). In Part I \cite{a1}, such solutions were not
distinguished; therefore, the main purpose of this work (Part II) is to
determine whether or not such solutions exist, and if so, to find their
asymptotics.

As shown in \cite{a1}, Equation~\eqref{eq1.1} can be presented as a
Hamiltonian system, with the Hamiltonian function, $\mathcal{H}(\tau)$, given
by
\begin{equation}
\mathcal{H}(\tau) \! := \! \left(a \! - \! \dfrac{\mi}{2} \right)
\dfrac{b}{u(\tau)} \! + \! \dfrac{1}{2 \tau} \left(a \! - \! \dfrac{\mi}{2}
\right)^{2} \! + \! \dfrac{\tau}{4u^{2}(\tau)}((u^{\prime}(\tau))^2 \! + \!
b^{2}) \! + \! 4 \epsilon u(\tau). \label{eqh1}
\end{equation}
In \cite{a1} we introduced the auxiliary function
\begin{equation}
f(\tau) \! := \! \dfrac{\tau (u^{\prime}(\tau) \! - \! \mi b)}{4u(\tau)} \!
- \! \dfrac{1}{2} \left(\mi a \! + \! \dfrac{1}{2} \right), \label{eqf}
\end{equation}
which defines B\"{a}cklund transformations for Equation~\eqref{eq1.1}:
\begin{equation*}
\dfrac{\md (\tau \mathcal{H}(\tau))}{\md \tau} \! = \! 4 \epsilon q(\tau) \!
+ \! \mi bp(\tau), \quad 2f(\tau) \! = \! q(\tau)p(\tau), \quad -\tau^{2}
\dfrac{\md \mathcal{H}(\tau)}{\md \tau} \! = \! \left(2f(\tau) \! + \! \mi
a \! + \! \dfrac{1}{2} \right)^{2} \! - \! \dfrac{1}{2} \left(\mi a \! + \!
\dfrac{1}{2} \right)^{2},
\end{equation*}
where $q(\tau) \! = \! u(\tau)$ and $p(\tau)$ is a B\"{a}cklund transformation
of $u(\tau)$. In the Hamiltonian setting, the functions $q(\tau)$ and
$p(\tau)$ are, respectively, the generalized coordinate and impulse. (More
detailed information about the functions $\mathcal{H}(\tau)$ and $f(\tau)$,
in particular, the corresponding ODEs they satisfy, can be found in
Proposition~1.3 and Remark~1.3 of \cite{a1}.) In this work, these functions
play an important role in the study of the zeroes and poles of $u(\tau)$.

Section~1 of \cite{a1} contains a review of the literature on the theory and
applications of Equation~\eqref{eq1.1}; so here we mention only those papers
that are related to Equation~\eqref{eq1.1} and which have appeared since the
publication of \cite{a1}. According to the algebro-geometric classification
scheme given in \cite{a3}, the space of initial conditions of
Equation~\eqref{eq1.1} can be characterized by the extended Dynkin diagram of
type $D_{7}$. There is another case of the degenerate third Painlev\'{e}
equation which can be obtained by a similarity reduction of the well-known
Sine-Gordon equation: in the classification scheme of \cite{a3}, this equation
is of type $D_{8}$. The latter equation is known to be equivalent, via a
quadratic transformation, to a special case of the ``complete" third
Painlev\'{e} equation (type $D_{6}$ in the classification of \cite{a3});
therefore, we use the term ``degenerate'' to specify only
Equation~\eqref{eq1.1}, or its equivalent forms. Asymptotics of solutions of
the $D_{8}$ equation were studied via the Isomonodromy Deformation Method in
\cite{a2}. Recently, asymptotics of the so-called instanton solutions of the
$D_{8}$ equation were obtained in \cite{a4} via the exact WKB analysis. We
also mention the recent work \cite{DP}, where a class of semi-flat Calabi-Yau
metrics is obtained in terms of real solutions of Equation~\eqref{eq1.1} with
$a=0$.

In this work we apply the Isomonodromy Deformation Method: the reader can
find some basic ideas and references concerning this method in Sections~1
and~2 of \cite{a1}. We also mention the new monograph~\cite{a5}, which
reflects some recent developments of this method and of a closely related
technique based on a steepest-descent-type analysis of the associated
Riemann-Hilbert problem \cite{a6}. Although Equation~\eqref{eq1.1} resembles
one of the standard representatives of the list of the Painlev\'{e} equations,
its asymptotic study via the Isomonodromy Deformation Method contains, in
contrast to the other Painlev\'{e} equations, additional technical
difficulties: the problem is that the corresponding Fuchs-Garnier (or Lax)
pair is degenerate (see \cite{a1} for details); therefore, in contrast to the
non-degenerate versions of the Painlev\'{e} equations, its associated WKB
analysis requires a much more careful evaluation of the correction terms. In
fact, this is one of the reasons why the present work took $6$ years to
complete since the appearance of \cite{a1}.

In order to make the presentation as self-contained as possible, we now embark
on succinctly reminding the reader about some of the basic facts introduced
in Sections~1 and~2 of \cite{a1}.

Our study of Equation~\eqref{eq1.1} is based on the following Fuchs-Garnier
(or Lax) pair (see Proposition~2.1 of \cite{a1}):
\begin{equation}
\partial_{\mu} \Psi (\mu,\tau) \! = \! \widetilde{\mathscr{U}}(\mu,\tau)
\Psi (\mu,\tau), \quad \quad \partial_{\tau} \Psi (\mu,\tau) \! = \!
\widetilde{\mathscr{V}}(\mu,\tau) \Psi (\mu,\tau), \label{eqFGmain}
\end{equation}
where
\begin{align*}
\widetilde{\mathscr{U}}(\mu,\tau) =& \, -2 \mi \tau \mu \sigma_{3} \! + \! 2
\tau
\begin{pmatrix}
0 & \frac{2 \mi A(\tau)}{\sqrt{-A(\tau)B(\tau)}} \\
-D(\tau) & 0
\end{pmatrix}
\! - \! \dfrac{1}{\mu} \left(\mi a \! + \!
\dfrac{2 \tau A(\tau)D(\tau)}{\sqrt{-A(\tau)B(\tau)}} \! + \! \dfrac{1}{2}
\right) \sigma_{3} \\
+& \, \dfrac{1}{\mu^{2}}
\begin{pmatrix}
0 & \widetilde{\alpha}(\tau) \\
\mi \tau B(\tau) & 0
\end{pmatrix}, \\
\widetilde{\mathscr{V}}(\mu,\tau) =& \, -\mi \mu^{2} \sigma_{3} \! + \! \mu
\begin{pmatrix}
0 & \frac{2 \mi A(\tau)}{\sqrt{-A(\tau)B(\tau)}} \\
-D(\tau) & 0
\end{pmatrix}
\! + \! \left(\dfrac{\mi a}{2 \tau} \! - \! \dfrac{A(\tau)D(\tau)}{
\sqrt{-A(\tau)B(\tau)}} \right) \sigma_{3} \\
-& \, \frac{1}{2 \tau \mu}
\begin{pmatrix}
0 & \widetilde{\alpha}(\tau) \\
\mi \tau B(\tau) & 0
\end{pmatrix},
\end{align*}
with $\sigma_{3} \! = \!
\left(
\begin{smallmatrix}
1 & 0 \\
0 & -1
\end{smallmatrix}
\right)$, and
\begin{equation*}
\widetilde{\alpha}(\tau) \! := \! -\dfrac{2}{B(\tau)} \left(\mi a \sqrt{-A
(\tau)B(\tau)} + \! \tau (A(\tau)D(\tau) \! + \! B(\tau)C(\tau)) \right).
\end{equation*}
\begin{bbbbb} \label{prop1.3}
The Frobenius compatibility condition of System~\eqref{eqFGmain} for arbitrary
values of $\mu \! \in \! \mathbb{C}$ and for differentiable, scalar-valued
functions $A(\tau)$, $B(\tau)$, $C(\tau)$, and $D(\tau)$ is that these
functions satisfy the following system of isomonodromy deformations:
\begin{equation} \label{eq1.4}
\begin{gathered}
A^{\prime}(\tau) \! = \! 4C(\tau) \sqrt{-A(\tau)B(\tau)}, \qquad \qquad
B^{\prime}(\tau) \! = \! -4D(\tau) \sqrt{-A(\tau)B(\tau)}, \\
(\tau C(\tau))^{\prime} \! = \! 2 \mi a C(\tau) \! - \! 2 \tau A(\tau),
\qquad \qquad (\tau D(\tau))^{\prime} \! = \! -2 \mi a D(\tau) \! + \!
2 \tau B(\tau), \\
(\sqrt{-A(\tau)B(\tau)} \,)^{\prime} \! = \! 2(A(\tau)D(\tau) \! - \!
B(\tau)C(\tau)).
\end{gathered}
\end{equation}
\end{bbbbb}
\begin{eeeee} \label{rem1.1}
\textsl{Hereafter, all explicit $\tau$ dependencies are suppressed, except
where confusion may arise}.
\end{eeeee}
A relation between the Fuchs-Garnier pair~\eqref{eqFGmain} and the degenerate
third Painlev\'{e} equation~\eqref{eq1.1} is given by
\begin{bbbbb}[{\rm \cite{a1}}] \label{prop1.2}
Let $u \! = \! u(\tau)$ and $\varphi \! = \! \varphi (\tau)$ solve the system
\begin{equation} \label{eq1.5}
\begin{gathered}
u^{\prime \prime} \! = \! \dfrac{(u^{\prime})^{2}}{u} \! - \!
\dfrac{u^{\prime}}{\tau} \! + \! \dfrac{1}{\tau}(-8 \epsilon u^{2} \! + \!
2ab) \! + \! \dfrac{b^{2}}{u} \,, \\
\varphi^{\prime} \! = \! \dfrac{2a} \tau \! + \! \dfrac{b}{u} \,,
\end{gathered}
\end{equation}
where $\epsilon \! = \! \pm 1$, and $a,b \! \in \! \mathbb{C}$ are independent
of $\tau$. Then,
\begin{equation}
A(\tau) \! := \! \dfrac{u(\tau)}{\tau} \me^{\mi \varphi (\tau)}, \quad
B(\tau) \! := \! -\dfrac{u(\tau)}{\tau} \me^{-\mi \varphi (\tau)}, \quad
C(\tau) \! := \! \dfrac{\epsilon \tau A^{\prime}(\tau)}{4u(\tau)}, \quad
D(\tau) \! := \! -\dfrac{\epsilon \tau B^{\prime}(\tau)}{4u(\tau)},
\label{eq:ABCD}
\end{equation}
solve System~\eqref{eq1.4}. Conversely, let $A(\tau) \! \not\equiv \! 0$,
$B(\tau) \! \not\equiv \! 0$, $C(\tau)$, and $D(\tau)$ solve
System~\eqref{eq1.4}, and define
\begin{equation*}
u(\tau) \! := \! \epsilon \tau \sqrt{-A(\tau)B(\tau)} \,, \qquad \, \varphi
(\tau) \! := \! -\dfrac{\mi}{2} \ln \! \left(-\dfrac{A(\tau)}{B(\tau)}
\right), \quad \, \text{and} \quad \, b \! := \! u(\tau) \! \left(
\varphi^{\prime}(\tau) \! - \! \dfrac{2a}{\tau} \right).
\end{equation*}
Then $b$ is independent of $\tau$, and $u(\tau)$ and $\varphi (\tau)$ solve
System~\eqref{eq1.5}.
\end{bbbbb}

In this work asymptotics (as $\tau \! \to \! \pm \infty$ and $\tau \! \to \!
\pm \mi \infty)$ of solutions of Equation~\eqref{eq1.1} are parametrized in
terms of the monodromy data of System~\eqref{eqFGmain}. This parametrization
is equivalent to finding the corresponding connection formulae; indeed, given
asymptotics of some solution as $\tau \! \to \! +\infty$, say, one uses it
to determine the corresponding monodromy data, and therefore to obtain
asymptotics of the same solution as $\tau \! \to \! -\infty$ or $\tau \! \to
\! \pm \mi \infty$. Furthermore, employing results {}from \cite{a1}, one
arrives at asymptotics of the same solution as $\tau \! \to \! \pm 0$ or
$\tau \! \to \! \pm \mi 0$. Therefore, it is important to remind the reader
about the definition of the monodromy data of System~\eqref{eqFGmain} given
in \cite{a1}.

On the complex $\mu$-plane, System~\eqref{eqFGmain} has two irregular singular
points, $\mu \! = \! \infty$ and $\mu \! = \! 0$. For $\delta_{\infty},
\delta_{0} \! > \! 0$ and $k \! \in \! \mathbb{Z}$, define the (sectorial)
neighborhoods $\Omega_{k}^{\infty}$ and $\Omega_{k}^{0}$, respectively, of
these (singular) points:
\begin{align*}
\Omega_{k}^{\infty} &:= \, \left\{\mathstrut \mu \! \in \! \mathbb{C} \,
\colon \, \vert \mu \vert \! > \! \delta_{\infty}^{-1}, \, -\dfrac{\pi}{2}
\! + \! \dfrac{\pi k}{2} \! < \! \arg \mu \! + \! \dfrac{1}{2} \arg \tau \!
< \! \dfrac{\pi}{2} \! + \! \dfrac{\pi k}{2} \right\}, \\
\Omega_{k}^{0} &:= \, \left\{\mathstrut \mu \! \in \! \mathbb{C} \, \colon \,
\vert \mu \vert \! < \! \delta_{0}, \, -\pi \! + \! \pi k \! < \! \arg \mu \!
- \! \dfrac{1}{2} \arg \tau \! - \! \dfrac{1}{2} \arg (\epsilon b) \! < \!
\pi \! + \! \pi k \right\}.
\end{align*}
The following Proposition is a direct consequence of general asymptotic
results for linear ODEs \cite{W,F}.
\begin{bbbbb}[{\rm \cite{a1}}] \label{prop1.4}
For $k \! \in \! \mathbb{Z}$, there exist solutions $Y_{k}^{\infty}(\mu) \! =
\! Y_{k}^{\infty}(\mu,\tau)$ and $X_{k}^{0}(\mu) \! = \! X_{k}^{0}(\mu,\tau)$
of System~\eqref{eqFGmain} which are uniquely defined by the following
asymptotic expansions:
\begin{align*}
Y_{k}^{\infty}(\mu) \underset{\underset{\scriptstyle \mu \in \Omega_{k}^{
\infty}}{\mu \to \infty}}{:=}& \, \left(\mathrm{I} \! + \! \dfrac{1}{\mu}
\Psi^{(1)} \! + \! \dfrac{1}{\mu^{2}} \Psi^{(2)} \! + \! \dotsb \right)
\exp \! \left(-\mi \left(\tau \mu^{2} \! + \! \left(a \! - \! \dfrac{\mi}{2}
\right) \ln \mu \right) \sigma_{3} \right), \\
X_{k}^{0}(\mu) \underset{\underset{\scriptstyle \mu \in \Omega_{k}^{0}}{\mu
\to 0}}{:=}& \, \Psi_{0} \left(\mathrm{I} \! + \! \mathcal{Z}_{1} \mu \! +
\! \dotsb \right) \exp \! \left(-\frac{\mi \sqrt{\tau \epsilon b}}{\mu} \,
\sigma_{3} \right),
\end{align*}
where $\ln \mu \! := \! \ln \vert \mu \vert \! + \! \mi \arg \mu$,
\begin{gather*}
\Psi^{(1)} \! = \!
\begin{pmatrix}
0 & \frac{A}{\sqrt{-AB}} \\
\frac{D}{2 \mi} & 0
\end{pmatrix}, \qquad \qquad \Psi^{(2)} \! = \!
\begin{pmatrix}
\psi^{(2)}_{11} & 0 \\
0 & \psi^{(2)}_{22}
\end{pmatrix}, \\
\psi^{(2)}_{11} \! = \! -\dfrac{\mi}{2} \left(\tau \sqrt{-AB}+ \! \tau DC \!
+ \! \dfrac{AD}{\sqrt{-AB}} \right), \qquad \qquad \psi^{(2)}_{22} \! = \!
\dfrac{\mi \tau}{2} \left(\sqrt{-AB} + \! CD \right), \\
\Psi_{0} \! = \! \dfrac{\mi}{\sqrt{2}} \left(\dfrac{(\epsilon b)^{1/4}}{
\tau^{1/4} \sqrt{B}} \right)^{\sigma_{3}} \left(\sigma_{1} \! + \! \sigma_{3}
\right), \qquad \qquad \mathcal{Z}_{1} \! = \!
\begin{pmatrix}
z_{1}^{(11)} & z_{1}^{(12)} \\
-z_{1}^{(12)} & - z_{1}^{(11)}
\end{pmatrix}, \\
z_{1}^{(11)} \! = \dfrac{\left(\mi a \! + \! \frac{1}{2} \! + \! \frac{2 \tau
AD}{\sqrt{-AB}} \right)^{2}}{2 \mi \sqrt{\tau \epsilon b}} \! - \! \dfrac{2
\mi \tau^{3/2} \sqrt{-AB}}{\sqrt{\epsilon b}} \! - \! \dfrac{D \sqrt{\tau
\epsilon b}}{B}, \qquad \qquad z_{1}^{(12)} \! = \! \frac{\left(\mi a \! + \!
\frac{1}{2} \! + \! \frac{2 \tau AD}{\sqrt{-AB}} \right)}{2 \mi \sqrt{\tau
\epsilon b}},
\end{gather*}
with $\mathrm{I} \! = \!
\left(
\begin{smallmatrix}
1 & 0 \\
0 & 1
\end{smallmatrix}
\right)$, and $\sigma_{1} \! = \!
\left(
\begin{smallmatrix}
0 & 1 \\
1 & 0
\end{smallmatrix}
\right)$.
\end{bbbbb}
\begin{eeeee} \label{newrem12}
\textsl{The canonical solutions $X_{k}^{0}(\mu)$ are defined uniquely,
provided the branch of $\sqrt{B(\tau)}$ is fixed. Hereafter, the branch of
$\sqrt{B(\tau)}$ is not fixed; therefore, the set of canonical solutions
$\lbrace X_{k}^{0}(\mu) \rbrace_{k \in \mathbb{Z}}$ is defined up to a sign.
This ambiguity doesn't affect the definition of the Stokes multipliers $($see
Equations~{\rm \eqref{eqnstokmult}} below$)$, but results in an ambiguity
of sign in the definition of the connection matrix, $G$ $($see
Equation~{\rm \eqref{eqdefg}} below$)$.}
\end{eeeee}

The \emph{canonical solutions}, $Y_{k}^{\infty}(\mu)$ and $X_{k}^{0}(\mu)$,
enable one to define the \emph{Stokes matrices}, $S_{k}^{\infty}$ and
$S_{k}^{0}$:
\begin{equation}
Y_{k+1}^{\infty}(\mu) \! = \! Y_{k}^{\infty}(\mu)S_{k}^{\infty}, \qquad
X_{k+1}^{0}(\mu) \! = \! X_{k}^{0}(\mu)S_{k}^{0}. \label{eqnstokmult}
\end{equation}
The Stokes matrices are independent of the parameters $\mu$ and $\tau$, and
have the following structures:
\begin{equation*}
S_{2k}^{\infty} \! = \!
\begin{pmatrix}
1 & 0 \\
s_{2k}^{\infty} & 1
\end{pmatrix}, \qquad S_{2k+1}^{\infty} \! = \!
\begin{pmatrix}
1 & s_{2k+1}^{\infty} \\
0 & 1
\end{pmatrix}, \qquad S_{2k}^{0} \! = \!
\begin{pmatrix}
1 & s_{2k}^{0} \\
0 & 1
\end{pmatrix}, \qquad S_{2k+1}^{0} \! = \!
\begin{pmatrix}
1 & 0 \\
s_{2k+1}^{0} & 1
\end{pmatrix}.
\end{equation*}
The parameters $s_{n}^{\infty}$ and $s_{n}^{0}$, $n \! \in \! \mathbb{Z}$,
are called the \emph{Stokes multipliers}. One can show that
\begin{equation}
S_{k+4}^{\infty} \! = \! \me^{-2 \pi (a-\frac{\mi}{2}) \sigma_{3}}
S_{k}^{\infty} \me^{2 \pi (a-\frac{\mi}{2}) \sigma_{3}}, \qquad \quad
S_{k+2}^{0} \! = \! S_{k}^{0}. \label{eq1.8}
\end{equation}
Equations~\eqref{eq1.8} show that the number of independent Stokes multipliers
does not exceed six; for example, $s_{0}^{0}$, $s_{1}^0$, $s_{0}^{\infty}$,
$s_{1}^{\infty}$, $s_{2}^{\infty}$, and $s_{3}^{\infty}$. Furthermore, due
to the special structure of System~\eqref{eqFGmain}, that is, the coefficient
matrices of odd (resp., even) powers of $\mu$ in $\widetilde{\mathscr{U}}
(\mu,\tau)$ are diagonal (resp., off-diagonal) and \emph{vice-versa} for
$\widetilde{\mathscr{V}}(\mu,\tau)$, one can deduce the following relations
for the Stokes matrices (multipliers):
\begin{equation}
S_{k+2}^{\infty} \! = \! \sigma_{3} \me^{-\pi (a-\frac{\mi}{2}) \sigma_{3}}
S_{k}^{\infty} \me^{\pi (a-\frac{\mi}{2}) \sigma_{3}} \sigma_{3}, \qquad
\quad S_{k}^{0} = \! \sigma_{1}S_{k+1}^{0} \sigma_{1}. \label{eq1.9}
\end{equation}
Equations~\eqref{eq1.9} reduce the number of independent Stokes multipliers by
a factor of $2$, that is, all Stokes multipliers can be expressed in terms of
$s_{0}^{0}$, $s_{0}^{\infty}$, $s_{1}^{\infty}$, and the parameter of formal
monodromy, $a$. There is one more relation between the Stokes multipliers,
which follows {}from the so-called cyclic relation (see below). Define the
monodromy matrices at infinity, $M^{\infty}$, and at zero, $M^{0}$, by the
following relations:
\begin{gather*}
Y_{0}^{\infty}(\mu \me^{-2 \pi \mi}) \! := \! Y_{0}^{\infty}(\mu)M^{\infty},
\qquad X_{0}^{0}(\mu \me^{-2 \pi \mi}) \! := \! X_{0}^{0}(\mu)M^{0}.
\end{gather*}
Since $Y_{0}^{\infty}(\mu)$ and $X_{0}^{0}(\mu)$ are solutions of
System~\eqref{eqFGmain}, they differ by a right-hand (matrix) factor $G$:
\begin{equation}
Y_{0}^{\infty}(\mu) \! := \! X_{0}^{0}(\mu)G, \label{eqdefg}
\end{equation}
where $G$ is called the \emph{connection matrix}. As matrices relating
fundamental solutions of System~\eqref{eqFGmain}, the monodromy, connection,
and Stokes matrices are independent of $\mu$ and $\tau$; furthermore, since
$\mathrm{tr}(\widetilde{\mathscr{U}}(\mu,\tau)) \! = \! \mathrm{tr}
(\widetilde{\mathscr{V}}(\mu,\tau)) \! = \! 0$, it follows that
\begin{equation*}
\det (M^{0}) \! = \! \det (M^{\infty}) \! = \! \det (G) \! = \! 1.
\end{equation*}
{}From the definition of the monodromy and connection matrices, one deduces
the following \emph{cyclic relation}:
\begin{equation*}
GM^{\infty} \! = \! M^{0}G.
\end{equation*}
The monodromy matrices can be expressed in terms of the Stokes matrices:
\begin{equation*}
M^{\infty} \! = \! S_{0}^{\infty}S_{1}^{\infty}S_{2}^{\infty}S_{3}^{\infty}
\me^{-2 \pi (a-\frac{\mi}{2}) \sigma_{3}}, \qquad \quad M^{0} \! = \!
S_{0}^{0}S_{1}^{0}.
\end{equation*}
The Stokes multipliers, $s_{0}^{0}$, $s_{0}^{\infty}$, and $s_{1}^{\infty}$,
the elements of the connection matrix, $(G)_{ij} \! =: \! g_{ij}$, $i,j \!
= \! 1,2$, and the parameter of formal monodromy, $a$, are called the
\emph{monodromy data}. Consider $\mathbb{C}^{8}$ with co-ordinates $(a,
s_{0}^{0},s_{0}^{\infty},s_{1}^{\infty},g_{11},g_{12},g_{21},g_{22})$. The
algebraic variety defined by $\det (G) \! = \! 1$ and the \emph{semi-cyclic
relation}
\begin{equation*}
G^{-1}S_{0}^{0} \sigma_{1}G \! = \! S_{0}^{\infty}S_{1}^{\infty} \sigma_{3}
\me^{-\pi (a-\frac{\mi}{2}) \sigma_{3}}
\end{equation*}
is called the \emph{manifold of monodromy data}, $\mathscr{M}$. Since only
three of the four equations in the semi-cyclic relation are independent, it is
clear that $\mathrm{dim}_{\mathbb{C}}(\mathscr{M}) \! = \! 4$; more precisely,
the equations defining $\mathscr{M}$ are:
\begin{gather*}
s_{0}^{\infty}s_{1}^{\infty} \! = \! -1 \! - \! \me^{-2 \pi a} \! - \! \mi
s_{0}^{0} \me^{-\pi a}, \qquad \qquad g_{22}g_{21} \! - \! g_{11} g_{12} \!
+ \! s_{0}^{0}g_{11}g_{22} \! = \! \mi \me^{-\pi a}, \\
g_{11}^{2} \! - \! g_{21}^{2} \! - \! s_{0}^{0} g_{11} g_{21} \! = \! \mi
s_{0}^{\infty} \me^{-\pi a}, \qquad g_{22}^{2} \! - \! g_{12}^{2} \! + \!
s_{0}^{0} g_{12} g_{22} \! = \! \mi s_{1}^{\infty} \me^{\pi a}, \qquad
g_{11}g_{22} \! - \! g_{12} g_{21} \! = \! 1.
\end{gather*}
\begin{eeeee} \label{newrem13}
\textsl{To achieve a one-to-one correspondence between the coefficients
of System~{\rm \eqref{eqFGmain}} and the points on $\mathscr{M}$, one
has to factorize $\mathscr{M}$ by the involution $G \! \to \! -G$ $($cf.
Remark~{\rm \ref{newrem12}}$)$.}
\end{eeeee}

We now describe the contents of this paper.

In Section~\ref{sec2} the main asymptotic results for $u(\tau)$, $\mathcal{H}
(\tau)$, and $f(\tau)$ as $\tau \! \to \! \pm \infty$ are presented. Solutions
$u(\tau)$ having infinite sequences of zeroes and poles accumulating at the
point at infinity of the real and imaginary axes are also specified, and the
asymptotic distribution of these sequences are obtained.

In Section~\ref{sec3} (rather technical in nature) the asymptotic solution of
the direct monodromy problem for the $\mu$-part of System~\eqref{eqFGmain},
under certain restrictions on its coefficients, is presented for positive real
$\tau$. This asymptotic solution is based on matching WKB-asymptotics of the
fundamental solution of System~\eqref{eqFGmain} with its approximation in
terms of parabolic-cylinder functions near the double-turning point.

In Section~\ref{finalsec} the results of Section~\ref{sec3} are inverted
in order to solve the inverse monodromy problem for the $\mu$-part of
System~\eqref{eqFGmain}. At this stage, explicit asymptotics for the
coefficients of the $\mu$-part of System~\eqref{eqFGmain} are parametrized
in terms of the monodromy data. Under the assumption that the monodromy
data are constant and satisfy certain conditions, one finds that the
asymptotics obtained satisfy all of the conditions that were imposed in
Section~\ref{sec3}. According to the justification scheme presented in
\cite{a20}, it follows that there exists exact isomonodromy deformations
corresponding to these monodromy data, that is, solutions of the system of
isomonodromy deformations~\eqref{eq1.4}, whose asymptotics coincide with
the ones obtained in this section.

Appendix~A contains the Laurent expansions at zeroes and poles of the function
$u(\tau)$ together with the corresponding expansions of the associated
functions $\mathcal{H}(\tau)$ and $f(\tau)$. These results are used in
Section~\ref{finalsec} to complete the proof of the asymptotic distribution
of zeroes and poles of $u(\tau)$ for positive real $\tau$.

In Appendix~B the asymptotic results for $\tau \! \to \! \infty$ obtained in
this paper are compared with the corresponding asymptotitcs in \cite{a1},
and misprints {}from \cite{a1} are corrected. This comparison is used in
Section~\ref{finalsec} to resolve a sign ambiguity in the solution of the
inverse monodromy problem.

The main body of this paper is devoted to asymptotics as $\tau \! \to \!
+\infty$. To extend the results to negative and pure imaginary $\tau$, one
applies the action of the Lie-point symmetries $\tau \! \to \! -\tau$ and
$\tau \! \to \! \pm \mi \tau$ on the manifold of monodromy data that were
derived in Subsection~6.2 of \cite{a1}. Asymptotics for real $\tau$ are
presented in Section~\ref{sec2} (as mentioned above), and asymptotics for
imaginary $\tau$ are given in Appendix~C.

We plan to devote the third part of our studies of the degenerate third
Painlev\'{e} Equation~\eqref{eq1.1} to the extension of our results to
complex values of the parameter $b$, and the discussion of the behavior of
some special solutions on the real and imaginary axes together with the
comparison of asymptotic and numerical results.
\section{Summary of Results} \label{sec2}
In this work the detailed analysis for asymptotics of $u(\tau)$ is presented
for the case $\tau \! \to \! +\infty$ and $\epsilon b \! > \! 0$. The
analytic continuation of the function $u(\tau)$ {}from positive values of the
parameters $\tau$ and $\epsilon b$ to negative values of these parameters
is not single valued; therefore, in order to reflect this fact, write $\tau
\! = \! \lvert \tau \rvert \me^{\mi \pi \varepsilon_{1}}$ and $\epsilon b
\! = \! \lvert \epsilon b \rvert \me^{\mi \pi \varepsilon_{2}}$, where
$\varepsilon_{1},\varepsilon_{2} \! = \! 0,\pm 1$. The corresponding
asymptotics, for both positive and negative values of these parameters,
of $u(\tau)$, and the associated functions $\mathcal{H}(\tau)$ and
$f(\tau)$, are convenient to express in terms of the auxiliary mapping
$\mathscr{F}_{\varepsilon_{1},\varepsilon_{2}}$ (see below), which is an
isomorphism of the manifold of monodromy data\footnote{There is a misprint
on page~1173 (Section~3) of \cite{a1}: for items~(2), (3), (5),
(6), (8) and~(9) in the definition of the auxiliary mapping
$\mathscr{F}_{\varepsilon_{1},\varepsilon_{2}}$, the change $a \! \to \! -a$
should be made everywhere.}, $\mathscr{M}$. The following definition of
$\mathscr{F}_{\varepsilon_{1},\varepsilon_{2}}$ is based on Section~6 of
\cite{a1}, that is, transformation~6.2.1 changing\footnote{In
transformation~6.2.1, $\epsilon b \! \to \! \epsilon b$ and $a \! \to \! a$,
that is, $\epsilon_{n}b_{n} \! = \! \epsilon_{o}b_{o}$ and $a_{n} \! = \!
a_{o}$.} $\tau \! \to \! -\tau$ and transformation~6.2.2 changing\footnote{In
transformation~6.2.2, $\tau \! \to \! \tau$, that is, $\tau_{n} \! = \!
\tau_{o}$.} $a \! \to \! -a$.

Define\footnote{$s^{0}_{0}(\varepsilon_{1},\varepsilon_{2}) \! = \!
s^{0}_{0}$.} $\mathscr{F}_{\varepsilon_{1},\varepsilon_{2}} \colon (a,
s^{0}_{0},s^{\infty}_{0},s^{\infty}_{1},g_{11},g_{12},g_{21},g_{22}) \!
\to \! ((-1)^{\varepsilon_{2}}a,s^{0}_{0},s^{\infty}_{0}(\varepsilon_{1},
\varepsilon_{2}),s^{\infty}_{1}(\varepsilon_{1},\varepsilon_{2}),
g_{11}(\varepsilon_{1},\varepsilon_{2}),\linebreak[4]
g_{12}(\varepsilon_{1},\varepsilon_{2}),g_{21}(\varepsilon_{1},
\varepsilon_{2}),g_{22}(\varepsilon_{1},\varepsilon_{2}))$,
$\varepsilon_{1},\varepsilon_{2} \! = \! 0,\pm 1$:
\begin{enumerate}
\item[(1)] $\mathscr{F}_{0,0}$ is the identity mapping: $s^{\infty}_{0}(0,0)
\! = \! s^{\infty}_{0}$, $s^{\infty}_{1}(0,0) \! = \! s^{\infty}_{1}$, and
$g_{ij}(0,0) \! = \! g_{ij}$, $i,j \! = \! 1,2$;
\item[(2)] $\mathscr{F}_{0,-1}$: $s^{\infty}_{0}(0,-1) \! = \! s^{\infty}_{1}
\me^{\pi a}$, $s^{\infty}_{1}(0,-1) \! = \! s^{\infty}_{0} \me^{\pi a}$,
$g_{11}(0,-1) \! = \! -g_{22} \me^{\frac{\pi a}{2}}$, $g_{12}(0,-1) \! = \!
-(g_{21} \! + \! s^{\infty}_{0}g_{22}) \me^{-\frac{\pi a}{2}}$, $g_{21}(0,-1)
\! = \! -(g_{12} \! - \! s^{0}_{0}g_{22}) \me^{\frac{\pi a}{2}}$, and $g_{22}
(0,-1) \! = \! -(g_{11} \! - \! s^{0}_{0}g_{21} \! + \! (g_{12} \! - \!
s^{0}_{0}g_{22})s^{\infty}_{0}) \me^{-\frac{\pi a}{2}}$;
\item[(3)] $\mathscr{F}_{0,1}$: $s^{\infty}_{0}(0,1) \! = \! s^{\infty}_{1}
\me^{\pi a}$, $s^{\infty}_{1}(0,1) \! = \! s^{\infty}_{0} \me^{\pi a}$,
$g_{11}(0,1) \! = \! -\mi g_{12} \me^{\frac{\pi a}{2}}$, $g_{12}(0,1) \! = \!
-\mi (g_{11} \! + \! s^{\infty}_{0}g_{12}) \me^{-\frac{\pi a}{2}}$, $g_{21}
(0,1) \! = \! -\mi g_{22} \me^{\frac{\pi a}{2}}$, and $g_{22}(0,1) \! = \!
-\mi (g_{21} \! + \! s^{\infty}_{0}g_{22}) \me^{-\frac{\pi a}{2}}$;
\item[(4)] $\mathscr{F}_{-1,0}$: $s^{\infty}_{0}(-1,0) \! = \! -s^{\infty}_{0}
\me^{-\pi a}$, $s^{\infty}_{1}(-1,0) \! = \! -s^{\infty}_{1} \me^{\pi a}$,
$g_{11}(-1,0) \! = \! g_{21} \me^{-\frac{\pi a}{2}}$, $g_{12}(-1,0) \! = \!
-g_{22} \me^{\frac{\pi a}{2}}$, $g_{21}(-1,0) \! = \! (g_{11} \! - \!
s^{0}_{0}g_{21}) \me^{-\frac{\pi a}{2}}$, and $g_{22}(-1,0) \! = \! -(g_{12}
\! - \! s^{0}_{0}g_{22}) \me^{\frac{\pi a}{2}}$;
\item[(5)] $\mathscr{F}_{-1,-1}$: $s^{\infty}_{0}(-1,-1) \! = \!
-s^{\infty}_{1}$, $s^{\infty}_{1}(-1,-1) \! = \! -s^{\infty}_{0} \me^{2 \pi
a}$, $g_{11}(-1,-1) \! = \! g_{12} \! - \! s^{0}_{0}g_{22}$, $g_{12}(-1,-1)
\! = \! -g_{11} \! + \! s^{0}_{0}g_{21} \! - \! (g_{12} \! - \! s^{0}_{0}
g_{22})s^{\infty}_{0}$, $g_{21}(-1,-1) \! = \! g_{22} \! - \! (g_{12} \! - \!
s^{0}_{0}g_{22})s^{0}_{0}$, and $g_{22}(-1,-1) \! = \! -g_{21} \! + \! (g_{11}
\! - \! s^{0}_{0}g_{21})s^{0}_{0} \! - \! (g_{22} \! - \! (g_{12} \! - \!
s^{0}_{0}g_{22})s^{0}_{0})s^{\infty}_{0}$;
\item[(6)] $\mathscr{F}_{-1,1}$: $s^{\infty}_{0}(-1,1) \! = \!
-s^{\infty}_{1}$, $s^{\infty}_{1}(-1,1) \! = \! -s^{\infty}_{0} \me^{2 \pi
a}$, $g_{11}(-1,1) \! = \! \mi g_{22}$, $g_{12}(-1,1) \! = \! -\mi (g_{21} \!
+ \! s^{\infty}_{0}g_{22})$, $g_{21}(-1,1) \! = \! \mi (g_{12} \! - \!
s^{0}_{0}g_{22})$, and $g_{22}(-1,1) \! = \! -\mi (g_{11} \! - \! s^{0}_{0}
g_{21} \! + \! (g_{12} \! - \! s^{0}_{0}g_{22})s^{\infty}_{0})$;
\item[(7)] $\mathscr{F}_{1,0}$: $s^{\infty}_{0}(1,0) \! = \! -s^{\infty}_{0}
\me^{\pi a}$, $s^{\infty}_{1}(1,0) \! = \! -s^{\infty}_{1} \me^{-\pi a}$,
$g_{11}(1,0) \! = \! (g_{21} \!+ \! s^{0}_{0}g_{11}) \me^{\frac{\pi a}{2}}$,
$g_{12}(1,0) \! = \! -(g_{22} \! + \! s^{0}_{0}g_{12})
\me^{-\frac{\pi a}{2}}$, $g_{21}(1,0) \! = \! g_{11} \me^{\frac{\pi a}{2}}$,
and $g_{22}(1,0) \! = \! -g_{12} \me^{-\frac{\pi a}{2}}$;
\item[(8)] $\mathscr{F}_{1,-1}$: $s^{\infty}_{0}(1,-1) \! = \! -s^{\infty}_{1}
\me^{2 \pi a}$, $s^{\infty}_{1}(1,-1) \! = \! -s^{\infty}_{0}$, $g_{11}
(1,-1) \! = \! g_{12} \me^{\pi a}$, $g_{12}(1,-1) \! = \! -(g_{11} \! + \!
s^{\infty}_{0}g_{12}) \me^{-\pi a}$, $g_{21}(1,-1) \! = \! g_{22} \me^{\pi
a}$, and $g_{22}(1,-1) \! = \! -(g_{21} \! + \! s^{\infty}_{0}g_{22})
\me^{-\pi a}$;
\item[(9)] $\mathscr{F}_{1,1}$: $s^{\infty}_{0}(1,1) \! = \! -s^{\infty}_{1}
\me^{2 \pi a}$, $s^{\infty}_{1}(1,1) \! = \! -s^{\infty}_{0}$, $g_{11}(1,1)
\! = \! \mi (g_{22} \! + \! s^{0}_{0}g_{12}) \me^{\pi a}$, $g_{12}(1,1) \!
= \! -\mi (g_{21} \! + \! s^{0}_{0}g_{11} \! + \! (g_{22} \! + \! s^{0}_{0}
g_{12})s^{\infty}_{0}) \me^{-\pi a}$, and $g_{22}(1,1) \! = \! -\mi (g_{11}
\! + \! s^{\infty}_{0}g_{12}) \me^{-\pi a}$.
\end{enumerate}
\begin{eeeee} \label{rem2.1}
\textsl{The roots of positive quantities are assumed positive, whilst the
branches of the roots of complex quantities can be taken arbitrarily, unless
stated otherwise. Furthermore, it is assumed that, for negative real $z$,
the following branches are always taken: $z^{1/3} \! := \! -\lvert z
\rvert^{1/3}$ and $z^{2/3} \! := \! (z^{1/3})^{2}$.}
\end{eeeee}
\begin{ddddd} \label{theor2.1}
Let $\varepsilon_{1},\varepsilon_{2} \! = \! 0,\pm 1$, $\epsilon b \! = \!
\vert \epsilon b \vert \me^{\mi \pi \varepsilon_{2}}$, and $u(\tau)$ be a
solution of Equation~\eqref{eq1.1} corresponding to the monodromy data $(a,
s^{0}_{0},s^{\infty}_{0},s^{\infty}_{1},g_{11},g_{12},g_{21},g_{22})$. Suppose
that
\begin{equation}
g_{11}(\varepsilon_{1},\varepsilon_{2})g_{12}(\varepsilon_{1},\varepsilon_{2})
g_{21}(\varepsilon_{1},\varepsilon_{2})g_{22}(\varepsilon_{1},\varepsilon_{2})
\! \neq \! 0, \quad \Re (\widetilde{\nu}(\varepsilon_{1},\varepsilon_{2}) \!
+ \! 1) \! \in \! (0,1) \setminus \left\lbrace \tfrac{1}{2} \right\rbrace,
\label{eqth1.1}
\end{equation}
where
\begin{equation}
\widetilde{\nu}(\varepsilon_{1},\varepsilon_{2}) \! + \! 1 \! := \!
\dfrac{\mi}{2 \pi} \ln (g_{11}(\varepsilon_{1},\varepsilon_{2})
g_{22}(\varepsilon_{1},\varepsilon_{2})). \label{eqth1.5}
\end{equation}
Then there exist $\delta_{G}$ satisfying, for $0 \! < \! \Re (\widetilde{\nu}
(\varepsilon_{1},\varepsilon_{2}) \! + \! 1) \! < \! \tfrac{1}{2}$, the
inequality
\begin{equation*}
0 \! < \! \delta_{G} \! < \! \dfrac{1}{3} \! \left(\dfrac{1 \! + \! 2 \Re
(\widetilde{\nu}(\varepsilon_{1},\varepsilon_{2}) \! + \! 1)}{7 \! + \! 6
\Re (\widetilde{\nu}(\varepsilon_{1},\varepsilon_{2}) \! + \! 1)} \right),
\end{equation*}
and, for $\tfrac{1}{2} \! < \! \Re (\widetilde{\nu}(\varepsilon_{1},
\varepsilon_{2}) \! + \! 1) \! < \! 1$, the inequality
\begin{equation*}
0 \! < \! \delta_{G} \! < \! \dfrac{1}{3} \! \left(\dfrac{3 \! - \! 2 \Re
(\widetilde{\nu}(\varepsilon_{1},\varepsilon_{2}) \! + \! 1)}{9 \! + \! 2
\Re (\widetilde{\nu}(\varepsilon_{1},\varepsilon_{2}) \! + \! 1)} \right),
\end{equation*}
such that $u(\tau)$ has the asymptotic expansion
\begin{align}
u(\tau) \underset{\tau \to \infty \me^{\mi \pi \varepsilon_{1}}}{=}& \,
\dfrac{\epsilon (\epsilon b)^{2/3}}{2} \tau^{1/3} \left(1 \! - \! \dfrac{3}{2
\sin^{2}(\tfrac{1}{2} \vartheta (\varepsilon_{1},\varepsilon_{2},\tau))}
\right) \label{eqth1.2} \\
\underset{\tau \to \infty \me^{\mi \pi \varepsilon_{1}}}{=}& \, \dfrac{
\epsilon (\epsilon b)^{2/3}}{2} \tau^{1/3} \dfrac{\sin (\tfrac{1}{2}
\vartheta (\varepsilon_{1},\varepsilon_{2},\tau) \! - \! \vartheta_{0})
\sin (\tfrac{1}{2} \vartheta (\varepsilon_{1},\varepsilon_{2},\tau) \! +
\! \vartheta_{0})}{\sin^{2}(\tfrac{1}{2} \vartheta (\varepsilon_{1},
\varepsilon_{2},\tau))}, \label{eq:u-asympt-mult}
\end{align}
where
\begin{align}
\vartheta (\varepsilon_{1},\varepsilon_{2},\tau) :=& \, \phi (\tau) \! - \!
\mi \left((\widetilde{\nu}(\varepsilon_{1},\varepsilon_{2}) \! + \! 1) \! -
\! \dfrac{1}{2} \right) \ln \phi (\tau) \! - \! \mi \left((\widetilde{\nu}
(\varepsilon_{1},\varepsilon_{2}) \! + \! 1) \! - \! \dfrac{1}{2} \right)
\ln 12 \! + \! (-1)^{\varepsilon_{2}}a \ln (2 \! + \! \sqrt{3}) \nonumber \\
+& \, \dfrac{\pi}{4} \! - \! \dfrac{3 \pi}{2}(\widetilde{\nu}(\varepsilon_{1},
\varepsilon_{2}) \! + \! 1) \! + \! \mi \ln \left(\dfrac{g_{11}
(\varepsilon_{1},\varepsilon_{2})g_{12}(\varepsilon_{1},\varepsilon_{2})
\Gamma (\widetilde{\nu}(\varepsilon_{1},\varepsilon_{2}) \! + \! 1)}{
\sqrt{2 \pi}} \right) \! + \! \mathcal{O}(\tau^{-\delta_{G}} \ln \tau),
\label{eqth1.3}
\end{align}
with
\begin{gather}
\phi (\tau) \! = \! 3 \sqrt{3} \, (-1)^{\varepsilon_{2}}(\epsilon b)^{1/3}
\tau^{2/3}, \label{eqth1.4} \\
\vartheta_{0} \! = \! -\dfrac{\pi}{2} \! + \! \dfrac{\mi}{2} \ln (2 \! + \!
\sqrt{3}), \label{eq:theta0-def}
\end{gather}
and $\Gamma (\pmb{\cdot})$ is the Euler gamma function {\rm \cite{a24}}.

Let $\mathcal{H}(\tau)$ be the Hamiltonian function defined by
Equation~\eqref{eqh1} corresponding to the function $u(\tau)$ given above.
Then $\mathcal{H}(\tau)$ has the asymptotic expansion
\begin{align}
\mathcal{H}(\tau) \underset{\tau \to \infty \me^{\mi \pi \varepsilon_{1}}}{=}&
\, 3(\epsilon b)^{2/3} \tau^{1/3} \! - \! \mi (-1)^{\varepsilon_{2}}4 \sqrt{3}
\, (\epsilon b)^{1/3} \tau^{-1/3} \left((\widetilde{\nu}(\varepsilon_{1},
\varepsilon_{2})+1) \! - \! \dfrac{1}{2} \! + \! \dfrac{1}{2 \sqrt{3}} \left(
\mi (-1)^{\varepsilon_{2}}a \! + \! \dfrac{1}{2} \right) \right. \nonumber \\
+&\left. \, \dfrac{\mi}{4} \cot (\tfrac{1}{2} \vartheta(\varepsilon_{1},
\varepsilon_{2},\tau)) \! + \! \dfrac{\mi}{4} \cot (\tfrac{1}{2} \vartheta
(\varepsilon_{1},\varepsilon_{2},\tau) \! - \! \vartheta_{0}) \! + \!
\mathcal{O}(\tau^{-\delta_{G}}) \right). \label{eqth1.6}
\end{align}
The function $f(\tau)$ defined by Equation~\eqref{eqf} has the following
asymptotics:
\begin{equation}
f(\tau) \underset{\tau \to \infty \me^{\mi \pi \varepsilon_{1}}}{=}
-\dfrac{(-1)^{\varepsilon_{2}}(\epsilon b)^{1/3}}{2} \tau^{2/3} \left(\mi \!
+ \! \dfrac{3}{\sqrt{2} \, \sin (\tfrac{1}{2} \vartheta (\varepsilon_{1},
\varepsilon_{2},\tau)) \sin (\tfrac{1}{2} \vartheta (\varepsilon_{1},
\varepsilon_{2},\tau) \! - \! \vartheta_{0})} \right). \label{asympt-f-main}
\end{equation}
\end{ddddd}
\begin{eeeee} \label{defofd}
\textsl{Define the strip $($in the $\phi$-plane$)$
\begin{equation}
\mathcal{D} \! := \! \left\{\mathstrut \tau \! \in \! \mathbb{C} \, \colon
\, \Re (\phi (\tau)) \! > \! c_{1}, \, \lvert \Im (\phi (\tau)) \rvert \!
< \! c_{2} \right\}, \label{eqdofd}
\end{equation}
where $\phi (\tau)$ is given in Equation~\eqref{eqth1.4}, and $c_{1},c_{2} \!
> \! 0$ are parameters. In terms of the original variable $\tau$, the strip
$\mathcal{D}$ is a simply-connected domain with convex boundary of increasing
width proportional to $\lvert \tau \rvert^{1/3}$. The asymptotics of
$u(\tau)$, $\mathcal{H}(\tau)$, and $f(\tau)$ presented in
Theorem~{\rm \ref{theor2.1}} are actually valid in the strip domain
$\mathcal{D}$.}
\end{eeeee}
\begin{ddddd} \label{theorem2.2}
Let $\varepsilon_{1},\varepsilon_{2} \! = \! 0,\pm 1$, $\epsilon b \! = \!
\vert \epsilon b \vert \me^{\mi \pi \varepsilon_{2}}$, and $u(\tau)$ be a
solution of Equation~\eqref{eq1.1} corresponding to the monodromy data
$(a,s^{0}_{0},s^{\infty}_{0},s^{\infty}_{1},g_{11},g_{12},g_{21},g_{22})$.
Suppose that
\begin{equation}
g_{11}(\varepsilon_{1},\varepsilon_{2})g_{12}(\varepsilon_{1},\varepsilon_{2})
g_{21}(\varepsilon_{1},\varepsilon_{2})g_{22}(\varepsilon_{1},\varepsilon_{2})
\! \neq \! 0, \quad \Re \! \left(\dfrac{\mi}{2 \pi} \ln (g_{11}
(\varepsilon_{1},\varepsilon_{2})g_{22}(\varepsilon_{1},\varepsilon_{2}))
\right) \! = \! \dfrac{1}{2}. \label{eqtheorem2.21}
\end{equation}
Let the branch of the function $\ln (\pmb{\cdot})$ be chosen\footnote{The
second condition of Equations~\eqref{eqtheorem2.21} suggests that this branch
of $\ln (\pmb{\cdot})$ exists.} such that $\Im (\ln (-g_{11}(\varepsilon_{1},
\varepsilon_{2})g_{22}(\varepsilon_{1},\varepsilon_{2}))) \! = \! 0$. Define
\begin{equation}
\varrho_{1}(\varepsilon_{1},\varepsilon_{2}) \! := \! \dfrac{1}{2 \pi}
\ln (-g_{11}(\varepsilon_{1},\varepsilon_{2})g_{22}(\varepsilon_{1},
\varepsilon_{2})) \quad (\in \! \mathbb{R}). \label{eqtheorem2.22}
\end{equation}
Then $\exists$ $\delta \! \in \! (0,1/39)$ such that the function $u(\tau)$
has, for all large enough $m \! \in \! \mathbb{N}$, second-order poles,
$\tau_{m}^{\infty}$, accumulating at the point at infinity,
\begin{equation}
\tau_{m}^{\infty} \underset{m \to \infty}{=} \me^{\mi \pi \varepsilon_{1}}
\left(\dfrac{2 \pi (-1)^{\varepsilon_{2}}m}{3 \sqrt{3} \, (\epsilon b)^{1/3}}
\right)^{3/2} \left(1 \! - \! \dfrac{3 \varrho_{1}(\varepsilon_{1},
\varepsilon_{2})}{4 \pi} \dfrac{\ln m}{m} \! - \! \dfrac{3 \varrho_{2}
(\varepsilon_{1},\varepsilon_{2})}{4 \pi} \dfrac{1}{m} \right) \! + \!
\mathcal{O} \! \left(m^{\frac{1}{2}-\frac{3 \delta}{2}} \right),
\label{eqtheorem2.23}
\end{equation}
where
\begin{align}
\varrho_{2}(\varepsilon_{1},\varepsilon_{2}) :=& \, \varrho_{1}
(\varepsilon_{1},\varepsilon_{2}) \ln (24 \pi) \! + \! (-1)^{\varepsilon_{2}}
\Re (a) \ln (2 \! + \! \sqrt{3}) \! + \! \dfrac{\pi}{2} \! - \! \dfrac{1}{2}
\arg \! \left(\dfrac{g_{11}(\varepsilon_{1},\varepsilon_{2})g_{12}
(\varepsilon_{1},\varepsilon_{2})}{g_{21}(\varepsilon_{1},\varepsilon_{2})
g_{22}(\varepsilon_{1},\varepsilon_{2})} \right) \nonumber \\
-& \, \arg \! \left(\Gamma \! \left(\tfrac{1}{2} \! + \! \mi \varrho_{1}
(\varepsilon_{1},\varepsilon_{2}) \right) \right) \! + \! \mi \! \left(
(-1)^{\varepsilon_{2}} \Im (a) \ln (2 \! + \! \sqrt{3}) \! + \! \dfrac{1}{2}
\ln \! \left(\left\vert \dfrac{g_{11}(\varepsilon_{1},\varepsilon_{2})g_{12}
(\varepsilon_{1},\varepsilon_{2})}{g_{21}(\varepsilon_{1},\varepsilon_{2})
g_{22}(\varepsilon_{1},\varepsilon_{2})} \right\vert \right) \! \right);
\label{eqtheorem2.25}
\end{align}
furthermore, the function $u(\tau)$ has, for all large enough $m \! \in \!
\mathbb{N}$, a pair of first-order zeroes, $\tau_{m}^{\pm}$, accumulating
at the point at infinity,
\begin{equation}
\tau_{m}^{\pm} \underset{m \to \infty}{=} \me^{\mi \pi \varepsilon_{1}}
\left(\dfrac{2 \pi (-1)^{\varepsilon_{2}}m}{3 \sqrt{3} \, (\epsilon b)^{1/3}}
\right)^{3/2} \left(1 \! - \! \dfrac{3 \varrho_{1}(\varepsilon_{1},
\varepsilon_{2})}{4 \pi} \dfrac{\ln m}{m} \! - \! \dfrac{3}{4 \pi} \left(
\varrho_{2}(\varepsilon_{1},\varepsilon_{2}) \! \pm \! 2 \vartheta_{0}
\right) \dfrac{1}{m} \right) \! + \! \mathcal{O} \! \left(m^{\frac{1}{2}-
\frac{3 \delta}{2}} \right), \label{eqtheorem2.26}
\end{equation}
where $\vartheta_{0}$ is given in Equation~\eqref{eq:theta0-def}.
\end{ddddd}
\begin{eeeee} \label{remark2.2}
\textsl{Further information concerning the Laurent expansions of the functions
$u(\tau)$, $\mathcal{H}(\tau)$, and $f(\tau)$ at poles and zeroes can be found
in Appendix~{\rm A}.}
\end{eeeee}
\begin{eeeee} \label{remark2.3}
\textsl{To present asymptotics of $u(\tau)$, $\mathcal{H}(\tau)$, and
$f(\tau)$ outside of neighborhoods of poles and zeroes, introduce the
technical notion of the cheese-like domain, $\mathcal{D}_{u}$, for a
solution $u(\tau):$
\begin{equation*}
\mathcal{D}_{u} \! := \! \left\{\mathstrut \tau \! \in \! \mathcal{D} \,
\colon \, \lvert \phi (\tau) \! - \! \phi (\tau_{m}^{\kappa}) \rvert \!
\geqslant \! C \lvert \tau_{m}^{\kappa} \rvert^{-\delta} \right\},
\end{equation*}
where the strip domain $\mathcal{D}$ is defined by Equation~\eqref{eqdofd},
$\phi (\tau)$ is given in Equation~\eqref{eqth1.4}, $C \! > \! 0$ is a
parameter, $\kappa \! = \! \infty,\pm$ $(\tau_{m}^{\kappa}$ are the poles and
zeroes introduced in Theorem~{\rm \ref{theorem2.2}}$)$, and $0 \! < \! \delta
\! < \! 1/39$. In terms of the variable $\phi$, the cheese-like domain
$\mathcal{D}_{u}$ is a multiply-connected domain which resembles $\mathcal{D}$
with circular ``cheese-holes'' centered at $\phi (\tau_{m}^{\kappa})$ of
shrinking radii, whilst in terms of the original variable $\tau$ the
diameter of the cheese-holes are increasing, that is, $\lvert \tau \! - \!
\tau_{m}^{\kappa} \rvert \! \sim \! \lvert \tau_{m}^{\kappa}
\rvert^{\frac{1}{3}- \delta}$.}
\end{eeeee}
\begin{ddddd} \label{theorem2.3}
Let $\varepsilon_{1},\varepsilon_{2} \! = \! 0,\pm 1$, $\epsilon b \! = \!
\vert \epsilon b \vert \me^{\mi \pi \varepsilon_{2}}$, and $u(\tau)$ be a
solution of Equation~\eqref{eq1.1} corresponding to the monodromy data
$(a,s^{0}_{0},s^{\infty}_{0},s^{\infty}_{1},g_{11},g_{12},g_{21},g_{22})$.
Suppose that conditions~\eqref{eqtheorem2.21} are valid, the branch of
$\ln (\pmb{\cdot})$ is chosen as in Theorem~{\rm \ref{theorem2.2}}, and
$\varrho_{1}(\varepsilon_{1},\varepsilon_{2})$ is defined by
Equation~\eqref{eqtheorem2.22}. Then there exist $\delta,\delta_{G} \!
\in \! \mathbb{R}_{+}$ satisfying the inequalities
\begin{equation*}
0 \! < \! \delta \! < \! \dfrac{1}{39}, \quad \quad 0 \! < \! \delta
\! < \! \delta_{G} \! < \! \dfrac{1}{15} \! - \! \dfrac{8 \delta}{5},
\end{equation*}
such that $u(\tau)$ has the asymptotic expansion
\begin{equation}
u(\tau) \underset{\underset{\scriptstyle \tau \in \mathcal{D}_{u}}{\tau
\to \infty \me^{\mi \pi \varepsilon_{1}}}}{=} \dfrac{\epsilon (\epsilon
b)^{2/3}}{2} \tau^{1/3} \dfrac{\sin (\tfrac{1}{2} \theta (\varepsilon_{1},
\varepsilon_{2},\tau) \! - \! \vartheta_{0}) \sin (\tfrac{1}{2} \theta
(\varepsilon_{1},\varepsilon_{2},\tau) \! + \! \vartheta_{0})}{\sin^{2}
(\tfrac{1}{2} \theta (\varepsilon_{1},\varepsilon_{2},\tau))},
\label{eqtheorem2.31}
\end{equation}
where
\begin{align}
\theta (\varepsilon_{1},\varepsilon_{2},\tau) :=& \, \phi (\tau) \! + \!
\varrho_{1}(\varepsilon_{1},\varepsilon_{2}) \ln \phi (\tau) \! + \!
\varrho_{1}(\varepsilon_{1},\varepsilon_{2}) \ln 12 \! + \!
(-1)^{\varepsilon_{2}} \Re (a) \ln (2 \! + \! \sqrt{3}) \! + \! \dfrac{\pi}{2}
\nonumber \\
-& \, \dfrac{1}{2} \arg \! \left(\dfrac{g_{11}(\varepsilon_{1},
\varepsilon_{2})g_{12}(\varepsilon_{1},\varepsilon_{2})}{g_{21}
(\varepsilon_{1},\varepsilon_{2})g_{22}(\varepsilon_{1},\varepsilon_{2})}
\right) \! - \! \arg \! \left(\Gamma \! \left(\tfrac{1}{2} \! + \! \mi
\varrho_{1}(\varepsilon_{1},\varepsilon_{2}) \right) \right) \! + \! \mi
\! \left((-1)^{\varepsilon_{2}} \Im (a) \ln (2 \! + \! \sqrt{3}) \right.
\nonumber \\
+&\left. \, \dfrac{1}{2} \ln \! \left(\left\vert \dfrac{g_{11}(\varepsilon_{1},
\varepsilon_{2})g_{12}(\varepsilon_{1},\varepsilon_{2})}{g_{21}
(\varepsilon_{1},\varepsilon_{2})g_{22}(\varepsilon_{1},\varepsilon_{2})}
\right\vert \right) \! \right) \! + \! \mathcal{O}(\tau^{-\delta_{G}}
\ln \tau), \label{eqtheorem2.32}
\end{align}
with $\phi (\tau)$ and $\vartheta_{0}$ given, respectively, in
Equations~\eqref{eqth1.4} and~\eqref{eq:theta0-def}.

Let $\mathcal{H}(\tau)$ be the Hamiltonian function defined by
Equation~\eqref{eqh1} corresponding to the function $u(\tau)$ given above.
Then $\mathcal{H}(\tau)$ has the asymptotic expansion
\begin{align}
\mathcal{H}(\tau) \underset{\underset{\scriptstyle \tau \in
\mathcal{D}_{u}}{\tau \to \infty \me^{\mi \pi \varepsilon_{1}}}}{=}& \, 3
(\epsilon b)^{2/3} \tau^{1/3} \! + \! (-1)^{\varepsilon_{2}}4 \sqrt{3} \,
(\epsilon b)^{1/3} \tau^{-1/3} \left(\varrho_{1}(\varepsilon_{1},
\varepsilon_{2}) \! + \! \dfrac{1}{2 \sqrt{3}} \left((-1)^{\varepsilon_{2}}
a \! - \! \dfrac{\mi}{2} \right) \right. \nonumber \\
+&\left. \, \dfrac{1}{4} \cot (\tfrac{1}{2} \theta(\varepsilon_{1},
\varepsilon_{2},\tau)) \! + \! \dfrac{1}{4} \cot (\tfrac{1}{2} \theta
(\varepsilon_{1},\varepsilon_{2},\tau) \! - \! \vartheta_{0}) \! + \!
\mathcal{O}(\tau^{-\delta_{G}}) \right). \label{eqtheorem2.33}
\end{align}
The function $f(\tau)$ defined by Equation~\eqref{eqf} has the following
asymptotics:
\begin{equation}
f(\tau) \underset{\underset{\scriptstyle \tau \in \mathcal{D}_{u}}{\tau \to
\infty \me^{\mi \pi \varepsilon_{1}}}}{=} -\dfrac{(-1)^{\varepsilon_{2}}
(\epsilon b)^{1/3}}{2} \tau^{2/3} \left(\mi \! + \! \dfrac{3}{\sqrt{2} \,
\sin (\tfrac{1}{2} \theta (\varepsilon_{1},\varepsilon_{2},\tau)) \sin
(\tfrac{1}{2} \theta (\varepsilon_{1},\varepsilon_{2},\tau) \! - \!
\vartheta_{0})} \right). \label{eqtheorem2.34}
\end{equation}
\end{ddddd}
\begin{eeeee} \label{remark2.4}
\textsl{For real, non-zero values of $b$, singular real solutions $u(\tau)$
of Equation~\eqref{eq1.1} are specified by the following ``singular real
reduction'' for the monodromy data\footnote{There exist regular real
solutions $($cf. Part~I~\cite{a1}$)$ which are specified by a different
real reduction: we plan to discuss this in a subsequent work.}$:$
\begin{equation}
\begin{gathered}
s^{0}_{0} \! = \! -\overline{s^{0}_{0}}, \qquad s^{\infty}_{0}(\varepsilon_{1},
\varepsilon_{2}) \! = \! -\overline{s^{\infty}_{1}(\varepsilon_{1},
\varepsilon_{2})} \, \me^{2 \pi a}, \qquad g_{11}(\varepsilon_{1},
\varepsilon_{2}) \! = \! -\overline{g_{22}(\varepsilon_{1},\varepsilon_{2})},
\\
g_{12}(\varepsilon_{1},\varepsilon_{2}) \! = \! -\overline{g_{21}
(\varepsilon_{1},\varepsilon_{2})}, \qquad \Im (a) \! = \! 0.
\label{eqremark2.41}
\end{gathered}
\end{equation}
In this case, asymptotics of $\tau_{m}^{\infty}$, $\tau_{m}^{\pm}$,
$u(\tau)$, $\mathcal{H}(\tau)$, and $f(\tau)$ are as given in
Equations~\eqref{eqtheorem2.23}, \eqref{eqtheorem2.26}, \eqref{eqtheorem2.31},
\eqref{eqtheorem2.33}, and~\eqref{eqtheorem2.34}, respectively, but with the
changes $\varrho_{1}(\varepsilon_{1},\varepsilon_{2}) \! \to \! \varrho_{0}
(\varepsilon_{1},\varepsilon_{2})$, $\varrho_{2}(\varepsilon_{1},
\varepsilon_{2}) \! \to \! \varrho_{0}^{\sharp}(\varepsilon_{1},
\varepsilon_{2})$, and $\theta (\varepsilon_{1},\varepsilon_{2},\tau)
\! \to \! \Theta_{0}(\varepsilon_{1},\varepsilon_{2},\tau)$, where
\begin{align}
\varrho_{0}(\varepsilon_{1},\varepsilon_{2}) :=& \, \dfrac{1}{\pi} \ln
\lvert g_{11}(\varepsilon_{1},\varepsilon_{2}) \rvert, \label{eqremark2.42} \\
\varrho_{0}^{\sharp}(\varepsilon_{1},\varepsilon_{2}) :=& \, \varrho_{0}
(\varepsilon_{1},\varepsilon_{2}) \ln (24 \pi) \! + \! (-1)^{\varepsilon_{2}}
\Re (a) \ln (2 \! + \! \sqrt{3}) \! - \! \dfrac{\pi}{2} \nonumber \\
-&\, \arg \! \left(g_{11}(\varepsilon_{1},\varepsilon_{2})g_{12}
(\varepsilon_{1},\varepsilon_{2}) \Gamma \! \left(\tfrac{1}{2} \! +
\! \mi \varrho_{0}(\varepsilon_{1},\varepsilon_{2}) \right) \right),
\label{eqremark2.44} \\
\Theta_{0}(\varepsilon_{1},\varepsilon_{2},\tau) :=& \, \phi (\tau) \!
+ \! \varrho_{0}(\varepsilon_{1},\varepsilon_{2}) \ln \phi (\tau) \!
+ \! \varrho_{0}(\varepsilon_{1},\varepsilon_{2}) \ln 12 \! + \!
(-1)^{\varepsilon_{2}} \Re (a) \ln (2 \! + \! \sqrt{3}) \nonumber \\
-& \, \dfrac{\pi}{2} \! - \! \arg \! \left(g_{11}(\varepsilon_{1},
\varepsilon_{2})g_{12}(\varepsilon_{1},\varepsilon_{2}) \Gamma \! \left(
\tfrac{1}{2} \! + \! \mi \varrho_{0}(\varepsilon_{1},\varepsilon_{2}) \right)
\right) \! + \! \mathcal{O}(\tau^{-\delta_{G}} \ln \tau). \label{eqremark2.45}
\end{align}}
\end{eeeee}
\section{Asymptotic Solution of the Direct Monodromy Problem}
\label{sec3}
\subsection{Notation} \label{notat}
In this section the monodromy data introduced in Section~\ref{sec1} is
calculated as $\tau \! \to \! +\infty$ for $\epsilon b \! > \! 0$
(corresponding to $\varepsilon_{1} \! = \! \varepsilon_{2} \! = \! 0)$:
this constitutes the first step towards the proof of the results stated in
Section~\ref{sec2}. This calculation consists of three components: (i) the
WKB analysis of the $\mu$-part of System~\eqref{eqFGmain}, that is,
\begin{equation} \label{eq3.1}
\partial_{\mu} \Psi (\mu) \! = \! \widetilde{\mathscr{U}}(\mu,\tau) \Psi
(\mu),\end{equation}
where $\Psi (\mu) \! = \! \Psi (\mu,\tau)$; (ii) the approximation of $\Psi
(\mu)$ in the neighborhoods of the turning points; and (iii) the matching of
these asymptotics.

Some solutions $u(\tau)$ of Equation~\eqref{eq1.1} might (and actually do)
have poles and zeroes located on the positive real semi-axis. In order to
be able to study such solutions, one must consider a slightly more general
complex domain $\widetilde{\mathcal{D}}_{u}$: this coincides with the
cheese-like domain $\mathcal{D}_{u}$ introduced in Remark~\ref{remark2.3};
however, since, \emph{a priori}, one does not know the solutions $u(\tau)$
that have such poles and zeroes, nor the exact locations of these poles and
zeroes (cf. Equations~\eqref{eqtheorem2.23} and~\eqref{eqtheorem2.26}), it is
necessary to introduce a formal definition for $\widetilde{\mathcal{D}}_{u}$.
Denote by $\mathcal{P}_{u}$ and $\mathcal{Z}_{u}$, respectively, the countable
sets of poles and zeroes of $u(\tau)$: as follows {}from the Painlev\'{e}
property, these sets might have accumulation points at $0$ and $\infty$.
Define neighborhoods of $\mathcal{P}_{u}$ and $\mathcal{Z}_{u}$,
respectively: for $\delta \! > \! 0$,
\begin{gather}
\mathcal{P}_{u}^{\delta} \! := \! \left\{\mathstrut \tau \! \in \! \mathbb{C}
\, \colon \, \lvert \phi (\tau) \! - \! \phi (\tau_{p}) \rvert \! > \! C
\lvert \tau_{p} \rvert^{-\delta}, \, \tau_{p} \! \in \! \mathcal{P}_{u}
\right\}, \label{eqnotat2} \\
\mathcal{Z}_{u}^{\delta} \! := \! \left\{\mathstrut \tau \! \in \! \mathbb{C}
\, \colon \, \lvert \phi (\tau) \! - \! \phi (\tau_{z}) \rvert \! > \! C
\lvert \tau_{z} \rvert^{-\delta}, \, \tau_{z} \! \in \! \mathcal{Z}_{u}
\right\}. \label{eqnotat3}
\end{gather}
Now, the cheese-like domain $\widetilde{\mathcal{D}}_{u}$ is defined:
\begin{equation}
\widetilde{\mathcal{D}}_{u} \! := \! \mathcal{D} \setminus
(\mathcal{P}_{u}^{\delta} \cup \mathcal{Z}_{u}^{\delta}), \label{eqnotat4}
\end{equation}
where $\mathcal{D}$ is given in Equation~\eqref{eqdofd}.
\begin{eeee} \label{remarknotat}
\textsl{Throughout this section, and for brevity of notation, the following
convention is adopted: in asymptotics of all expressions and formulae
depending on $u(\tau)$, the ``notation'' $\tau \! \to \! +\infty$ means $\Re
(\tau) \! \to \! +\infty$ and $\tau \! \in \! \widetilde{\mathcal{D}}_{u}$.}
\end{eeee}
Further notation used throughout this section is now summarized:
\begin{enumerate}
\item[(1)] $\mathrm{I} \! = \!
\left(
\begin{smallmatrix}
1 & 0 \\
0 & 1
\end{smallmatrix}
\right)$ is the $2 \times 2$ identity matrix, $\sigma_{1} \! = \!
\left(
\begin{smallmatrix}
0 & 1 \\
1 & 0
\end{smallmatrix}
\right)$, $\sigma_{2} \! = \!
\left(
\begin{smallmatrix}
0 & -\mi \\
\mi & 0
\end{smallmatrix}
\right)$, and $\sigma_{3} \! = \!
\left(
\begin{smallmatrix}
1 & 0 \\
0 & -1
\end{smallmatrix}
\right)$ are the Pauli matrices, $\sigma_{\pm} \! := \! \tfrac{1}{2}
(\sigma_{1} \! \pm \! \mi \sigma_{2})$ are the raising $(+)$ and lowering
$(-)$ matrices, and $\mathbb{R}_{\pm} \! := \! \lbrace \mathstrut x \! \in
\! \mathbb{R} \colon \pm x \! > \! 0 \rbrace$;
\item[(2)] for $(\varsigma_{1},\varsigma_{2}) \! \in \! \mathbb{R} \times
\mathbb{R}$, the function $(z \! - \! \varsigma_{1})^{\mi \varsigma_{2}}
\colon \mathbb{C} \setminus (-\infty,\varsigma_{1}) \! \to \! \mathbb{C}$, $z
\! \mapsto \! \exp (\mi \varsigma_{2} \ln (z \! - \! \varsigma_{1}))$, with
the branch cut taken along $(-\infty,\varsigma_{1})$ and the principal branch
of the logarithm chosen;
\item[(3)] for a scalar $\omega_{o}$ and a $2 \times 2$ matrix $\hat{
\Upsilon}$, $\omega_{o}^{\ad (\sigma_{3})} \hat{\Upsilon }\! := \!
\omega_{o}^{\sigma_{3}} \hat{\Upsilon} \omega_{o}^{-\sigma_{3}}$;
\item[(4)] for a $2 \times 2$ matrix-valued function $\mathfrak{Y}(z)$,
$\mathfrak{Y}(z) \! =_{z \to z_{0}} \! \mathcal{O}(\pmb{\ast})$ (resp.,
$o(\pmb{\ast}))$ means $\mathfrak{Y}_{ij}(z) \! =_{z \to z_{0}} \! \mathcal{O}
(\pmb{\ast}_{ij})$ (resp., $o(\pmb{\ast}_{ij}))$, $i,j \! = \! 1,2$;
\item[(5)] for $\mathfrak{B}(\pmb{\cdot}) \! \in \! \mathrm{M}_{2}
(\mathbb{C})$, $\lvert \lvert \mathfrak{B}(\pmb{\cdot}) \rvert \rvert \! := \!
(\sum_{i,j=1}^{2} \mathfrak{B}_{ij}(\pmb{\cdot}) \overline{\mathfrak{B}_{ij}
(\pmb{\cdot})})^{1/2}$ denotes the Hilbert-Schmidt norm, where
$\overline{\pmb{\star}}$ denotes complex conjugation of $\pmb{\star}$;
\item[(6)] for some $\delta_{\ast} \! > \! 0$ and sufficiently small,
$\mathscr{O}_{\delta_{\ast}}(p)$ denotes the $\delta_{\ast}$-neighborhood of
the point $p$, that is, $\mathscr{O}_{\delta_{\ast}}(p) \! := \! \lbrace
\mathstrut z \! \in \! \mathbb{C} \colon \, \lvert z \! - \! p \rvert \! <
\! \delta_{\ast} \rbrace$.
\end{enumerate}
\subsection{WKB Analysis} \label{subsec3.1}
This subsection is devoted to the WKB analysis of Equation~\eqref{eq3.1} as
$\tau \! \to \! +\infty$ (and $\epsilon b \! > \! 0)$. In order to transform
Equation~\eqref{eq3.1} into a form amenable to WKB analysis, one uses the
result of Proposition~4.1.1 in \cite{a1}, which is summarized here for the
reader's convenience.
\begin{bbbb}[\textrm{\cite{a1}}] \label{prop3.1.1}
Let
\begin{equation} \label{eq3.2}
\begin{gathered}
A(\tau) \! = \! a(\tau) \tau^{-2/3}, \quad B(\tau) \! = \! b(\tau)
\tau^{-2/3}, \quad C(\tau) \! = \! c(\tau) \tau^{-1/3}, \quad D(\tau) \!
= \! d(\tau) \tau^{-1/3}, \\
\widetilde{\mu} \! = \! \mu \tau^{1/6}, \quad \quad \widetilde{\Psi}
(\widetilde{\mu}) \! := \! \tau^{-(1/12) \sigma_{3}} \Psi(\widetilde{\mu}
\tau^{-1/6}),
\end{gathered}
\end{equation}
where $\widetilde{\Psi}(\widetilde{\mu}) \! = \! \widetilde{\Psi}
(\widetilde{\mu},\tau)$. Then
\begin{equation} \label{eq3.3}
\partial_{\widetilde{\mu}} \widetilde{\Psi}(\widetilde{\mu}) \! = \!
\tau^{2/3} \mathcal{A}(\widetilde{\mu},\tau) \widetilde{\Psi}(\widetilde{\mu}),
\end{equation}
where
\begin{equation} \label{eq3.4}
\mathcal{A}(\widetilde{\mu},\tau) \! := \! -2 \mi \widetilde{\mu} \sigma_{3}
\! + \!
\begin{pmatrix}
0 & -\frac{4 \mi \sqrt{-a(\tau)b(\tau)}}{b(\tau)} \\
-2d(\tau) & 0
\end{pmatrix} \! - \! \dfrac{\mi r(\tau)(\epsilon b)^{1/3}}{2 \widetilde{\mu}}
\sigma_{3} \! + \! \dfrac{1}{\widetilde{\mu}^{2}}
\begin{pmatrix}
0 & \frac{\mi \epsilon b}{b(\tau)} \\
\mi b(\tau) & 0
\end{pmatrix},
\end{equation}
with
\begin{equation} \label{eq3.5}
\dfrac{\mi r(\tau)(\epsilon b)^{1/3}}{2} \! = \! \left(
\mi a \! + \! \dfrac{1}{2} \right) \tau^{-2/3} \! + \!
\dfrac{2a(\tau)d(\tau)}{\sqrt{-a(\tau)b(\tau)}}.
\end{equation}
\end{bbbb}

It is now important to state under what conditions the WKB analysis of
Equation~\eqref{eq3.3} is undertaken. Define the functions $\hat{r}_{0}
(\tau)$, $\hat{u}_{0}(\tau)$, and $h_{0}(\tau)$ via the following
equations\footnote{Equations~(63)--(66) in \cite{a1} are expressed
in terms of the functions $r_{0}(\tau)$ and $u_{0}(\tau)$, whereas
Equations~\eqref{eq3.8}--\eqref{eq3.10} and conditions~\eqref{eq3.11} are
expressed in terms of the functions $\hat{r}_{0}(\tau) \! := \! r_{0}(\tau)
\tau^{-1/3}$ and $\hat{u}_{0}(\tau) \! := \! u_{0}(\tau) \tau^{-1/3}$.}:
\begin{gather}
\sqrt{-a(\tau)b(\tau)} + \! c(\tau)d(\tau) \! + \! \dfrac{a(\tau)d(\tau)
\tau^{-2/3}}{2 \sqrt{-a(\tau)b(\tau)}} \! - \! \dfrac{1}{4} \left(a \! -
\! \dfrac{\mi}{2} \right)^{2} \tau^{-4/3} \! = \! \dfrac{3}{4}(\epsilon
b)^{2/3} \! - \! h_{0}(\tau) \tau^{-2/3}, \label{eq3.8} \\
r(\tau) \! = \! -2 \! + \! \hat{r}_{0}(\tau), \label{eq3.9} \\
\sqrt{-a(\tau)b(\tau)} = \! \dfrac{(\epsilon b)^{2/3}}{2}(1 \! + \!
\hat{u}_{0}(\tau)), \label{eq3.10}
\end{gather}
and assume that the functions $\hat{r}_{0}(\tau)$, $\hat{u}_{0}(\tau)$,
and $h_{0}(\tau)$, which are holomorphic in the domain
$\widetilde{\mathcal{D}}_{u}$, satisfy the conditions
\begin{equation}
\begin{gathered}
\mathcal{O}(\tau^{-\frac{1}{3}+ \delta_{1}}) \underset{\tau \to +\infty}{
\leqslant} \lvert \hat{r}_{0}(\tau) \rvert \underset{\tau \to +\infty}{
\leqslant} \mathcal{O}(\tau^{\delta}), \qquad \mathcal{O}(\tau^{-\frac{1}{3}+
\delta_{1}}) \underset{\tau \to +\infty}{\leqslant} \lvert \hat{u}_{0}(\tau)
\rvert \underset{\tau \to +\infty}{\leqslant} \mathcal{O}(\tau^{2 \delta}), \\
\dfrac{1}{\lvert 1 \! + \! \hat{u}_{0}(\tau) \rvert} \underset{\tau \to
+\infty}{\leqslant} \mathcal{O}(\tau^{\delta}), \qquad \lvert h_{0}(\tau)
\rvert \underset{\tau \to +\infty}{\leqslant} \mathcal{O}(\tau^{\delta}),
\qquad 0 \! < \! \delta_{1} \! < \! \dfrac{1}{3}, \qquad 0 \! \leqslant \!
\delta \! < \! \dfrac{1}{39}, \label{eq3.11}
\end{gathered}
\end{equation}
where the parameter $\delta$ is the one that was introduced in the
definition of the cheese-like domain $\widetilde{\mathcal{D}}_{u}$ (cf.
Equations~\eqref{eqnotat2}, \eqref{eqnotat3}, and~\eqref{eqnotat4}).
\begin{eeee} \label{rem3.1.1}
\textsl{Even though, at this juncture, the upper bound $1/39$ for the growth
exponent, $\delta$, given in conditions~\eqref{eq3.11} might seem artificial,
it must be noted that it arises naturally during the course of the ensuing
asymptotic analysis.}
\end{eeee}
It is also assumed that the functions $a(\tau)$, $b(\tau)$, $c(\tau)$, and
$d(\tau)$ are related via the ``integral of motion'' associated with the
underlying Hamiltonian structure of Equation~\eqref{eq1.1} (see \cite{a1},
Lemma~2.1):
\begin{equation}
a(\tau)d(\tau) \! + \! b(\tau)c(\tau) \! + \! \mi a \sqrt{-a(\tau)b(\tau)} \,
\tau^{-2/3} \! = \! -\dfrac{\mi \epsilon b}{2}, \quad \epsilon \! = \! \pm 1.
\label{eq3.7}
\end{equation}
\begin{eeee}
\textsl{It is worth noting that Equations~\eqref{eq3.8}--\eqref{eq3.10}
and~\eqref{eq3.7} are self-consistent; in fact, they are equivalent to}
\begin{align}
a(\tau)d(\tau) =& \, \dfrac{(\epsilon b)^{2/3}}{2}(1 \! + \! \hat{u}_{0}
(\tau)) \! \left(-\dfrac{\mi (\epsilon b)^{1/3}}{2} \! + \! \dfrac{\mi
(\epsilon b)^{1/3} \hat{r}_{0}(\tau)}{4} \! - \! \dfrac{\mi}{2} \left(a \!
- \! \dfrac{\mi}{2} \right) \tau^{-2/3} \right), \label{eq3.12} \\
b(\tau)c(\tau) =& \, \dfrac{(\epsilon b)^{2/3}}{2}(1 \! + \! \hat{u}_{0}
(\tau)) \! \left(-\dfrac{\mi (\epsilon b)^{1/3}}{2} \! + \! \mi (\epsilon
b)^{1/3} \left(\dfrac{\hat{u}_{0}(\tau)}{1 \! + \! \hat{u}_{0}(\tau)} \!
- \! \dfrac{\hat{r}_{0}(\tau)}{4} \right) \! - \! \dfrac{\mi}{2} \left(a
\! + \! \dfrac{\mi}{2} \right) \tau^{-2/3} \right), \label{eq3.13} \\
-h_{0}(\tau) \tau^{-2/3}=& \, \dfrac{(\epsilon b)^{2/3}}{2} \! \left(\dfrac{
(\hat{u}_{0}(\tau))^{2} \! + \! \frac{1}{2} \hat{u}_{0}(\tau) \hat{r}_{0}
(\tau)}{1 \! + \! \hat{u}_{0}(\tau)} \! - \! \dfrac{(\hat{r}_{0}(\tau))^{
2}}{8} \right) \! + \! \dfrac{(\epsilon b)^{1/3}}{2} \left(a \! - \!
\dfrac{\mi}{2} \right) \dfrac{\tau^{-2/3}}{(1 \! + \! \hat{u}_{0}(\tau))};
\label{eq3.14}
\end{align}
\textsl{furthermore, via Equations~\eqref{eq3.10}, \eqref{eq3.12},
and~\eqref{eq3.13}},
\begin{align}
c(\tau)d(\tau) =& \, -\left(\dfrac{\mi (\epsilon b)^{1/3}}{2} \! - \!
\dfrac{\mi (\epsilon b)^{1/3} \hat{r}_{0}(\tau)}{4} \! + \! \dfrac{\mi}{2}
\left(a \! - \! \dfrac{\mi}{2} \right) \tau^{-2/3} \right) \nonumber \\
\times& \, \left(\dfrac{\mi (\epsilon b)^{1/3}}{2} \! - \! \mi (\epsilon
b)^{1/3} \left(\dfrac{\hat{u}_{0}(\tau)}{1 \! + \! \hat{u}_{0}(\tau)} \! - \!
\dfrac{\hat{r}_{0}(\tau)}{4} \right) \! + \! \dfrac{\mi}{2} \left(a \! + \!
\dfrac{\mi}{2} \right) \tau^{-2/3} \right). \label{eq3.15}
\end{align}
\end{eeee}

In certain domains of the complex $\widetilde{\mu}$-plane (see the discussion
below), the leading term of asymptotics (as $\tau \! \to \! +\infty)$ of
a fundamental solution, $\widetilde{\Psi}(\widetilde{\mu})$, of
Equation~\eqref{eq3.3} is given by the following WKB formula\footnote{For
simplicity of notation, the $\tau$ dependencies of
$\widetilde{\Psi}_{\text{WKB}}(\widetilde{\mu})$, $l(\widetilde{\mu})$,
and $T(\widetilde{\mu})$ are suppressed.} (see, for example, Chapter~5 of
\cite{F}),
\begin{equation}
\widetilde{\Psi}_{\text{WKB}}(\widetilde{\mu}) \! = \! T(\widetilde{\mu})
\exp \! \left(-\sigma_{3} \mi \tau^{2/3} \int_{}^{\widetilde{\mu}}l(\xi)
\, \md \xi \! - \! \int_{}^{\widetilde{\mu}} \text{diag}((T(\xi))^{-1}
\partial_{\xi} T(\xi)) \, \md \xi \right), \label{eq3.16}
\end{equation}
where
\begin{equation}
l(\widetilde{\mu}) \! := \! \sqrt{\det (\mathcal{A}(\widetilde{\mu}))},
\label{eq3.17}
\end{equation}
and $T(\widetilde{\mu})$, which diagonalizes $\mathcal{A}(\widetilde{\mu})$,
that is, $(T(\widetilde{\mu}))^{-1} \mathcal{A}(\widetilde{\mu})
T(\widetilde{\mu}) \! = \! -\mi l(\widetilde{\mu}) \sigma_{3}$, is given by
\begin{equation}
T(\widetilde{\mu}) \! = \! \dfrac{\mi}{\sqrt{2 \mi l(\widetilde{\mu})
(\mathcal{A}_{11}(\widetilde{\mu}) \! - \! \mi l(\widetilde{\mu}))}} \left(
\mathcal{A}(\widetilde{\mu}) \! - \! \mi l(\widetilde{\mu}) \sigma_{3}
\right) \sigma_{3}, \label{eq3.18}
\end{equation}
where $\mathcal{A}_{11}(\widetilde{\mu})$ is the $(1 \, 1)$-element of the
matrix $\mathcal{A}(\widetilde{\mu})$ (cf. Equation~\eqref{eq3.4}).
\begin{bbbb} \label{prop3.1.2}
Let $T(\widetilde{\mu})$ be given in Equation~\eqref{eq3.18}, with
$\mathcal{A}(\widetilde{\mu})$ and $l(\widetilde{\mu})$ defined by
Equations~\eqref{eq3.4} and~\eqref{eq3.17}, respectively. Then $\det
(T(\widetilde{\mu})) \! = \! 1$, and $\tr ((T(\widetilde{\mu}))^{-1}
\partial_{\widetilde{\mu}}T(\widetilde{\mu})) \! = \! 0;$ furthermore,
\begin{equation}
\diag ((T(\widetilde{\mu}))^{-1} \partial_{\widetilde{\mu}}T(\widetilde{\mu}))
\! = \! -\dfrac{1}{2} \left(\dfrac{\mathcal{A}_{12}(\widetilde{\mu})
\partial_{\widetilde{\mu}} \mathcal{A}_{21}(\widetilde{\mu}) \! - \!
\mathcal{A}_{21}(\widetilde{\mu}) \partial_{\widetilde{\mu}} \mathcal{A}_{12}
(\widetilde{\mu})}{2 \mi l(\widetilde{\mu}) \mathcal{A}_{11}(\widetilde{\mu})
\! + \! 2l^{2}(\widetilde{\mu})} \right) \sigma_{3}. \label{eq3.23}
\end{equation}
\end{bbbb}

\emph{Proof}. Set $T(\widetilde{\mu}) \! = \! \left(
\begin{smallmatrix}
T_{11}(\widetilde{\mu}) & T_{12}(\widetilde{\mu}) \\
T_{21}(\widetilde{\mu}) & T_{22}(\widetilde{\mu})
\end{smallmatrix}
\right)$. {}From the formula for $T(\widetilde{\mu})$ given in
Equation~\eqref{eq3.18}, with $\mathcal{A}(\widetilde{\mu})$ and
$l(\widetilde{\mu})$ defined by Equations~\eqref{eq3.4} and~\eqref{eq3.17},
respectively, one shows at
\begin{equation}
\begin{gathered}
T_{11}(\widetilde{\mu}) \! = \! T_{22}(\widetilde{\mu}) \! = \! \dfrac{\mi
(\mathcal{A}_{11}(\widetilde{\mu}) \! - \! \mi l(\widetilde{\mu}))}{\sqrt{2
\mi l(\mu)(\mathcal{A}_{11}(\widetilde{\mu}) \! - \! \mi
l(\widetilde{\mu}))}}, \qquad \qquad T_{12}(\widetilde{\mu}) \! = \! -\dfrac{
\mi \mathcal{A}_{12}(\widetilde{\mu})}{\sqrt{2 \mi l(\mu)(\mathcal{A}_{11}
(\widetilde{\mu}) \! - \! \mi l(\widetilde{\mu}))}}, \\
T_{21}(\widetilde{\mu}) \! = \! \dfrac{\mi \mathcal{A}_{21}(\widetilde{\mu})}{
\sqrt{2 \mi l(\mu)(\mathcal{A}_{11}(\widetilde{\mu}) \! - \! \mi l
(\widetilde{\mu}))}}, \label{eqintertee}
\end{gathered}
\end{equation}
whence (as $\det (\mathcal{A}(\widetilde{\mu})) \! = \!
l^{2}(\widetilde{\mu})$ and $\tr (\mathcal{A}(\widetilde{\mu})) \! = \! 0)$
$\det (T(\widetilde{\mu})) \! = \! (T_{11}(\widetilde{\mu}))^{2} \! - \!
T_{12}(\widetilde{\mu})T_{21}(\widetilde{\mu}) \! = \! 1$. Since $\det
(T(\widetilde{\mu})) \! = \! 1$, it follows that
\begin{equation*}
(T(\widetilde{\mu}))^{-1} \partial_{\widetilde{\mu}}T(\widetilde{\mu}) \! = \!
\begin{pmatrix}
T_{11}(\widetilde{\mu}) \partial_{\widetilde{\mu}}T_{11}(\widetilde{\mu})-
T_{12}(\widetilde{\mu}) \partial_{\widetilde{\mu}}T_{21}(\widetilde{\mu}) &
T_{11}(\widetilde{\mu}) \partial_{\widetilde{\mu}}T_{12}(\widetilde{\mu})-
T_{12}(\widetilde{\mu}) \partial_{\widetilde{\mu}}T_{11}(\widetilde{\mu}) \\
T_{11}(\widetilde{\mu}) \partial_{\widetilde{\mu}}T_{21}(\widetilde{\mu})-
T_{21}(\widetilde{\mu}) \partial_{\widetilde{\mu}}T_{11}(\widetilde{\mu}) &
T_{11}(\widetilde{\mu}) \partial_{\widetilde{\mu}}T_{11}(\widetilde{\mu})-
T_{21}(\widetilde{\mu}) \partial_{\widetilde{\mu}}T_{12}(\widetilde{\mu})
\end{pmatrix},
\end{equation*}
whence
\begin{equation*}
\diag ((T(\widetilde{\mu}))^{-1} \partial_{\widetilde{\mu}}T(\widetilde{\mu}))
\! = \!
\begin{pmatrix}
\frac{1}{2} \partial_{\widetilde{\mu}}(T_{11}(\widetilde{\mu}))^{2}-T_{12}
(\widetilde{\mu}) \partial_{\widetilde{\mu}}T_{21}(\widetilde{\mu}) & 0 \\
0 & \frac{1}{2} \partial_{\widetilde{\mu}}(T_{11}(\widetilde{\mu}))^{2}-
T_{21}(\widetilde{\mu}) \partial_{\widetilde{\mu}}T_{12}(\widetilde{\mu})
\end{pmatrix},
\end{equation*}
which implies that $\tr ((T(\widetilde{\mu}))^{-1} \partial_{\widetilde{\mu}}
T(\widetilde{\mu})) \! = \! \partial_{\widetilde{\mu}}((T_{11}(\widetilde{\mu}
))^{2} \! - \! T_{12}(\widetilde{\mu})T_{21}(\widetilde{\mu})) \! = \!
\partial_{\widetilde{\mu}}(1) \! = \! 0$. Writing out explicitly $\diag
((T(\widetilde{\mu}))^{-1} \partial_{\widetilde{\mu}}T(\widetilde{\mu}))$,
one shows that
\begin{align*}
((T(\widetilde{\mu}))^{-1} \partial_{\widetilde{\mu}}T(\widetilde{\mu}))_{11}
=& \, \dfrac{\mi}{4} \! \left(\dfrac{\partial_{\widetilde{\mu}}
\mathcal{A}_{11}(\widetilde{\mu})}{l(\widetilde{\mu})} \! - \! \dfrac{
\mathcal{A}_{11}(\widetilde{\mu}) \partial_{\widetilde{\mu}}l(\widetilde{
\mu})}{l^{2}(\widetilde{\mu})} \right) \! - \! \dfrac{\mathcal{A}_{12}
(\widetilde{\mu}) \partial_{\widetilde{\mu}} \mathcal{A}_{21}(\widetilde{
\mu})}{(2 \mi l(\widetilde{\mu}) \mathcal{A}_{11}(\widetilde{\mu}) \! + \!
2l^{2}(\widetilde{\mu}))} \\
+& \, \dfrac{\mathcal{A}_{12}(\widetilde{\mu}) \mathcal{A}_{21}(\widetilde{
\mu}) \left(\mi (l(\widetilde{\mu}) \partial_{\widetilde{\mu}} \mathcal{A}_{
11}(\widetilde{\mu}) \! + \! \mathcal{A}_{11}(\widetilde{\mu}) \partial_{
\widetilde{\mu}}l(\widetilde{\mu})) \! + \! 2l(\widetilde{\mu}) \partial_{
\widetilde{\mu}}l(\widetilde{\mu}) \right)}{(2 \mi l(\widetilde{\mu})
\mathcal{A}_{11}(\widetilde{\mu}) \! + \! 2l^{2}(\widetilde{\mu}))^{2}}, \\
((T(\widetilde{\mu}))^{-1} \partial_{\widetilde{\mu}}T(\widetilde{\mu}))_{22}
=& \, \dfrac{\mi}{4} \! \left(\dfrac{\partial_{\widetilde{\mu}} \mathcal{A}_{
11}(\widetilde{\mu})}{l(\widetilde{\mu})} \! - \! \dfrac{\mathcal{A}_{11}
(\widetilde{\mu}) \partial_{\widetilde{\mu}}l(\widetilde{\mu})}{l^{2}
(\widetilde{\mu})} \right) \! - \! \dfrac{\mathcal{A}_{21}(\widetilde{\mu})
\partial_{\widetilde{\mu}} \mathcal{A}_{12}(\widetilde{\mu})}{(2 \mi
l(\widetilde{\mu}) \mathcal{A}_{11}(\widetilde{\mu}) \! + \! 2l^{2}
(\widetilde{\mu}))} \\
+& \, \dfrac{\mathcal{A}_{12}(\widetilde{\mu}) \mathcal{A}_{21}(\widetilde{
\mu}) \left(\mi (l(\widetilde{\mu}) \partial_{\widetilde{\mu}} \mathcal{A}_{
11}(\widetilde{\mu}) \! + \! \mathcal{A}_{11}(\widetilde{\mu}) \partial_{
\widetilde{\mu}}l(\widetilde{\mu})) \! + \! 2l(\widetilde{\mu}) \partial_{
\widetilde{\mu}}l(\widetilde{\mu}) \right)}{(2 \mi l(\widetilde{\mu})
\mathcal{A}_{11}(\widetilde{\mu}) \! + \! 2l^{2}(\widetilde{\mu}))^{2}}:
\end{align*}
now, using the fact that $\tr ((T(\widetilde{\mu}))^{-1} \partial_{
\widetilde{\mu}}T(\widetilde{\mu})) \! = \! 0$, one arrives at
Equation~\eqref{eq3.23}. \hfill $\qed$
\begin{ffff} \label{cor3.1.1}
Let $\widetilde{\Psi}_{\mathrm{WKB}}(\widetilde{\mu})$ be given in
Equation~\eqref{eq3.16}, with $l(\widetilde{\mu})$ defined by
Equation~\eqref{eq3.17} and $T(\widetilde{\mu})$ given in
Equation~\eqref{eq3.18}. Then $\det (\widetilde{\Psi}_{\mathrm{WKB}}
(\widetilde{\mu})) \! = \! 1$.
\end{ffff}

The domains in the complex $\widetilde{\mu}$-plane where
Equation~\eqref{eq3.16} gives the asymptotic approximation of solutions to
Equation~\eqref{eq3.3} are defined in terms of the \emph{Stokes graph} (see,
for example, \cite{W,F}). The vertices of the Stokes graph are the singular
points of Equation~\eqref{eq3.3}, that is, $\widetilde{\mu} \! = \! 0$ and
$\widetilde{\mu} \! = \! \infty$, and the \emph{turning points}, which are
the roots of the equation $l^{2}(\widetilde{\mu}) \! = \! 0$. The edges
of the Stokes graph are the \emph{Stokes curves}, defined as $\Re
(\smallint_{\widetilde{\mu}_{\mathrm{TP}}}^{\widetilde{\mu}}l(\xi) \, \md
\xi) \! = \! 0$, where $\widetilde{\mu}_{\mathrm{TP}}$ denotes a turning
point. \emph{Canonical domains} are those domains in the complex
$\widetilde{\mu}$-plane containing one, and only one, Stokes curve and bounded
by two adjacent Stokes curves. Note that the restriction of any branch of
$l(\widetilde{\mu})$ to a canonical domain is a single-valued function. In
each canonical domain, for any choice of the branch of $l(\widetilde{\mu})$,
there exists a fundamental solution of Equation ~\eqref{eq3.3} which has
asymptotics whose leading term is given by Equation~\eqref{eq3.16}. {}From
the definition of $l(\widetilde{\mu})$ given by Equation~\eqref{eq3.17}, one
arrives at
\begin{equation}
l^{2}(\widetilde{\mu}) \! = \! \dfrac{4}{\widetilde{\mu}^{4}} \!
\left(\left(\widetilde{\mu}^{2} \! - \! \alpha^{2} \right)^{2} \left(
\widetilde{\mu}^{2} \! + \! 2 \alpha^{2} \right) \! + \! \left(
\widetilde{\mu}^{4} \left(a \! - \! \dfrac{\mi}{2} \right) \! + \!
\widetilde{\mu}^{2}h_{0}(\tau) \right) \tau^{-2/3} \right), \label{eq3.19}
\end{equation}
where
\begin{equation*}
\alpha \! := \! \dfrac{(\epsilon b)^{1/6}}{\sqrt{2}} \! > \! 0.
\end{equation*}
It follows {}from Equation~\eqref{eq3.19} that there are six turning points:
two turning points coalesce (as $\tau \! \to \! +\infty$) at $\alpha \! + \!
\mathcal{O}(\tau^{-\frac{1}{3}+\frac{\delta}{2}})$, another pair coalesce at
$-\alpha \! + \! \mathcal{O}(\tau^{-\frac{1}{3}+\frac{\delta}{2}})$, and the
remaining two turning points approach (as $\tau \! \to \! +\infty$) $\pm
\mi \sqrt{2} \, \alpha \! + \! \mathcal{O}(\tau^{-\frac{2}{3}+\delta})$,
respectively. Denote by $\widetilde{\mu}_{1}$ the turning point in the first
quadrant of the complex $\widetilde{\mu}$-plane which approaches $\alpha
\! + \! \mathcal{O}(\tau^{-\frac{1}{3}+\frac{\delta}{2}})$, and by
$\widetilde{\mu}_{2}$ the pure imaginary turning point which approaches $\mi
\sqrt{2} \, \alpha \! + \! \mathcal{O}(\tau^{-\frac{2}{3}+\delta})$. Denote
by $\mathscr{G}_{1}$ the part of the Stokes graph in the first quadrant of the
complex $\widetilde{\mu}$-plane which consists of the vertices $0$, $\infty$,
$\widetilde{\mu}_{1}$, and $\widetilde{\mu}_{2}$, and the edges $(+\mi \infty,
\widetilde{\mu}_{2})$, $(0,\widetilde{\mu}_{1})$, $(\widetilde{\mu}_{2},
\widetilde{\mu}_{1})$, and $(\widetilde{\mu}_{1},+\infty)$: the complete
Stokes graph is the union of the mirror images of $\mathscr{G}_{1}$ with
respect to the real and imaginary axes in the complex $\widetilde{\mu}$-plane.
\begin{bbbb} \label{prop3.1.3}
Let $l^{2}(\widetilde{\mu})$ be given in Equation~\eqref{eq3.19}. Then
\begin{equation}
\int_{\widetilde{\mu}_{0}}^{\widetilde{\mu}}l(\xi) \, \md \xi \underset{\tau
\to +\infty}{=} \varUpsilon_{\tau}(\widetilde{\mu}) \! - \! \varUpsilon_{\tau}
(\widetilde{\mu}_{0}) \! + \! \mathcal{O}(\mathscr{E}_{l}(\widetilde{\mu})),
\label{eq3.21}
\end{equation}
where $\widetilde{\mu}_{0} \! \in \! \mathbb{C} \setminus
(\mathscr{O}_{\tau^{-\frac{1}{3}+ \frac{\delta}{2}}}(\pm \alpha) \cup
\mathscr{O}_{\tau^{-\frac{2}{3}+ \delta}}(\pm \mi \sqrt{2} \, \alpha)
\cup \lbrace 0,\infty \rbrace)$ and the path of integration lies in
the corresponding\footnote{The same canonical domain to which
$\widetilde{\mu}_{0}$ belongs.} canonical domain,
\begin{align}
\varUpsilon_{\tau}(\xi) :=& \, \left(\xi \! + \! \dfrac{2 \alpha^{2}}{\xi}
\right) \sqrt{\xi^{2} \! + \! 2 \alpha^{2}} + \! \tau^{-2/3} \left(a \! - \!
\dfrac{\mi}{2} \right) \ln \left(\xi \! + \! \sqrt{\xi^{2} \! + \! 2
\alpha^{2}} \, \right) \nonumber \\
+& \, \dfrac{\tau^{-2/3}}{2 \sqrt{3}} \left(\left(a \! - \! \dfrac{\mi}{2}
\right) \! + \! \dfrac{h_{0}(\tau)}{\alpha^{2}} \right) \ln \left(\left(
\dfrac{\sqrt{3} \, \sqrt{\xi^{2} \! + \! 2 \alpha^{2}} - \! \xi \! + \! 2
\alpha}{\sqrt{3} \, \sqrt{\xi^{2} \! + \! 2 \alpha^{2}} + \! \xi \! + \! 2
\alpha} \right) \! \left(\dfrac{\xi \! - \! \alpha}{\xi \! + \! \alpha}
\right) \right), \label{equpsi}
\end{align}
and $\mathscr{E}_{l}(\widetilde{\mu}) \! = \! \mathscr{E}_{l}^{\natural}
(\widetilde{\mu}) \! - \! \mathscr{E}_{l}^{\natural}(\widetilde{\mu}_{0})$,
with
\begin{equation}
\mathscr{E}_{l}^{\natural}(\xi) \! := \!
\begin{cases}
\dfrac{(c_{1}(\delta_{0})(h_{0}(\tau))^{2} \! + \! c_{2}(\delta_{0}))
\tau^{-4/3}}{(\xi \! \mp \! \alpha)^{2}}, &\text{$\xi \! \in \!
\mathscr{O}_{\delta_{0}}(\pm \alpha)$,} \\
\dfrac{(c_{3}(\delta_{0})(h_{0}(\tau))^{2} \! + \! c_{4}(\delta_{0}))
\tau^{-4/3}}{(\xi \! \mp \! \mi \sqrt{2} \, \alpha)^{1/2}}, &\text{$\xi \!
\in \! \mathscr{O}_{\delta_{0}}(\pm \mi \sqrt{2} \, \alpha)$,} \\
(c_{5}(\delta_{0})(h_{0}(\tau))^{2} \! + \! c_{6}(\delta_{0})) \tau^{-4/3},
&\text{$\xi \! \in \! \mathscr{O}_{\delta_{0}}(0)$,} \\
(c_{7}(\delta_{0})(h_{0}(\tau))^{2} \! + \! c_{8}(\delta_{0})) \tau^{-4/3},
&\text{$\xi \! \in \! \mathscr{O}_{\delta_{0}}(\infty)$,}
\end{cases} \label{eqeel}
\end{equation}
where $c_{k}(z)$, $k \! = \! 1,\dotsc,8$, are holomorphic functions of $z$ in
a neighborhood of $z \! = \! 0$ with $c_{k}(0) \! \neq \! 0$, and $\delta_{0}
\! > \! 0$ and sufficiently small.
\end{bbbb}

\emph{Proof}. Recall the expression for $l^{2}(\widetilde{\mu})$ $(= \! l^{2}
(\widetilde{\mu},\tau))$ given in Equation~\eqref{eq3.19}. Set
\begin{equation}
l^{2}_{\infty}(\widetilde{\mu}) \! := \! l^{2}(\widetilde{\mu},+\infty)
\! = \! 4 \widetilde{\mu}^{-4}(\widetilde{\mu}^{2} \! - \! \alpha^{2})^{2}
(\widetilde{\mu}^{2} \! + \! 2 \alpha^{2}). \label{eqlsquared}
\end{equation}
Define
\begin{equation}
\Delta_{\tau}(\widetilde{\mu}) \! := \! \dfrac{l^{2}(\widetilde{\mu}) \! -
\! l^{2}_{\infty}(\widetilde{\mu})}{l^{2}_{\infty}(\widetilde{\mu})} \! =
\! \dfrac{\widetilde{\mu}^{2}(h_{0}(\tau) \! + \! \widetilde{\mu}^{2}(a \!
- \! \frac{\mi}{2})) \tau^{-2/3}}{(\widetilde{\mu}^{2} \! - \! \alpha^{2})^{2}
(\widetilde{\mu}^{2} \! + \! 2 \alpha^{2})}. \label{eqdellt}
\end{equation}
Now, via conditions~\eqref{eq3.11}, and writing
\begin{equation}
l(\widetilde{\mu}) \! = \! l_{\infty}(\widetilde{\mu})(1 \! + \! \Delta_{\tau}
(\widetilde{\mu}))^{1/2} \underset{\tau \to +\infty}{=} l_{\infty}
(\widetilde{\mu}) \! \left(1 \! + \! \dfrac{1}{2} \Delta_{\tau}
(\widetilde{\mu}) \! + \! \mathcal{O}((\Delta_{\tau}(\widetilde{\mu}))^{2})
\right), \label{eqjustell}
\end{equation}
one arrives at
\begin{equation}
l(\widetilde{\mu}) \underset{\tau \to +\infty}{=} 2 \! \left(1 \! - \!
\dfrac{\alpha^{2}}{\widetilde{\mu}^{2}} \right) \sqrt{\widetilde{\mu}^{2}
\! + \! 2 \alpha^{2}}+ \dfrac{(h_{0}(\tau) \! + \! \widetilde{\mu}^{2}(a \!
- \! \frac{\mi}{2})) \tau^{-2/3}}{(\widetilde{\mu}^{2} \! - \! \alpha^{2})
\sqrt{\widetilde{\mu}^{2} \! + \! 2 \alpha^{2}}} \! + \! \mathcal{O} \! \left(
\dfrac{\widetilde{\mu}^{2}(h_{0}(\tau) \! + \! \widetilde{\mu}^{2}(a \! - \!
\frac{\mi}{2}))^{2} \tau^{-4/3}}{(\widetilde{\mu}^{2} \! - \! \alpha^{2})^{3}
(\widetilde{\mu}^{2} \! + \! 2 \alpha^{2})^{3/2}} \right). \label{eqforl}
\end{equation}
In order to obtain the leading term in Equation~\eqref{eq3.21}, one integrates
the first two terms in Equation~\eqref{eqforl}. Analogously, integrating the
error term in Equation~\eqref{eqforl}, one finds an explicit expression for
the function $\mathscr{E}_{l}^{\natural}(\widetilde{\mu})$: since this latter
expression is quite cumbersome, only its asymptotics at the singular and
turning points are presented in Equation~\eqref{eqeel}. \hfill $\qed$
\begin{ffff} \label{cor3.1.2}
Set $\widetilde{\mu}_{0} \! = \! \alpha \! + \! \tau^{-1/3}
\widetilde{\Lambda}$, where $\widetilde{\Lambda} \! =_{\tau \to +\infty} \!
\mathcal{O}(\tau^{\varepsilon})$, $0 \! < \! \delta \! < \! \varepsilon \!
< \! 1/9$. Then
\begin{equation}
\int_{\widetilde{\mu}_{0}}^{\widetilde{\mu}}l(\xi) \, \md \xi \underset{\tau
\to +\infty}{=} \varUpsilon_{\tau}(\widetilde{\mu}) \! + \! \varUpsilon_{
\tau}^{\sharp} \! + \! \mathcal{O}(\mathscr{E}_{l}^{\natural}(\widetilde{
\mu})) \! + \! \mathcal{O}(\tau^{-\frac{2}{3}-2(\varepsilon - \delta)}) \! +
\! \mathcal{O}(\tau^{-1+3 \varepsilon}), \label{eq3.22}
\end{equation}
where $\varUpsilon_{\tau}(\widetilde{\mu})$ and $\mathscr{E}_{l}^{\natural}
(\widetilde{\mu})$ are defined in Proposition~{\rm \ref{prop3.1.3}},
\begin{align}
\varUpsilon_{\tau}^{\sharp} :=& \, \mp 3 \sqrt{3} \, \alpha^{2} \! \mp \! 2
\sqrt{3} \, \tau^{-2/3} \widetilde{\Lambda}^{2} \! - \! \tau^{-2/3} \left(
a \! - \! \dfrac{\mi}{2} \right) \ln \left((\sqrt{3} \pm \! 1) \alpha
\me^{\mi \pi \mathfrak{s}(\pm)} \right) \nonumber \\
\mp& \, \dfrac{\tau^{-2/3}}{2 \sqrt{3}} \left(\left(a \! - \! \dfrac{\mi}{2}
\right) \! + \! \dfrac{h_{0}(\tau)}{\alpha^{2}} \right) \! \left(\ln
\widetilde{\Lambda} \! - \! \dfrac{1}{3} \ln \tau \! - \! \ln (3 \alpha)
\right), \label{equpsisharp}
\end{align}
and $\mathfrak{s}(\pm) \! = \! (1 \! \mp \! 1)/2$, with the upper (resp.,
lower) signs taken if the positive $(+)$ (resp., negative $(-))$ branch of
the square-root function $\sqrt{\xi^{2} \! + \! 2 \alpha^{2}}$ is chosen.
\end{ffff}

\emph{Proof}. Substituting $\widetilde{\mu}_{0}$ as the argument of the
functions $\varUpsilon_{\tau}(\xi)$ and $\mathscr{E}_{l}^{\natural}(\xi)$
(cf. Equation~\eqref{equpsi} and the first line of Equation~\eqref{eqeel},
respectively) and expanding with respect to (the small parameter) $\tau^{-1/3}
\widetilde{\Lambda}$, one arrives at the following estimates:
\begin{equation*}
\varUpsilon_{\tau}(\widetilde{\mu}_{0}) \underset{\tau \to +\infty}{=}
\varUpsilon_{\tau}^{\sharp} \! + \! \mathcal{O}(\tau^{-1}
\widetilde{\Lambda}^{3}) \quad \text{and} \quad \mathscr{E}_{l}^{\natural}
(\widetilde{\mu}_{0}) \underset{\tau \to +\infty}{=} \mathcal{O} \! \left(
\tau^{-2/3}(c_{1}(\delta_{0})(h_{0}(\tau))^{2} \! + \! c_{2}(\delta_{0}))
\widetilde{\Lambda}^{-2} \right),
\end{equation*}
where $\varUpsilon_{\tau}^{\sharp}$ is defined by
Equation~\eqref{equpsisharp}. The inequality $\delta \! < \! \varepsilon \!
< \! 1/9$ (see conditions~\eqref{eq3.11} for the function $h_{0}(\tau))$
is introduced in order to guarantee that the last two error estimates in
Equation~\eqref{eq3.22} are decaying after multiplication by the large
parameter $\tau^{2/3}$ $($cf. Equation~\eqref{eq3.16}$)$. \hfill $\qed$
\begin{bbbb} \label{prop3.1.4}
Let $T(\widetilde{\mu})$ be given in Equation~\eqref{eq3.18}, with
$\mathcal{A}(\widetilde{\mu})$ defined by Equation~\eqref{eq3.4}, and
$l^{2}(\widetilde{\mu})$ given in Equation~\eqref{eq3.19}. Then
\begin{equation}
\int_{\widetilde{\mu}_{0}}^{\widetilde{\mu}} \diag ((T(\xi))^{-1}
\partial_{\xi}T(\xi)) \, \md \xi \underset{\tau \to +\infty}{=} \left(
\mathcal{I}_{\tau}(\widetilde{\mu}) \! + \! \mathcal{O}(\mathcal{E}_{T}
(\widetilde{\mu})) \right) \sigma_{3}, \label{eqintTminus1T}
\end{equation}
where $\widetilde{\mu}_{0} \! \in \! \mathbb{C} \setminus
(\mathscr{O}_{\tau^{-\frac{1}{3}+ \frac{\delta}{2}}}(\pm \alpha) \cup
\mathscr{O}_{\tau^{-\frac{2}{3}+ \delta}}(\pm \mi \sqrt{2} \, \alpha)
\cup \lbrace 0,\infty \rbrace)$ and the path of integration lies in the
corresponding canonical domain,
\begin{equation}
\mathcal{I}_{\tau}(\widetilde{\mu}) \! = \! \dfrac{\mathfrak{p}_{\tau}}{8 \mi}
\left(\digamma_{\tau}(\widetilde{\mu}) \! - \! \digamma_{\tau}
(\widetilde{\mu}_{0}) \right), \label{eqItee}
\end{equation}
with
\begin{equation}
\mathfrak{p}_{\tau} \! := \! \dfrac{4 \mi}{\hat{r}_{0}(\tau)} \! \left(2(1
\! + \! \hat{u}_{0}(\tau)) \! + \! \dfrac{(-2 \! + \! \hat{r}_{0}(\tau))}{1
\! + \! \hat{u}_{0}(\tau)} \right), \label{eqpeetee}
\end{equation}
\begin{align}
\digamma_{\tau}(\xi) :=& \, \dfrac{1}{\sqrt{3}} \ln \! \left(\left(\dfrac{
\sqrt{3} \, \sqrt{\xi^{2} \! + \! 2 \alpha^{2}}-\xi \! + \! 2 \alpha}{\sqrt{3}
\, \sqrt{\xi^{2} \! + \! 2 \alpha^{2}}+\xi \! + \! 2 \alpha} \right) \! \left(
\dfrac{\xi \! - \! \alpha}{\xi \! + \! \alpha} \right) \right) \! + \!
\dfrac{4}{\sqrt{32 \! + \! (\hat{r}_{0}(\tau) \! - \! 4)^{2}}} \nonumber \\
\times& \, \left(\operatorname{arctanh} \left(\dfrac{4 \sqrt{\xi^{2} \! + \!
2 \alpha^{2}}}{\xi (\hat{r}_{0}(\tau) \! + \! \sqrt{32 \! + \! (\hat{r}_{0}
(\tau) \! - \! 4)^{2}})} \right) \! - \! \operatorname{arctanh} \left(
\dfrac{4 \sqrt{\xi^{2} \! + \! 2 \alpha^{2}}}{\xi (\hat{r}_{0}(\tau) \! -
\! \sqrt{32 \! + \! (\hat{r}_{0}(\tau) \! - \! 4)^{2}})} \right) \right.
\nonumber \\
+&\left. \, \operatorname{arctanh} \left(\dfrac{8(-2 \! + \! \hat{r}_{0}
(\tau)) \xi^{2} \! + \! \alpha^{2}(-2 \! + \! \hat{r}_{0}(\tau))^{2} \! + \!
12 \alpha^{2}}{\alpha^{2} \hat{r}_{0}(\tau) \sqrt{32 \! + \! (\hat{r}_{0}
(\tau) \! - \! 4)^{2}}} \right) \! \right), \label{eqFtee}
\end{align}
and $\mathcal{E}_{T}(\widetilde{\mu}) \! = \! \mathcal{E}_{T}^{\natural}
(\widetilde{\mu}) \! - \! \mathcal{E}_{T}^{\natural}(\widetilde{\mu}_{0})$,
where
\begin{equation}
\mathcal{E}_{T}^{\natural}(\xi) \! := \!
\begin{cases}
\frac{\tau^{-2/3} \ln (\xi \mp \alpha)}{(1+ \hat{u}_{0}(\tau))(C_{1}
(\tau)+ \hat{r}_{0}(\tau))} \! + \! \frac{(\hat{u}_{0}(\tau)+ \frac{
\hat{u}_{0}(\tau)+ \frac{1}{2} \hat{r}_{0}(\tau)}{1+ \hat{u}_{0}(\tau)})
(C_{2}(\tau)+h_{0}(\tau))}{\hat{r}_{0}(\tau)(\xi \mp \alpha)^{2} \tau^{2/3}},
&\!\!\!\text{$\lvert \xi \! \mp \! \alpha \rvert \! < \! \tau^{-\frac{1}{3}+
\delta_{1}} \! < \! \lvert \hat{r}_{0}(\tau) \rvert$,} \\
\frac{\tau^{-2/3}}{(1+ \hat{u}_{0}(\tau)) \hat{r}_{0}(\tau)} \! + \!
\frac{(\hat{u}_{0}(\tau)+ \frac{\hat{u}_{0}(\tau)+ \frac{1}{2} \hat{r}_{0}
(\tau)}{1+ \hat{u}_{0}(\tau)})(C_{3}(\tau)+h_{0}(\tau))}{(C_{4}(\tau)+
\hat{r}_{0}(\tau))(\xi \mp \mi \sqrt{2} \, \alpha)^{1/2} \tau^{2/3}},
\qquad\hat{r}_{0}(\tau) \! \neq \! 6,
&\!\!\!\text{$ \xi \! \in \!
\mathscr{O}_{\delta_{0}}(\pm \mi \sqrt{2} \, \alpha)$,} \\
\frac{\tau^{-2/3}}{(1+ \hat{u}_{0}(\tau)) \hat{r}_{0}(\tau)} \! + \!
(\hat{u}_{0}(\tau) \! + \! \tfrac{\hat{u}_{0}(\tau)+ \frac{1}{2} \hat{r}_{0}
(\tau)}{1+ \hat{u}_{0}(\tau)}) \frac{(\frac{C_{5}(\tau)}{\hat{r}_{0}(\tau)}+
\frac{C_{6}(\tau)}{(\hat{r}_{0}(\tau))^{3}})(C_{7}(\tau)+h_{0}(\tau))}{\tau^{2/3}}, 
&\!\!\!\text{$\xi \! \in \! \mathscr{O}_{\delta_{0}}(0)$,} \\
\frac{\tau^{-2/3}}{(1+ \hat{u}_{0}(\tau)) \hat{r}_{0}(\tau)} \! + \!
(\hat{u}_{0}(\tau) \! + \! \tfrac{\hat{u}_{0}(\tau)+ \frac{1}{2} \hat{r}_{0}
(\tau)}{1+ \hat{u}_{0}(\tau)}) \frac{(\frac{C_{8}(\tau)}{\hat{r}_{0}(\tau)}+
\frac{C_{9}(\tau)}{(\hat{r}_{0}(\tau))^{3}})(C_{10}(\tau)+h_{0}(\tau))}{\tau^{2/3}}, 
&\!\!\!\text{$\xi \! \in \! \mathscr{O}_{\delta_{0}}(\infty)$,}
\end{cases} \label{eqee2}
\end{equation}
with\footnote{In case $\hat{r}_{0}(\tau) \! \to \! 6$, the factor $(\xi \!
\mp \! \mi \sqrt{2} \, \alpha)^{1/2}$ which appears in the second line of
Equation~\eqref{eqee2} should be changed to $\xi \! \mp \! \mi \sqrt{2} \,
\alpha$.} $C_{j}(\tau) \! =_{\tau \to +\infty} \! \mathcal{O}(1)$, $j \! =
\! 1,\dotsc,10$.
\end{bbbb}

\emph{Proof}. {}From Equation~\eqref{eq3.4} and
Equations~\eqref{eqlsquared}--\eqref{eqjustell}, one shows that, via
Equation~\eqref{eq3.9} (cf. Equation~\eqref{eq3.23}),
\begin{equation}
2 \mi l(\xi) \mathcal{A}_{11}(\xi) \! + \! 2l^{2}(\xi) \underset{\tau
\to +\infty}{=} \mathcal{P}_{\infty}(\xi) \! + \! \mathcal{P}_{1}(\xi)
\Delta_{\tau}(\xi) \! + \! \mathcal{O} \! \left(\xi^{-1}l_{\infty}(\xi)
\hat{r}_{0}(\tau)(\Delta_{\tau}(\xi))^{2} \right), \label{eqdenomTminus1T}
\end{equation}
where
\begin{align}
\mathcal{P}_{\infty}(\xi) :=& \, 2l_{\infty}^{2}(\xi) \! + \! 2l_{\infty}(\xi)
\! \left(2 \xi \! + \! \dfrac{\alpha^{2}(-2 \! + \! \hat{r}_{0}(\tau))}{\xi}
\right), \label{eq3.24} \\
\mathcal{P}_{1}(\xi) :=& \, 2l_{\infty}^{2}(\xi) \! + \! l_{\infty}(\xi) \!
\left(2 \xi \! + \! \dfrac{\alpha^{2}(-2 \! + \! \hat{r}_{0}(\tau))}{\xi}
\right), \label{eq3.25}
\end{align}
and, via Equations~\eqref{eq3.4}, \eqref{eq3.10}, and~\eqref{eq3.12} (cf.
Equation~\eqref{eq3.23}),
\begin{equation}
\mathcal{A}_{12}(\xi) \partial_{\xi} \mathcal{A}_{21}(\xi) \! - \!
\mathcal{A}_{21}(\xi) \partial_{\xi} \mathcal{A}_{12}(\xi) \underset{\tau
\to +\infty}{=} \dfrac{2 \mi \alpha^{4} \hat{r}_{0}(\tau) \mathfrak{p}_{
\tau}}{\xi^{3}} \! + \! \mathcal{O} \! \left(\dfrac{\tau^{-2/3}}{\xi^{3}
(1 \! + \! \hat{u}_{0}(\tau))} \right), \label{eqnumTminus1T}
\end{equation}
where $\mathfrak{p}_{\tau}$ is defined by Equation~\eqref{eqpeetee}.
Substituting Equations~\eqref{eqdenomTminus1T} and~\eqref{eqnumTminus1T}
into Equation~\eqref{eq3.23} and expanding $(2 \mi l(\xi) \mathcal{A}_{11}
(\xi) \! + \! 2l^{2}(\xi))^{-1}$ into a series of powers of $\Delta_{\tau}
(\xi)$, one arrives at
\begin{equation*}
\int_{\widetilde{\mu}_{0}}^{\widetilde{\mu}} \diag ((T(\xi))^{-1}
\partial_{\xi}T(\xi)) \, \md \xi \underset{\tau \to +\infty}{=} \left(
\mathcal{I}_{\tau}(\widetilde{\mu}) \! + \! \mathcal{O}(\mathcal{E}_{T}
(\widetilde{\mu})) \right) \sigma_{3},
\end{equation*}
where
\begin{equation}
\mathcal{I}_{\tau}(\widetilde{\mu}) \! := \! -\mi \alpha^{4} \hat{r}_{0}(\tau)
\mathfrak{p}_{\tau} \int_{\widetilde{\mu}_{0}}^{\widetilde{\mu}} \dfrac{1}{
\xi^{3} \mathcal{P}_{\infty}(\xi)} \, \md \xi, \label{eq3.26}
\end{equation}
and
\begin{equation}
\mathcal{E}_{T}(\widetilde{\mu}) \! := \! \mi \alpha^{4} \hat{r}_{0}(\tau)
\mathfrak{p}_{\tau} \int_{\widetilde{\mu}_{0}}^{\widetilde{\mu}} \dfrac{
\mathcal{P}_{1}(\xi) \Delta_{\tau}(\xi)}{\xi^{3}(\mathcal{P}_{\infty}
(\xi))^{2}} \, \md \xi \! + \! \dfrac{\tau^{-2/3}}{(1 \! + \! \hat{u}_{0}
(\tau))} \int_{\widetilde{\mu}_{0}}^{\widetilde{\mu}} \dfrac{1}{\xi^{3}
\mathcal{P}_{\infty}(\xi)} \, \md \xi. \label{eq3.27}
\end{equation}
Via Equation~\eqref{eqlsquared}, it follows {}from Equation~\eqref{eq3.24}
that
\begin{equation}
\dfrac{1}{\xi^{3} \mathcal{P}_{\infty}(\xi)} \! = \! \dfrac{\xi (\xi (4
\xi^{2} \! + \! 2 \alpha^{2}(-2 \! + \! \hat{r}_{0}(\tau))) \! - \! 4(\xi^{2}
\! - \! \alpha^{2})(\xi^{2} \! + \! 2 \alpha^{2})^{1/2})}{2(\xi^{2} \! - \!
\alpha^{2})(\xi^{2} \! + \! 2 \alpha^{2})^{1/2}((\xi (4 \xi^{2} \! + \! 2
\alpha^{2}(-2 \! + \! \hat{r}_{0}(\tau))))^{2} \! - \! 16(\xi^{2} \! - \!
\alpha^{2})^{2}(\xi^{2} \! + \! 2 \alpha^{2}))}:
\label{eq3.28}
\end{equation}
writing the factorization
\begin{align*}
&(\xi (4 \xi^{2} \! + \! 2 \alpha^{2}(-2 \! + \! \hat{r}_{0}(\tau))))^{2} \!
- \! 16(\xi^{2} \! - \! \alpha^{2})^{2}(\xi^{2} \! + \! 2 \alpha^{2}) \\
&= 16 \alpha^{2}(-2 \! + \! \hat{r}_{0}(\tau)) \left(\xi^{4} \! + \!
\dfrac{\xi^{2} \alpha^{2}}{(-2 \! + \! \hat{r}_{0}(\tau))} \left(3 \! +
\! \left(\dfrac{-2 \! + \! \hat{r}_{0}(\tau)}{2} \right)^{2} \right) \! -
\! \dfrac{2 \alpha^{4}}{(-2 \! + \! \hat{r}_{0}(\tau))} \right) \\
&= 16 \alpha^{2}(-2 \! + \! \hat{r}_{0}(\tau))(\xi^{2} \! + \!
\hat{z}_{+})(\xi^{2} \! + \! \hat{z}_{-}),
\end{align*}
where
\begin{equation}
\hat{z}_{\pm} \! := \! \dfrac{\alpha^{2}}{2(-2 \! + \! \hat{r}_{0}(\tau))} \!
\left(3 \! + \! \left(\dfrac{-2 \! + \! \hat{r}_{0}(\tau)}{2} \right)^{2} \!
\mp \! \dfrac{\hat{r}_{0}(\tau)}{4} \sqrt{32 \! + \! (\hat{r}_{0}(\tau) \! -
\! 4)^{2}} \right), \label{eqzplusminus}
\end{equation}
one shows that, for $\hat{z}_{+} \! \not\equiv \! \hat{z}_{-}$,
\begin{align}
\mathcal{I}_{\tau}(\widetilde{\mu}) =& \, -\dfrac{\mi \alpha^{2} \hat{r}_{0}
(\tau) \mathfrak{p}_{\tau}}{32(-2 \! + \! \hat{r}_{0}(\tau))} \left(
\int_{\widetilde{\mu}_{0}}^{\widetilde{\mu}} \dfrac{\xi^{2}(4 \xi^{2} \!
+ \! 2 \alpha^{2}(-2 \! + \! \hat{r}_{0}(\tau))) \sqrt{\xi^{2} \! + \! 2
\alpha^{2}}}{(\xi^{2} \! - \! \alpha^{2})(\xi^{2} \! + \! 2 \alpha^{2})
(\xi^{2} \! + \! \hat{z}_{+})(\xi^{2} \! + \! \hat{z}_{-})} \, \md \xi
\right. \nonumber \\
-&\left. \, 4 \int_{\widetilde{\mu}_{0}}^{\widetilde{\mu}}
\dfrac{\xi}{(\xi^{2} \! + \! \hat{z}_{+})(\xi^{2} \! + \hat{z}_{-})}
\, \md \xi \right). \label{eq3.29}
\end{align}
Write the partial fraction decomposition
\begin{align}
\dfrac{\xi^{2}(4 \xi^{2} \! + \! 2 \alpha^{2}(-2 \! + \! \hat{r}_{0}
(\tau)))}{(\xi^{2} \! - \! \alpha^{2})(\xi^{2} \! + \! 2 \alpha^{2})(\xi^{2}
\! + \! \hat{z}_{+})(\xi^{2} \! + \! \hat{z}_{-})} =& \, \dfrac{A_{0}
(\tau)}{\xi \! - \! \alpha} \! + \! \dfrac{B_{0}(\tau)}{\xi \! + \! \alpha}
\! + \! \dfrac{C_{0}(\tau) \xi \! + \! D_{0}(\tau)}{\xi^{2} \! + \! 2
\alpha^{2}} \! + \! \dfrac{E_{0}(\tau) \xi \! + \! F_{0}(\tau)}{\xi^{2} \!
+ \! \hat{z}_{+}} \nonumber \\
+& \, \dfrac{G_{0}(\tau) \xi \! + \! H_{0}(\tau)}{\xi^{2} \! + \!
\hat{z}_{-}}: \label{eq3.30}
\end{align}
to determine $A_{0}(\tau),\dotsc,H_{0}(\tau)$, one notes that the left-hand
side of Equation~\eqref{eq3.30} is symmetric $(\xi \! \to \! -\xi)$, whence
it follows that
\begin{equation}
A_{0}(\tau) \! = \! -B_{0}(\tau), \qquad C_{0}(\tau) \! = \! E_{0}(\tau) \!
= \! G_{0}(\tau) \! = \! 0, \label{eqcoeff1}
\end{equation}
and to determine the remaining coefficients, that is, $A_{0}(\tau)$, $D_{0}
(\tau)$, $F_{0}(\tau)$, and $H_{0}(\tau)$, one compares residues (at $\xi
\! = \! \alpha$, $\xi^{2} \! = \! -2 \alpha^{2}$, and $\xi^{2} \! = \!
-\hat{z}_{\pm})$ on the left- and right-hand sides of Equation~\eqref{eq3.30}
and uses the relations (cf. Equation~\eqref{eqzplusminus})
\begin{equation*}
\hat{z}_{+} \! + \! \hat{z}_{-} \! = \! \dfrac{\alpha^{2}}{(-2 \! + \!
\hat{r}_{0}(\tau))} \! \left(3 \! + \! \left(\dfrac{-2 \! + \! \hat{r}_{0}
(\tau)}{2} \right)^{2} \right) \qquad \text{and} \qquad \hat{z}_{+}
\hat{z}_{-} \! = \! -\dfrac{2 \alpha^{4}}{(-2 \! + \! \hat{r}_{0}(\tau))}
\end{equation*}
to show that
\begin{equation}
\begin{gathered}
A_{0}(\tau) \! = \! \dfrac{\alpha \hat{r}_{0}(\tau)}{3(\alpha^{2} \! + \!
\hat{z}_{+})(\alpha^{2} \! + \! \hat{z}_{-})} \! = \! \dfrac{4(-2 \! + \!
\hat{r}_{0}(\tau))}{3 \alpha^{3} \hat{r}_{0}(\tau)}, \quad D_{0}(\tau)
\! = \! \dfrac{4 \alpha^{2} \hat{r}_{0}(\tau) \! - \! 24 \alpha^{2}}{3(-2
\alpha^{2} \! + \! \hat{z}_{-})(-2 \alpha^{2} \! + \! \hat{z}_{+})} \! = \!
\dfrac{8(-2 \! + \! \hat{r}_{0}(\tau))}{3 \alpha^{2}(6 \! - \! \hat{r}_{0}
(\tau))}, \\
F_{0}(\tau) \! = \! \dfrac{2 \alpha^{2} \hat{z}_{+}(-2 \! + \! \hat{r}_{0}
(\tau)) \! - \! 4(\hat{z}_{+})^{2}}{(\alpha^{2} \! + \! \hat{z}_{+})(-2
\alpha^{2} \! + \! \hat{z}_{+})(\hat{z}_{+} \! - \! \hat{z}_{-})}, \quad
H_{0}(\tau) \! = \! \dfrac{4(\hat{z}_{-})^{2} \! - \! 2 \alpha^{2} \hat{z}_{-}
(-2 \! + \! \hat{r}_{0}(\tau))}{(\alpha^{2} \! + \! \hat{z}_{-})(-2 \alpha^{2}
\! + \! \hat{z}_{-})(\hat{z}_{+} \! - \! \hat{z}_{-})} \quad \Rightarrow \\
F_{0}(\tau) \! + \! H_{0}(\tau) \! = \! -\dfrac{16(-2 \! + \! \hat{r}_{0}
(\tau))}{\alpha^{2} \hat{r}_{0}(\tau)(6 \! - \! \hat{r}_{0}(\tau))}.
\label{eqcoeff2}
\end{gathered}
\end{equation}
(Note: it is the quantity $F_{0}(\tau) \! + \! H_{0}(\tau)$, and not $F_{0}
(\tau)$ and $H_{0}(\tau)$ individually, that is requisite for the ensuing
calculation.) Substituting Equations~\eqref{eqcoeff2}, \eqref{eqcoeff1},
and~\eqref{eq3.30} into Equation~\eqref{eq3.29}, one arrives at, after an
integration argument and neglecting $\xi$-independent terms (that is, with
abuse of nomenclature, ``constants of integration''),
\begin{align}
\mathcal{I}_{\tau}(\widetilde{\mu}) =& \, -\dfrac{\mi \alpha^{2} \hat{r}_{0}
(\tau) \mathfrak{p}_{\tau}}{32(-2 \! + \! \hat{r}_{0}(\tau))} \left. \! \left(
\vphantom{M^{M^{M^{M^{M^{M}}}}}} \! \left(2 \alpha A_{0}(\tau) \! + \! D_{0}
(\tau) \! + \! \dfrac{8(-2 \! + \! \hat{r}_{0}(\tau))}{\alpha^{2} \hat{r}_{0}
(\tau)(\hat{r}_{0}(\tau) \! - \! 6)} \right) \ln (\sqrt{\xi^{2} \! + \! 2
\alpha^{2}}+ \! \xi) \right. \right. \nonumber \\
-&\left. \left. \dfrac{8(-2 \! + \! \hat{r}_{0}(\tau))}{\alpha^{2}
\hat{r}_{0}(\tau)(\hat{r}_{0}(\tau) \! - \! 6)} \ln (\sqrt{\xi^{2} \! + \!
2 \alpha^{2}}- \! \xi) \! + \! \sqrt{3} \, \alpha A_{0}(\tau) \ln \! \left(
\left(\dfrac{\sqrt{3} \, \sqrt{\xi^{2} \! + \! 2 \alpha^{2}}- \! \xi \! + \!
2 \alpha}{\sqrt{3} \, \sqrt{\xi^{2} \! + \! 2 \alpha^{2}}+ \! \xi \! + \! 2
\alpha} \right) \! \left(\dfrac{\xi \! - \! \alpha}{\xi \! + \! \alpha}
\right) \right) \right. \right. \nonumber \\
+&\left. \left. \dfrac{16(-2 \! + \! \hat{r}_{0}(\tau))}{\alpha^{2}
\hat{r}_{0}(\tau) \sqrt{32 \! + \! (\hat{r}_{0}(\tau) \! - \! 4)^{2}}} \!
\left(\operatorname{arctanh} \! \left(\dfrac{8(-2 \! + \! \hat{r}_{0}(\tau))
\xi^{2} \! + \! \alpha^{2}(-2 \! + \! \hat{r}_{0}(\tau))^{2} \! + \! 12
\alpha^{2}}{\alpha^{2} \hat{r}_{0}(\tau) \sqrt{32 \! + \! (\hat{r}_{0}(\tau)
\! - \! 4)^{2}}} \right) \right. \right. \right. \nonumber \\
+&\left. \left. \left. \operatorname{arctanh} \! \left(\dfrac{4 \sqrt{\xi^{2}
\! + \! 2 \alpha^{2}}}{\xi (\hat{r}_{0}(\tau) \! + \! \sqrt{32 \! + \!
(\hat{r}_{0}(\tau) \! - \! 4)^{2}})} \right) \! - \! \operatorname{arctanh}
\! \left(\dfrac{4 \sqrt{\xi^{2} \! + \! 2 \alpha^{2}}}{\xi (\hat{r}_{0}(\tau)
\! - \! \sqrt{32 \! + \! (\hat{r}_{0}(\tau) \! - \! 4)^{2}})} \! \right) \!
\! \right) \! \! \right) \! \right\vert_{\widetilde{\mu}_{0}}^{
\widetilde{\mu}}. \label{eq3.34}
\end{align}
{}From Equations~\eqref{eqcoeff2}, one notes that
\begin{equation*}
2 \alpha A_{0}(\tau) \! + \! D_{0}(\tau) \! + \! \dfrac{16(-2 \! + \!
\hat{r}_{0}(\tau))}{\alpha^{2} \hat{r}_{0}(\tau)(\hat{r}_{0}(\tau) \! - \!
6)} \! = \! 0;
\end{equation*}
whence (cf. Equation~\eqref{eq3.34})
\begin{align*}
\left(2 \alpha A_{0}(\tau) \! + \! D_{0}(\tau) \! + \! \dfrac{8(-2 \! + \!
\hat{r}_{0}(\tau))}{\alpha^{2} \hat{r}_{0}(\tau)(\hat{r}_{0}(\tau) \! - \!
6)} \right)& \ln (\sqrt{\xi^{2} \! + \! 2 \alpha^{2}}+ \! \xi) \! - \!
\dfrac{8(-2 \! + \! \hat{r}_{0}(\tau))}{\alpha^{2} \hat{r}_{0}(\tau)
(\hat{r}_{0}(\tau) \! - \! 6)} \ln (\sqrt{\xi^{2} \! + \! 2 \alpha^{2}}- \!
\xi) \\
=& -\dfrac{8(-2 \! + \! \hat{r}_{0}(\tau))}{\alpha^{2} \hat{r}_{0}(\tau)
(\hat{r}_{0}(\tau) \! - \! 6)} \ln (2 \alpha^{2}),
\end{align*}
which is a $\xi$-independent term (a ``constant of integration''):
neglecting this term, one arrives at (after simplification)
Equations~\eqref{eqItee}--\eqref{eqFtee}.

One now studies the error term $\mathcal{O}(\mathcal{E}_{T}(\widetilde{\mu}))$
defined by Equation~\eqref{eq3.27}. Recall the expression for $(\xi^{3}
\mathcal{P}_{\infty}(\xi))^{-1}$ given in Equation~\eqref{eq3.28}:
{}from this latter expression, and those for $\Delta_{\tau}(\xi)$
(cf. Equation~\eqref{eqdellt}) and $\mathcal{P}_{1}(\xi)$ (cf.
Equation~\eqref{eq3.25}), one shows that
\begin{equation}
\dfrac{\mathcal{P}_{1}(\xi) \Delta_{\tau}(\xi)}{\xi^{3}(\mathcal{P}_{\infty}
(\xi))^{2}} \! = \! \dfrac{\xi^{2}(h_{0}(\tau) \! + \! \xi^{2}(a \! - \!
\frac{\mi}{2}))(2 \xi l_{\infty}(\xi) \! + \! 2 \xi^{2} \! + \! \alpha^{2}
(-2 \! + \! \hat{r}_{0}(\tau))) \tau^{-2/3}}{8(\xi^{2} \! - \! \alpha^{2})^{3}
(\xi^{2} \! + \! 2 \alpha^{2})^{3/2}(\xi l_{\infty}(\xi) \! + \! 2 \xi^{2}
\! + \! \alpha^{2}(-2 \! + \! \hat{r}_{0}(\tau)))^{2}}. \label{eq3.35}
\end{equation}
Using Equation~\eqref{eq3.35}, one evaluates the first integral in
Equation~\eqref{eq3.27} explicitly (as done above, for the second integral
in Equation~\eqref{eq3.27}): the resulting expression for $\mathcal{E}_{T}
(\widetilde{\mu})$ is quite cumbersome; therefore, only its asymptotics at
the singular and turning points are defined by Equation~\eqref{eqee2}.
\hfill $\qed$
\begin{ffff} \label{cor3.1.3}
Let $\widetilde{\mu}_{0}$ be defined as in Corollary~{\rm \ref{cor3.1.2}},
that is, $\widetilde{\mu}_{0} \! = \! \alpha \! + \! \tau^{-1/3}
\widetilde{\Lambda}$, where $\widetilde{\Lambda} \! =_{\tau \to +\infty} \!
\mathcal{O}(\tau^{\varepsilon})$, and $\delta,\delta_{1}$ be defined in
conditions~\eqref{eq3.11}. Then, for $0 \! < \! \delta \! < \! \varepsilon
\! < \! 1/9$ and $\varepsilon \! < \! \delta_{1} \! < \! 1/3$,
\begin{equation}
\int_{\widetilde{\mu}_{0}}^{\widetilde{\mu}} \diag ((T(\xi))^{-1}
\partial_{\xi}T(\xi)) \, \md \xi \underset{\tau \to +\infty}{=} \left(
\dfrac{\mathfrak{p}_{\tau}}{8 \mi}(\digamma_{\tau}(\widetilde{\mu}) \! - \!
\digamma_{\tau}^{\sharp}) \! + \! \mathcal{O}(\mathcal{E}_{T}^{\natural}
(\widetilde{\mu})) \! + \! \mathcal{O}(\tau^{-2(\varepsilon - \delta)}) \! +
\! \mathcal{O}(\tau^{-\frac{1}{3}+ \varepsilon + \delta}) \right) \sigma_{3},
\label{eq3.36}
\end{equation}
where $\mathfrak{p}_{\tau}$, $\digamma_{\tau}(\widetilde{\mu})$,
and $\mathcal{E}_{T}^{\natural}(\widetilde{\mu})$ are defined by
Equations~\eqref{eqpeetee}, \eqref{eqFtee}, and~\eqref{eqee2},
respectively, and
\begin{align}
\digamma_{\tau}^{\sharp} :=& \, \mp \dfrac{1}{\sqrt{3}} \ln (3 \alpha) \! +
\! \dfrac{2}{\sqrt{32 \! + \! (\hat{r}_{0}(\tau) \! - \! 4)^{2}}} \ln \!
\left(\dfrac{\hat{r}_{0}(\tau) \! + \! 4 \! + \! \sqrt{32 \! + \! (\hat{r}_{0}
(\tau) \! - \! 4)^{2}}}{-\hat{r}_{0}(\tau) \! - \! 4 \! + \! \sqrt{32 \! +
\! (\hat{r}_{0}(\tau) \! - \! 4)^{2}}} \right) \nonumber \\
\mp& \, \dfrac{2}{\sqrt{32 \! + \! (\hat{r}_{0}(\tau) \! - \! 4)^{2}}} \ln \!
\left(\! \left(\dfrac{\hat{r}_{0}(\tau) \! + \! 4 \sqrt{3}- \! \sqrt{32 \! +
\! (\hat{r}_{0}(\tau) \! - \! 4)^{2}}}{\hat{r}_{0}(\tau) \! + \! 4 \sqrt{3}
+ \! \sqrt{32 \! + \! (\hat{r}_{0}(\tau) \! - \! 4)^{2}}} \right) \right.
\nonumber \\
\times&\left. \, \left(\dfrac{\hat{r}_{0}(\tau) \! - \! 4 \sqrt{3}+ \!
\sqrt{32 \! + \! (\hat{r}_{0}(\tau) \! - \! 4)^{2}}}{\hat{r}_{0}(\tau) \!
- \! 4 \sqrt{3} -\! \sqrt{32 \! + \! (\hat{r}_{0}(\tau) \! - \! 4)^{2}}}
\right) \! \right) \! \pm \! \dfrac{1}{\sqrt{3}} \! \left(-\dfrac{1}{3}
\ln \tau \! + \! \ln \widetilde{\Lambda} \right), \label{eqfteesharp}
\end{align}
with the upper (resp., lower) signs taken if the positive $(+)$ (resp.,
negative $(-))$ branch of the square-root function $\sqrt{\xi^{2} \! + \!
2 \alpha^{2}}$ is chosen.
\end{ffff}

\emph{Proof}. Substituting $\widetilde{\mu}_{0}$ as the argument of the
functions $\digamma_{\tau}(\xi)$ and $\mathcal{E}_{T}^{\natural}(\xi)$
(cf. Equation~\eqref{eqFtee} and the first line of Equation~\eqref{eqee2},
respectively) and expanding with respect to (the small parameter) $\tau^{-1/3}
\widetilde{\Lambda}$, one arrives at the following estimates:
\begin{equation*}
\digamma_{\tau}(\widetilde{\mu}_{0}) \underset{\tau \to +\infty}{=}
\digamma_{\tau}^{\sharp} \! + \! \mathcal{O}(\tau^{-1/3} \widetilde{\Lambda}),
\end{equation*}
and
\begin{equation*}
\mathcal{E}_{T}^{\natural}(\widetilde{\mu}_{0}) \underset{\tau \to +\infty}{=}
\mathcal{O} \! \left(\dfrac{\tau^{-2/3} \ln \tau}{(1 \! + \! \hat{u}_{0}(\tau))
(C_{1}(\tau) \! + \! \hat{r}_{0}(\tau))} \right) \! + \! \mathcal{O} \! \left(
\dfrac{(2 \hat{u}_{0}(\tau) \! + \! \frac{1}{2} \hat{r}_{0}(\tau))(C_{2}(\tau)
\! + \! h_{0}(\tau))}{\hat{r}_{0}(\tau) \widetilde{\Lambda}^{2}} \right),
\end{equation*}
where $\digamma_{\tau}^{\sharp}$ is defined by Equation~\eqref{eqfteesharp},
and $C_{j}(\tau) \! =_{\tau \to +\infty} \! \mathcal{O}(1)$, $j \! = \!
1,2$. Demanding that the above error estimates are decaying, and using
the definition of the parameters $\delta,\delta_{1}$ given in
conditions~\eqref{eq3.11} for the functions $\hat{r}_{0}(\tau)$, $\hat{u}_{0}
(\tau)$, and $h_{0}(\tau)$, one arrives at the inequalities stated in the
Corollary. \hfill $\qed$
\begin{ffff} \label{cor3.1.4}
Under the conditions of Corollary~{\rm \ref{cor3.1.2}}, for the branch of
$l(\xi)$ that is positive for large and small positive $\xi$,
\begin{align}
-\mi \tau^{2/3} \int_{\widetilde{\mu}_{0}}^{\widetilde{\mu}}l(\xi) \, \md \xi
\underset{\underset{\scriptstyle \Re (\widetilde{\mu}) \to +\infty}{\tau \to
+\infty}}{=}& \, -\mi \left(\tau^{2/3} \widetilde{\mu}^{2} \! + \! \left(a \!
- \! \dfrac{\mi}{2} \right) \ln \widetilde{\mu} \right) \! + \! \mi \tau^{2/3}
3(\sqrt{3}- \! 1) \alpha^{2} \! + \! \mi 2 \sqrt{3} \, \widetilde{\Lambda}^{2}
\nonumber \\
-& \, \dfrac{\mi}{2 \sqrt{3}} \left(\left(a \! - \! \dfrac{\mi}{2} \right)
\! + \! \dfrac{h_{0}(\tau)}{\alpha^{2}} \right) \! \left(\dfrac{1}{3} \ln
\tau \! - \! \ln \widetilde{\Lambda} \right) \! + \! C^{\mathrm{WKB}}_{\infty}
\! + \! \mathcal{O}(\tau^{-2(\varepsilon - \delta)}) \nonumber \\
+& \, \mathcal{O}(\tau^{-\frac{1}{3}+3 \varepsilon}), \label{eq3.37}
\end{align}
where
\begin{equation}
C^{\mathrm{WKB}}_{\infty} \! := \! \mi \left(a \! - \! \dfrac{\mi}{2} \right)
\ln \left(\dfrac{(\sqrt{3}+ \! 1) \alpha}{2} \right) \! - \! \dfrac{\mi}{2
\sqrt{3}} \left(\left(a \! - \! \dfrac{\mi}{2} \right) \! + \!
\dfrac{h_{0}(\tau)}{\alpha^{2}} \right) \ln \left(\dfrac{6 \alpha}{(\sqrt{3}
+ \! 1)^{2}} \right), \label{eq3.38}
\end{equation}
and
\begin{align}
-\mi \tau^{2/3} \int_{\widetilde{\mu}_{0}}^{\widetilde{\mu}}l(\xi) \, \md \xi
\underset{\underset{\scriptstyle \Re (\widetilde{\mu}) \to +0}{\tau \to
+\infty}}{=}& \, \dfrac{\mi \tau^{2/3}2 \sqrt{2} \, \alpha^{3}}{\widetilde{
\mu}} \! - \! \mi \tau^{2/3} 3 \sqrt{3} \, \alpha^{2} \! - \! \mi 2 \sqrt{3}
\, \widetilde{\Lambda}^{2} \! + \! \dfrac{\mi}{2 \sqrt{3}} \left(\left(a \!
- \! \dfrac{\mi}{2} \right) \! + \! \dfrac{h_{0}(\tau)}{\alpha^{2}} \right)
\nonumber \\
\times& \, \left(\dfrac{1}{3} \ln \tau \! - \! \ln \widetilde{\Lambda} \right)
\! + \! C^{\mathrm{WKB}}_{0} \! + \! \mathcal{O}(\tau^{-2(\varepsilon -
\delta)}) \! + \! \mathcal{O}(\tau^{-\frac{1}{3}+3 \varepsilon}),
\label{eq3.39}
\end{align}
where
\begin{equation}
C^{\mathrm{WKB}}_{0} \! := \! -\mi \left(a \! - \! \dfrac{\mi}{2} \right) \ln
\left(\dfrac{\sqrt{3}+ \! 1}{\sqrt{2}} \right) \! + \! \dfrac{\mi}{2 \sqrt{3}}
\left(\left(a \! - \! \dfrac{\mi}{2} \right) \! + \! \dfrac{h_{0}(\tau)}{
\alpha^{2}} \right) \ln (3 \alpha \me^{-\mi \pi}). \label{eq3.40}
\end{equation}
\end{ffff}

\emph{Proof}. Follows {}from Corollary~\ref{cor3.1.2},
Equation~\eqref{eq3.22}, by choosing the corresponding branches in
Equations~\eqref{equpsi} and~\eqref{equpsisharp} and taking the limits $\Re
(\widetilde{\mu}) \! \to \! +\infty$ and $\Re (\widetilde{\mu}) \! \to \!
+0$: the error estimate $\mathcal{O}(\mathscr{E}_{l}^{\natural}(\xi))$ in
Equation~\eqref{eq3.22} is defined by Equation~\eqref{eqeel}. \hfill $\qed$
\begin{ffff} \label{cor3.1.5}
Under the conditions of Corollary~{\rm \ref{cor3.1.3}}, for the branch of
$l(\xi)$ that is positive for large and small positive $\xi$,
\begin{align}
\int_{\widetilde{\mu}_{0}}^{\widetilde{\mu}} \diag ((T(\xi))^{-1} \partial_{
\xi}T(\xi)) \, \md \xi \underset{\underset{\scriptstyle \Re (\widetilde{\mu})
\to +\infty}{\tau \to +\infty}}{=}& \, \left(-\dfrac{\sqrt{32 \! + \!
(\hat{r}_{0}(\tau) \! - \! 4)^{2}}}{8 \sqrt{3}} \! \left(-\dfrac{1}{3}
\ln \tau \! + \! \ln \widetilde{\Lambda} \right) \! + \!
\mathcal{I}_{\infty}^{\sharp}(\tau) \right. \nonumber \\
+&\left. \, \mathcal{O} \! \left(\dfrac{\hat{r}_{0}(\tau)}{\widetilde{
\mu}^{2}} \right) \! + \! \mathcal{O}(\tau^{\delta -2 \delta_{1}} \ln \tau)
\! + \! \mathcal{O}(\tau^{-2(\varepsilon - \delta)}) \right. \nonumber \\
+&\left. \, \mathcal{O}(\tau^{-\frac{2}{3}+2 \delta}) \! + \! \mathcal{O}
(\tau^{-\frac{1}{3}+ \varepsilon + \delta}) \right) \sigma_{3}, \label{eq3.41}
\end{align}
where
\begin{align}
\mathcal{I}_{\infty}^{\sharp}(\tau) :=& \, \dfrac{1}{4} \ln \! \left(\! \left(
\dfrac{\hat{r}_{0}(\tau) \! + \! 4 \sqrt{3}- \! \sqrt{32 \! + \! (\hat{r}_{0}
(\tau) \! - \! 4)^{2}}}{\hat{r}_{0}(\tau) \! + \! 4 \sqrt{3}+ \! \sqrt{32 \!
+ \! (\hat{r}_{0}(\tau) \! - \! 4)^{2}}} \right) \! \left(\dfrac{\hat{r}_{0}
(\tau) \! - \! 4 \sqrt{3}+ \! \sqrt{32 \! + \! (\hat{r}_{0}(\tau) \! - \!
4)^{2}}}{\hat{r}_{0}(\tau) \! - \! 4 \sqrt{3}- \! \sqrt{32 \! + \!
(\hat{r}_{0}(\tau) \! - \! 4)^{2}}} \right) \! \right) \nonumber \\
+& \, \dfrac{\sqrt{32 \! + \! (\hat{r}_{0}(\tau) \! - \! 4)^{2}}}{8 \sqrt{3}}
\ln \! \left(\dfrac{6 \alpha}{(\sqrt{3}+ \! 1)^{2}} \right) \! - \!
\dfrac{1}{2} \ln \! \left(\dfrac{\hat{u}_{0}(\tau)}{\hat{r}_{0}(\tau)}
\right) \! - \! \dfrac{3}{4} \ln 2 \! - \! \dfrac{\mi \pi}{4}, \label{eq3.87}
\end{align}
and
\begin{align}
\int_{\widetilde{\mu}_{0}}^{\widetilde{\mu}} \diag ((T(\xi))^{-1} \partial_{
\xi}T(\xi)) \, \md \xi \underset{\underset{\scriptstyle \Re (\widetilde{\mu})
\to +0}{\tau \to +\infty}}{=}& \, \left(\dfrac{\sqrt{32 \! + \! (\hat{r}_{0}
(\tau) \! - \! 4)^{2}}}{8 \sqrt{3}} \! \left(-\dfrac{1}{3} \ln \tau \! + \!
\ln \widetilde{\Lambda} \right) \! + \! \mathcal{I}_{0}^{\sharp}(\tau)
\right. \nonumber \\
+&\left. \, \mathcal{O}(\hat{r}_{0}(\tau) \widetilde{\mu}) \! + \!
\mathcal{O}(\tau^{\delta -2 \delta_{1}} \ln \tau) \! + \! \mathcal{O}
(\tau^{-2(\varepsilon - \delta)}) \right. \nonumber \\
+&\left. \, \mathcal{O}(\tau^{-\frac{2}{3}+2 \delta}) \! + \! \mathcal{O}
(\tau^{-\frac{1}{3} + \varepsilon + \delta}) \! + \! \mathcal{O}
(\tau^{-\frac{1}{3}+ \delta - \delta_{1}}) \right) \sigma_{3}, \label{eq3.43}
\end{align}
where
\begin{align}
\mathcal{I}_{0}^{\sharp}(\tau) :=& \, -\dfrac{1}{4} \ln \! \left(\! \left(
\dfrac{\hat{r}_{0}(\tau) \! + \! 4 \sqrt{3}- \! \sqrt{32 \! + \! (\hat{r}_{0}
(\tau) \! - \! 4)^{2}}}{\hat{r}_{0}(\tau) \! + \! 4 \sqrt{3}+ \! \sqrt{32 \!
+ \! (\hat{r}_{0}(\tau) \! - \! 4)^{2}}} \right) \! \left(\dfrac{\hat{r}_{0}
(\tau) \! - \! 4 \sqrt{3}+ \! \sqrt{32 \! + \! (\hat{r}_{0}(\tau) \! - \!
4)^{2}}}{\hat{r}_{0}(\tau) \! - \! 4 \sqrt{3}- \! \sqrt{32 \! + \!
(\hat{r}_{0}(\tau) \! - \! 4)^{2}}} \right) \! \right) \nonumber \\
-& \, \dfrac{\sqrt{32 \! + \! (\hat{r}_{0}(\tau) \! - \! 4)^{2}}}{8 \sqrt{3}}
\ln (3 \alpha \me^{-\mi \pi}) \! + \! \dfrac{1}{2} \ln \! \left(\dfrac{
\hat{u}_{0}(\tau)}{\hat{r}_{0}(\tau)} \right) \! + \! \dfrac{3}{4} \ln 2
\! + \! \dfrac{\mi \pi}{4}. \label{eq3.89}
\end{align}
\end{ffff}

\emph{Proof}. Taking the limit $\Re (\widetilde{\mu}) \! \to \! +\infty$ in
Equation~\eqref{eq3.36}, with $\digamma_{\tau}(\xi)$ and $\mathcal{E}_{T}^{
\natural}(\xi)$ defined by Equations~\eqref{eqFtee} and~\eqref{eqee2},
respectively, and using conditions~\eqref{eq3.11}, one arrives at
\begin{align}
\digamma_{\tau}(\widetilde{\mu}) \underset{\underset{\scriptstyle \Re
(\widetilde{\mu}) \to +\infty}{\tau \to +\infty}}{=}& \, \dfrac{1}{\sqrt{3}}
\ln \left(\dfrac{\sqrt{3}- \! 1}{\sqrt{3}+ \! 1} \right) \! - \! \dfrac{2}{
\sqrt{32 \! + \! (\hat{r}_{0}(\tau) \! - \! 4)^{2}}} \ln \! \left(\! \left(
\dfrac{\hat{r}_{0}(\tau) \! - \! 4 \! + \! \sqrt{32 \! + \! (\hat{r}_{0}(\tau)
\! - \! 4)^{2}}}{\hat{r}_{0}(\tau) \! - \! 4 \! - \! \sqrt{32 \! + \!
(\hat{r}_{0}(\tau) \! - \! 4)^{2}}} \right) \right. \nonumber \\
\times&\left. \, \left(\dfrac{-\hat{r}_{0}(\tau) \! - \! 4 \! + \! \sqrt{32
\! + \! (\hat{r}_{0}(\tau) \! - \! 4)^{2}}}{\hat{r}_{0}(\tau) \! + \! 4 \! +
\! \sqrt{32 \! + \! (\hat{r}_{0}(\tau) \! - \! 4)^{2}}} \right) \! \right)
\! + \! \mathcal{O}(\widetilde{\mu}^{-2}), \label{eqF1}
\end{align}
and
\begin{equation}
\mathcal{E}_{T}^{\natural}(\widetilde{\mu}) \underset{\underset{\scriptstyle
\Re (\widetilde{\mu}) \to +\infty}{\tau \to +\infty}}{=} \mathcal{O}(\tau^{-2
\delta_{1}}) \! + \! \mathcal{O}(\tau^{-\frac{2}{3}+ 2 \delta}). \label{eqE1}
\end{equation}
One now computes asymptotics of $\mathfrak{p}_{\tau}$ (cf.
Equation~\eqref{eqpeetee}). Using conditions~\eqref{eq3.11} and Equation
\eqref{eq3.14}, one shows that
\begin{equation}
\dfrac{-2 \! + \! \hat{r}_{0}(\tau)}{1 \! + \! \hat{u}_{0}(\tau)} \underset{
\tau \to +\infty}{=} \dfrac{\hat{r}_{0}(\tau)(4 \! - \! \hat{r}_{0}(\tau))}{4}
\! + \! 2(1 \! + \! \hat{u}_{0}(\tau)) \! - \! 4 \! + \! \mathcal{O}
(\tau^{-2/3} \varkappa_{0}^{2}(\tau)), \label{eq3.46}
\end{equation}
where
\begin{equation}
\varkappa_{0}^{2}(\tau) \! := \! \dfrac{4}{\alpha^{2}} \left(h_{0}(\tau) \! +
\! \dfrac{\alpha^{2}(a \! - \! \frac{\mi}{2})}{1 \! + \! \hat{u}_{0}(\tau)}
\right), \label{eq3.58}
\end{equation}
whence, substituting Equation~\eqref{eq3.46} into Equation~\eqref{eqpeetee},
one obtains
\begin{equation}
\mathfrak{p}_{\tau} \underset{\tau \to +\infty}{=} \dfrac{4 \mi}{\hat{r}_{0}
(\tau)} \left(4 \hat{u}_{0}(\tau) \! + \! \dfrac{\hat{r}_{0}(\tau)(4 \! - \!
\hat{r}_{0}(\tau))}{4} \! + \! \mathcal{O}(\tau^{-2/3} \varkappa_{0}^{2}
(\tau)) \right). \label{eq3.47}
\end{equation}
Solving Equation~\eqref{eq3.14} for $\hat{u}_{0}(\tau)$ and taking into
account conditions~\eqref{eq3.11}, one deduces that
\begin{equation}
16 \hat{u}_{0}(\tau) \underset{\tau \to +\infty}{=} -\hat{r}_{0}(\tau)(4 \! -
\! \hat{r}_{0}(\tau)) \! + \! \hat{r}_{0}(\tau) \sqrt{32 \! + \! (\hat{r}_{0}
(\tau) \! - \! 4)^{2}} \left(1 \! + \! \mathcal{O} \left(\dfrac{\tau^{-2/3}
\varkappa_{0}^{2}(\tau)}{(\hat{r}_{0}(\tau))^{2}} \right) \right).
\label{eq3.48}
\end{equation}
Now, substituting Equation~\eqref{eq3.48} into Equation~\eqref{eq3.47}, one
arrives at
\begin{equation}
\mathfrak{p}_{\tau} \underset{\tau \to +\infty}{=} \mi \sqrt{32 \! + \!
(\hat{r}_{0}(\tau) \! - \! 4)^{2}} \left(1 \! + \! \mathcal{O} \left(
\dfrac{\tau^{-2/3} \varkappa_{0}^{2}(\tau)}{(\hat{r}_{0}(\tau))^{2}} \right)
\right). \label{eq3.49}
\end{equation}
Substituting the expansions~\eqref{eqF1}, \eqref{eqE1}, and~\eqref{eq3.49}
into Equation~\eqref{eq3.36} (with the upper signs in the formula for
$\digamma_{\tau}^{\sharp}$ taken (cf. Equation~\eqref{eqfteesharp})), and
taking into account conditions~\eqref{eq3.11}, one shows that
\begin{align}
\int_{\widetilde{\mu}_{0}}^{\widetilde{\mu}} \diag ((T(\xi))^{-1}
\partial_{\xi}T(\xi)) \, \md \xi \underset{\underset{\scriptstyle
\Re (\widetilde{\mu}) \to +\infty}{\tau \to +\infty}}{=}& \, \left(
\mathcal{I}_{\infty}(\tau) \! + \! \mathcal{O}(\tau^{\delta -2 \delta_{1}}
\ln \tau) \! + \! \mathcal{O}(\tau^{-2(\varepsilon -\delta)}) \! + \!
\mathcal{O}(\tau^{-\frac{1}{3}+ \varepsilon+ \delta}) \right. \nonumber \\
+&\left. \, \mathcal{O}(\tau^{-\frac{2}{3}+2 \delta}) \! + \! \mathcal{O}
(\hat{r}_{0}(\tau) \widetilde{\mu}^{-2}) \right) \sigma_{3}, \label{eqT1}
\end{align}
where
\begin{equation}
\mathcal{I}_{\infty}(\tau) \! = \! -\dfrac{\sqrt{32 \! + \! (\hat{r}_{0}
(\tau) \! - \! 4)^{2}}}{8 \sqrt{3}} \left(-\dfrac{1}{3} \ln \tau \! + \! \ln
\widetilde{\Lambda} \right) \! + \! \widehat{\mathcal{I}}_{\infty}(\tau),
\label{eq3.90}
\end{equation}
with
\begin{align}
\widehat{\mathcal{I}}_{\infty}(\tau) :=& \, \dfrac{\sqrt{32 \! + \!
(\hat{r}_{0}(\tau) \! - \! 4)^{2}}}{8 \sqrt{3}} \ln \left(\dfrac{6 \alpha}{
(\sqrt{3}+ \! 1)^{2}} \right) \! - \! \dfrac{1}{4} \ln \left(\dfrac{
\hat{r}_{0}(\tau) \! - \! 4 \! + \! \sqrt{32 \! + \! (\hat{r}_{0}(\tau) \!
- \! 4)^{2}}}{\hat{r}_{0}(\tau) \! - \! 4 \! - \! \sqrt{32 \! + \!
(\hat{r}_{0}(\tau) \! - \! 4)^{2}}} \right) \nonumber \\
+& \, \dfrac{1}{4} \ln \! \left(\! \left(\dfrac{\hat{r}_{0}(\tau) \! + \! 4
\sqrt{3}- \! \sqrt{32 \! + \! (\hat{r}_{0}(\tau) \! - \! 4)^{2}}}{\hat{r}_{0}
(\tau) \! + \! 4 \sqrt{3}+ \! \sqrt{32 \! + \! (\hat{r}_{0}(\tau) \! - \!
4)^{2}}} \right) \! \left(\dfrac{\hat{r}_{0}(\tau) \! - \! 4 \sqrt{3}+ \!
\sqrt{32 \! + \! (\hat{r}_{0}(\tau) \! - \! 4)^{2}}}{\hat{r}_{0}(\tau) \!
- \! 4 \sqrt{3}- \! \sqrt{32 \! + \! (\hat{r}_{0}(\tau) \! - \! 4)^{2}}}
\right) \! \right).
\label{eq3.92}
\end{align}
Taking the limit $\Re (\widetilde{\mu}) \! \to \! +0$ in
Equation~\eqref{eq3.36}, with $\digamma_{\tau}(\xi)$ and $\mathcal{E}_{T}^{
\natural}(\xi)$ defined by Equations~\eqref{eqFtee} and~\eqref{eqee2},
respectively, and using conditions~\eqref{eq3.11}, one obtains
\begin{align}
\digamma_{\tau}(\widetilde{\mu}) \underset{\underset{\scriptstyle \Re
(\widetilde{\mu}) \to +0}{\tau \to +\infty}}{=}& \, \dfrac{2}{\sqrt{32 \! +
\! (\hat{r}_{0}(\tau) \! - \! 4)^{2}}} \ln \left(\dfrac{\hat{r}_{0}(\tau)
\sqrt{32 \! + \! (\hat{r}_{0}(\tau) \! - \! 4)^{2}}+ \! ((\hat{r}_{0}
(\tau))^{2} \! - \! 4 \hat{r}_{0}(\tau) \! + \! 16)}{\hat{r}_{0}(\tau)
\sqrt{32 \! + \! (\hat{r}_{0}(\tau) \! - \! 4)^{2}}- \! ((\hat{r}_{0}
(\tau))^{2} \! - \! 4 \hat{r}_{0}(\tau) \! + \! 16)} \right) \nonumber \\
+& \, \dfrac{1}{\sqrt{3}} \ln (\me^{\mi \pi}) \! + \! \mathcal{O}
(\widetilde{\mu}), \label{eqF2}
\end{align}
and
\begin{equation}
\mathcal{E}_{T}^{\natural}(\widetilde{\mu}) \underset{\underset{\scriptstyle
\Re (\widetilde{\mu}) \to +0}{\tau \to +\infty}}{=} \mathcal{O}(\tau^{-2
\delta_{1}}) \! + \! \mathcal{O}(\tau^{-\frac{2}{3}+ 2 \delta}), \label{eqE2}
\end{equation}
whence, substituting the expansions~\eqref{eqF2}, \eqref{eqE2},
and~\eqref{eq3.49} into Equation~\eqref{eq3.36} (with the lower
signs in the formula for $\digamma_{\tau}^{\sharp}$ taken (cf.
Equation~\eqref{eqfteesharp})), and taking into account
conditions~\eqref{eq3.11}, one arrives at
\begin{align}
\int_{\widetilde{\mu}_{0}}^{\widetilde{\mu}} \diag ((T(\xi))^{-1}
\partial_{\xi}T(\xi)) \, \md \xi \underset{\underset{\scriptstyle
\Re (\widetilde{\mu}) \to +0}{\tau \to +\infty}}{=}& \, \left(\mathcal{I}_{0}
(\tau) \! + \! \mathcal{O}(\tau^{\delta -2 \delta_{1}} \ln \tau) \! + \!
\mathcal{O}(\tau^{-2(\varepsilon - \delta)}) \! + \! \mathcal{O}
(\tau^{-\frac{1}{3}+ \varepsilon + \delta}) \right. \nonumber \\
+&\left. \, \mathcal{O}(\tau^{-\frac{2}{3}+ 2 \delta}) \! + \! \mathcal{O}
(\hat{r}_{0}(\tau) \widetilde{\mu}) \right) \sigma_{3}, \label{eqT2}
\end{align}
where
\begin{equation}
\mathcal{I}_{0}(\tau) \! = \! \dfrac{\sqrt{32 \! + \! (\hat{r}_{0}(\tau) \! -
\! 4)^{2}}}{8 \sqrt{3}} \left(-\dfrac{1}{3} \ln \tau \! + \! \ln \widetilde{
\Lambda} \right) \! + \! \widehat{\mathcal{I}}_{0}(\tau), \label{eq3.91}
\end{equation}
with
\begin{align}
\widehat{\mathcal{I}}_{0}(\tau) :=& \, -\dfrac{\sqrt{32 \! + \! (\hat{r}_{0}
(\tau) \! - \! 4)^{2}}}{8 \sqrt{3}} \ln (3 \alpha \me^{-\mi \pi}) \! - \!
\dfrac{1}{4} \ln \left(\dfrac{\hat{r}_{0}(\tau) \! + \! 4 \! + \! \sqrt{32 \!
+ \! (\hat{r}_{0}(\tau) \! - \! 4)^{2}}}{-\hat{r}_{0}(\tau) \! - \! 4 \! + \!
\sqrt{32 \! + \! (\hat{r}_{0}(\tau) \! - \! 4)^{2}}} \right) \nonumber \\
-& \, \dfrac{1}{4} \ln \! \left(\! \left(\dfrac{\hat{r}_{0}(\tau) \! + \! 4
\sqrt{3}- \! \sqrt{32 \! + \! (\hat{r}_{0}(\tau) \! - \! 4)^{2}}}{\hat{r}_{0}
(\tau) \! + \! 4 \sqrt{3}+ \! \sqrt{32 \! + \! (\hat{r}_{0}(\tau) \! - \!
4)^{2}}} \right) \! \left(\dfrac{\hat{r}_{0}(\tau) \! - \! 4 \sqrt{3}+ \!
\sqrt{32 \! + \! (\hat{r}_{0}(\tau) \! - \! 4)^{2}}}{\hat{r}_{0}(\tau) \!
- \! 4 \sqrt{3}- \! \sqrt{32 \! + \! (\hat{r}_{0}(\tau) \! - \! 4)^{2}}}
\right) \! \right) \nonumber \\
+& \, \dfrac{1}{4} \ln \left(\dfrac{\hat{r}_{0}(\tau) \sqrt{32 \! + \!
(\hat{r}_{0}(\tau) \! - \! 4)^{2}}+ \! ((\hat{r}_{0}(\tau))^{2} \! - \! 4
\hat{r}_{0}(\tau) \! + \! 16)}{\hat{r}_{0}(\tau) \sqrt{32 \! + \! (\hat{r}_{0}
(\tau) \! - \! 4)^{2}}- \! ((\hat{r}_{0}(\tau))^{2} \! - \! 4 \hat{r}_{0}
(\tau) \! + \! 16)} \right). \label{eq3.93}
\end{align}
One now simplifies Equations~\eqref{eq3.92} and~\eqref{eq3.93}; in order to do
so, however, several estimates are necessary. Rewrite Equation~\eqref{eq3.14}
as follows:
\begin{equation}
\dfrac{(\hat{u}_{0}(\tau))^{2}}{1 \! + \! \hat{u}_{0}(\tau)} \! + \!
\dfrac{\hat{r}_{0}(\tau) \hat{u}_{0}(\tau)}{2(1 \! + \! \hat{u}_{0}(\tau))}
\! - \! \dfrac{(\hat{r}_{0}(\tau))^{2}}{8} \! = \! -\dfrac{\tau^{-2/3}
\varkappa_{0}^{2}(\tau)}{8 \alpha^{2}}. \label{eq3.65}
\end{equation}
Via Equations~\eqref{eqpeetee} and~\eqref{eq3.49}, using
Equation~\eqref{eq3.65} to eliminate $(\hat{u}_{0}(\tau))^{2}$, and taking
into account conditions~\eqref{eq3.11}, one arrives at, after simplification,
\begin{align}
\hat{r}_{0}(\tau) \sqrt{32 \! + \! (\hat{r}_{0}(\tau) \! - \! 4)^{2}}- \!
((\hat{r}_{0}(\tau))^{2} \! - \! 4 \hat{r}_{0}(\tau) \! + \! 16) \underset{
\tau \to +\infty}{=}& \, \dfrac{8(-2 \! + \! \hat{r}_{0}(\tau))}{1 \! + \!
\hat{u}_{0}(\tau)} \! \left(1 \! + \! \mathcal{O}(\tau^{-\frac{2}{3}+
\delta}) \right. \nonumber \\
+&\left. \, \mathcal{O}(\tau^{-\frac{1}{3}+ \delta -\delta_{1}}) \right).
\label{eq3.94}
\end{align}
The estimates~\eqref{eq3.95}--\eqref{eq3.99} below are derived analogously:
\begin{align}
\hat{r}_{0}(\tau) \sqrt{32 \! + \! (\hat{r}_{0}(\tau) \! - \! 4)^{2}}+ \!
((\hat{r}_{0}(\tau))^{2} \! - \! 4 \hat{r}_{0}(\tau) \! + \! 16) \underset{
\tau \to +\infty}{=}& \, 16(1 \! + \! \hat{u}_{0}(\tau)) \! \left(1 \!
+ \! \mathcal{O}(\tau^{-\frac{2}{3}+ \delta}) \right. \nonumber \\
+&\left. \, \mathcal{O}(\tau^{-\frac{1}{3}+ \delta -\delta_{1}}) \right),
\label{eq3.95}
\end{align}
\begin{gather}
\hat{r}_{0}(\tau) \! + \! 4 \! + \! \sqrt{32 \! + \! (\hat{r}_{0}(\tau) \! -
\! 4)^{2}} \underset{\tau \to +\infty}{=} \dfrac{2 \hat{r}_{0}(\tau)(1 \! +
\! \hat{u}_{0}(\tau))}{\hat{u}_{0}(\tau)} \! \left(1 \! + \! \mathcal{O}
(\tau^{-\frac{2}{3}- \delta}) \! + \! \mathcal{O}(\tau^{\delta -2 \delta_{1}})
\right), \label{eq3.96} \\
\hat{r}_{0}(\tau) \! + \! 4 \! - \! \sqrt{32 \! + \! (\hat{r}_{0}(\tau) \! -
\! 4)^{2}} \underset{\tau \to +\infty}{=} \dfrac{8 \hat{u}_{0}(\tau)(-2 \! +
\! \hat{r}_{0}(\tau))}{\hat{r}_{0}(\tau)(1 \! + \! \hat{u}_{0}(\tau))} \!
\left(1 \! + \! \mathcal{O}(\tau^{-\frac{2}{3}- \delta}) \! + \! \mathcal{O}
(\tau^{\delta -2 \delta_{1}}) \right), \label{eq3.97} \\
\hat{r}_{0}(\tau) \! - \! 4 \! + \! \sqrt{32 \! + \! (\hat{r}_{0}(\tau) \!
- \! 4)^{2}} \underset{\tau \to +\infty}{=} \dfrac{16 \hat{u}_{0}(\tau)}{
\hat{r}_{0}(\tau)} \! \left(1 \! + \! \mathcal{O}(\tau^{-\frac{2}{3}- \delta})
\! + \! \mathcal{O}(\tau^{\delta -2 \delta_{1}}) \right), \label{eq3.98} \\
\hat{r}_{0}(\tau) \! - \! 4 \! - \! \sqrt{32 \! + \! (\hat{r}_{0}(\tau) \!
- \! 4)^{2}} \underset{\tau \to +\infty}{=} -\dfrac{2 \hat{r}_{0}(\tau)}{
\hat{u}_{0}(\tau)} \! \left(1 \! + \! \mathcal{O}(\tau^{-\frac{2}{3}- \delta})
\! + \! \mathcal{O}(\tau^{\delta -2 \delta_{1}}) \right). \label{eq3.99}
\end{gather}
Substituting the estimates~\eqref{eq3.98} and~\eqref{eq3.99} into
Equation~\eqref{eq3.92}, one proves that
\begin{equation}
\widehat{\mathcal{I}}_{\infty}(\tau) \underset{\tau \to +\infty}{=}
\mathcal{I}_{\infty}^{\sharp}(\tau) \! + \! \mathcal{O}(\tau^{-\frac{2}{3}-
\delta}) \! + \! \mathcal{O}(\tau^{\delta -2 \delta_{1}}), \label{eq3.100}
\end{equation}
where $\mathcal{I}_{\infty}^{\sharp}(\tau)$ is defined by
Equation~\eqref{eq3.87}, and, substituting the
estimates~\eqref{eq3.94}--\eqref{eq3.97} into Equation~\eqref{eq3.93}, one
proves that
\begin{equation}
\widehat{\mathcal{I}}_{0}(\tau) \underset{\tau \to +\infty}{=} \mathcal{I}_{
0}^{\sharp}(\tau) \! + \! \mathcal{O}(\tau^{-\frac{2}{3}+ \delta}) \! + \!
\mathcal{O}(\tau^{-\frac{1}{3}+ \delta -\delta_{1}}) \! + \! \mathcal{O}
(\tau^{\delta -2 \delta_{1}}), \label{eq3.101}
\end{equation}
where $\mathcal{I}_{0}^{\sharp}(\tau)$ is defined by Equation~\eqref{eq3.89}.
Hence, via expansions~\eqref{eq3.100} and~\eqref{eq3.101}, and
Equations \eqref{eqT1}, \eqref{eq3.90}, \eqref{eqT2}, and~\eqref{eq3.91},
one obtains the results stated in the Proposition. \hfill $\qed$
\begin{bbbb} \label{prop3.1.5}
Let $T(\widetilde{\mu})$ be given in Equation~\eqref{eq3.18}, with
$\mathcal{A}(\widetilde{\mu})$ defined by Equation~\eqref{eq3.4}, and
$l(\widetilde{\mu})$ given in Equation~\eqref{eq3.19} with the branches
defined as in Corollary~{\rm \ref{cor3.1.4}}. Then
\begin{align}
T(\widetilde{\mu}) \underset{\underset{\scriptstyle \Re (\widetilde{\mu})
\to +\infty}{\tau \to +\infty}}{=}& \, (b(\tau))^{-\frac{1}{2} \ad
(\sigma_{3})} \left(\mathrm{I} \! + \! \dfrac{1}{\widetilde{\mu}}
\begin{pmatrix}
0 & -2 \alpha^{4}(1 \! + \! \hat{u}_{0}(\tau)) \\
\frac{(2-\hat{r}_{0}(\tau)) \alpha^{2}+(a-\frac{\mi}{2}) \tau^{-2/3}}{8
\alpha^{4}(1+ \hat{u}_{0}(\tau))} & 0
\end{pmatrix} \right. \nonumber \\
+&\left. \, \mathcal{O} \! \left(\dfrac{(c_{1}(\tau)(\hat{r}_{0}(\tau))^{2}
\! + \! c_{2}(\tau))}{\widetilde{\mu}^{2}} \right) \right), \label{eq3.51}
\end{align}
and
\begin{align}
T(\widetilde{\mu}) \underset{\underset{\scriptstyle \Re (\widetilde{\mu}) \to
+0}{\tau \to +\infty}}{=}& \, \dfrac{1}{\sqrt{2}} \! \left(\dfrac{b(\tau)}{2
\sqrt{2} \, \alpha^{3}} \right)^{-\frac{1}{2} \ad (\sigma_{3})} \left(
\begin{pmatrix}
1 & 1 \\
-1 & 1
\end{pmatrix} \! + \! \widetilde{\mu} \, \dfrac{(-2 \! + \! \hat{r}_{0}
(\tau))}{4 \sqrt{2} \, \alpha}
\begin{pmatrix}
1 & -1 \\
1 & 1
\end{pmatrix} \right. \nonumber \\
+&\left. \, \mathcal{O} \! \left(\widetilde{\mu}^{2}(c_{3}(\tau)(\hat{r}_{0}
(\tau))^{2} \! + \! c_{4}(\tau) \hat{u}_{0}(\tau) \! + \! c_{5}(\tau))
\right) \vphantom{M^{M^{M^{M^{M}}}}} \! \right),
\label{eq3.52}
\end{align}
where $c_{j}(\tau) \! =_{\tau \to +\infty} \! \mathcal{O}(1)$, $j \! = \! 1,
\dotsc,5$.
\end{bbbb}

\emph{Proof}. Here, the proof of estimate~\eqref{eq3.51} is presented: the
estimate~\eqref{eq3.52} is proven analogously. Via Equations~\eqref{eq3.9},
\eqref{eq3.10}, and~\eqref{eq3.12}, and conditions~\eqref{eq3.11}, one shows
that
\begin{align*}
l(\widetilde{\mu}) \underset{\underset{\scriptstyle \Re (\widetilde{\mu}) \to
+\infty}{\tau \to +\infty}}{=}& \, \, 2 \widetilde{\mu} \! + \! \dfrac{1}{
\widetilde{\mu}} \left(a \! - \! \dfrac{\mi}{2} \right) \tau^{-2/3} \! + \!
\mathcal{O}(\widetilde{\mu}^{-3}), \\
\mi (\mathcal{A}(\widetilde{\mu}) \! - \! \mi l(\widetilde{\mu}) \sigma_{3})
\sigma_{3} \underset{\underset{\scriptstyle \Re (\widetilde{\mu}) \to
+\infty}{\tau \to +\infty}}{=}& \, 4 \mathrm{I} \widetilde{\mu} \! + \!
\begin{pmatrix}
0 & -\frac{4 \sqrt{-a(\tau)b(\tau)}}{b(\tau)} \\
-2 \mi d(\tau) & 0
\end{pmatrix} \! + \! \dfrac{1}{\widetilde{\mu}} \left((-2 \! + \! \hat{r}_{0}
(\tau)) \alpha^{2} \right. \\
+&\left. \, \left(a \! - \! \dfrac{\mi}{2} \right) \tau^{-2/3} \right)
\mathrm{I} \! + \! \dfrac{1}{\widetilde{\mu}^{2}}
\begin{pmatrix}
0 & \frac{8 \alpha^{6}}{b(\tau)} \\
-b(\tau) & 0
\end{pmatrix} \! + \! \mathcal{O}(\widetilde{\mu}^{-3} \mathrm{I}), \\
\dfrac{1}{\sqrt{2 \mi l(\widetilde{\mu})(\mathcal{A}_{11}(\widetilde{\mu})
\! - \! \mi l(\widetilde{\mu}))}} \underset{\underset{\scriptstyle \Re
(\widetilde{\mu}) \to +\infty}{\tau \to +\infty}}{=}& \, \, \dfrac{1}{4
\widetilde{\mu}} \left(1 \! - \! \dfrac{1}{8 \widetilde{\mu}^{2}} \left(
(-2 \! + \! \hat{r}_{0}(\tau)) \alpha^{2} \! + \! 3 \left(a \! - \!
\dfrac{\mi}{2} \right) \tau^{-2/3} \right) \right. \\
+&\left. \, \mathcal{O} \! \left(\dfrac{(C_{1}(\tau)(\hat{r}_{0}(\tau))^{2}
\! + \! C_{2}(\tau))}{\widetilde{\mu}^{4}} \right) \right),
\end{align*}
where $C_{j}(\tau) \! =_{\tau \to +\infty} \! \mathcal{O}(1)$, $j \!
= \! 1,2$; hence, via Equation~\eqref{eq3.18} and conditions~\eqref{eq3.11},
one arrives at the estimate~\eqref{eq3.51}. \hfill $\qed$
\begin{bbbb} \label{prop3.1.6}
Let $T(\widetilde{\mu})$ be given in Equation~\eqref{eq3.18}, with
$\mathcal{A}(\widetilde{\mu})$ defined by Equation~\eqref{eq3.4}, and $l^{2}
(\widetilde{\mu})$ given in Equation~\eqref{eq3.19}. Under the conditions of
Corollary~{\rm \ref{cor3.1.3}}$:$
\begin{align}
T_{11}(\widetilde{\mu}) \! = \! T_{22}(\widetilde{\mu}) \underset{\underset{
\scriptstyle \tau \to +\infty}{\widetilde{\mu}=\widetilde{\mu}_{0}}}{=}& \,
\dfrac{\alpha \hat{r}_{0}(\tau) \tau^{1/6}}{\sqrt{8 \sqrt{3} \, \alpha
\hat{r}_{0}(\tau) \varpi \widetilde{\Lambda}}} \left(1 \! - \! \dfrac{h_{0}
(\tau)}{48 \alpha^{2} \widetilde{\Lambda}^{2}} \! - \! \dfrac{(a \! - \!
\frac{\mi}{2})}{48 \widetilde{\Lambda}^{2}} \! + \! \dfrac{1}{2 \alpha}
\left(\dfrac{4 \! - \! \hat{r}_{0}(\tau) \! + \! 4 \sqrt{3} \,
\varpi}{\hat{r}_{0}(\tau)} \! + \! \dfrac{7}{6} \right) \right. \nonumber \\
\times&\left. \, \tau^{-1/3} \widetilde{\Lambda} \! + \! \mathcal{O}
(\tau^{2 \delta -4 \varepsilon}) \! + \! \mathcal{O}(\tau^{-\varepsilon
-\delta_{1}}) \! + \! \mathcal{O}(\tau^{2 \varepsilon -2 \delta_{1}})
\vphantom{M^{M^{M^{M^{M^{M}}}}}} \right), \label{eq3.53}
\end{align}
\begin{align}
T_{12}(\widetilde{\mu}) \underset{\underset{\scriptstyle \tau \to +\infty}{
\widetilde{\mu}=\widetilde{\mu}_{0}}}{=}& \, -\dfrac{8 \alpha^{4} \hat{u}_{0}
(\tau) \tau^{1/6}}{b(\tau) \sqrt{8 \sqrt{3} \, \alpha \hat{r}_{0}(\tau) \varpi
\widetilde{\Lambda}}} \left(1 \! - \! \dfrac{h_{0}(\tau)}{48 \alpha^{2}
\widetilde{\Lambda}^{2}} \! - \! \dfrac{(a \! - \! \frac{\mi}{2})}{48
\widetilde{\Lambda}^{2}} \! - \! \left(\dfrac{1}{2 \alpha} \! \left(\dfrac{4
\! - \! \hat{r}_{0}(\tau) \! + \! 4 \sqrt{3} \, \varpi}{\hat{r}_{0}(\tau)}
\! - \! \dfrac{7}{6} \right) \right. \right. \nonumber \\
-&\left. \left. \, \dfrac{2}{\alpha \hat{u}_{0}(\tau)} \right) \tau^{-1/3}
\widetilde{\Lambda} \! + \! \mathcal{O}(\tau^{2 \delta -4 \varepsilon}) \! +
\! \mathcal{O}(\tau^{-\varepsilon -\delta_{1}}) \right), \label{eq3.54} \\
T_{21}(\widetilde{\mu}) \underset{\underset{\scriptstyle \tau \to +\infty}{
\widetilde{\mu}=\widetilde{\mu}_{0}}}{=}& \, -\dfrac{b(\tau) \tau^{1/6}}{2
\alpha^{2} \sqrt{8 \sqrt{3} \, \alpha \hat{r}_{0}(\tau) \varpi \widetilde{
\Lambda}}} \! \left(\dfrac{\hat{r}_{0}(\tau) \! + \! 2 \hat{u}_{0}(\tau) \!
- \! \alpha^{-2}(a \! - \! \frac{\mi}{2}) \tau^{-2/3}}{1 \! + \! \hat{u}_{0}
(\tau)} \right) \! \left(1 \! - \! \dfrac{h_{0}(\tau)}{48 \alpha^{2}
\widetilde{\Lambda}^{2}} \! - \! \dfrac{(a \! - \! \frac{\mi}{2})}{48
\widetilde{\Lambda}^{2}} \right. \nonumber \\
-&\left. \, \left(\dfrac{1}{2 \alpha} \! \left(\dfrac{4 \! - \! \hat{r}_{0}
(\tau) \! + \! 4 \sqrt{3} \, \varpi}{\hat{r}_{0}(\tau)} \! - \! \dfrac{7}{6}
\right) \! + \! \dfrac{4}{\alpha} \! \left(\dfrac{\hat{r}_{0}(\tau) \! + \! 2
\hat{u}_{0}(\tau) \! - \! \alpha^{-2}(a \! - \! \frac{\mi}{2}) \tau^{-2/3}}{1
\! + \! \hat{u}_{0}(\tau)} \right)^{-1} \right) \! \tau^{-1/3} \widetilde{
\Lambda} \right. \nonumber \\
+&\left. \, \mathcal{O}(\tau^{2 \delta -4 \varepsilon}) \! + \! \mathcal{O}
(\tau^{-\varepsilon -\delta_{1}}) \vphantom{M^{M^{M^{M^{M^{M}}}}}} \right),
\label{eq3.55}
\end{align}
where
\begin{equation*}
\varpi \! := \!
\begin{cases}
+1, &\text{$\arg (\widetilde{\Lambda}) \! = \! 0$,} \\
-1, &\text{$\arg (\widetilde{\Lambda}) \! = \! \pi$.}
\end{cases}
\end{equation*}
\end{bbbb}

\emph{Proof}. Using the definition of $\mathcal{A}(\widetilde{\mu})$ (cf.
Equation~\eqref{eq3.4}), together with Equations~\eqref{eq3.9},
\eqref{eq3.10}, and~\eqref{eq3.12}, and conditions~\eqref{eq3.11} and those
in Corollary~\ref{cor3.1.3}, one derives the following expansions:
\begin{align*}
\dfrac{1}{\sqrt{2 \mi l(\widetilde{\mu})(\mathcal{A}_{11}(\widetilde{\mu})
\! - \! \mi l(\widetilde{\mu}))}} \underset{\underset{\scriptstyle \tau
\to +\infty}{\widetilde{\mu}=\widetilde{\mu}_{0}}}{=}& \, \dfrac{\tau^{1/6}}{
\sqrt{8 \sqrt{3} \, \alpha \hat{r}_{0}(\tau) \varpi \widetilde{\Lambda}}}
\left(1 \! - \! \dfrac{h_{0}(\tau)}{48 \alpha^{2} \widetilde{\Lambda}^{2}}
\! - \! \dfrac{(a \! - \! \frac{\mi}{2})}{48 \widetilde{\Lambda}^{2}} \!
- \! \dfrac{1}{2 \alpha} \left(\dfrac{4 \! - \! \hat{r}_{0}(\tau) \! + \!
4 \sqrt{3} \, \varpi}{\hat{r}_{0}(\tau)} \right. \right. \\
-&\left. \left. \, \dfrac{7}{6} \right) \tau^{-1/3} \widetilde{\Lambda} \!
+ \! \mathcal{O} \! \left(\dfrac{(h_{0}(\tau))^{2}}{\widetilde{\Lambda}^{4}}
\right) \! + \! \mathcal{O} \! \left(\dfrac{\tau^{-1/3}h_{0}(\tau)}{\hat{r}_{0}
(\tau) \widetilde{\Lambda}} \right) \! \right),
\end{align*}
\begin{align*}
\mi \mathcal{A}_{11}(\widetilde{\mu}) \! + \! l(\widetilde{\mu}) \underset{
\underset{\scriptstyle \tau \to +\infty}{\widetilde{\mu}=\widetilde{\mu}_{
0}}}{=}& \, \, \alpha \hat{r}_{0}(\tau) \! \left(1 \! + \! \left(
\dfrac{4 \! - \! \hat{r}_{0}(\tau) \! + \! 4 \sqrt{3} \, \varpi}{\alpha
\hat{r}_{0}(\tau)} \! + \! \dfrac{4 \sqrt{3} \, \varpi}{\alpha \hat{r}_{0}
(\tau)} \left(\dfrac{h_{0}(\tau)}{24 \alpha^{2} \widetilde{\Lambda}^{2}} \!
+ \! \dfrac{(a \! - \! \frac{\mi}{2})}{24 \widetilde{\Lambda}^{2}} \right)
\right) \tau^{-1/3} \widetilde{\Lambda} \right. \\
+&\left. \, \left(\dfrac{-2 \! + \! \hat{r}_{0}(\tau) \! - \! \frac{14}{
\sqrt{3}} \varpi}{\alpha^{2} \hat{r}_{0}(\tau)} \right) \tau^{-2/3}
\widetilde{\Lambda}^{2} \! + \! \mathcal{O} \! \left(\dfrac{\tau^{-1/3}(h_{0}
(\tau))^{2}}{\hat{r}_{0}(\tau) \widetilde{\Lambda}^{3}} \right) \right), \\
-\mi \mathcal{A}_{12}(\widetilde{\mu}) \underset{\underset{\scriptstyle \tau
\to +\infty}{\widetilde{\mu}=\widetilde{\mu}_{0}}}{=}& \, \, \dfrac{1}{b
(\tau)} \! \left(-8 \alpha^{4} \hat{u}_{0}(\tau) \! - \! 16 \alpha^{3}
\tau^{-1/3} \widetilde{\Lambda} \! + \! 24 \alpha^{2} \tau^{-2/3} \widetilde{
\Lambda}^{2} \! + \! \mathcal{O}(\tau^{-1} \widetilde{\Lambda}^{3}) \right), \\
\mi \mathcal{A}_{21}(\widetilde{\mu}) \underset{\underset{\scriptstyle
\tau \to +\infty}{\widetilde{\mu}=\widetilde{\mu}_{0}}}{=}& \, \, b(\tau)
\! \left(-\dfrac{(-2 \! + \! \hat{r}_{0}(\tau))}{2 \alpha^{2}(1 \! + \!
\hat{u}_{0}(\tau))} \! - \! \dfrac{1}{\alpha^{2}} \! + \! \dfrac{(a \! - \!
\frac{\mi}{2}) \tau^{-2/3}}{2 \alpha^{4}(1 \! + \! \hat{u}_{0}(\tau))} \!
+ \! \dfrac{2}{\alpha^{3}} \tau^{-1/3} \widetilde{\Lambda} \! - \!
\dfrac{3}{\alpha^{4}} \tau^{-2/3} \widetilde{\Lambda}^{2} \right. \\
+&\left. \, \mathcal{O}(\tau^{-1} \widetilde{\Lambda}^{3}) \right).
\end{align*}
According to the choice of the branch of $l(\xi)$ in
Corollary~\ref{cor3.1.4}, $l(\xi) \! > \! 0$ for positive $\xi$ outside the
neighborhood of $\xi \! = \! \alpha$ containing the double-turning points;
therefore, one arrives at the definition of $\varpi$ stated in the
Proposition. The asymptotics~\eqref{eq3.53}--\eqref{eq3.55} are obtained by
substituting the above expansions into Equations~\eqref{eqintertee}. \hfill
$\qed$
\subsection{The Model Problem and Asymptotics Near the Turning Points}
\label{sec3.2}
For the calculation of the monodromy data, one needs an approximation that is
more accurate than that given by the WKB formula (cf. Equation~\eqref{eq3.16})
for the solution of Equation~\eqref{eq3.3} in proper neighborhoods of the
turning points. There are two simple turning points approaching $\pm \mi
\sqrt{2} \, \alpha$: the approximate solution of Equation~\eqref{eq3.3} in
neighborhoods of these turning points is given in terms of Airy functions.
There are also two pairs of turning points, one pair coalescing at $-\alpha$
and another pair coalescing at $+\alpha$ (double-turning points): it is well
known that, in neighborhoods of $\pm \alpha$, the approximate solution of
Equation~\eqref{eq3.3} is expressed in terms of parabolic-cylinder functions
(see, for example, \cite{W,F}). In order to obtain asymptotics of $u(\tau)$,
and the associated functions $\mathcal{H}(\tau)$ and $f(\tau)$, it is
sufficient to study a subset of the complete set of the monodromy data,
which can be calculated via the approximation of the general solution
of Equation~\eqref{eq3.3} in a neighborhood of the double-turning point
$+\alpha$. For the asymptotic conditions~\eqref{eq3.11} on the functions
$\hat{r}_{0}(\tau)$, $\hat{u}_{0}(\tau)$, and $h_{0}(\tau)$, this
approximation is not straightforward, and is given in Lemma~\ref{lem3.2.1}
below.
\begin{cccc} \label{lem3.2.1}
Let $\widetilde{\mu} \! = \! \alpha \! + \! \tau^{-1/3} \widetilde{\Lambda}$,
where $\widetilde{\Lambda} \! =_{\tau \to +\infty} \! \mathcal{O}(\tau^{
\varepsilon})$, $0 \! < \! \varepsilon \! < \! 1/9$, and
\begin{equation}
\nu \! + \! 1 \! := \! -\dfrac{\mi}{2 \sqrt{3}} \! \left(\dfrac{h_{0}(\tau)}{
\alpha^{2}} \! + \! \left(a \! - \! \dfrac{\mi}{2} \right) \! - \!
\dfrac{\mathfrak{p}_{\tau}}{4} \right) \! + \! \dfrac{1}{2}, \label{eq3.63}
\end{equation}
where $\mathfrak{p}_{\tau}$ is defined by Equation~\eqref{eqpeetee}. In
conjunction with Equations~\eqref{eq3.8}--\eqref{eq3.10} and~\eqref{eq3.7},
and conditions~\eqref{eq3.11}, impose the following restrictions:
\begin{gather}
0 \underset{\tau \to +\infty}{<} \Re (\nu \! + \! 1) \underset{\tau \to
+\infty}{<} 1, \quad \quad \Im (\nu \! + \! 1) \underset{\tau \to +\infty}{=}
\mathcal{O}(1), \label{newconds1} \\
0 \! < \! \delta \! < \! \varepsilon \! < \! \dfrac{1}{9}, \qquad 6
\varepsilon \! + \! 2 \varepsilon \Re (\nu \! + \! 1) \! + \!
\dfrac{\delta}{2} \! < \! \delta_{1} \! < \! \dfrac{1}{3}. \label{newcconds2}
\end{gather}
Then there exists a fundamental solution of Equation~\eqref{eq3.3} with
asymptotic representation
\begin{equation}
\widetilde{\Psi}(\widetilde{\mu},\tau) \underset{\tau \to +\infty}{=}
\mathcal{F}_{\tau}(\widetilde{\Lambda}) \! \left(\mathrm{I} \! + \!
\mathcal{O} \! \left(\tau^{6 \varepsilon +2 \varepsilon \Re (\nu +1)+
\frac{\delta}{2} -\delta_{1}} \right) \right) \! \psi_{0}
(\widetilde{\Lambda}), \label{eq3.56}
\end{equation}
where
\begin{equation}
\mathcal{F}_{\tau}(\widetilde{\Lambda}) \! := \!
\begin{pmatrix}
\frac{\mi \tau^{1/6} \sqrt{-8 \alpha^{4} \hat{u}_{0}(\tau)}}{\sqrt{b(\tau)}
\, \sqrt{\varkappa_{0}(\tau)}} & 0 \\
\frac{\mi \sqrt{b(\tau)} \, \left(\alpha \hat{r}_{0}(\tau) \tau^{1/6}- \mi
\tau^{-1/6} \ell_{\tau} \widetilde{\Lambda} \right)}{\sqrt{\varkappa_{0}
(\tau)} \, \sqrt{-8 \alpha^{4} \hat{u}_{0}(\tau)}} & -\frac{\mi \sqrt{
\varkappa_{0}(\tau)} \, \sqrt{b(\tau)} \, \tau^{-1/6}}{\sqrt{-8 \alpha^{4}
\hat{u}_{0}(\tau)}}
\end{pmatrix}, \label{eq3.57}
\end{equation}
$\varkappa_{0}^{2}(\tau)$ is defined by Equation~\eqref{eq3.58},
\begin{equation}
\ell_{\tau} \! = \! -\mathfrak{p}_{\tau} \! + \! \sqrt{\mathfrak{p}_{\tau}^{
2} \! - \! \hat{r}_{0}(\tau)(8 \! - \! \hat{r}_{0}(\tau))} \, \underset{\tau
\to +\infty}{=} \ell_{\tau}^{\infty} \! + \! \mathcal{O} \! \left(\dfrac{
\tau^{-2/3} \varkappa_{0}^{2}(\tau)(\mathcal{O}(1) \! + \! (\hat{r}_{0}
(\tau))^{2})}{(\hat{r}_{0}(\tau))^{2}} \right), \label{eq3.59}
\end{equation}
with
\begin{equation}
\ell_{\tau}^{\infty} \! := \! -\mathfrak{p}_{\tau} \! + \! \mi 4 \sqrt{3},
\label{eqlteeinf}
\end{equation}
and $\psi_{0}(\widetilde{\Lambda})$ is a fundamental solution of
\begin{equation}
\dfrac{\partial \psi_{0}(\widetilde{\Lambda})}{\partial \widetilde{\Lambda}}
\! = \! \left(\mi 4 \sqrt{3} \, \widetilde{\Lambda} \sigma_{3} \! + \!
\widetilde{q}(\tau) \sigma_{-} \! + \! \widetilde{p}(\tau) \sigma_{+}
\right) \! \psi_{0}(\widetilde{\Lambda}), \label{eq3.60}
\end{equation}
where
\begin{equation}
\widetilde{p}(\tau) \! = \! -\mi \varkappa_{0}(\tau), \qquad \quad
\widetilde{q}(\tau) \! = \! -\dfrac{\mi}{\varkappa_{0}(\tau)} \! \left(
\varkappa_{0}^{2}(\tau) \! + \! \ell_{\tau}^{\infty} \! + \! \dfrac{4(a \! -
\! \frac{\mi}{2}) \hat{u}_{0}(\tau)}{1 \! + \! \hat{u}_{0}(\tau)} \right).
\label{eq3.61}
\end{equation}
The function $\psi_{0}(\widetilde{\Lambda})$ can be explicitly presented in
the following form:
\begin{equation}
\psi_{0}(\widetilde{\Lambda}) \! = \!
\begin{pmatrix}
D_{-1-\nu}(\mi \me^{\frac{\mi \pi}{4}}2^{3/2}3^{1/4} \widetilde{\Lambda}) &
D_{\nu}(\me^{\frac{\mi \pi}{4}}2^{3/2}3^{1/4} \widetilde{\Lambda}) \\
\hat{\partial}_{\widetilde{\Lambda}}D_{-1-\nu}(\mi \me^{\frac{\mi \pi}{4}}
2^{3/2}3^{1/4} \widetilde{\Lambda}) & \hat{\partial}_{\widetilde{\Lambda}}
D_{\nu}(\me^{\frac{\mi \pi}{4}}2^{3/2}3^{1/4} \widetilde{\Lambda})
\end{pmatrix}, \label{eq3.62}
\end{equation}
where $D_{\pmb{\ast}}(\pmb{\cdot})$ is the parabolic-cylinder function
{\rm \cite{a24}}, and $\hat{\partial}_{\widetilde{\Lambda}} \! := \!
(\widetilde{p}(\tau))^{-1}(\tfrac{\partial}{\partial \widetilde{\Lambda}} \!
- \! \mi 4 \sqrt{3} \, \widetilde{\Lambda})$.
\end{cccc}

\emph{Proof}. The derivation of approximation~\eqref{eq3.56} consists
of the following sequence of invertible linear transformations, $F_{j}
\colon \operatorname{SL}_{2}(\mathbb{C}) \! \to \! \operatorname{SL}_{2}
(\mathbb{C})$, $j \! = \! 1,\dotsc,7$:
\begin{align*}
\text{\pmb{(i)}} \quad F_{1} \colon& \widetilde{\Psi}(\widetilde{\mu}) \!
\mapsto \! \Phi (\widetilde{\Lambda}) \! := \! \widetilde{\Psi}(\alpha \! +
\! \tau^{-1/3} \widetilde{\Lambda}), \\
\text{\pmb{(ii)}} \quad F_{2} \colon& \Phi (\widetilde{\Lambda}) \! \mapsto
\! \widehat{\Phi}(\widetilde{\Lambda}) \! := \! (b(\tau))^{\frac{1}{2}
\sigma_{3}} \Phi (\widetilde{\Lambda}), \\
\text{\pmb{(iii)}} \quad F_{3} \colon& \widehat{\Phi}(\widetilde{\Lambda})
\! \mapsto \! \phi (\widetilde{\Lambda}) \! := \! (\mathcal{N}(\tau))^{-1}
\widehat{\Phi}(\widetilde{\Lambda}), \\
\text{\pmb{(iv)}} \quad F_{4} \colon& \phi (\widetilde{\Lambda}) \! \mapsto
\! \widetilde{\Phi}(\widetilde{\Lambda}) \! := \! \tau^{-\frac{1}{6}
\sigma_{3}} \phi (\widetilde{\Lambda}), \\
\text{\pmb{(v)}} \quad F_{5} \colon& \widetilde{\Phi}(\widetilde{\Lambda})
\! \mapsto \! \widehat{\phi}(\widetilde{\Lambda}) \! := \!
\left(
\begin{smallmatrix}
1 & 0 \\
-\ell_{\tau} \widetilde{\Lambda} & 1
\end{smallmatrix}
\right) \widetilde{\Phi}(\widetilde{\Lambda}), \\
\text{\pmb{(vi)}} \quad F_{6} \colon& \widehat{\phi}(\widetilde{\Lambda})
\! \mapsto \! \psi (\widetilde{\Lambda}) \! := \! (\mathcal{G}(\tau))^{-1}
\widehat{\phi}(\widetilde{\Lambda}), \\
\text{\pmb{(vii)}} \quad F_{7} \colon& \psi (\widetilde{\Lambda}) \! \mapsto
\! \psi_{0}(\widetilde{\Lambda}) \! := \! (\chi_{0}(\widetilde{\Lambda}))^{-1}
\psi (\widetilde{\Lambda}), \\
\end{align*}
where the unimodular, matrix-valued functions $\mathcal{N}(\tau)$,
$\mathcal{G}(\tau)$, and $\chi_{0}(\widetilde{\Lambda})$, respectively,
are described in steps~\pmb{(iii)}, \pmb{(vi)}, and~\pmb{(vii)} below.

\pmb{(i)} Let $\widetilde{\Psi}(\widetilde{\mu})$ solve
Equation~\eqref{eq3.3}. Then applying the transformation $F_{1}$, using
Equations~\eqref{eq3.9}, \eqref{eq3.10}, and~\eqref{eq3.12}, and
conditions~\eqref{eq3.11}, one shows that
\begin{equation}
\dfrac{\partial \Phi (\widetilde{\Lambda})}{\partial \widetilde{\Lambda}}
\underset{\tau \to +\infty}{=} \left((\tau^{1/3} \mathcal{Q}_{+} \! + \!
\tau^{-1/3} \mathcal{Q}_{-}) \! + \! \widetilde{\Lambda} \mathcal{Q}_{1} \!
+ \! \tau^{-1/3} \widetilde{\Lambda}^{2} \mathcal{Q}_{2} \! + \! \mathcal{O}
(\tau^{-2/3} \widetilde{\Lambda}^{3} \mathcal{Q}_{3}) \right) \! \Phi
(\widetilde{\Lambda}), \label{eq3.64}
\end{equation}
where
\begin{equation*}
\mathcal{Q}_{\pm} \! = \! (b(\tau))^{-\frac{1}{2} \ad (\sigma_{3})}
\hat{\mathcal{P}}_{\pm}, \qquad \qquad \mathcal{Q}_{j} \! = \! (b(\tau))^{-
\frac{1}{2} \ad (\sigma_{3})} \hat{\mathcal{P}}_{j}, \quad j \! = \! 1,2,3,
\end{equation*}
with
\begin{gather*}
\hat{\mathcal{P}}_{+} \! := \!
\begin{pmatrix}
-\mi \alpha \hat{r}_{0}(\tau) & -8 \mi \alpha^{4} \hat{u}_{0}(\tau) \\
\frac{\mi (\hat{r}_{0}(\tau)+2 \hat{u}_{0}(\tau))}{2 \alpha^{2}(1+\hat{u}_{0}
(\tau))} & \mi \alpha \hat{r}_{0}(\tau)
\end{pmatrix}, \qquad \hat{\mathcal{P}}_{-} \! := \!
\begin{pmatrix}
0 & 0 \\
-\frac{\mi (a-\frac{\mi}{2})}{2 \alpha^{4}(1+\hat{u}_{0}(\tau))} & 0
\end{pmatrix}, \\
\hat{\mathcal{P}}_{1} \! := \!
\begin{pmatrix}
\mi (-4 \! + \! \hat{r}_{0}(\tau)) & -16 \mi \alpha^{3} \\
-\frac{2 \mi}{\alpha^{3}} & -\mi (-4 \! + \! \hat{r}_{0}(\tau))
\end{pmatrix}, \qquad \hat{\mathcal{P}}_{2} \! := \!
\begin{pmatrix}
-\frac{\mi (-2+\hat{r}_{0}(\tau))}{\alpha} & 24 \mi \alpha^{2} \\
\frac{3 \mi}{\alpha^{4}} & \frac{\mi (-2+\hat{r}_{0}(\tau))}{\alpha}
\end{pmatrix}, \\
\hat{\mathcal{P}}_{3} \! := \!
\begin{pmatrix}
\mathcal{O}(1) \! + \! \mathcal{O}(\hat{r}_{0}(\tau)) & \mathcal{O}(1) \\
\mathcal{O}(1) & \mathcal{O}(1) \! + \! \mathcal{O}(\hat{r}_{0}(\tau))
\end{pmatrix}.
\end{gather*}
Note that $\tr (\mathcal{Q}_{\pm}) \! = \! \tr (\hat{\mathcal{P}}_{\pm}) \! =
\! \tr (\mathcal{Q}_{j}) \! = \! \tr (\hat{\mathcal{P}}_{j}) \! = \! 0$, $j
\! = \! 1,2,3$; furthermore, Equation~\eqref{eq3.65} implies
\begin{equation}
\det (\hat{\mathcal{P}}_{+}) \! = \! -8 \alpha^{2} \! \left(\dfrac{(\hat{u}_{0}
(\tau))^{2}}{1 \! + \! \hat{u}_{0}(\tau)} \! + \! \dfrac{\hat{r}_{0}(\tau)
\hat{u}_{0}(\tau)}{2(1 \! + \! \hat{u}_{0}(\tau))} \! - \! \dfrac{(\hat{r}_{0}
(\tau))^{2}}{8} \right) \! = \! \varkappa_{0}^{2}(\tau) \tau^{-2/3},
\label{eq3.66}
\end{equation}
where $\varkappa_{0}^{2}(\tau)$ is defined by Equation~\eqref{eq3.58}.

\pmb{(ii)} Let $\Phi (\widetilde{\Lambda})$ solve Equation~\eqref{eq3.64}.
Then applying the transformation $F_{2}$, one obtains
\begin{equation}
\dfrac{\partial \widehat{\Phi}(\widetilde{\Lambda})}{\partial \widetilde{
\Lambda}} \underset{\tau \to +\infty}{=} \left((\tau^{1/3}
\hat{\mathcal{P}}_{+} \! + \! \tau^{-1/3} \hat{\mathcal{P}}_{-}) \! + \!
\widetilde{\Lambda} \hat{\mathcal{P}}_{1} \! + \! \tau^{-1/3} \widetilde{
\Lambda}^{2} \hat{\mathcal{P}}_{2} \! + \! \mathcal{O}(\tau^{-2/3}
\widetilde{\Lambda}^{3} \hat{\mathcal{P}}_{3}) \right) \widehat{\Phi}
(\widetilde{\Lambda}). \label{eq3.67}
\end{equation}

\pmb{(iii)} The idea behind the following transformation for
Equation~\eqref{eq3.67} is to put the matrix $\hat{\mathcal{P}}_{+}$ into
Jordan canonical form, that is, to find a function $\mathcal{N}(\tau)$ such
that
\begin{equation*}
(\mathcal{N}(\tau))^{-1} \hat{\mathcal{P}}_{+} \mathcal{N}(\tau) \! = \!
\mi \varkappa_{0}(\tau) \tau^{-1/3} \sigma_{3} \! + \! \sigma_{+}.
\end{equation*}
The following solution for $\mathcal{N}(\tau)$ is chosen:
\begin{equation*}
\mathcal{N}(\tau) \! = \!
\begin{pmatrix}
\sqrt{-8 \mi \alpha^{4} \hat{u}_{0}(\tau)} & 0 \\
\frac{\mi \left(\alpha \hat{r}_{0}(\tau)+ \varkappa_{0}(\tau) \tau^{-1/3}
\right)}{\sqrt{-8 \mi \alpha^{4} \hat{u}_{0}(\tau)}} & \frac{1}{\sqrt{-8
\mi \alpha^{4} \hat{u}_{0}(\tau)}}
\end{pmatrix}.
\end{equation*}
Let $\widehat{\Phi}(\widetilde{\Lambda})$ solve Equation~\eqref{eq3.67}.
Then applying the transformation $F_{3}$, using Equations~\eqref{eq3.9}
and~\eqref{eq3.10}, conditions~\eqref{eq3.11}, and Equations~\eqref{eq3.65},
\eqref{eq3.58}, and~\eqref{eq3.66}, one shows that
\begin{equation}
\dfrac{\partial \phi (\widetilde{\Lambda})}{\partial \widetilde{\Lambda}}
\underset{\tau \to +\infty}{=} \left((\tau^{1/3} \mathcal{P}_{+}^{\sharp}
\! + \! \tau^{-1/3} \mathcal{P}_{-}^{\sharp}) \! + \! \widetilde{\Lambda}
\mathcal{P}_{1}^{\sharp} \! + \! \tau^{-1/3} \widetilde{\Lambda}^{2}
\mathcal{P}_{2}^{\sharp} \! + \! \mathcal{O}(\tau^{-2/3} \widetilde{
\Lambda}^{3}(\mathcal{N}(\tau))^{-1} \hat{\mathcal{P}}_{3} \mathcal{N}
(\tau)) \right) \! \phi (\widetilde{\Lambda}), \label{eq3.69}
\end{equation}
where $\mathcal{P}_{\pm}^{\sharp} \! := \! (\mathcal{N}(\tau))^{-1}
\mathcal{P}_{\pm} \mathcal{N}(\tau)$, $\mathcal{P}_{j}^{\sharp} \! := \!
(\mathcal{N}(\tau))^{-1} \hat{\mathcal{P}}_{j} \mathcal{N}(\tau)$, $j \!
= \! 1,2$, with
\begin{gather*}
\mathcal{P}^{\sharp}_{+} \! = \! \mi \varkappa_{0}(\tau) \tau^{-1/3}
\sigma_{3} \! + \! \sigma_{+}, \quad \quad \mathcal{P}_{-}^{\sharp} \! = \!
-\dfrac{4(a \! - \! \frac{\mi}{2}) \hat{u}_{0}(\tau)}{(1 \! + \! \hat{u}_{0}
(\tau))} \sigma_{-}, \\
(\mathcal{P}^{\sharp}_{1})_{11} \! = \! -(\mathcal{P}^{\sharp}_{1})_{22}
\underset{\tau \to +\infty}{=} \mathfrak{p}_{\tau} \! + \! \dfrac{2 \mi
\varkappa_{0}(\tau) \tau^{-1/3}}{\alpha \hat{u}_{0}(\tau)} \! + \!
\mathcal{O} \! \left(\dfrac{\varkappa_{0}^{2}(\tau) \tau^{-2/3}}{\hat{r}_{0}
(\tau) \hat{u}_{0}(\tau)} \right) \! + \! \mathcal{O} \! \left(
\dfrac{\varkappa_{0}^{2}(\tau) \tau^{-2/3}}{\hat{r}_{0}(\tau)} \right), \\
(\mathcal{P}^{\sharp}_{1})_{12} \! = \! \dfrac{2}{\alpha \hat{u}_{0}(\tau)},
\qquad (\mathcal{P}^{\sharp}_{1})_{21} \underset{\tau \to +\infty}{=} -2 \mi
\varkappa_{0}(\tau) \mathfrak{p}_{\tau} \tau^{-1/3} \! + \! \mathcal{O} \!
\left(\dfrac{\varkappa_{0}^{2}(\tau) \tau^{-2/3}}{\hat{u}_{0}(\tau)} \right)
\! + \! \mathcal{O}(\varkappa_{0}^{2}(\tau) \tau^{-2/3}), \\
(\mathcal{P}^{\sharp}_{2})_{11} \! = \! -(\mathcal{P}^{\sharp}_{2})_{22} \! =
\! -\dfrac{\mi (6 \! + \! (-2 \! + \! \hat{r}_{0}(\tau))(3 \! + \! \hat{u}_{0}
(\tau)))}{\alpha \hat{u}_{0}(\tau)} \! - \! \dfrac{3 \mi \varkappa_{0}(\tau)
\tau^{-1/3}}{\alpha^{2} \hat{u}_{0}(\tau)}, \qquad
(\mathcal{P}^{\sharp}_{2})_{12} \! = \! -\dfrac{3}{\alpha^{2} \hat{u}_{0}
(\tau)},
\end{gather*}
\begin{align*}
(\mathcal{P}^{\sharp}_{2})_{21} \underset{\tau \to +\infty}{=}& \, -\hat{r}_{0}
(\tau)(8 \! - \! \hat{r}_{0}(\tau)) \! - \! \dfrac{2 \varkappa_{0}(\tau)}{
\alpha \hat{u}_{0}(\tau)} \left(6 \! + \! (-2 \! + \! \hat{r}_{0}(\tau))(3
\! + \! \hat{u}_{0}(\tau)) \right) \tau^{-1/3} \\
+& \, \mathcal{O} \! \left(\dfrac{\varkappa_{0}^{2}(\tau) \tau^{-2/3}}{
\hat{u}_{0}(\tau)} \right) \! + \! \mathcal{O}(\varkappa_{0}^{2}(\tau)
\tau^{-2/3}),
\end{align*}
where $\mathfrak{p}_{\tau}$ is defined by Equation~\eqref{eqpeetee}.

\pmb{(iv)} Let $\phi (\widetilde{\Lambda})$ solve Equation~\eqref{eq3.69}.
Then applying the transformation $F_{4}$, one proves that
\begin{align}
\dfrac{\partial \widetilde{\Phi}(\widetilde{\Lambda})}{\partial
\widetilde{\Lambda}} \underset{\tau \to +\infty}{=}& \, \left(
\begin{pmatrix}
\mi \varkappa_{0}(\tau) & 1 \\
-\frac{4(a-\frac{\mi}{2}) \hat{u}_{0}(\tau)}{1+ \hat{u}_{0}(\tau)} &
-\mi \varkappa_{0}(\tau)
\end{pmatrix} \! + \! \widetilde{\Lambda}
\begin{pmatrix}
\mathfrak{p}_{\tau} & 0 \\
-2 \mi \varkappa_{0}(\tau) \mathfrak{p}_{\tau} & -\mathfrak{p}_{\tau}
\end{pmatrix} \right. \nonumber \\
+&\left. \, \widetilde{\Lambda}^{2}
\begin{pmatrix}
0 & 0 \\
-\hat{r}_{0}(\tau)(8 \! - \! \hat{r}_{0}(\tau)) & 0
\end{pmatrix} \! + \! \mathcal{O}(\tau^{-1/3}E_{\tau}(\widetilde{\Lambda}))
\right) \! \widetilde{\Phi}(\widetilde{\Lambda}), \label{eq3.70}
\end{align}
where, with the help of conditions~\eqref{eq3.11},
\begin{align}
E_{\tau}(\widetilde{\Lambda}) &\underset{\tau \to +\infty}{=} \,
\widetilde{\Lambda}
\begin{pmatrix}
\mathcal{O} \! \left(\frac{\varkappa_{0}(\tau)}{\hat{u}_{0}(\tau)} \right) +
\mathcal{O} \! \left(\frac{\varkappa_{0}^{2}(\tau) \tau^{-1/3}}{\hat{r}_{0}
(\tau) \hat{u}_{0}(\tau)} \right) & \mathcal{O}((\hat{u}_{0}(\tau))^{-1}) \\
\mathcal{O} \! \left(\frac{\varkappa_{0}^{2}(\tau)}{\hat{u}_{0}(\tau)}
\right) + \mathcal{O}(\varkappa_{0}^{2}(\tau)) & \mathcal{O} \! \left(
\frac{\varkappa_{0}(\tau)}{\hat{u}_{0}(\tau)} \right) + \mathcal{O} \!
\left(\frac{\varkappa_{0}^{2}(\tau) \tau^{-1/3}}{\hat{r}_{0}(\tau)
\hat{u}_{0}(\tau)} \right)
\end{pmatrix} \nonumber \\
&+ \, \widetilde{\Lambda}^{2}
\begin{pmatrix}
\mathcal{O}(1) + \mathcal{O}(\hat{r}_{0}(\tau)) & \mathcal{O}(\tau^{-1/3}
(\hat{u}_{0}(\tau))^{-1}) \\
\mathcal{O}(\varkappa_{0}(\tau)) + \mathcal{O}(\varkappa_{0}(\tau) \hat{r}_{0}
(\tau)) & \mathcal{O}(1) + \mathcal{O}(\hat{r}_{0}(\tau))
\end{pmatrix} \nonumber \\
&+ \, \widetilde{\Lambda}^{3}
\begin{pmatrix}
\mathcal{O}(\tau^{-1/3}) + \mathcal{O}(\tau^{-1/3} \hat{r}_{0}(\tau)) &
\mathcal{O}(\tau^{-2/3}(\hat{u}_{0}(\tau))^{-1}) \\
\mathcal{O}(\hat{r}_{0}(\tau)) + \mathcal{O}((\hat{r}_{0}(\tau))^{2}) &
\mathcal{O}(\tau^{-1/3}) + \mathcal{O}(\tau^{-1/3} \hat{r}_{0}(\tau))
\end{pmatrix} \nonumber \\
\underset{\tau \to +\infty}{=}& \,
\begin{pmatrix}
\mathcal{O} \! \left(\frac{\widetilde{\Lambda} \varkappa_{0}^{2}(\tau)
\tau^{-1/3}}{\hat{r}_{0}(\tau) \hat{u}_{0}(\tau)} \right) + \mathcal{O}
(\widetilde{\Lambda}^{2}) + \mathcal{O}(\widetilde{\Lambda}^{2} \hat{r}_{0}
(\tau)) & \mathcal{O}(\widetilde{\Lambda}(\hat{u}_{0}(\tau))^{-1}) \\
\mathcal{O} \! \left(\frac{\widetilde{\Lambda} \varkappa_{0}^{2}(\tau)}{
\hat{u}_{0}(\tau)} \right) + \mathcal{O}(\widetilde{\Lambda}^{2} \varkappa_{0}
(\tau)) + \mathcal{O}(\widetilde{\Lambda}^{3}(\hat{r}_{0}(\tau))^{2}) &
\mathcal{O} \! \left(\frac{\widetilde{\Lambda} \varkappa_{0}^{2}(\tau)
\tau^{-1/3}}{\hat{r}_{0}(\tau) \hat{u}_{0}(\tau)} \right) + \mathcal{O}
(\widetilde{\Lambda}^{2}) + \mathcal{O}(\widetilde{\Lambda}^{2} \hat{r}_{0}
(\tau))
\end{pmatrix}. \label{eqEsubtau}
\end{align}

\pmb{(v)} One now proceeds to eliminate the $\widetilde{\Lambda}^{2}$ term
in Equation~\eqref{eq3.70}. Applying the transformation $F_{5}$, with a
$\tau$-dependent parameter $\ell_{\tau}$, one obtains
\begin{align}
\dfrac{\partial \widehat{\phi}(\widetilde{\Lambda})}{\partial
\widetilde{\Lambda}} \underset{\tau \to +\infty}{=}& \, \left(
\begin{pmatrix}
\mi \varkappa_{0}(\tau) & 1 \\
-\ell_{\tau} \! - \! \frac{4(a-\frac{\mi}{2}) \hat{u}_{0}(\tau)}{1+
\hat{u}_{0}(\tau)} & -\mi \varkappa_{0}(\tau)
\end{pmatrix} \! + \! \widetilde{\Lambda}
\begin{pmatrix}
\mathfrak{p}_{\tau} \! + \! \ell_{\tau} & 0 \\
-2 \mi \varkappa_{0}(\tau)(\ell_{\tau} \! + \! \mathfrak{p}_{\tau}) &
-(\mathfrak{p}_{\tau} \! + \! \ell_{\tau})
\end{pmatrix} \right. \nonumber \\
+&\left. \, \widetilde{\Lambda}^{2}
\begin{pmatrix}
0 & 0 \\
-\hat{r}_{0}(\tau)(8 \! - \! \hat{r}_{0}(\tau)) \! - \! 2 \mathfrak{p}_{\tau}
\ell_{\tau} \! - \! \ell^{2}_{\tau} & 0
\end{pmatrix} \! + \! \mathcal{O}(\tau^{-1/3} \widehat{E}_{\tau}
(\widetilde{\Lambda})) \right) \! \widehat{\phi}(\widetilde{\Lambda}),
\label{eq3.71}
\end{align}
where
\begin{equation}
\widehat{E}_{\tau}(\widetilde{\Lambda}) \! := \!
\begin{pmatrix}
1 & 0 \\
-\ell_{\tau} \widetilde{\Lambda} & 1
\end{pmatrix} E_{\tau}(\widetilde{\Lambda})
\begin{pmatrix}
1 & 0 \\
\ell_{\tau} \widetilde{\Lambda} & 1
\end{pmatrix}. \label{eqEhatsubtau}
\end{equation}
One now chooses $\ell_{\tau}$ so that the quadratic term in
Equation~\eqref{eq3.71} is annihilated:
\begin{equation*}
\ell^{2}_{\tau} \! + \! 2 \mathfrak{p}_{\tau} \ell_{\tau} \! + \! \hat{r}_{0}
(\tau)(8 \! - \! \hat{r}_{0}(\tau)) \! = \! 0;
\end{equation*}
thus,
\begin{equation}
\ell_{\tau} \! = \! -\mathfrak{p}_{\tau} \! + \! \sqrt{\mathfrak{p}_{\tau}^{2}
\! - \! \hat{r}_{0}(\tau)(8 \! - \! \hat{r}_{0}(\tau))}. \label{eq3.72}
\end{equation}
A straightforward calculation shows that (cf. Equation~\eqref{eq3.49})
\begin{equation*}
\mathfrak{p}_{\tau}^{2} \! - \! \hat{r}_{0}(\tau)(8 \! - \! \hat{r}_{0}(\tau))
\underset{\tau \to +\infty}{=} -48 \! + \! \mathcal{O} \! \left(\dfrac{\tau^{
-2/3} \varkappa_{0}^{2}(\tau)(\mathcal{O}(1) \! + \! (\hat{r}_{0}(\tau))^{2}
)}{(\hat{r}_{0}(\tau))^{2}} \right);
\end{equation*}
hence, making the choice $\sqrt{-1} \! = \! \mi$,
\begin{equation}
\ell_{\tau} \underset{\tau \to +\infty}{=} \ell_{\tau}^{\infty} \! + \!
\mathcal{O} \! \left(\dfrac{\tau^{-2/3} \varkappa_{0}^{2}(\tau)(\mathcal{O}
(1) \! + \! (\hat{r}_{0}(\tau))^{2})}{(\hat{r}_{0}(\tau))^{2}} \right),
\label{eq3.73}
\end{equation}
where $\ell_{\tau}^{\infty}$ is defined by Equation~\eqref{eqlteeinf}. Now,
with the help of Equations~\eqref{eqEsubtau}, \eqref{eqEhatsubtau},
and~\eqref{eq3.73}, and conditions~\eqref{eq3.11}, one re-writes
Equation~\eqref{eq3.71} as
\begin{equation}
\dfrac{\partial \widehat{\phi}(\widetilde{\Lambda})}{\partial
\widehat{\Lambda}} \underset{\tau \to +\infty}{=} \left(
\begin{pmatrix}
\mi \varkappa_{0}(\tau) & 1 \\
-\ell_{\tau}^{\infty} \! - \! \frac{4(a-\frac{\mi}{2}) \hat{u}_{0}(\tau)}{1+
\hat{u}_{0}(\tau)} & -\mi \varkappa_{0}(\tau)
\end{pmatrix} \! + \! \widetilde{\Lambda}
\begin{pmatrix}
\mi 4 \sqrt{3} & 0 \\
8 \sqrt{3} \, \varkappa_{0}(\tau) & -\mi 4 \sqrt{3}
\end{pmatrix} \! + \! \mathcal{O}(\widetilde{E}_{\tau}(\widetilde{\Lambda}))
\right) \! \widehat{\phi}(\widetilde{\Lambda}), \label{eq3.74}
\end{equation}
where
\begin{equation}
\widetilde{E}_{\tau}(\widetilde{\Lambda}) \! := \!
\begin{pmatrix}
(\widetilde{E}_{\tau}(\widetilde{\Lambda}))_{11} & (\widetilde{E}_{\tau}
(\widetilde{\Lambda}))_{12} \\
(\widetilde{E}_{\tau}(\widetilde{\Lambda}))_{21} & (\widetilde{E}_{\tau}
(\widetilde{\Lambda}))_{22}
\end{pmatrix}, \label{eqEwavetau}
\end{equation}
with
\begin{align*}
(\widetilde{E}_{\tau}(\widetilde{\Lambda}))_{11} \! = \! -(\widetilde{E}_{\tau}
(\widetilde{\Lambda}))_{22} \underset{\tau \to +\infty}{=}& \, \mathcal{O}
\! \left(\dfrac{\widetilde{\Lambda} \tau^{-2/3} \varkappa_{0}^{2}(\tau)}{
(\hat{r}_{0}(\tau))^{2}} \right) \! + \! \mathcal{O} \! \left(\widetilde{
\Lambda} \tau^{-2/3} \varkappa_{0}^{2}(\tau) \right) \\
+& \, \mathcal{O} \! \left(\widetilde{\Lambda}^{2} \tau^{-1/3} \hat{r}_{0}
(\tau) \right) \! + \! \mathcal{O} \! \left(\dfrac{\widetilde{\Lambda}^{2}
\tau^{-1/3} \ell_{\tau}^{\infty}}{\hat{u}_{0}(\tau)} \right), \\
(\widetilde{E}_{\tau}(\widetilde{\Lambda}))_{12} \underset{\tau \to
+\infty}{=}& \, \mathcal{O} \! \left(\dfrac{\widetilde{\Lambda} \tau^{-1/3}}{
\hat{u}_{0}(\tau)} \right),
\end{align*}
\begin{align*}
(\widetilde{E}_{\tau}(\widetilde{\Lambda}))_{21} \underset{\tau \to
+\infty}{=}& \, \mathcal{O} \! \left(\dfrac{\tau^{-2/3} \varkappa_{0}^{2}
(\tau)}{(\hat{r}_{0}(\tau))^{2}} \right) \! + \! \mathcal{O} \! \left(
\tau^{-2/3} \varkappa_{0}^{2}(\tau) \right) \! + \! \mathcal{O} \! \left(
\dfrac{\widetilde{\Lambda} \tau^{-2/3} \varkappa_{0}^{3}(\tau)}{(\hat{r}_{0}
(\tau))^{2}} \right) \! + \! \mathcal{O} \! \left(\widetilde{\Lambda}
\tau^{-2/3} \varkappa_{0}^{3}(\tau) \right) \\
+& \, \mathcal{O} \! \left(\dfrac{\widetilde{\Lambda} \tau^{-1/3}
\varkappa_{0}^{2}(\tau)}{\hat{u}_{0}(\tau)} \right) \! + \! \mathcal{O} \!
\left(\widetilde{\Lambda}^{2} \tau^{-1/3} \varkappa_{0}(\tau) \right) \!
+ \! \mathcal{O} \! \left(\dfrac{\widetilde{\Lambda}^{2} \tau^{-2/3}
\ell_{\tau}^{\infty} \varkappa_{0}^{2}(\tau)}{\hat{r}_{0}(\tau) \hat{u}_{0}
(\tau)} \right) \\
+& \, \mathcal{O} \! \left(\widetilde{\Lambda}^{3} \tau^{-1/3}(\hat{r}_{0}
(\tau))^{2} \right) \! + \! \mathcal{O} \! \left(\widetilde{\Lambda}^{3}
\tau^{-1/3} \ell_{\tau}^{\infty} \right) \! + \! \mathcal{O} \! \left(
\widetilde{\Lambda}^{3} \tau^{-1/3} \ell_{\tau}^{\infty} \hat{r}_{0}(\tau)
\right) \\
+& \, \mathcal{O} \! \left(\dfrac{\widetilde{\Lambda}^{3} \tau^{-1/3}
(\ell_{\tau}^{\infty})^{2}}{\hat{u}_{0}(\tau)} \right).
\end{align*}

\pmb{(vi)} One now proceeds to diagonalize the $\widetilde{\Lambda}$ term in
Equation~\eqref{eq3.74}. Set
\begin{equation}
\mathcal{G}(\tau) \! := \! \me^{-\frac{\mi \pi}{4}} \sqrt{\varkappa_{0}(\tau)}
\begin{pmatrix}
\mi (\varkappa_{0}(\tau))^{-1} & 0 \\
1 & 1
\end{pmatrix}. \label{eq3.75}
\end{equation}
Then applying the transformation $F_{6}$, one shows that
\begin{equation}
\dfrac{\partial \psi (\widetilde{\Lambda})}{\partial \widetilde{\Lambda}}
\! = \! \left(B_{0}(\widetilde{\Lambda}) \! + \! R_{0}(\widetilde{\Lambda})
\right) \! \psi (\widetilde{\Lambda}), \label{eq3.76}
\end{equation}
where
\begin{equation}
B_{0}(\widetilde{\Lambda}) \! := \! \mi 4 \sqrt{3} \, \widetilde{\Lambda}
\sigma_{3} \! + \! \widetilde{q}(\tau) \sigma_{-} \! + \! \widetilde{p}(\tau)
\sigma_{+}, \label{eqBee}
\end{equation}
with $\widetilde{p}(\tau)$, $\widetilde{q}(\tau)$  given in
Equations~\eqref{eq3.61}, and
\begin{equation}
R_{0}(\widetilde{\Lambda}) \underset{\tau \to +\infty}{:=} \mathcal{O} \!
\left((\mathcal{G}(\tau))^{-1} \widetilde{E}_{\tau}(\widetilde{\Lambda})
\mathcal{G}(\tau) \right), \label{eqRee}
\end{equation}
where $\widetilde{E}_{\tau}(\widetilde{\Lambda})$ is defined by
Equation~\eqref{eqEwavetau}.

\pmb{(vii)} Let $\psi_{0}(\widetilde{\Lambda})$ be a fundamental solution of
$\tfrac{\partial \psi_{0}(\widetilde{\Lambda})}{\partial \widetilde{\Lambda}}
\! = \! B_{0}(\widetilde{\Lambda}) \psi_{0}(\widetilde{\Lambda})$, which
coincides with Equation~\eqref{eq3.60}. Changing variables in
Equation~\eqref{eq3.60} according to $\mathcal{D}(x) \! := \! \psi_{0}
(\widetilde{\Lambda})$, where $\widetilde{\Lambda} \! = \! x_{0}x$, with
$x_{0} \! = \! \me^{-\frac{\mi \pi}{4}}2^{-3/2}3^{-1/4}$, one proves that
$\mathcal{D}(x)$ solves the standard equation
\begin{equation*}
\partial_{x} \mathcal{D}(x) \! = \! \left(\dfrac{x}{2} \sigma_{3} \! +
\! q^{\ast}(\tau) \sigma_{-} \! + \! p^{\ast}(\tau) \sigma_{+} \right)
\mathcal{D}(x),
\end{equation*}
where $p^{\ast}(\tau) \! := \! x_{0} \widetilde{p}(\tau)$ and $q^{\ast}(\tau)
\! := \! x_{0} \widetilde{q}(\tau)$, with fundamental solution given in terms
of the parabolic-cylinder function, $D_{\pmb{\ast}}(\pmb{\cdot})$ (see, for
example, \cite{a2,a5,a18}):
\begin{equation}
\mathcal{D}(x) \! = \!
\begin{pmatrix}
D_{-1-\nu}(\mi x) & D_{\nu}(x) \\
\dot{D}_{-1-\nu}(\mi x) & \dot{D}_{\nu}(x)
\end{pmatrix}, \label{eqDeecal}
\end{equation}
where $\dot{D}_{\pmb{\ast}}(z) \! := \! (p^{\ast}(\tau))^{-1}(\partial_{z}
D_{\pmb{\ast}}(z) \! - \! \tfrac{z}{2}D_{\pmb{\ast}}(z))$, and
\begin{equation}
\nu \! + \! 1 \! := \! -p^{\ast}(\tau)q^{\ast}(\tau) \! = \! -\dfrac{\mi}{8
\sqrt{3}} \! \left(\varkappa_{0}^{2}(\tau) \! + \! \ell_{\tau}^{\infty} \! +
\! \dfrac{4(a \! - \! \frac{\mi}{2}) \hat{u}_{0}(\tau)}{1 \! + \! \hat{u}_{0}
(\tau)} \right). \label{eq3.77}
\end{equation}
Now, applying definitions~\eqref{eq3.58} and~\eqref{eqlteeinf}, one re-writes
Equation~\eqref{eq3.77} as definition~\eqref{eq3.63}; thus, the representation
for $\psi_{0}(\widetilde{\Lambda})$ given in Equations~\eqref{eq3.63}
and~\eqref{eq3.60}--\eqref{eq3.62} is obtained.

Finally, in order to prove the error estimate in Equation~\eqref{eq3.56}, one
has to estimate the function $\chi_{0}(\widetilde{\Lambda})$ defined in the
transformation $F_{7}$. Applying the transformation $F_{7}$, one re-writes
Equation~\eqref{eq3.76} as follows:
\begin{equation}
\dfrac{\partial \chi_{0}(\widetilde{\Lambda})}{\partial \widetilde{\Lambda}}
\underset{\tau \to +\infty}{=} R_{0}(\widetilde{\Lambda}) \chi_{0}
(\widetilde{\Lambda}) \! + \! \left[B_{0}(\widetilde{\Lambda}),\chi_{0}
(\widetilde{\Lambda}) \right], \label{eq3.78}
\end{equation}
where $B_{0}(\widetilde{\Lambda})$ is defined by Equation~\eqref{eqBee},
$[B_{0}(\widetilde{\Lambda}),\chi_{0}(\widetilde{\Lambda})] \! := \! B_{0}
(\widetilde{\Lambda}) \chi_{0}(\widetilde{\Lambda}) \! - \! \chi_{0}
(\widetilde{\Lambda})B_{0}(\widetilde{\Lambda})$ is the matrix commutator,
and $R_{0}(\widetilde{\Lambda})$ is defined by Equation~\eqref{eqRee}. The
normalized solution $(\chi_{0}(0) \! = \! \mathrm{I})$ of
Equation~\eqref{eq3.78} is given by
\begin{equation}
\chi_{0}(\widetilde{\Lambda}) \! = \! \mathrm{I} \! + \! \int_{0}^{
\widetilde{\Lambda}} \psi_{0}(\widetilde{\Lambda})(\psi_{0}(\xi))^{-1}R_{0}
(\xi) \chi_{0}(\xi) \psi_{0}(\xi)(\psi_{0}(\widetilde{\Lambda}))^{-1} \,
\md \xi. \label{eq3.80}
\end{equation}
To prove the required estimate for $\chi_{0}(\widetilde{\Lambda})$, one uses
the method of successive approximations,
\begin{equation*}
\chi^{(n)}_{0}(\widetilde{\Lambda}) \! := \! \mathrm{I} \! + \! \int_{0}^{
\widetilde{\Lambda}} \psi_{0}(\widetilde{\Lambda})(\psi_{0}(\xi))^{-1}
R_{0}(\xi) \chi_{0}^{(n-1)}(\xi) \psi_{0}(\xi)(\psi_{0}
(\widetilde{\Lambda}))^{-1} \, \md \xi, \quad n \! \in \! \mathbb{N},
\end{equation*}
with $\chi^{(0)}_{0}(\widetilde{\Lambda}) \! = \! \mathrm{I}$, to construct
a Neumann series solution for $\chi_{0}(\widetilde{\Lambda})$, that is,
$\chi_{0}(\widetilde{\Lambda}) \! := \! \lim_{n \to \infty} \chi^{(n)}_{0}
(\widetilde{\Lambda})$. In this case, however, it suffices to estimate the
norm of the associated resolvent kernel. Via the above iteration argument,
it follows that
\begin{equation}
\lvert \lvert \chi_{0}(\widetilde{\Lambda}) \! - \! \mathrm{I} \rvert \rvert
\! \leqslant \! \exp \! \left(\int_{0}^{\widetilde{\Lambda}} \lvert \lvert
\psi_{0}(\widetilde{\Lambda}) \rvert \rvert \lvert \lvert (\psi_{0}(\xi))^{-1}
\rvert \rvert \lvert \lvert R_{0}(\xi) \rvert \rvert \lvert \lvert \psi_{0}
(\xi) \rvert \rvert \lvert \lvert (\psi_{0}(\widetilde{\Lambda}))^{-1} \rvert
\rvert \, \lvert \md \xi \rvert \right) \! - \! 1,
\label{eqchizeronorm}
\end{equation}
where $\lvert \md \xi \rvert$ denotes integration with respect to arc length.
One now estimates the norms appearing in Equation~\eqref{eqchizeronorm}.
Using Equations~\eqref{eqEwavetau}, \eqref{eq3.75}, and~\eqref{eqRee},
conditions~\eqref{eq3.11}, and the conditions of Corollary~\ref{cor3.1.3}, one
shows that
\begin{equation*}
R_{0}(\xi) \underset{\tau \to +\infty}{=}
\begin{pmatrix}
\mathcal{O}(\tau^{2 \varepsilon +\delta -\delta_{1}}) & \mathcal{O}
(\tau^{\varepsilon + \frac{\delta}{2} -\delta_{1}}) \\
\mathcal{O}(\tau^{3 \varepsilon + \frac{3 \delta}{2} -\delta_{1}}) &
\mathcal{O}(\tau^{2 \varepsilon +\delta -\delta_{1}})
\end{pmatrix};
\end{equation*}
recalling the definition of the matrix norm (cf. Subsection~\ref{notat}),
one arrives at
\begin{equation}
\lvert \lvert R_{0}(\xi) \rvert \rvert \underset{\tau \to +\infty}{=}
\mathcal{O} \! \left(\tau^{3 \varepsilon + \frac{3 \delta}{2} -\delta_{1}}
\right). \label{eqRzeronorm}
\end{equation}
For the function $\psi_{0}(\xi)$, one has to find a uniform approximation
for $\xi \! \in \! \mathbb{R} \cup \mi \mathbb{R}$. To do so, one uses the
following integral representation for the parabolic-cylinder function
\cite{EMOT}:
\begin{equation}
D_{\nu}(x) \! = \! \dfrac{2^{\frac{\nu}{2}} \me^{-\frac{x^{2}}{4}}}{\Gamma
(-\frac{\nu}{2})} \int_{0}^{+\infty} \me^{-\frac{tx^{2}}{2}}t^{-\frac{\nu}{2}
-1}(1 \! + \! t)^{\frac{\nu -1}{2}} \, \md t, \quad \Re (\nu) \! < \! 0, \quad
\lvert \arg (x) \rvert \! \leqslant \! \pi/4. \label{eqintegrepdsubnu}
\end{equation}
As this integral representation will be applied to the entries of the
matrix-valued function $\psi_{0}(\xi)$ (cf. Equation~\eqref{eq3.62}), it
implies the following restriction on $\nu$:
\begin{equation*}
0 \underset{\tau \to +\infty}{<} \Re (\nu \! + \! 1) \underset{\tau \to
+\infty}{<} 1.
\end{equation*}
Since, in the sector $\lvert \arg (x) \rvert \! \leqslant \! \pi/4$, the
exponents in the integral representation~\eqref{eqintegrepdsubnu} are less
than or equal to one, the following estimate for $D_{\nu}(x)$ can be deduced:
\begin{equation*}
\lvert D_{\nu}(x) \rvert \! \leqslant \! \dfrac{\sqrt{\pi} \,
2^{\Re (\nu)/2}}{\Gamma (\frac{1-\Re (\nu)}{2})}, \quad \lvert \arg (x)
\rvert \! \leqslant \! \pi/4.
\end{equation*}
{}From the above inequality, one derives, for the elements of the second
column of $\psi_{0}(\xi)$, the following estimates:
\begin{equation*}
\lvert (\psi_{0}(\xi))_{12} \rvert \underset{\tau \to +\infty}{=} \mathcal{O}
(1), \qquad \lvert (\psi_{0}(\xi))_{22} \rvert \underset{\tau \to +\infty}{=}
\mathcal{O}(1) \! + \! \mathcal{O}(\xi \lvert \varkappa_{0}(\tau)
\rvert^{-1}), \quad \arg (\xi) \! \in \! (-\pi/2,0).
\end{equation*}
In order to apply, simultaneously, the same integral
representation~\eqref{eqintegrepdsubnu} for the elements of the first column
of $\psi_{0}(\xi)$, one has to restrict $\arg (\xi)$ to $-\pi/2$; hence, one
arrives at the following, analogous estimates:
\begin{equation*}
\lvert (\psi_{0}(\xi))_{11} \rvert \underset{\tau \to +\infty}{=} \mathcal{O}
(1), \qquad \lvert (\psi_{0}(\xi))_{21} \rvert \underset{\tau \to +\infty}{=}
\mathcal{O}(1) \! + \! \mathcal{O}(\xi \lvert \varkappa_{0}(\tau)
\rvert^{-1}), \quad \arg (\xi) \! = \! -\pi/2.
\end{equation*}
Thus, for $\arg (\xi) \! = \! -\pi/2$, one arrives at the following estimate
for $\lvert \lvert \psi_{0}(\xi) \rvert \rvert$:
\begin{equation}
\lvert \lvert \psi_{0}(\xi) \rvert \rvert \underset{\tau \to +\infty}{=}
\mathcal{O}(1) \! + \! \mathcal{O}(\xi \lvert \varkappa_{0}(\tau)
\rvert^{-1}). \label{eqnormpseye}
\end{equation}
In order to find estimates for $\psi_{0}(\xi)$ on the other Stokes rays $\arg
(\xi) \! = \! 0,\pi/2,\pm \pi,\dotsc$, one has to use the linear relations for
the parabolic-cylinder functions relating any three of the four functions
$D_{\nu}(\pm x)$ and $D_{-\nu-1}(\pm \mi x)$ (see, for example, \cite{a24},
pg.~1094, Equations~\pmb{9.248} 1.--3.), and impose the additional restriction
$\Im (\nu \! + \! 1) \! =_{\tau \to +\infty} \! \mathcal{O}(1)$, in which
case, for all the Stokes rays $\arg (\xi) \! = \! 0,\pm \pi/2,\pm \pi,\dotsc$,
$0 \! \leqslant \! \lvert \xi \rvert \! < \! +\infty$, one verifies an
estimate similar to Equation~\eqref{eqnormpseye}. Noting that $\det (\psi_{0}
(\xi)) \! = \! -(p^{\ast}(\tau))^{-1} \exp (-\tfrac{\mi \pi}{2}(\nu \! + \!
1))$, one obtains
\begin{equation}
\lvert \lvert (\psi_{0}(\xi))^{-1} \rvert \rvert \underset{\tau \to
+\infty}{=} \lvert \varkappa_{0}(\tau) \rvert \me^{-\frac{\pi}{2} \Im (\nu
+1)} \lvert \lvert \psi_{0}(\xi) \rvert \rvert. \label{eqnormpseyeminus1}
\end{equation}
Using the asymptotic expansions for the parabolic-cylinder functions (see
Remark~\ref{rem3.2.2} below), one shows that
\begin{equation}
\lvert \lvert \psi_{0}(\widetilde{\Lambda}) \rvert \rvert \underset{\tau \to
+\infty}{=} \mathcal{O}(\lvert \varkappa_{0}(\tau) \rvert^{-1} \widetilde{
\Lambda}^{\Re (\nu +1)}), \qquad \lvert \lvert (\psi_{0}(\widetilde{\Lambda}
))^{-1} \rvert \rvert \underset{\tau \to +\infty}{=} \mathcal{O}(\widetilde{
\Lambda}^{\Re (\nu +1)}). \label{eqpseyezerolamdanorm}
\end{equation}
Combining the estimates~\eqref{eqRzeronorm}, \eqref{eqnormpseye},
\eqref{eqnormpseyeminus1}, and~\eqref{eqpseyezerolamdanorm}, and assuming
that\begin{equation*}
6 \varepsilon \! + \! 2 \varepsilon \Re (\nu \! + \! 1) \! + \!
\dfrac{\delta}{2} \! < \! \delta_{1},
\end{equation*}
one deduces {}from Equation~\eqref{eqchizeronorm} that
\begin{equation*}
\lvert \lvert \chi_{0}(\widetilde{\Lambda}) \! - \! \mathrm{I} \rvert
\rvert \underset{\tau \to +\infty}{\leqslant} \mathcal{O} \! \left(\tau^{6
\varepsilon +2 \varepsilon \Re (\nu +1)+ \frac{\delta}{2} -\delta_{1}} \right).
\end{equation*}
Hence, forming the composition of the invertible linear transformations
$F_{1},\dotsc,F_{7}$, that is (the ``symbol'' $\circ$ denotes
``composition''),
\begin{align*}
\widetilde{\Psi}(\widetilde{\mu}) =& \, \left(F_{1}^{-1} \circ F_{2}^{-1}
\circ F_{3}^{-1} \circ F_{4}^{-1} \circ F_{5}^{-1} \circ F_{6}^{-1} \circ
F_{7}^{-1} \right) \psi_{0}(\widetilde{\Lambda}) \\
=& \, \underbrace{(b(\tau))^{-\frac{1}{2} \sigma_{3}} \mathcal{N}(\tau)
\tau^{\frac{1}{6} \sigma_{3}}
\begin{pmatrix}
1 & 0 \\
\ell_{\tau} \widetilde{\Lambda} & 1
\end{pmatrix} \mathcal{G}(\tau)}_{=: \, \mathcal{F}_{\tau}
(\widetilde{\Lambda})} \chi_{0}(\widetilde{\Lambda}) \psi_{0}
(\widetilde{\Lambda}),
\end{align*}
one arrives at the asymptotic representation for $\widetilde{\Psi}
(\widetilde{\mu})$ given in Equation~\eqref{eq3.56}. \hfill $\qed$
\begin{eeee} \label{rem3.2.2}
\textsl{In Lemma~{\rm \ref{lem3.2.1}} and hereafter, the matrix-valued
function $\psi_{0}(\widetilde{\Lambda})$ $($cf. Equation \eqref{eq3.62}$)$
plays a pivotal role; therefore, for the reader's convenience, its
asymptotics are presented here:}
\begin{equation*}
\psi_{0}(\widetilde{\Lambda}) \underset{\underset{\scriptstyle
\arg (\widetilde{\Lambda})= \frac{k \pi}{2}}{\widetilde{\Lambda}
\to \infty}}{=} \left(\mathrm{I} \! + \! \sum_{j=1}^{\infty} \psi_{j}(\tau)
\widetilde{\Lambda}^{-j} \right) \! \exp \! \left(\left(\mi 2 \sqrt{3} \,
\widetilde{\Lambda}^{2} \! - \! (\nu \! + \! 1) \ln (\me^{\frac{\mi \pi}{4}}
2^{\frac{3}{2}}3^{\frac{1}{4}} \widetilde{\Lambda}) \right) \! \sigma_{3}
\right) \! \mathcal{R}_{k}, \, \, k \! = \! -1,0,1,2,
\end{equation*}
\textsl{where $\psi_{j}(\tau)$ are off-diagonal $($resp., diagonal$)$
matrices for $j$ odd $($resp., $j$ even$)$,}
\begin{gather*}
\mathcal{R}_{-1} \! := \!
\begin{pmatrix}
\me^{-\frac{\pi \mi}{2}(\nu +1)} & 0 \\
0 & -(p^{\ast}(\tau))^{-1}
\end{pmatrix}, \qquad \qquad \mathcal{R}_{0} \! := \!
\begin{pmatrix}
\me^{-\frac{\pi \mi}{2}(\nu +1)} & 0 \\
-\frac{\mi \sqrt{2 \pi}}{p^{\ast}(\tau) \Gamma (\nu +1)}
\me^{-\frac{\pi \mi}{2}(\nu +1)} & -(p^{\ast}(\tau))^{-1}
\end{pmatrix}, \\
\mathcal{R}_{1} \! := \!
\begin{pmatrix}
\me^{\frac{3 \pi \mi}{2}(\nu +1)} & \frac{\sqrt{2 \pi}}{\Gamma (-\nu)}
\me^{\pi \mi (\nu +1)} \\
-\frac{\mi \sqrt{2 \pi}}{p^{\ast}(\tau) \Gamma (\nu +1)}
\me^{-\frac{\pi \mi}{2}(\nu +1)} & -(p^{\ast}(\tau))^{-1}
\end{pmatrix}, \qquad \mathcal{R}_{2} \! := \!
\begin{pmatrix}
\me^{\frac{3 \pi \mi}{2}(\nu +1)} & \frac{\sqrt{2 \pi}}{\Gamma (-\nu)}
\me^{\pi \mi (\nu +1)} \\
0 &  -(p^{\ast}(\tau))^{-1} \me^{-2 \pi \mi (\nu +1)}
\end{pmatrix},
\end{gather*}
\textsl{and $\Gamma (\pmb{\cdot})$ is the (Euler) gamma function
{\rm \cite{a24}}. The above asymptotic expansion can be deduced {}from
the asymptotics of the parabolic-cylinder functions $($see, for example,
{\rm \cite{EMOT}}$)$. The diagonal/off-diagonal structure of the matrices
$\psi_{j}(\tau)$ is a consequence of the ``$\sigma_{3}$-reduction'' for
Equation~\eqref{eq3.60}$:$ $\psi_{0}(\widetilde{\Lambda}) \! \to \!
\sigma_{3} \psi_{0}(-\widetilde{\Lambda})$. In Lemmata~{\rm \ref{lem3.3.1}}
and~{\rm \ref{lem3.3.2}} $($see Subsection~{\rm \ref{sec3.3}} below$)$,
explicit knowledge of the following matrices is essential:}
\begin{gather*}
\psi_{1}(\tau) \! = \!
\begin{pmatrix}
0 & \frac{\varkappa_{0}(\tau)}{8 \sqrt{3}} \\
-\frac{\mi (\nu +1)}{\varkappa_{0}(\tau)} & 0
\end{pmatrix}, \quad \qquad \psi_{2}(\tau) \! = \!
\begin{pmatrix}
-\frac{\mi (\nu +1)(\nu +2)}{16 \sqrt{3}} & 0 \\
0 & \frac{\mi \nu (\nu +1)}{16 \sqrt{3}}
\end{pmatrix}, \\
\psi_{3}(\tau) \! = \!
\begin{pmatrix}
0 & \frac{\mi \nu (\nu -1) \varkappa_{0}(\tau)}{384} \\
-\frac{(\nu +1)(\nu +2)(\nu +3)}{16 \sqrt{3} \, \varkappa_{0}(\tau)} & 0
\end{pmatrix}.
\end{gather*}
\end{eeee}
\subsection{Matching of Asymptotics} \label{sec3.3}
In this subsection the connection matrix (cf. Equation~\eqref{eqdefg}), $G$,
is calculated asymptotically in terms of the matrix elements of the function
$\mathcal{A}(\widetilde{\mu},\tau)$ defined by Equation~\eqref{eq3.4}, that
is, the functions $\hat{r}_{0}(\tau)$, $\hat{u}_{0}(\tau)$, $h_{0}(\tau)$,
and $b(\tau)$; thus, under conditions~\eqref{eq3.11}, the direct monodromy
problem for Equation~\eqref{eq3.3} is solved asymptotically.
\begin{bbbb} \label{prop3.3.1}
Let $\widetilde{\mu} \! = \! \widetilde{\mu}_{0} \! = \! \alpha \! + \!
\tau^{-1/3} \widetilde{\Lambda}$, where $\widetilde{\Lambda} \! =_{\tau \to
+\infty} \! \mathcal{O}(\tau^{\varepsilon})$, $0 \! < \! \varepsilon \! < \!
1/9$. Then under conditions~\eqref{eq3.11}, with $0 \! < \! \delta \! < \!
\varepsilon \! < \! 1/9$ and $\varepsilon \! < \! \delta_{1} \! < \! 1/3$,
\begin{align}
(\mathcal{F}_{\tau}(\widetilde{\Lambda}))^{-1}T(\widetilde{\mu}) \underset{
\underset{\scriptstyle \tau \to +\infty}{\widetilde{\mu}=\widetilde{\mu}_{
0}}}{=}& \, \left(\dfrac{1}{3^{1/4}} \sqrt{\dfrac{\varkappa_{0}(\tau)}{\varpi
\widetilde{\Lambda}}} \, \right)^{\sigma_{3}}
\begin{pmatrix}
-\frac{1}{8}(1 \! + \! \mathcal{O}(\tau^{\delta -2 \varepsilon})) & 1 \!
+ \! \mathcal{O}(\tau^{\delta -2 \varepsilon}) \\
-\frac{1}{2}(\varpi \! + \! 1)(1 \! + \! \mathcal{O}(\tau^{2 \delta -2
\varepsilon})) & 4(\varpi \! - \! 1)(1 \! + \! \mathcal{O}(\tau^{2 \delta
-2 \varepsilon}))
\end{pmatrix} \nonumber \\
\times& \, \left(\dfrac{1}{\alpha^{3/2}} \sqrt{\dfrac{b(\tau) \hat{r}_{0}
(\tau)}{\hat{u}_{0}(\tau)}} \, \right)^{\sigma_{3}}, \label{eq3.81}
\end{align}
where $\varpi$ is defined in Proposition~{\rm \ref{prop3.1.6}}.
\end{bbbb}

\emph{Proof}. The proof is a straightforward, though algebraically
tedious, consequence of Proposition~\ref{prop3.1.6} (cf.
Equations~\eqref{eq3.53}--\eqref{eq3.55}) and Lemma~\ref{lem3.2.1}
(cf. Equations~\eqref{eq3.57}--\eqref{eqlteeinf}). More precisely, using
conditions~\eqref{eq3.11}, the conditions of Corollary~\ref{cor3.1.3},
Equations~\eqref{eqpeetee} and~\eqref{eq3.58}, and repeated application
of Equation~\eqref{eq3.65}, one shows that
\begin{align*}
\left((\mathcal{F}_{\tau}(\widetilde{\Lambda}))^{-1}T(\alpha \! + \!
\tau^{-1/3} \widetilde{\Lambda}) \right)_{11} \underset{\tau \to +\infty}{=}&
\, -\dfrac{1}{8 \pmb{\cdot} 3^{1/4} \alpha \sqrt{\alpha}} \sqrt{\dfrac{\hat{
r}_{0}(\tau)b(\tau)}{\hat{u}_{0}(\tau)}} \, \sqrt{\dfrac{\varkappa_{0}(\tau)}{
\varpi \widetilde{\Lambda}}} \left(1 \! + \! \mathcal{O}(h_{0}(\tau)
\widetilde{\Lambda}^{-2}) \right), \\
\left((\mathcal{F}_{\tau}(\widetilde{\Lambda}))^{-1}T(\alpha \! + \!
\tau^{-1/3} \widetilde{\Lambda}) \right)_{12} \underset{\tau \to +\infty}{=}&
\, \dfrac{\alpha \sqrt{\alpha}}{3^{1/4}} \sqrt{\dfrac{\hat{u}_{0}(\tau)}{
\hat{r}_{0}(\tau)b(\tau)}} \, \sqrt{\dfrac{\varkappa_{0}(\tau)}{\varpi
\widetilde{\Lambda}}} \left(1 \! + \! \mathcal{O}(h_{0}(\tau)
\widetilde{\Lambda}^{-2}) \right), \\
\left((\mathcal{F}_{\tau}(\widetilde{\Lambda}))^{-1}T(\alpha \! + \!
\tau^{-1/3} \widetilde{\Lambda}) \right)_{21} \underset{\tau \to +\infty}{=}&
\, -\dfrac{\sqrt{3} \, (\varpi \! + \! 1)}{2 \pmb{\cdot} 3^{1/4} \alpha
\sqrt{\alpha}} \sqrt{\dfrac{b(\tau) \hat{r}_{0}(\tau)}{\hat{u}_{0}(\tau)}}
\, \sqrt{\dfrac{\varpi \widetilde{\Lambda}}{\varkappa_{0}(\tau)}} \\
\times& \, \left(1 \! + \! \mathcal{O} \! \left(\dfrac{(4 \hat{u}_{0}(\tau)
\! + \! \hat{r}_{0}(\tau) \! + \! 2(\hat{u}_{0}(\tau))^{2})h_{0}(\tau)}{
\hat{r}_{0}(\tau)(1 \! + \! \hat{u}_{0}(\tau)) \widetilde{\Lambda}^{2}}
\right) \right), \\
\left((\mathcal{F}_{\tau}(\widetilde{\Lambda}))^{-1}T(\alpha \! + \!
\tau^{-1/3} \widetilde{\Lambda}) \right)_{22} \underset{\tau \to +\infty}{=}&
\, \dfrac{4 \sqrt{3} \, (\varpi \! - \! 1) \alpha \sqrt{\alpha}}{3^{1/4}}
\sqrt{\dfrac{\varpi \widetilde{\Lambda}}{\varkappa_{0}(\tau)}} \,
\sqrt{\dfrac{\hat{u}_{0}(\tau)}{b(\tau) \hat{r}_{0}(\tau)}} \\
\times& \, \left(1 \! + \! \mathcal{O} \! \left(\dfrac{(4 \hat{u}_{0}(\tau)
\! + \! \hat{r}_{0}(\tau) \! + \! 2(\hat{u}_{0}(\tau))^{2})h_{0}(\tau)}{
\hat{r}_{0}(\tau)(1 \! + \! \hat{u}_{0}(\tau)) \widetilde{\Lambda}^{2}}
\right) \right);
\end{align*}
hence, via conditions~\eqref{eq3.11} and the conditions of
Corollary~\ref{cor3.1.3}, one arrives at, after simplification,
Equation~\eqref{eq3.81}. \hfill $\qed$
\begin{cccc} \label{lem3.3.1}
Let $\widetilde{\Psi}(\widetilde{\mu},\tau)$ be the fundamental solution of
Equation~\eqref{eq3.3} with asymptotics given in Lemma~{\rm \ref{lem3.2.1}},
and $Y^{\infty}_{0}(\widetilde{\mu},\tau)$ be the canonical solution of
Equation~\eqref{eq3.1}. Define\footnote{Note that $\tau^{-\frac{1}{12}
\sigma_{3}}Y^{\infty}_{0}(\tau^{-\frac{1}{6}} \widetilde{\mu},\tau)$ is also
a fundamental solution of Equation~\eqref{eq3.3}; therefore, $L_{\infty}
(\tau)$ is independent of $\widetilde{\mu}$.}
\begin{equation*}
L_{\infty}(\tau) \! := \! \left(\widetilde{\Psi}(\widetilde{\mu},\tau)
\right)^{-1} \tau^{-\frac{1}{12} \sigma_{3}}Y^{\infty}_{0}
(\tau^{-\frac{1}{6}} \widetilde{\mu},\tau).
\end{equation*}
Assume that the parameters $\varepsilon$, $\delta$, $\delta_{1}$, and $\nu
\! + \! 1$ satisfy the conditions~\eqref{newconds1} and~\eqref{newcconds2}$;$
furthermore, let the function $b(\tau)$ satisfy the following conditions:
\begin{gather}
\left(\dfrac{\tau^{\frac{1}{3} \Im (a)}(\hat{u}_{0}(\tau))^{2}}{b(\tau)
(\hat{r}_{0}(\tau))^{2}} \right) \! \tau^{6 \varepsilon + \delta - \delta_{1}}
\underset{\tau \to +\infty}{=} \mathcal{O}(\tau^{-\hat{\delta}_{1}}), \,
\left(\dfrac{b(\tau)(\hat{r}_{0}(\tau))^{2}}{\tau^{\frac{1}{3} \Im (a)}
(\hat{u}_{0}(\tau))^{2}} \right) \! \tau^{6 \varepsilon +4 \varepsilon
\Re (\nu +1)- \delta_{1}} \underset{\tau \to +\infty}{=} \mathcal{O}
(\tau^{-\hat{\delta}_{2}}), \label{eq3.103} \\
\left(\dfrac{\tau^{\frac{1}{3} \Im (a)}(\hat{u}_{0}(\tau))^{2}}{b(\tau)
(\hat{r}_{0}(\tau))^{2}} \!\right) \! \tau^{-\varepsilon -2 \varepsilon
\Re (\nu +1)} \underset{\tau \to +\infty}{=} \mathcal{O}
(\tau^{-\hat{\delta}_{3}}), \, \left(\dfrac{b(\tau)(\hat{r}_{0}
(\tau))^2}{\tau^{\frac{1}{3} \Im (a)}(\hat{u}_{0}(\tau))^2}\! \right) \!
\tau^{2 \delta+2 \varepsilon \Re (\nu +1)-3 \varepsilon}\! \underset{\tau
\to +\infty}{=}\! \mathcal{O}(\tau^{-\hat{\delta}_{4}}), \label{eq3.104}
\end{gather}
where $\hat{\delta}_{k} \! > \! 0$, $k \! = \! 1,2,3,4$. Then
\begin{align}
L_{\infty}(\tau) \underset{\tau \to +\infty}{=}& \, \mi \,
\mathcal{R}_{0}^{-1} \exp \left(\left(\mathfrak{t}_{\infty} \! + \!
\mathfrak{c}_{\infty} \! - \! \mathcal{I}_{\infty}^{\sharp}(\tau) \right)
\sigma_{3} \right) \left(\sqrt{\dfrac{\varkappa_{0}(\tau) \hat{u}_{0}
(\tau)}{b(\tau) \hat{r}_{0}(\tau)}} \, \right)^{\sigma_{3}} \, \sigma_{2}
\nonumber \\
\times& \, \left(\mathrm{I} \! + \!
\begin{pmatrix}
\mathcal{O}(\tau^{-\overset{\infty}{\delta_{11}}}) &
\mathcal{O}(\tau^{-\overset{\infty}{\delta_{12}}}) \\
\mathcal{O}(\tau^{-\overset{\infty}{\delta_{21}}}) &
\mathcal{O}(\tau^{-\overset{\infty}{\delta_{22}}})
\end{pmatrix} \right), \label{eq3.107}
\end{align}
where $\mathcal{R}_{0}$ is defined in Remark~{\rm \ref{rem3.2.2}},
$\mathcal{I}_{\infty}^{\sharp}(\tau)$ is defined by Equation~\eqref{eq3.87},
\begin{gather}
\mathfrak{t}_{\infty} \! := \! \mi 3(\sqrt{3}- \! 1) \alpha^{2} \tau^{2/3} \!
+ \! \dfrac{1}{3} \! \left(\nu \! +  \! \dfrac{1}{2} \! - \! \dfrac{\mi a}{2}
\right) \ln \tau, \label{eq3.108} \\
\mathfrak{c}_{\infty} \! := \! \dfrac{3}{2} \ln \alpha \! - \! \dfrac{1}{4}
\ln 3 \! + \! (\nu \! + \! 1) \ln (\me^{\frac{\mi \pi}{4}}2^{3/2}3^{1/4})
\! + \! C^{\mathrm{WKB}}_{\infty}, \label{eq3.109}
\end{gather}
with $C^{\mathrm{WKB}}_{\infty}$ defined by Equation~\eqref{eq3.38}, and
$\overset{\infty}{\delta_{11}} \! = \! \overset{\infty}{\delta_{22}} \! :=
\! \min \lbrace \delta_{1} \! - \! \tfrac{\delta}{2} \! - \! 6 \varepsilon
\! - \! 2 \varepsilon \Re (\nu \! + \! 1),2(\varepsilon \! - \! \delta),
\tfrac{1}{3} \! - \! 3 \varepsilon \rbrace$, $\overset{\infty}{\delta_{12}}
\! := \! \min \lbrace \hat{\delta}_{1},\hat{\delta}_{3} \rbrace$, and
$\overset{\infty}{\delta_{21}} \! := \! \min \lbrace \hat{\delta}_{2},
\hat{\delta}_{4} \rbrace$.
\end{cccc}

\emph{Proof}. Denote by $\widetilde{\Psi}_{\mathrm{WKB}}(\widetilde{\mu},
\tau)$ the solution of Equation~\eqref{eq3.3} which has $\mathrm{WKB}$
asymptotics given by Equations~\eqref{eq3.16}--\eqref{eq3.18} in the canonical
domain containing the Stokes curve approaching the positive real $\widetilde{
\mu}$-axis {}from above as $\widetilde{\mu} \! \to \! +\infty$. Now, re-write
$L_{\infty}(\tau)$ in the following, equivalent form:
\begin{equation*}
L_{\infty}(\tau) \! = \! \left((\widetilde{\Psi}(\widetilde{\mu},\tau))^{-1}
\widetilde{\Psi}_{\mathrm{WKB}}(\widetilde{\mu},\tau) \right) \!
\left((\widetilde{\Psi}_{\mathrm{WKB}}(\widetilde{\mu},\tau))^{-1}
\tau^{-\frac{1}{12} \sigma_{3}}Y^{\infty}_{0}(\tau^{-\frac{1}{6}}
\widetilde{\mu},\tau) \right).
\end{equation*}
Noting that $\widetilde{\Psi}(\widetilde{\mu},\tau)$, $\widetilde{\Psi}_{
\mathrm{WKB}}(\widetilde{\mu},\tau)$, and $\tau^{-\frac{1}{12} \sigma_{3}}
Y^{\infty}_{0}(\tau^{-\frac{1}{6}} \widetilde{\mu},\tau)$ are all solutions
of Equation~\eqref{eq3.3}, it follows that they differ by right-hand,
$\widetilde{\mu}$-independent matrix factors. Using this fact, one evaluates
asymptotically (as $\tau \! \to \! +\infty)$ each pair by taking separate
limits, that is, $\widetilde{\mu} \! \to \! \alpha$ and $\widetilde{\mu} \!
\to \! +\infty$, respectively; more precisely,
\begin{gather*}
(\widetilde{\Psi}(\widetilde{\mu},\tau))^{-1} \widetilde{\Psi}_{\mathrm{WKB}}
(\widetilde{\mu},\tau) \underset{\tau \to +\infty}{=} \underbrace{(\mathcal{
F}_{\tau}(\widetilde{\Lambda})(\mathrm{I} \! + \! \mathcal{O}(\tau^{6
\varepsilon +2 \varepsilon \Re (\nu +1)+ \frac{\delta}{2}-\delta_{1}}))
\psi_{0}(\widetilde{\Lambda}))^{-1}T(\widetilde{\mu}) \me^{\mathrm{WKB}
(\widetilde{\mu},\tau)}}_{\widetilde{\mu}=\widetilde{\mu}_{0}= \alpha +
\tau^{-1/3} \widetilde{\Lambda}, \quad \widetilde{\Lambda} \sim \tau^{
\varepsilon}, \, \, 0< \varepsilon <1/9, \quad \arg (\widetilde{\Lambda})=0}
\, , \\
(\widetilde{\Psi}_{\mathrm{WKB}}(\widetilde{\mu},\tau))^{-1} \tau^{-
\frac{1}{12} \sigma_{3}}Y^{\infty}_{0}(\tau^{-\frac{1}{6}} \widetilde{\mu},
\tau) \underset{\tau \to +\infty}{=} \underbrace{(T(\widetilde{\mu})
\me^{\mathrm{WKB}(\widetilde{\mu},\tau)})^{-1} \tau^{-\frac{1}{12} \sigma_{3}}
Y^{\infty}_{0}(\tau^{-\frac{1}{6}} \widetilde{\mu},\tau)}_{\widetilde{\mu}
\to \infty, \quad \arg (\widetilde{\mu})=0} \, ,
\end{gather*}
where $\mathrm{WKB}(\widetilde{\mu},\tau) \! := \! -\sigma_{3} \mi \tau^{2/3}
\int_{\widetilde{\mu}_{0}}^{\widetilde{\mu}}l(\xi) \, \md \xi \! - \! \int_{
\widetilde{\mu}_{0}}^{\widetilde{\mu}} \diag ((T(\xi))^{-1} \partial_{\xi}
T(\xi)) \, \md \xi$.

To calculate asymptotics of $(\widetilde{\Psi}(\widetilde{\mu},\tau))^{-1}
\widetilde{\Psi}_{\mathrm{WKB}}(\widetilde{\mu},\tau)$, one uses
Equation~\eqref{eq3.81}, with $\varpi \! = \! +1$, and Equation~\eqref{eq3.62}
(in conjunction with Remark~\ref{rem3.2.2}) to show that
\begin{align}
(\widetilde{\Psi}(\widetilde{\mu},\tau))^{-1} \widetilde{\Psi}_{\mathrm{WKB}}
(\widetilde{\mu},\tau) &\underset{\tau \to +\infty}{=} \, \left(\mathrm{I} \!
+ \! (\psi_{0}(\widetilde{\Lambda}))^{-1} \mathcal{O}(\tau^{6 \varepsilon
+2 \varepsilon \Re (\nu +1)+ \frac{\delta}{2}- \delta_{1}}) \psi_{0}
(\widetilde{\Lambda}) \right) \! (\psi_{0}(\widetilde{\Lambda}))^{-1}
(\mathcal{F}_{\tau}(\widetilde{\Lambda}))^{-1}T(\widetilde{\mu}) \nonumber \\
\underset{\tau \to +\infty}{=}& \, \left(\mathrm{I} \! + \! (\psi_{0}
(\widetilde{\Lambda}))^{-1} \mathcal{O}(\tau^{6 \varepsilon +2 \varepsilon
\Re (\nu +1)+ \frac{\delta}{2}- \delta_{1}}) \psi_{0}(\widetilde{\Lambda})
\right) \mathcal{R}_{0}^{-1} \nonumber \\
\times& \, \exp \! \left(- \! \left(\mi 2 \sqrt{3} \, \widetilde{\Lambda}^{2}
\! - \! (\nu \! + \! 1) \ln (\me^{\frac{\mi \pi}{4}}2^{3/2}3^{1/4}
\widetilde{\Lambda}) \right) \sigma_{3} \right) \nonumber \\
\times& \, \left(
\begin{pmatrix}
1 & -\frac{\varkappa_{0}(\tau)}{8 \sqrt{3} \, \widetilde{\Lambda}} \\
0 & 1
\end{pmatrix} \! + \! \dfrac{1}{\widetilde{\Lambda}}
\begin{pmatrix}
0 & 0 \\
\frac{\mi (\nu +1)}{\varkappa_{0}(\tau)} & 0
\end{pmatrix} \! + \! \dfrac{1}{\widetilde{\Lambda}^{2}}
\begin{pmatrix}
\frac{\mi \nu (\nu +1)}{16 \sqrt{3}} & 0 \\
0 & -\frac{\mi (\nu +1)(\nu +2)}{16 \sqrt{3}}
\end{pmatrix} \right. \nonumber \\
+&\left. \, \dfrac{1}{\widetilde{\Lambda}^{3}}
\begin{pmatrix}
0 & -\frac{\mi \nu (\nu -1) \varkappa_{0}(\tau)}{384} \\
\frac{(\nu +1)(\nu +2)(\nu +3)}{16 \sqrt{3} \, \varkappa_{0}(\tau)} & 0
\end{pmatrix} \! + \! \mathcal{O} \left(\dfrac{1}{\widetilde{\Lambda}^{4}}
\begin{pmatrix}
\mathcal{O}(1) & 0 \\
0 & \mathcal{O}(1)
\end{pmatrix} \right) \right) \nonumber \\
\times& \, \left(\dfrac{1}{3^{1/4}} \sqrt{\dfrac{\varkappa_{0}(\tau)}{
\widetilde{\Lambda}}} \, \right)^{\sigma_{3}}
\left(\begin{pmatrix}
-\frac{1}{8} & 1 \\
-1 & 0
\end{pmatrix} \! + \!
\begin{pmatrix}
\mathcal{O}(\tau^{\delta -2 \varepsilon}) & \mathcal{O}(\tau^{\delta -2
\varepsilon}) \\
\mathcal{O}(\tau^{2 \delta -2 \varepsilon}) & 0
\end{pmatrix} \right) \nonumber \\
\times& \, \left(\dfrac{1}{\alpha^{3/2}} \sqrt{\dfrac{b(\tau) \hat{r}_{0}
(\tau)}{\hat{u}_{0}(\tau)}} \, \right)^{\sigma_{3}} \nonumber \\
\underset{\tau \to +\infty}{=}& \, \left(\mathrm{I} \! + \! (\psi_{0}
(\widetilde{\Lambda}))^{-1} \mathcal{O}(\tau^{6 \varepsilon +2 \varepsilon
\Re (\nu +1)+ \frac{\delta}{2}- \delta_{1}}) \psi_{0}(\widetilde{\Lambda})
\right) \mathcal{R}_{0}^{-1} \nonumber \\
\times& \, \left(\dfrac{\alpha^{3/2}}{3^{1/4}} \sqrt{\dfrac{\varkappa_{0}
(\tau) \hat{u}_{0}(\tau)}{b(\tau) \hat{r}_{0}(\tau) \widetilde{\Lambda}}}
\, \right)^{\sigma_{3}} \exp \! \left(- \! \left(\mi 2 \sqrt{3} \,
\widetilde{\Lambda}^{2} \! - \! (\nu \! + \! 1) \ln
(\me^{\frac{\mi \pi}{4}}2^{3/2}3^{1/4} \widetilde{\Lambda}) \right)
\sigma_{3} \right) \mi \sigma_{2} \nonumber \\
\times& \, \left(\mathrm{I} \! + \! \tau^{-2 \varepsilon}
\begin{pmatrix}
\mathcal{O}(\tau^{2 \delta}) & \mathcal{O} \! \left(\frac{\hat{u}_{0}(\tau)}{b
(\tau) \hat{r}_{0}(\tau)} \right) \\
\mathcal{O} \! \left(\frac{b(\tau) \hat{r}_{0}(\tau) \tau^{2 \delta}}{
\hat{u}_{0}(\tau)} \right) & \mathcal{O}(\tau^{\delta})
\end{pmatrix} \right). \label{eq3.113}
\end{align}
For the calculation of asymptotics for $(\widetilde{\Psi}_{\mathrm{WKB}}
(\widetilde{\mu},\tau))^{-1} \tau^{-\frac{1}{12} \sigma_{3}}Y^{\infty}_{0}
(\tau^{-\frac{1}{6}} \widetilde{\mu},\tau)$, one uses the asymptotics of
$Y^{\infty}_{0}(\tau^{-\frac{1}{6}} \widetilde{\mu},\tau)$ given in
Proposition~\ref{prop1.4}, Equation~\eqref{eq3.51},
Equations~\eqref{eq3.37} and~\eqref{eq3.38}, Equations~\eqref{eq3.41}
and~\eqref{eq3.87}, and conditions~\eqref{newcconds2}, to show that
\begin{equation}
(\widetilde{\Psi}_{\mathrm{WKB}}(\widetilde{\mu},\tau))^{-1}
\tau^{-\frac{1}{12} \sigma_{3}}Y^{\infty}_{0}(\tau^{-\frac{1}{6}}
\widetilde{\mu},\tau) \underset{\tau \to +\infty}{=} \me^{\widetilde{
\mathcal{W}}_{\infty}(\tau) \sigma_{3}} \! \left(\mathrm{I} \! + \!
\mathcal{O}(\tau^{-2(\varepsilon - \delta)}) \! + \! \mathcal{O}
(\tau^{-\frac{1}{3}+3 \varepsilon}) \right), \label{eq3.111}
\end{equation}
where
\begin{align}
\widetilde{\mathcal{W}}_{\infty}(\tau) \underset{\tau \to +\infty}{:=}& \,
-\mi 3(\sqrt{3}- \! 1) \alpha^{2} \tau^{2/3} \! - \! \mi 2 \sqrt{3} \,
\widetilde{\Lambda}^{2} \! + \! \dfrac{\mi a}{6} \ln \tau \! + \!
\mathcal{I}_{\infty}^{\sharp}(\tau) \! - \! C^{\mathrm{WKB}}_{\infty}
\nonumber \\
+& \, \dfrac{\mi}{2 \sqrt{3}} \left(\left(a \! - \! \dfrac{\mi}{2} \right) \!
+ \! \dfrac{h_{0}(\tau)}{\alpha^{2}} \! - \! \dfrac{\mathfrak{p}_{\tau}}{4}
\right) \! \left(\dfrac{1}{3} \ln \tau \! - \! \ln \widetilde{\Lambda}
\right). \label{eq3.112}
\end{align}
Hence, via expansions~\eqref{eq3.113} and~\eqref{eq3.111}, upon taking
into account definition~\eqref{eq3.63}, and conditions~\eqref{newconds1}
and~\eqref{newcconds2}, one arrives at
\begin{equation*}
L_{\infty}(\tau) \underset{\tau \to +\infty}{=} \mi \, \mathcal{R}_{0}^{-1}
\exp \left(\left(\mathfrak{t}_{\infty} \! + \! \mathfrak{c}_{\infty} \! -
\! \mathcal{I}_{\infty}^{\sharp}(\tau) \right) \sigma_{3} \right) \! \left(
\sqrt{\dfrac{\varkappa_{0}(\tau) \hat{u}_{0}(\tau)}{b(\tau) \hat{r}_{0}
(\tau)}} \, \right)^{\sigma_{3}} \sigma_{2} \left(\mathrm{I} \! + \!
E_{\infty}(\tau) \right),
\end{equation*}
where $\mathfrak{t}_{\infty}$ and $\mathfrak{c}_{\infty}$ are defined by
Equations~\eqref{eq3.108} and~\eqref{eq3.109}, respectively, and
\begin{align*}
E_{\infty}(\tau) \underset{\tau \to +\infty}{=}& \, \tau^{6 \varepsilon
+2 \varepsilon \Re (\nu +1)+ \frac{\delta}{2}- \delta_{1}}
\begin{pmatrix}
\mathcal{O}(1) & \mathcal{O} \! \left(\frac{\varkappa_{0}(\tau)(\hat{u}_{0}
(\tau))^{2} \me^{\eta_{1}(\tau)}}{b(\tau)(\hat{r}_{0}(\tau))^{2}} \right) \\
\mathcal{O} \! \left(\frac{b(\tau)(\hat{r}_{0}(\tau))^{2} \me^{-\eta_{1}
(\tau)}}{\varkappa_{0}(\tau)(\hat{u}_{0}(\tau))^{2}} \right) & \mathcal{O}(1)
\end{pmatrix} \\
+& \, \tau^{-2 \varepsilon}
\begin{pmatrix}
\mathcal{O}(\tau^{2 \delta}) & \mathcal{O} \! \left(\frac{(\hat{u}_{0}
(\tau))^{2} \tau^{\varepsilon} \me^{\eta_{1}(\tau)}}{b(\tau)(\hat{r}_{0}
(\tau))^{2}} \right) \\
\mathcal{O} \! \left(\frac{(\hat{r}_{0}(\tau))^{2}b(\tau) \tau^{2 \delta}
\tau^{-\varepsilon} \me^{-\eta_{1}(\tau)}}{(\hat{u}_{0}(\tau))^{2}} \right)
& \mathcal{O}(\tau^{\delta})
\end{pmatrix} \\
+& \, \left(\mathcal{O}(\tau^{-2(\varepsilon - \delta)}) \! + \!
\mathcal{O}(\tau^{-\frac{1}{3}+3 \varepsilon}) \right) \!
\begin{pmatrix}
\mathcal{O}(1) & o(1) \\
o(1) & \mathcal{O}(1)
\end{pmatrix},
\end{align*}
where $\eta_{1}(\tau) \! := \! \tfrac{1}{3} \Im (a) \ln \tau \! - \! 2
\varepsilon \Re (\nu \! + \! 1) \ln \tau$. Invoking conditions~\eqref{eq3.103}
and~\eqref{eq3.104}, and conditions~\eqref{newconds1} and~\eqref{newcconds2},
one arrives at Equation~\eqref{eq3.107}. \hfill $\qed$
\begin{cccc} \label{lem3.3.2}
Let $\widetilde{\Psi}(\widetilde{\mu},\tau)$ be the fundamental solution of
Equation~\eqref{eq3.3} with asymptotics given in Lemma~{\rm \ref{lem3.2.1}},
and $X^{0}_{0}(\widetilde{\mu},\tau)$ be the canonical solution of
Equation~\eqref{eq3.1}. Define\footnote{Note that $\tau^{-\frac{1}{12}
\sigma_{3}}X^{0}_{0}(\tau^{-\frac{1}{6}} \widetilde{\mu},\tau)$ is also a
fundamental solution of Equation~\eqref{eq3.3}; therefore, $L_{0}(\tau)$ is
independent of $\widetilde{\mu}$.}
\begin{equation*}
L_{0}(\tau) \! := \! \left(\widetilde{\Psi}(\widetilde{\mu},\tau) \right)^{-1}
\tau^{-\frac{1}{12} \sigma_{3}}X^{0}_{0}(\tau^{-\frac{1}{6}} \widetilde{\mu},
\tau).
\end{equation*}
Assume that the parameters $\varepsilon$, $\delta$, $\delta_{1}$, and $\nu
\! + \! 1$ satisfy the conditions~\eqref{newconds1} and~\eqref{newcconds2}$;$
furthermore, let
\begin{gather}
6 \varepsilon \! + \! \delta \! - \! \delta_{1} \! < \! 0, \qquad
6 \varepsilon \! + \! 4 \varepsilon \Re (\nu  \! + \! 1) \! - \! \delta_{1}
\! < \! 0, \qquad -3 \varepsilon \! + \! 2 \varepsilon \Re (\nu \! + \! 1)
\! + \! 2 \delta \! < \! 0. \label{newconds3}
\end{gather}
Then
\begin{align}
L_{0}(\tau) \underset{\tau \to +\infty}{=}& \, -\mi \, \mathcal{R}_{2}^{-1}
\exp \left(\left(\mathfrak{t}_{0} \! + \! \mathfrak{c}_{0} \! + \!
\mathcal{I}_{0}^{\sharp}(\tau) \right) \sigma_{3} \right) \left(
\sqrt{\dfrac{\varkappa_{0}(\tau) \hat{r}_{0}(\tau)}{\hat{u}_{0}(\tau)}} \,
\right)^{\sigma_{3}} \, \sigma_{2} \nonumber \\
\times& \, \left(\mathrm{I} \! + \!
\begin{pmatrix}
\mathcal{O}(\tau^{-\overset{0}{\delta_{11}}}) & \mathcal{O}
(\tau^{-\overset{0}{\delta_{12}}}) \\
\mathcal{O}(\tau^{-\overset{0}{\delta_{21}}}) & \mathcal{O}
(\tau^{-\overset{0}{\delta_{22}}})
\end{pmatrix} \right), \label{eq3.115}
\end{align}
where $\mathcal{R}_{2}$ is defined in Remark~{\rm \ref{rem3.2.2}},
$\mathcal{I}_{0}^{\sharp}(\tau)$ is defined by Equation~\eqref{eq3.89},
\begin{gather}
\mathfrak{t}_{0} \! := \! \mi 3 \sqrt{3} \, \alpha^{2} \tau^{2/3} \! +
\! \dfrac{1}{3} \! \left(\nu \! +  \! \dfrac{1}{2} \right) \ln \tau,
\label{eq3.116} \\
\mathfrak{c}_{0} \! := \! -\dfrac{9}{4} \ln 2 \! - \! \dfrac{1}{4} \ln 3 \!
+ \! (\nu \! + \! 1) \ln (\me^{\frac{\mi \pi}{4}}2^{3/2}3^{1/4}) \! - \!
C^{\mathrm{WKB}}_{0}, \label{eq3.117}
\end{gather}
with $C^{\mathrm{WKB}}_{0}$ defined by Equation~\eqref{eq3.40}, and
$\overset{0}{\delta_{11}} \! = \! \overset{0}{\delta_{22}} \! := \! \min
\lbrace \delta_{1} \! - \! \tfrac{\delta}{2} \! - \! 6 \varepsilon \! - \! 2
\varepsilon \Re (\nu \! + \! 1),2(\varepsilon \! - \! \delta),\tfrac{1}{3} \!
- \! 3 \varepsilon \rbrace$, $\overset{0}{\delta_{12}} \! := \! \min \lbrace
\delta_{1} \! - \! \delta \! - \! 6 \varepsilon,\varepsilon \! + \! 2
\varepsilon \Re (\nu \! + \! 1) \rbrace$, and $\overset{0}{\delta_{21}} \! :=
\! \min \lbrace \delta_{1} \! - \! 6 \varepsilon \! - \! 4 \varepsilon \Re
(\nu \! + \! 1),3 \varepsilon \! - \! 2 \delta \! - \! 2 \varepsilon \Re (\nu
\! + \! 1) \rbrace$.
\end{cccc}

\emph{Proof}. Denote by $\widetilde{\Psi}_{\mathrm{WKB}}(\widetilde{\mu},
\tau)$ the solution of Equation~\eqref{eq3.3} which has $\mathrm{WKB}$
asymptotics given by Equations~\eqref{eq3.16}--\eqref{eq3.18} in the canonical
domain containing the Stokes curve approaching the positive real $\widetilde{
\mu}$-axis {}from above as $\widetilde{\mu} \! \to \! +0$. Now, re-write
$L_{0}(\tau)$ in the following, equivalent form:
\begin{equation*}
L_{0}(\tau) \! = \! \left((\widetilde{\Psi}(\widetilde{\mu},\tau))^{-1}
\widetilde{\Psi}_{\mathrm{WKB}}(\widetilde{\mu},\tau) \right) \!
\left((\widetilde{\Psi}_{\mathrm{WKB}}(\widetilde{\mu},\tau))^{-1}
\tau^{-\frac{1}{12} \sigma_{3}}X^{0}_{0}(\tau^{-\frac{1}{6}} \widetilde{\mu},
\tau) \right).
\end{equation*}
Noting that $\widetilde{\Psi}(\widetilde{\mu},\tau)$, $\widetilde{\Psi}_{
\mathrm{WKB}}(\widetilde{\mu},\tau)$, and $\tau^{-\frac{1}{12} \sigma_{3}}
X^{0}_{0}(\tau^{-\frac{1}{6}} \widetilde{\mu},\tau)$ are all solutions of
Equation~\eqref{eq3.3}, it follows that they differ by right-hand,
$\widetilde{\mu}$-independent matrix factors. Using this fact, one evaluates
asymptotically (as $\tau \! \to \! +\infty)$ each pair by taking separate
limits, that is, $\widetilde{\mu} \! \to \! \alpha$ and $\widetilde{\mu} \!
\to \! +0$, respectively; more precisely,
\begin{gather*}
(\widetilde{\Psi}(\widetilde{\mu},\tau))^{-1} \widetilde{\Psi}_{\mathrm{WKB}}
(\widetilde{\mu},\tau) \underset{\tau \to +\infty}{=} \underbrace{(\mathcal{
F}_{\tau}(\widetilde{\Lambda})(\mathrm{I} \! + \! \mathcal{O}(\tau^{6
\varepsilon +2 \varepsilon \Re (\nu +1)+ \frac{\delta}{2}-\delta_{1}}))
\psi_{0}(\widetilde{\Lambda}))^{-1}T(\widetilde{\mu}) \me^{\mathrm{WKB}
(\widetilde{\mu},\tau)}}_{\widetilde{\mu}=\widetilde{\mu}_{0}= \alpha +
\tau^{-1/3} \widetilde{\Lambda}, \quad \widetilde{\Lambda} \sim \tau^{
\varepsilon}, \, \, 0< \varepsilon <1/9, \quad \arg (\widetilde{\Lambda})=
\pi} \, , \\
(\widetilde{\Psi}_{\mathrm{WKB}}(\widetilde{\mu},\tau))^{-1} \tau^{-
\frac{1}{12} \sigma_{3}}X^{0}_{0}(\tau^{-\frac{1}{6}} \widetilde{\mu},\tau)
\underset{\tau \to +\infty}{=} \underbrace{(T(\widetilde{\mu})
\me^{\mathrm{WKB}(\widetilde{\mu},\tau)})^{-1} \tau^{-\frac{1}{12} \sigma_{3}}
X^{0}_{0}(\tau^{-\frac{1}{6}} \widetilde{\mu},\tau)}_{\widetilde{\mu} \to 0,
\quad \arg (\widetilde{\mu})=0} \, ,
\end{gather*}
where $\mathrm{WKB}(\widetilde{\mu},\tau) \! := \! -\sigma_{3} \mi \tau^{2/3}
\int_{\widetilde{\mu}_{0}}^{\widetilde{\mu}}l(\xi) \, \md \xi \! - \! \int_{
\widetilde{\mu}_{0}}^{\widetilde{\mu}} \diag ((T(\xi))^{-1} \partial_{\xi}
T(\xi)) \, \md \xi$.

To calculate asymptotics of $(\widetilde{\Psi}(\widetilde{\mu},\tau))^{-1}
\widetilde{\Psi}_{\mathrm{WKB}}(\widetilde{\mu},\tau)$, one uses
Equation~\eqref{eq3.81}, with $\varpi \! = \! -1$, and Equation~\eqref{eq3.62}
(in conjunction with Remark~\ref{rem3.2.2}) to show that
\begin{align}
(\widetilde{\Psi}(\widetilde{\mu},\tau))^{-1} \widetilde{\Psi}_{\mathrm{WKB}}
(\widetilde{\mu},\tau) \underset{\tau \to +\infty}{=}& \, \left(\mathrm{I} \!
+ \! (\psi_{0}(\widetilde{\Lambda}))^{-1} \mathcal{O}(\tau^{6 \varepsilon +2
\varepsilon \Re (\nu +1)+ \frac{\delta}{2}- \delta_{1}}) \psi_{0}
(\widetilde{\Lambda}) \right) \! (\psi_{0}(\widetilde{\Lambda}))^{-1}
(\mathcal{F}_{\tau}(\widetilde{\Lambda}))^{-1}T(\widetilde{\mu}) \nonumber \\
\underset{\tau \to +\infty}{=}& \, \left(\mathrm{I} \! + \! (\psi_{0}
(\widetilde{\Lambda}))^{-1} \mathcal{O}(\tau^{6 \varepsilon +2 \varepsilon
\Re (\nu +1)+ \frac{\delta}{2}- \delta_{1}}) \psi_{0}(\widetilde{\Lambda})
\right) \mathcal{R}_{2}^{-1} \nonumber \\
\times& \, \exp \! \left(- \! \left(\mi 2 \sqrt{3} \, \widetilde{\Lambda}^{2}
\! - \! (\nu \! + \! 1) \ln (\me^{\frac{\mi \pi}{4}}2^{3/2}3^{1/4}
\widetilde{\Lambda}) \right) \sigma_{3} \right) \nonumber \\
\times& \, \left(
\begin{pmatrix}
1 & -\frac{\varkappa_{0}(\tau)}{8 \sqrt{3} \, \widetilde{\Lambda}} \\
0 & 1
\end{pmatrix} \! + \! \dfrac{1}{\widetilde{\Lambda}}
\begin{pmatrix}
0 & 0 \\
\frac{\mi (\nu +1)}{\varkappa_{0}(\tau)} & 0
\end{pmatrix} \! + \! \dfrac{1}{\widetilde{\Lambda}^{2}}
\begin{pmatrix}
\frac{\mi \nu (\nu +1)}{16 \sqrt{3}} & 0 \\
0 & -\frac{\mi (\nu +1)(\nu +2)}{16 \sqrt{3}}
\end{pmatrix} \right. \nonumber \\
+&\left. \, \dfrac{1}{\widetilde{\Lambda}^{3}}
\begin{pmatrix}
0 & -\frac{\mi \nu (\nu -1) \varkappa_{0}(\tau)}{384} \\
\frac{(\nu +1)(\nu +2)(\nu +3)}{16 \sqrt{3} \, \varkappa_{0}(\tau)} & 0
\end{pmatrix} \! + \! \mathcal{O} \left(\dfrac{1}{\widetilde{\Lambda}^{4}}
\begin{pmatrix}
\mathcal{O}(1) & 0 \\
0 & \mathcal{O}(1)
\end{pmatrix} \right) \right) \nonumber \\
\times& \, \left(\dfrac{\me^{-\frac{\mi \pi}{2}}}{3^{1/4}} \sqrt{\dfrac{
\varkappa_{0}(\tau)}{\widetilde{\Lambda}}} \, \right)^{\sigma_{3}}
\left(\begin{pmatrix}
-\frac{1}{8} & 1 \\
0 & -8
\end{pmatrix} \! + \!
\begin{pmatrix}
\mathcal{O}(\tau^{\delta -2 \varepsilon}) & \mathcal{O}(\tau^{\delta -2
\varepsilon}) \\
0 & \mathcal{O}(\tau^{2 \delta -2 \varepsilon})
\end{pmatrix} \right) \nonumber \\
\times& \, \left(\dfrac{1}{\alpha^{3/2}} \sqrt{\dfrac{b(\tau) \hat{r}_{0}
(\tau)}{\hat{u}_{0}(\tau)}} \, \right)^{\sigma_{3}} \nonumber \\
\underset{\tau \to +\infty}{=}& \, -\left(\mathrm{I} \! + \! (\psi_{0}
(\widetilde{\Lambda}))^{-1} \mathcal{O}(\tau^{6 \varepsilon +2 \varepsilon
\Re (\nu +1)+ \frac{\delta}{2}- \delta_{1}}) \psi_{0}(\widetilde{\Lambda})
\right) \nonumber \\
\times& \, \mathcal{R}_{2}^{-1} \left(\dfrac{-\mi}{8 \pmb{\cdot} 3^{1/4}
\alpha^{3/2}} \sqrt{\dfrac{b(\tau) \varkappa_{0}(\tau) \hat{r}_{0}(\tau)}{
\hat{u}_{0}(\tau) \widetilde{\Lambda}}} \, \right)^{\sigma_{3}} \nonumber \\
\times& \, \exp \! \left(- \! \left(\mi 2 \sqrt{3} \, \widetilde{\Lambda}^{2}
\! - \! (\nu \! + \! 1) \ln (\me^{\frac{\mi \pi}{4}}2^{3/2}3^{1/4}
\widetilde{\Lambda}) \right) \sigma_{3} \right) \nonumber \\
\times& \, \left(\mathrm{I} \! + \! \tau^{-2 \varepsilon}
\begin{pmatrix}
\mathcal{O}(\tau^{\delta}) & \mathcal{O} \! \left(\frac{\hat{u}_{0}
(\tau) \tau^{2 \delta}}{b(\tau) \hat{r}_{0}(\tau)} \right) \\
\mathcal{O} \! \left(\frac{b(\tau) \hat{r}_{0}(\tau)}{\hat{u}_{0}(\tau)}
\right) & \mathcal{O}(\tau^{2 \delta})
\end{pmatrix} \right). \label{eq3.121}
\end{align}
For the calculation of asymptotics for $(\widetilde{\Psi}_{\mathrm{WKB}}
(\widetilde{\mu},\tau))^{-1} \tau^{-\frac{1}{12} \sigma_{3}}X^{0}_{0}
(\tau^{-\frac{1}{6}} \widetilde{\mu},\tau)$, one uses the asymptotics
of $X^{0}_{0}(\tau^{-\frac{1}{6}} \widetilde{\mu},\tau)$ given in
Proposition~\ref{prop1.4}, Equation~\eqref{eq3.52}, Equations~\eqref{eq3.39}
and~\eqref{eq3.40}, Equations~\eqref{eq3.43} and~\eqref{eq3.89}, and
conditions~\eqref{newcconds2}, to show that
\begin{align}
(\widetilde{\Psi}_{\mathrm{WKB}}(\widetilde{\mu},\tau))^{-1}
\tau^{-\frac{1}{12} \sigma_{3}}X^{0}_{0}(\tau^{-\frac{1}{6}} \widetilde{\mu},
\tau) \underset{\tau \to +\infty}{=}& \, \me^{\widetilde{\mathcal{W}}_{0}
(\tau) \sigma_{3}} \!
\begin{pmatrix}
0 & \frac{\mi 2^{3/4} \alpha^{3/2}}{\sqrt{b(\tau)}} \\
\frac{\mi \sqrt{b(\tau)}}{2^{3/4} \alpha^{3/2}} & 0
\end{pmatrix} \nonumber \\
\times& \, \left(\mathrm{I} \! + \! \mathcal{O}(\tau^{-2(\varepsilon
-\delta)}) \! + \! \mathcal{O}(\tau^{-\frac{1}{3}+3 \varepsilon}) \right),
\label{eq3.119}
\end{align}
where
\begin{align}
\widetilde{\mathcal{W}}_{0}(\tau) \underset{\tau \to +\infty}{:=}& \,
\mi 3 \sqrt{3} \, \alpha^{2} \tau^{2/3} \! + \! \mi 2 \sqrt{3} \,
\widetilde{\Lambda}^{2} \! + \! \mathcal{I}_{0}^{\sharp}(\tau) \! - \!
C^{\mathrm{WKB}}_{0} \nonumber \\
-& \, \dfrac{\mi}{2 \sqrt{3}} \left(\left(a \! - \! \dfrac{\mi}{2} \right) \!
+ \! \dfrac{h_{0}(\tau)}{\alpha^{2}} \! - \! \dfrac{\mathfrak{p}_{\tau}}{4}
\right) \! \left(\dfrac{1}{3} \ln \tau \! - \! \ln \widetilde{\Lambda}
\right). \label{eq3.120}
\end{align}
Hence, via expansions~\eqref{eq3.121} and~\eqref{eq3.119}, upon taking
into account definition~\eqref{eq3.63}, and conditions~\eqref{newconds1}
and~\eqref{newcconds2}, one arrives at
\begin{equation*}
L_{0}(\tau) \underset{\tau \to +\infty}{=} -\mi \, \mathcal{R}_{2}^{-1}
\exp \! \left(\left(\mathfrak{t}_{0} \! + \! \mathfrak{c}_{0} \! + \!
\mathcal{I}_{0}^{\sharp}(\tau) \right) \sigma_{3} \right) \! \left(
\sqrt{\dfrac{\varkappa_{0}(\tau) \hat{r}_{0}(\tau)}{\hat{u}_{0}(\tau)}} \,
\right)^{\sigma_{3}} \sigma_{2} \left(\mathrm{I} \! + \! E_{0}(\tau) \right),
\end{equation*}
where $\mathfrak{t}_{0}$ and $\mathfrak{c}_{0}$ are defined by
Equations~\eqref{eq3.116} and~\eqref{eq3.117}, respectively, and
\begin{align*}
E_{0}(\tau) \underset{\tau \to +\infty}{=}& \, \tau^{6 \varepsilon +2
\varepsilon \Re (\nu +1)+ \frac{\delta}{2}- \delta_{1}}
\begin{pmatrix}
\mathcal{O}(1) & \mathcal{O}(\varkappa_{0}(\tau) \tau^{-2 \varepsilon \Re
(\nu +1)}) \\
\mathcal{O}((\varkappa_{0}(\tau))^{-1} \tau^{2 \varepsilon \Re (\nu +1)}) &
\mathcal{O}(1)
\end{pmatrix} \\
+& \, \tau^{-2 \varepsilon}
\begin{pmatrix}
\mathcal{O}(\tau^{2 \delta}) & \mathcal{O}(\tau^{-2 \varepsilon \Re (\nu
+1/2)}) \\
\mathcal{O}(\tau^{2 \delta +2 \varepsilon \Re (\nu +1/2)}) & \mathcal{O}
(\tau^{\delta})
\end{pmatrix} \\
+& \, \left(\mathcal{O}(\tau^{-2(\varepsilon - \delta)}) \! + \! \mathcal{O}
(\tau^{-\frac{1}{3}+3 \varepsilon}) \right) \!
\begin{pmatrix}
\mathcal{O}(1) & o(1) \\
o(1) & \mathcal{O}(1)
\end{pmatrix}.
\end{align*}
Invoking conditions~\eqref{newconds3}, and conditions~\eqref{newconds1}
and~\eqref{newcconds2}, one arrives at Equation~\eqref{eq3.115}. \hfill
$\qed$
\begin{dddd} \label{theor3.3.1}
Assuming conditions~\eqref{eq3.11}, \eqref{newconds1}, \eqref{newcconds2},
\eqref{eq3.103}, \eqref{eq3.104}, and~\eqref{newconds3}, define $\delta^{0} \!
:= \! \min \lbrace \overset{0}{\delta_{ij}} \rbrace_{i,j=1,2}$. Furthermore,
suppose that
\begin{equation}
\left(\dfrac{\tau^{\frac{1}{3} \Im (a)}(\hat{u}_{0}(\tau))^{2}}{b(\tau)
(\hat{r}_{0}(\tau))^{2}} \right) \tau^{-\delta^{0}} \underset{\tau \to
+\infty}{=} \mathcal{O}(\tau^{-\delta^{G}_{12}}), \quad \quad \left(
\dfrac{b(\tau)(\hat{r}_{0}(\tau))^{2}}{\tau^{\frac{1}{3} \Im (a)}(\hat{u}_{0}
(\tau))^{2}} \right) \tau^{-\delta^{0}} \underset{\tau \to +\infty}{=}
\mathcal{O}(\tau^{-\delta^{G}_{21}}), \label{newconds4}
\end{equation}
where $\delta^{G}_{12},\delta^{G}_{21} \! > \! 0$. Let
\begin{equation}
\delta_{G} \! := \! \min \lbrace \delta^{0},\delta^{G}_{12},\delta^{G}_{21},
\overset{\infty}{\delta_{12}},\overset{\infty}{\delta_{21}} \rbrace.
\label{newconds5}
\end{equation}
Then the connection matrix $($cf. Equation~\eqref{eqdefg}$)$ has the following
asymptotics:
\begin{equation}
G \underset{\tau \to +\infty}{=} G_{\infty}(\mathrm{I} \! + \! \mathcal{O}
(\tau^{-\delta_{G}})), \label{eq3.122}
\end{equation}
where
\begin{equation}
G_{\infty} \! := \!
\begin{pmatrix}
-\frac{\hat{r}_{0}(\tau) \sqrt{b(\tau)}}{\hat{u}_{0}(\tau)}
\me^{(\mathfrak{z}_{0}(\tau)-\mathfrak{z}_{\infty}(\tau))} \me^{-2 \pi \mi
(\nu +1)} & -\frac{\me^{\frac{\mi \pi}{4}}2^{3/2}3^{1/4} \sqrt{2 \pi}}{
\sqrt{b(\tau)} \, \Gamma (\nu +1)} \me^{(\mathfrak{z}_{0}(\tau)+
\mathfrak{z}_{\infty}(\tau))} \me^{-2 \pi \mi (\nu +1)} \\
\frac{\me^{\frac{\mi \pi}{4}} \sqrt{2 \pi} \, \sqrt{b(\tau)}}{2^{3/2}3^{1/4}
\Gamma (-\nu)} \me^{-(\mathfrak{z}_{0}(\tau)+\mathfrak{z}_{\infty}(\tau))}
\me^{\mi \pi (\nu +1)} & -\frac{\hat{u}_{0}(\tau)}{\hat{r}_{0}(\tau)
\sqrt{b(\tau)}} \me^{-(\mathfrak{z}_{0}(\tau)-\mathfrak{z}_{\infty}(\tau))}
\end{pmatrix}, \label{eq3.123}
\end{equation}
with
\begin{gather}
\mathfrak{z}_{\infty}(\tau) \! = \! \mathfrak{t}_{\infty} \! + \!
\mathfrak{c}_{\infty} \! - \! \mathcal{I}_{\infty}^{\sharp}(\tau),
\label{eq3.124} \\
\mathfrak{z}_{0}(\tau) \! = \! \mathfrak{t}_{0} \! + \! \mathfrak{c}_{0}
\! + \! \mathcal{I}_{0}^{\sharp}(\tau), \label{eq3.125}
\end{gather}
where $\mathfrak{t}_{\infty}$, $\mathfrak{c}_{\infty}$, $\mathcal{I}_{
\infty}^{\sharp}(\tau)$, $\mathfrak{t}_{0}$, $\mathfrak{c}_{0}$, and
$\mathcal{I}_{0}^{\sharp}(\tau)$ are defined by Equations~\eqref{eq3.108},
\eqref{eq3.109}, \eqref{eq3.87}, \eqref{eq3.116}, \eqref{eq3.117},
and~\eqref{eq3.89}, respectively.
\end{dddd}

\emph{Proof}. Starting {}from the definition of the connection matrix
(cf. Equation~\eqref{eqdefg}), one has the following sequence of formulae:
$Y^{\infty}_{0}(\mu) \! := \! X^{0}_{0}(\mu)G$ $\Rightarrow$ (cf.
Equations~\eqref{eq3.2}) $Y^{\infty}_{0}(\tau^{-1/6} \widetilde{\mu}) \! = \!
X^{0}_{0}(\tau^{-1/6} \widetilde{\mu})G$ $\Rightarrow$ $\tau^{-\frac{1}{12}
\sigma_{3}}Y^{\infty}_{0}(\tau^{-1/6} \widetilde{\mu}) \! = \! \tau^{-
\frac{1}{12} \sigma_{3}}X^{0}_{0}(\tau^{-1/6} \widetilde{\mu})G$ $\Rightarrow$
(cf. Lemmata~\ref{lem3.3.1} and~\ref{lem3.3.2}) $\widetilde{\Psi}
(\widetilde{\mu},\tau)L_{\infty}(\tau) \! = \! \widetilde{\Psi}
(\widetilde{\mu},\tau)L_{0}(\tau)G$ $\Rightarrow$ $L_{\infty}(\tau) \! = \!
L_{0}(\tau)G$ $\Rightarrow$
\begin{equation}
G \! = \! (L_{0}(\tau))^{-1}L_{\infty}(\tau). \label{eq3.128}
\end{equation}
Using Equations~\eqref{eq3.107}--\eqref{eq3.109},
Equations~\eqref{eq3.115}--\eqref{eq3.117}, and the formulae for
$\mathcal{R}_{0}$ and $\mathcal{R}_{2}$ given in Remark~\ref{rem3.2.2}, one
arrives at
\begin{equation*}
G \underset{\tau \to +\infty}{=} G_{\infty}(\mathrm{I} \! + \! E_{G}(\tau)),
\end{equation*}
where $G_{\infty}$ is defined by Equations~\eqref{eq3.123}--\eqref{eq3.125},
and
\begin{equation*}
E_{G}(\tau) \underset{\tau \to +\infty}{:=}
\begin{pmatrix}
\mathcal{O}(\tau^{-\overset{\infty}{\delta_{11}}}) &
\mathcal{O}(\tau^{-\overset{\infty}{\delta_{12}}}) \\
\mathcal{O}(\tau^{-\overset{\infty}{\delta_{21}}}) &
\mathcal{O}(\tau^{-\overset{\infty}{\delta_{22}}})
\end{pmatrix} \! + \! (G_{\infty})^{-1}
\begin{pmatrix}
\mathcal{O}(\tau^{-\overset{0}{\delta_{11}}}) &
\mathcal{O}(\tau^{-\overset{0}{\delta_{12}}}) \\
\mathcal{O}(\tau^{-\overset{0}{\delta_{21}}}) &
\mathcal{O}(\tau^{-\overset{0}{\delta_{22}}})
\end{pmatrix} G_{\infty}.
\end{equation*}
The estimate for $E_{G}(\tau)$ (cf. Equations~\eqref{newconds5}
and~\eqref{eq3.122}) follows {}from the definitions of
$\overset{\infty}{\delta_{ij}}$ and $\overset{0}{\delta_{ij}}$, $i,j \! = \!
1,2$, given in Lemmata~\ref{lem3.3.1} and~\ref{lem3.3.2}, respectively, and
conditions~\eqref{newconds4}. \hfill $\qed$
\section{Asymptotic Solution of the Inverse Monodromy Problem} \label{finalsec}
In Section~\ref{sec3} the connection matrix, $G$, for Equation~\eqref{eq3.3}
is calculated asymptotically under certain conditions on its coefficient
functions\footnote{See Equations~\eqref{eq3.4}-\eqref{eq3.10}
and~\eqref{eq3.7}, and conditions~\eqref{eq3.11}, \eqref{newconds1},
\eqref{newcconds2}, \eqref{eq3.103}, \eqref{eq3.104}, \eqref{newconds3},
and~\eqref{newconds4}.}. Using these conditions, one can derive the
$\tau$-dependent class of functions to which $G$ belongs. We are not, however,
going to proceed in this, most general, way; rather, we will invoke the
isomonodromy condition on $G$, that is, its matrix elements, $g_{ij}$, $i,j
\! = \! 1,2$, are constants, and then invert the formula for $G$ to obtain
the coefficient functions of Equation~\eqref{eq3.3} and verify that they
satisfy all of the imposed conditions for this isomonodromy case. The latter
procedure produces explicit formulae for the coefficient functions of
Equation~\eqref{eq3.3}, which give rise to asymptotics of the solution
of the system of isomonodromy deformations~\eqref{eq1.4} (see, also,
Equations~\eqref{eq3.2}) and, in turn, define asymptotics of the solution
of the degenerate third Painlev\'{e} equation~\eqref{eq1.1} and the related
functions $\mathcal{H}(\tau)$ and $f(\tau)$ (cf. Equations~\eqref{eqh1}
and~\eqref{eqf}, respectively). Zeroes and poles of the leading term of
asymptotics of $u(\tau)$ define the centers of the cheese-holes of the
domain $\widetilde{\mathcal{D}}_{u}$ (cf. definitions~\eqref{eqdofd},
\eqref{eqnotat2}, \eqref{eqnotat3}, and~\eqref{eqnotat4}). We also prove
that, for large enough $\tau$, there is a one-to-one correspondence between
the zeroes and poles of the solution $u(\tau)$ and the zeroes and poles of
its leading term of asymptotics, that is, in each cheese-hole there is
located one, and only one, zero or pole of $u(\tau)$.
\begin{bbbbb} \label{cor3.3.2}
Let $g_{ij}$, $i,j \! = \! 1,2$, denote the matrix elements of $G$. Assume
that $g_{11}g_{12}g_{21}g_{22} \! \neq \! 0$ and the conditions of
Theorem~{\rm \ref{theor3.3.1}} are valid. Then the functions $\hat{r}_{0}
(\tau)$, $\hat{u}_{0}(\tau)$, $h_{0}(\tau)$, $\varkappa_{0}^{2}(\tau)$, and
$b(\tau)$ have the following asymptotics:
\begin{gather}
\hat{r}_{0}(\tau) \underset{\tau \to +\infty}{=} \dfrac{12}{1 \! - \! \mi
\sqrt{2} \, \cos (\widetilde{\vartheta}(\tau))} \underset{\tau \to +\infty}{=}
\dfrac{6 \cos \vartheta_{0}}{\cos (\tfrac{1}{2}(\widetilde{\vartheta}(\tau)
\! + \! \vartheta_{0})) \cos (\tfrac{1}{2}(\widetilde{\vartheta}(\tau) \! - \!
\vartheta_{0}))}, \label{eq3.129} \\
\hat{u}_{0}(\tau) \underset{\tau \to +\infty}{=} -\dfrac{3}{2 \cos^{2}
(\tfrac{1}{2}(\widetilde{\vartheta}(\tau) \! + \! \vartheta_{0}))} \! + \!
\mathcal{O}(\tau^{-\frac{1}{3}+ \delta -\delta_{1}}) \underset{\tau \to
+\infty}{=} -\dfrac{\sin^{2} \vartheta_{0}}{\cos^{2}(\tfrac{1}{2}
(\widetilde{\vartheta}(\tau) \! + \! \vartheta_{0}))} \! + \! \mathcal{O}
(\tau^{-\frac{1}{3}+ \delta -\delta_{1}}), \label{eq3.130}
\end{gather}
\begin{align}
h_{0}(\tau) \underset{\tau \to +\infty}{=}& \, \mi 2 \sqrt{3} \, \alpha^{2}
\left((\widetilde{\nu} \! + \! 1) \! - \! \dfrac{1}{2} \! + \! \dfrac{\mi}{
2 \sqrt{3}} \left(a \! - \! \dfrac{\mi}{2} \right) \! - \! \dfrac{\sin
(\widetilde{\vartheta}(\tau)) \cos \vartheta_{0}}{2 \sqrt{2} \, \cos
(\tfrac{1}{2}(\widetilde{\vartheta}(\tau) \! + \! \vartheta_{0})) \cos
(\tfrac{1}{2}(\widetilde{\vartheta}(\tau) \! - \! \vartheta_{0}))} \right.
\nonumber \\
+&\left. \, \mathcal{O}(\tau^{-\delta_{G}}) \! + \! \mathcal{O}
(\tau^{\delta -2 \delta_{1}}) \vphantom{M^{M^{M^{M^{M}}}}} \right),
\label{eq3.131}
\end{align}
\begin{align}
\varkappa_{0}^{2}(\tau) \underset{\tau \to +\infty}{=}& \, \mi 8 \sqrt{3}
\left((\widetilde{\nu} \! + \! 1) \! - \! \dfrac{1}{2} \! + \! \dfrac{\mi}{2
\sqrt{3}} \left(a \! - \! \dfrac{\mi}{2} \right) \! - \! \dfrac{\sin
(\widetilde{\vartheta}(\tau)) \cos \vartheta_{0}}{2 \sqrt{2} \, \cos
(\tfrac{1}{2}(\widetilde{\vartheta}(\tau) \! + \! \vartheta_{0})) \cos
(\tfrac{1}{2}(\widetilde{\vartheta}(\tau) \! - \! \vartheta_{0}))} \right)
\nonumber \\
+& \, \dfrac{4(a \! - \! \tfrac{\mi}{2}) \cos^{2}(\tfrac{1}{2}(\widetilde{
\vartheta}(\tau) \! + \! \vartheta_{0}))}{\cos (\tfrac{1}{2}(\widetilde{
\vartheta}(\tau) \! + \! \vartheta_{0}) \! + \! \vartheta_{0}) \cos
(\tfrac{1}{2}(\widetilde{\vartheta}(\tau) \! + \! \vartheta_{0}) \! - \!
\vartheta_{0})} \! + \! \mathcal{O}(\tau^{-\delta_{G}}) \! + \! \mathcal{O}
(\tau^{\delta -2 \delta_{1}}) \! + \! \mathcal{O}(\tau^{-\frac{2}{3}+3
\delta}), \label{eq3.132}
\end{align}
\begin{equation}
\sqrt{b(\tau)} \underset{\tau \to +\infty}{=} \left(\dfrac{\cos (\tfrac{1}{2}
(\widetilde{\vartheta}(\tau) \! - \! \vartheta_{0}))}{4 \cos (\vartheta_{0})
\cos (\tfrac{1}{2}(\widetilde{\vartheta}(\tau) \! + \! \vartheta_{0}))}
\! + \! \mathcal{O}(\tau^{-\frac{2}{3}}) \! + \! \mathcal{O}(\tau^{\delta
-2 \delta_{1}}) \right) \exp (\Phi (\tau)), \label{eq3.133}
\end{equation}
where
\begin{equation}
\widetilde{\nu} \! + \! 1 \! := \! \dfrac{\mi}{2 \pi} \ln (g_{11}g_{22}),
\qquad 0 \! < \! \Re (\widetilde{\nu} \! + \! 1) \! < \! 1, \label{eq3.134}
\end{equation}
\begin{align}
\widetilde{\vartheta}(\tau) :=& \, 6 \sqrt{3} \, \alpha^{2} \tau^{2/3} \! - \!
\mi \left((\widetilde{\nu} \! + \! 1) \! - \! \dfrac{1}{2} \right) \ln (6
\sqrt{3} \, \alpha^{2} \tau^{2/3}) \! - \! \mi \left((\widetilde{\nu} \! + \!
1) \! - \! \dfrac{1}{2} \right) \ln 12 \! + \! \left(a \! - \! \dfrac{\mi}{2}
\right) \ln (2 \! + \! \sqrt{3}) \nonumber \\
+& \, \mi \ln \left(\dfrac{g_{11}g_{12} \Gamma (\widetilde{\nu} \! + \! 1)}{
\sqrt{2 \pi}} \right) \! - \! \dfrac{3 \pi}{2}(\widetilde{\nu} \! + \! 1) \!
+ \! \dfrac{7 \pi}{4} \! + \! \mathcal{O}(\tau^{-\delta_{G}} \ln \tau),
\label{eq3.135}
\end{align}
\begin{equation}
\vartheta_{0} \! := \! -\dfrac{\pi}{2} \! + \! \dfrac{\mi}{2} \ln (2 \! + \!
\sqrt{3}), \qquad \cos \vartheta_{0} \! = \! \dfrac{\mi}{\sqrt{2}}, \qquad
\sin \vartheta_{0} \! = \! -\sqrt{\dfrac{3}{2}}, \label{eq3.136}
\end{equation}
\begin{align}
\Phi (\tau) :=& \, -\mi 3 \alpha^{2} \tau^{2/3} \! - \! \dfrac{\mi a}{6} \ln
\tau \! + \! \mi \pi \left((\widetilde{\nu} \! + \! 1) \! - \! \dfrac{1}{2}
\right) \! + \! \ln g_{11} \! + \! 2 \ln (2 \alpha) \! + \! \dfrac{\mi a}{2}
\ln \! \left(\dfrac{\alpha^{2}}{2} \right) \nonumber \\
-& \, \left((\widetilde{\nu} \! + \! 1) \! - \! \dfrac{1}{2} \right)
\ln (2 \! + \! \sqrt{3}) \! + \! \mathcal{O}(\tau^{-\delta_{G}}),
\label{eq3.137}
\end{align}
and $\delta_{G} \! > \! 0$ is defined by Equation~\eqref{newconds5}.
\end{bbbbb}

\emph{Proof}. Multiplying the diagonal elements of $G$ (cf.
Equations~\eqref{eq3.122}--\eqref{eq3.125}) and taking into account
conditions~\eqref{newconds1} and definition~\eqref{eq3.134}, one shows that
\begin{equation}
\nu \! + \! 1 \underset{\tau \to +\infty}{=} (\widetilde{\nu} \! + \! 1)
(1 \! + \! \mathcal{O}(\tau^{-\delta_{G}})), \label{eq3.139}
\end{equation}
where $\delta_{G} \! > \! 0$ is defined by Equation~\eqref{newconds5}.

Multiplying the elements of the first row of $G$, one shows, via
Equation~\eqref{eq3.139}, that
\begin{equation}
g_{11}g_{12} \underset{\tau \to +\infty}{=} \exp \! \left(\dfrac{\mi \pi}{4}
\! + \! \dfrac{1}{2} \ln 2 \pi \! + \! \ln (2^{\frac{3}{2}}3^{\frac{1}{4}}) \!
- \! \ln \Gamma (\widetilde{\nu} \! + \! 1) \! - \! 4 \pi \mi (\widetilde{\nu}
\! + \! 1) \! + \! \ln \! \left(\dfrac{\hat{r}_{0}(\tau)}{\hat{u}_{0}(\tau)}
\right) \! + \! 2 \mathfrak{z}_{0}(\tau) \! + \! \mathcal{O}
(\tau^{-\delta_{G}}) \right), \label{intermedg}
\end{equation}
where, via Equations~\eqref{eq3.125}, \eqref{eq3.117}, \eqref{eq3.89},
and~\eqref{eq3.139},
\begin{align}
2 \mathfrak{z}_{0}(\tau) \underset{\tau \to +\infty}{=}& \, 2 \mathfrak{t}_{0}
\! - \! 3 \ln 2 \! - \! \dfrac{1}{2} \ln 3 \! + \! \dfrac{3 \pi \mi}{2} \! +
\! 2(\widetilde{\nu} \! + \! 1) \ln (\me^{\frac{\mi \pi}{4}}2^{\frac{3}{2}}
3^{\frac{1}{4}}) \! + \! 2 \mi \! \left(a \! - \! \dfrac{\mi}{2} \right)
\ln \! \left(\dfrac{\sqrt{3}+ \! 1}{\sqrt{2}} \right) \nonumber \\
+& \, 2 \left((\widetilde{\nu} \! + \! 1) \! - \! \dfrac{1}{2} \right) \ln
3 \alpha \! + \! 2 \pi \mi (\widetilde{\nu} \! + \! 1) \! - \! \ln \! \left(
\dfrac{\hat{r}_{0}(\tau)}{\hat{u}_{0}(\tau)} \right) \! - \! \dfrac{1}{2} \ln
\mathfrak{Q}(\tau) \! + \! \mathcal{O}(\tau^{-\delta_{G}}), \label{eq3.140}
\end{align}
with $\mathfrak{t}_{0}$ defined by Equation~\eqref{eq3.116}, and
\begin{equation}
\mathfrak{Q}(\tau) \! := \! \left(\dfrac{\hat{r}_{0}(\tau) \! - \! 4 \sqrt{3}
+ \! \sqrt{32 \! + \! (\hat{r}_{0}(\tau) \! - \! 4)^{2}}}{\hat{r}_{0}(\tau)
\! - \! 4 \sqrt{3}- \! \sqrt{32 \! + \! (\hat{r}_{0}(\tau) \! - \! 4)^{2}}}
\right) \! \left(\dfrac{\hat{r}_{0}(\tau) \! + \! 4 \sqrt{3}- \! \sqrt{32 \!
+ \! (\hat{r}_{0}(\tau) \! - \! 4)^{2}}}{\hat{r}_{0}(\tau) \! + \! 4 \sqrt{3}
+ \! \sqrt{32 \! + \! (\hat{r}_{0}(\tau) \! - \! 4)^{2}}} \right).
\label{eq3.141}
\end{equation}
Substituting Equations~\eqref{eq3.140} and~\eqref{eq3.141} into
Equation~\eqref{intermedg} and solving for $\mathfrak{Q}(\tau)$, one
arrives at
\begin{equation}
\dfrac{12 \! - \! \hat{r}_{0}(\tau) \! - \! \sqrt{3} \, \sqrt{32 \! + \!
(\hat{r}_{0}(\tau) \! - \! 4)^{2}}}{12 \! - \! \hat{r}_{0}(\tau) \! + \!
\sqrt{3} \, \sqrt{32 \! + \! (\hat{r}_{0}(\tau) \! - \! 4)^{2}}}
\underset{\tau \to +\infty}{=} \exp (2 \mi \widetilde{\vartheta}(\tau)),
\label{eq3.142}
\end{equation}
where $\widetilde{\vartheta}(\tau)$ is defined by Equation~\eqref{eq3.135}.
There are two solutions of Equation~\eqref{eq3.142} for $\hat{r}_{0}(\tau)$.
Define the parameter $\vartheta_{0}$ by the first equation of \eqref{eq3.136};
then, one of the solutions for $\hat{r}_{0}(\tau)$ is given in
Equation~\eqref{eq3.129}, whilst the other is given by the same formula as
in Equation~\eqref{eq3.129} but with the change $\widetilde{\vartheta}(\tau)
\! \to \! \widetilde{\vartheta}(\tau) \! + \! \pi$. This ambiguity in the
solution of the inverse monodromy problem can be rectified by calculating one
additional parameter, namely, one of the Stokes multipliers\footnote{This
ambiguity is distinct {}from the one for the direct monodromy problem
mentioned in Remark~\ref{newrem13}.}; however, we choose a different
approach. There is a common domain of validity of the results obtained
in this paper and Part I \cite{a1}. In Appendix~B, this fact is used
in order to justify the choice for $\hat{r}_{0}(\tau)$ that is given
in Equation~\eqref{eq3.129}. Furthermore, as a by-product of
Equations~\eqref{eq3.142}, \eqref{eq3.136}, and~\eqref{eq3.129}, one shows
that
\begin{equation}
\sqrt{32 \! + \! (\hat{r}_{0}(\tau) \! - \! 4)^{2}} \underset{\tau \to
+\infty}{=} -\dfrac{4 \sqrt{6} \, \sin \widetilde{\vartheta}(\tau)}{1 \!
- \! \mi \sqrt{2} \, \cos \widetilde{\vartheta}(\tau)} \underset{\tau \to
+\infty}{=} \dfrac{2 \sin (\widetilde{\vartheta}(\tau)) \sin 2 \vartheta_{0}}{
\cos (\tfrac{1}{2}(\widetilde{\vartheta}(\tau) \! + \! \vartheta_{0}))
\cos (\tfrac{1}{2}(\widetilde{\vartheta}(\tau) \! - \! \vartheta_{0}))}.
\label{eq3.144}
\end{equation}
Via Equations~\eqref{eq3.65}, \eqref{eq3.129}, \eqref{eq3.136},
and~\eqref{eq3.144}, and conditions~\eqref{eq3.11}, one shows that
\begin{equation}
\hat{u}_{0}(\tau) \underset{\tau \to +\infty}{=} -\dfrac{3}{2} \dfrac{(\mi
\sqrt{2} \, (\cos (\vartheta_{0}) \! + \! \cos \widetilde{\vartheta}(\tau))
\! - \! \sqrt{6} \, (\sin (\vartheta_{0}) \! + \! \sin \widetilde{\vartheta}
(\tau)))}{(\cos (\vartheta_{0}) \! + \! \cos \widetilde{\vartheta}(\tau))^{2}}
\! + \! \mathcal{O}(\tau^{-\frac{1}{3}+ \delta -\delta_{1}}). \label{eq3.145}
\end{equation}
Now, {}from Equations~\eqref{eq3.136}, \eqref{eq3.139}, and~\eqref{eq3.145},
and the trigonometric identities $\sin (\vartheta_{0}) \! + \! \sin
\widetilde{\vartheta}(\tau) \! = \! 2 \sin (\tfrac{1}{2}(\vartheta_{0} \!
+ \! \widetilde{\vartheta}(\tau))) \cos (\tfrac{1}{2}(\vartheta_{0} \! -
\! \widetilde{\vartheta}(\tau)))$ and $\cos (\vartheta_{0}) \! + \! \cos
\widetilde{\vartheta}(\tau) \! = \! 2 \cos (\tfrac{1}{2}(\vartheta_{0} \! +
\! \widetilde{\vartheta}(\tau))) \cos (\tfrac{1}{2}(\vartheta_{0} \! - \!
\widetilde{\vartheta}(\tau)))$, one arrives at Equation~\eqref{eq3.130}.
To obtain $h_{0}(\tau)$ (cf. Equation~\eqref{eq3.131}), one uses
Equations~\eqref{eq3.49}, \eqref{eq3.63}, \eqref{eq3.136}, \eqref{eq3.139},
and~\eqref{eq3.144}, and conditions~\eqref{eq3.11}. Finally, using
Equations~\eqref{eq3.58}, \eqref{eq3.130}, and~\eqref{eq3.131}, one arrives
at Equation~\eqref{eq3.132} for $\varkappa_{0}^{2}(\tau)$.

It follows {}from the $(1 \, 2)$-element of $G$ and Equation~\eqref{eq3.139}
that
\begin{equation}
\sqrt{b(\tau)} \underset{\tau \to +\infty}{=} \exp \! \left(
\dfrac{5 \pi \mi}{4} \! + \! \ln (2^{\frac{3}{2}}3^{\frac{1}{4}})
\! - \! 2 \pi \mi (\widetilde{\nu} \! + \! 1) \! + \! \ln \! \left(
\dfrac{\sqrt{2 \pi}}{g_{12} \Gamma (\widetilde{\nu} \! + \! 1)} \right)
\! + \! \mathfrak{z}_{0}(\tau) \! + \! \mathfrak{z}_{\infty}(\tau) \! + \!
\mathcal{O}(\tau^{-\delta_{G}}) \right), \label{intermedg2}
\end{equation}
where $\mathfrak{z}_{0}(\tau)$ is given in Equations~\eqref{eq3.140}
and~\eqref{eq3.141}, and (cf. Equation~\eqref{eq3.124}) $\mathfrak{z}_{\infty}
(\tau) \! = \! \mathfrak{t}_{\infty} \! + \! \mathfrak{c}_{\infty} \! - \!
\mathcal{I}_{\infty}^{\sharp}(\tau)$, with $\mathfrak{t}_{\infty}$,
$\mathfrak{c}_{\infty}$, and $\mathcal{I}_{\infty}^{\sharp}(\tau)$ defined,
respectively, by Equations~\eqref{eq3.108}, \eqref{eq3.109},
and~\eqref{eq3.87}. Hence, using Equations~\eqref{eq3.142}
and~\eqref{intermedg2}, one shows that
\begin{equation}
\sqrt{b(\tau)} \underset{\tau \to +\infty}{=}
\dfrac{\hat{u}_{0}(\tau)}{\hat{r}_{0}(\tau)} \exp (\Phi (\tau)),
\label{intermedg3}
\end{equation}
where $\Phi (\tau)$ is defined by Equation~\eqref{eq3.137}, whence, via
Equations~\eqref{eq3.129} and~\eqref{eq3.130}, and conditions~\eqref{eq3.11},
one arrives at Equation~\eqref{eq3.133}. \hfill $\qed$
\begin{eeeee} \label{newrem14}
\textsl{It is important to note that, as a consequence of
Equation~{\rm \eqref{intermedg3}},
\begin{equation}
\dfrac{\tau^{\frac{1}{6} \Im (a)} \hat{u}_{0}(\tau)}{\sqrt{b(\tau)} \,
\hat{r}_{0}(\tau)} \underset{\tau \to +\infty}{=} \mathcal{O}(1);
\label{eqtawonesixth}
\end{equation}
therefore, conditions~{\rm \eqref{newconds4}} are satisfied; moreover, via
conditions~{\rm \eqref{newconds3}} and the
estimate~{\rm \eqref{eqtawonesixth}}, it follows that
conditions~{\rm \eqref{eq3.103}} and~{\rm \eqref{eq3.104}} are also satisfied.
Making the following choice for the parameters $\varepsilon$, $\delta$, and
$\delta_{1}$,
\begin{equation*}
0 \! < \! 20 \delta \! < \! 10 \varepsilon \! < \! \delta_{1} \! < \! 1/3,
\end{equation*}
one shows that conditions~{\rm \eqref{newcconds2}} and~{\rm \eqref{newconds3}}
are valid uniformly for $0 \! < \! \Re (\widetilde{\nu} \! + \! 1) \! <
\! 1$. For the functions $\hat{r}_{0}(\tau)$, $\hat{u}_{0}(\tau)$, and
$h_{0}(\tau)$ defined, respectively, by Equations~{\rm \eqref{eq3.129}},
{\rm \eqref{eq3.130}}, and~{\rm \eqref{eq3.131}},
conditions~{\rm \eqref{eq3.11}} are satisfied.}
\end{eeeee}
\begin{bbbbb} \label{prop3.3.2}
Under the conditions of Proposition~{\rm \ref{cor3.3.2}}, the functions
$a(\tau)$, $b(\tau)$, $c(\tau)$, and $d(\tau)$, defining, via
Equations~\eqref{eq3.2}, the solution of the system of isomonodromy
deformations~\eqref{eq1.4}, have the following asymptotic representations:
\begin{equation}
\sqrt{-a(\tau)b(\tau)} \underset{\tau \to +\infty}{=} \dfrac{2 \alpha^{4}
\cos (\kappa (\tau) \! + \! \vartheta_{0}) \cos (\kappa (\tau) \! - \!
\vartheta_{0})}{\cos^{2}(\kappa (\tau))} \! + \! \mathcal{O}
(\tau^{-\frac{1}{3}+ \delta -\delta_{1}}), \label{eq3.146}
\end{equation}
\begin{align}
a(\tau)d(\tau) \underset{\tau \to +\infty}{=}& \, \dfrac{\mi 4 \alpha^{6}
\cos (\vartheta_{0}) \cos (\kappa (\tau) \! + \! \vartheta_{0}) \cos (\kappa
(\tau) \! - \! \vartheta_{0})}{\cos^{2}(\kappa (\tau))} \left(\cos
\vartheta_{0} \! + \! \dfrac{\sin^{2} \vartheta_{0}}{\cos (\kappa (\tau))
\cos (\kappa (\tau) \! - \! \vartheta_{0})} \right) \nonumber \\
+& \, \mathcal{O}(\tau^{-\frac{2}{3}+2 \delta}) \! + \! \mathcal{O}
(\tau^{-\frac{1}{3}+ \delta -\delta_{1}}), \label{eq3.147}
\end{align}
\begin{align}
b(\tau)c(\tau) \underset{\tau \to +\infty}{=}& \, -\dfrac{\mi 4 \alpha^{6}
\cos (\kappa (\tau) \! + \! \vartheta_{0}) \cos (\kappa (\tau) \! - \!
\vartheta_{0})}{\cos^{2}(\kappa (\tau))} \left(\dfrac{\cos (\vartheta_{0})
\sin^{2} \vartheta_{0}}{\cos (\kappa (\tau)) \cos (\kappa (\tau) \! - \!
\vartheta_{0})} \right. \nonumber \\
+&\left. \, \dfrac{\sin^{2} \vartheta_{0}}{\cos (\kappa (\tau) \! + \!
\vartheta_{0}) \cos (\kappa (\tau) \! - \! \vartheta_{0})} \! - \! \cos^{2}
\vartheta_{0} \right) \! + \! \mathcal{O}(\tau^{-\frac{2}{3}+2 \delta}) \! +
\! \mathcal{O}(\tau^{-\frac{1}{3}+ \delta - \delta_{1}}), \label{eq3.148}
\end{align}
\begin{align}
c(\tau)d(\tau) \underset{\tau \to +\infty}{=}& \, -4 \alpha^{4} \left(
\cos^{2} \vartheta_{0} \! + \! \dfrac{\cos \vartheta_{0}}{\cos (\kappa
(\tau)) \cos (\kappa (\tau) \! - \! \vartheta_{0})} \right) \! \left(
\dfrac{\cos (\vartheta_{0}) \sin^{2} \vartheta_{0}}{\cos (\kappa (\tau))
\cos (\kappa (\tau) \! - \! \vartheta_{0})} \right. \nonumber \\
+&\left. \, \dfrac{\sin^{2} \vartheta_{0}}{\cos (\kappa (\tau) \! + \!
\vartheta_{0}) \cos (\kappa (\tau) \! - \! \vartheta_{0})} \! - \! \cos^{2}
\vartheta_{0} \right) \! + \! \mathcal{O}(\tau^{-\frac{1}{3}+ \delta
-\delta_{1}}), \label{eq3.149}
\end{align}
where $b(\tau)$ is given in Equation~\eqref{eq3.133} $($see, also,
Equation~\eqref{eq3.137}$)$, $\kappa (\tau) \! := \! (\widetilde{\vartheta}
(\tau) \! + \! \vartheta_{0})/2$, and $\widetilde{\vartheta}(\tau)$,
$\vartheta_{0}$, and $\delta_{G} \! > \! 0$ are defined, respectively, by
Equations~\eqref{eq3.135}, \eqref{eq3.136}, and~\eqref{newconds5}.
\end{bbbbb}

\emph{Proof}. If $g_{ij}$, $i,j \! = \! 1,2$, are $\tau$-dependent, then any
functions whose asymptotics are given by
Equations~\eqref{eq3.129}--\eqref{eq3.133} satisfy conditions~\eqref{eq3.11},
\eqref{newconds1}, \eqref{newcconds2}, \eqref{eq3.103}, \eqref{eq3.104},
\eqref{newconds3}, and~\eqref{newconds4}; therefore, one can now use the
justification scheme suggested in \cite{a20} (see, also, \cite{a22}).
Equations~\eqref{eq3.146}--\eqref{eq3.149} are obtained {}from
Equations~\eqref{eq3.10}, \eqref{eq3.12}, \eqref{eq3.13}, and~\eqref{eq3.15}
via the direct substitution of Equations~\eqref{eq3.129} and~\eqref{eq3.130}.
\hfill $\qed$

{}From Proposition~\ref{prop1.2}, Equations~\eqref{eq3.2}, and
Equation~\eqref{eq3.10}, one shows that $u(\tau)$, the solution of the
degenerate third Painlev\'{e} equation~\eqref{eq1.1}, is given by
\begin{equation}
u(\tau) \! = \! \epsilon \tau^{1/3} \sqrt{-a(\tau)b(\tau)} = \! \dfrac{
\epsilon (\epsilon b)^{2/3}}{2} \tau^{1/3}(1 \! + \! \hat{u}_{0}(\tau)),
\quad \epsilon \! = \! \pm 1. \label{eq3.150}
\end{equation}
It was shown in \cite{a1} that, in terms of the function $h_{0}(\tau)$, the
Hamiltonian (cf. Equation~\eqref{eqh1}), $\mathcal{H}(\tau)$, corresponding
to $u(\tau)$, is given by
\begin{equation}
\mathcal{H}(\tau) \! = \! 3(\epsilon b)^{2/3} \tau^{1/3} \! - \! 4h_{0}(\tau)
\tau^{-1/3} \! + \! \dfrac{1}{2 \tau} \! \left(a \! - \! \dfrac{\mi}{2}
\right)^{2}. \label{eq3.151}
\end{equation}
Finally, the function $f(\tau)$ (cf. Equation~\eqref{eqf}) can be presented
in terms of $\hat{r}_{0}(\tau)$ with the help of Equations~\eqref{eq3.9},
\eqref{eq3.5}, \eqref{eq3.2}, and~\eqref{eq:ABCD}:
\begin{equation}
f(\tau) \! = \! -\dfrac{\mi (\epsilon b)^{1/3}}{2} \tau^{2/3} \! \left(
1 \! - \! \dfrac{\hat{r}_{0}(\tau)}{2} \right) \! - \! \dfrac{1}{2} \!
\left(\mi a \! + \! \dfrac{1}{2} \right). \label{eq:f}
\end{equation}
Thus, Proposition~\ref{cor3.3.2} implies the corresponding asymptotic results
for the functions $\hat{r}_{0}(\tau)$, $\mathcal{H}(\tau)$, and $f(\tau)$
stated in Theorem~\ref{theor2.1}.

It is important to note that asymptotics of $\hat{r}_{0}(\tau)$ can also be
reformulated as asymptotics of the logarithmic derivative of $u(\tau)$.
\begin{bbbbb} \label{prop:Alog-u-derivative}
The leading term of asymptotics of $(\ln u(\tau))^{\prime}$ coincides
with the logarithmic derivative of the leading term of asymptotics of
the function $u(\tau)$ given in Theorem~{\rm \ref{theor2.1}},
Equation~\eqref{eq:u-asympt-mult}$:$
\begin{equation}
\dfrac{\md \ln u(\tau)}{\md \tau} \underset{\tau \to \infty \me^{\mi \pi
\varepsilon_{1}}}{=} \dfrac{3 \sqrt{3} \, (-1)^{\varepsilon_{2}}(\epsilon
b)^{1/3} \tau^{-1/3} \cos \! \left(\frac{\vartheta}{2} \right)}{\sin \!
\left(\frac{\vartheta}{2} \right) \sin \! \left(\frac{\vartheta}{2} \! - \!
\vartheta_{0} \right) \sin \! \left(\frac{\vartheta}{2} \! + \! \vartheta_{0}
\right)} \! + \! \mathcal{O}(\tau^{-\frac{2}{3}+ \delta -\delta_{1}}),
\label{eq:Alog-u-derivative}
\end{equation}
where $\vartheta \! = \! \vartheta (\varepsilon_{1},\varepsilon_{2},\tau)$ is
defined by Equations~\eqref{eqth1.5}, \eqref{eqth1.3}, and~\eqref{eqth1.4},
and $\vartheta_{0}$ is given in Equation~\eqref{eq:theta0-def}.
\end{bbbbb}

\emph{Proof}. As in the main body of the paper, the case $\varepsilon_{1}
\! = \! \varepsilon_{2} \! = \! 0$ is considered in detail; in particular,
\begin{equation*}
\vartheta \! = \! \vartheta (0,0,\tau) \! = \! \widetilde{\vartheta}(\tau)
\! + \! \vartheta_{0} \! - \! \pi:
\end{equation*}
the remaining cases are proved analogously upon applying the transformation
$\mathscr{F}_{\varepsilon_{1},\varepsilon_{2}}$ introduced at the beginning
of Section~\ref{sec2}. {}From Equations~\eqref{eq:f} and~\eqref{eqf}, one
shows that
\begin{equation}
\dfrac{u^{\prime}(\tau)}{u(\tau)} \! = \! \dfrac{\mi b}{u(\tau)} \! - \!
\dfrac{2 \mi (\epsilon b)^{1/3}}{\tau^{1/3}} \! \left(1 \! - \! \dfrac{
\hat{r}_{0}(\tau)}{2} \right). \label{eq:Aloguprime}
\end{equation}
Substituting into Equation~\eqref{eq:Aloguprime} the asymptotics for
$u(\tau)$ and $\hat{r}_{0}(\tau)$ given in Equations~\eqref{eq:u-asympt-mult}
and~\eqref{eq3.129}, respectively, and using the identity
\begin{equation*}
2 \sin \! \left(\dfrac{\vartheta}{2} \right) \sin \! \left(
\dfrac{\vartheta}{2} \! - \! \vartheta_{0} \right) \! = \! \cos
\vartheta_{0} \! - \! \cos(\vartheta \! - \! \vartheta_{0}),
\end{equation*}
one arrives at
\begin{equation*}
\dfrac{u^{\prime}(\tau)}{u(\tau)} \underset{\tau \to +\infty}{=} \dfrac{\mi
(\epsilon b)^{1/3} \left(2 \sin^{3} \! \left(\frac{\vartheta}{2} \right)
\! + \! 5 \cos (\vartheta_{0}) \sin \! \left(\frac{\vartheta}{2} \! + \!
\vartheta_{0} \right) \! + \! \sin \! \left(\frac{\vartheta}{2} \! + \!
\vartheta_{0} \right) \cos (\vartheta \! - \! \vartheta_{0}) \right)}{
\tau^{1/3} \sin \! \left(\frac{\vartheta}{2} \right) \sin \! \left(
\frac{\vartheta}{2} \! - \! \vartheta_{0} \right) \sin \! \left(
\frac{\vartheta}{2} \! + \! \vartheta_{0} \right)} \! + \! \mathcal{O}
(\tau^{-\frac{2}{3}+ \delta -\delta_{1}}).
\end{equation*}
Consider, now, the numerator of the above fraction. One uses the identities
\begin{align*}
\sin \! \left(\dfrac{\vartheta}{2} \! + \! \vartheta_{0} \right) &= \sin
\! \left(\dfrac{\vartheta}{2} \right) \cos \vartheta_{0} \! + \! \cos
\! \left(\dfrac{\vartheta}{2} \right) \sin \vartheta_{0}, \\
\sin \! \left(\dfrac{\vartheta}{2} \! + \! \vartheta_{0} \right) \cos
(\vartheta \! - \! \vartheta_{0}) &= \dfrac{1}{2} \sin \! \left(\dfrac{3
\vartheta}{2} \right) \! - \! \dfrac{1}{2} \sin \! \left(\dfrac{\vartheta}{2}
\! - \! 2 \vartheta_{0} \right) \\
&= \dfrac{1}{2} \sin \! \left(\dfrac{3 \vartheta}{2} \right) \! - \!
\dfrac{1}{2}(2 \cos^2 \vartheta_{0} \! - \! 1) \sin \! \left(
\dfrac{\vartheta}{2} \right) \! + \! \sin (\vartheta_{0}) \cos
(\vartheta_{0}) \cos \! \left(\frac{\vartheta}{2} \right),
\end{align*}
to transform the second and third terms, respectively, of the summand:
regrouping, now, the terms in the numerator, one transforms it into the form
\begin{equation*}
\dfrac{1}{2} \! \left(4 \sin^{3} \! \left(\dfrac{\vartheta}{2} \right)
\! + \! \sin \! \left(\dfrac{3 \vartheta}{2} \right) \! + \! (8 \cos^2
\vartheta_{0} \! + \! 1) \sin \! \left(\dfrac{\vartheta}{2} \right) \right)
\! + \! 6 \sin (\vartheta_{0}) \cos (\vartheta_{0}) \cos \! \left(
\dfrac{\vartheta}{2} \right) \! = \! -\mi 3 \sqrt{3} \, \cos \! \left(
\dfrac{\vartheta}{2} \right),
\end{equation*}
where, in order to get the last equality, one uses the values of $\sin
\vartheta_{0}$ and $\cos \vartheta_{0}$ given in Equation~\eqref{eq3.136}.
Thus, asymptotics~\eqref{eq:Alog-u-derivative} is proved.

To finish the proof, one has to confirm that the logarithmic derivative of
the leading term of asymptotics~\eqref{eq:u-asympt-mult} coincides, modulo
$\mathcal{O}(\tau^{-\frac{2}{3}+ \delta -\delta_{1}})$, with the leading
term of asymptotics~\eqref{eq:Alog-u-derivative}. Taking the logarithmic
derivative of asymptotics~\eqref{eq:u-asympt-mult}, neglecting terms that
are $\mathcal{O}(\tau^{-1})$, and assuming that the derivatives of the errors
can also be neglected, one shows that
\begin{align*}
&\dfrac{\sqrt{3} \, (\epsilon b)^{1/3}}{\tau^{1/3}} \! \left(\dfrac{\cos \!
\left(\frac{\vartheta}{2} \! - \! \vartheta_{0} \right)}{\sin \! \left(
\frac{\vartheta}{2} \! - \! \vartheta_{0} \right)} \! + \! \dfrac{\cos \!
\left(\frac{\vartheta}{2} \! + \! \vartheta_{0} \right)}{\sin \! \left(
\frac{\vartheta}{2} \! + \! \vartheta_{0} \right)} \! - \! \dfrac{2 \cos \!
\left(\frac{\vartheta}{2} \right)}{\sin \! \left(\frac{\vartheta}{2} \right)}
\right) \\
&= \dfrac{\sqrt{3} \, (\epsilon b)^{1/3} \cos \! \left(\frac{\vartheta}{2}
\right)}{\tau^{1/3}} \! \left(\frac{2 \sin \! \left(\frac{\vartheta}{2}
\right)}{\sin \! \left(\frac{\vartheta}{2} \! - \! \vartheta_{0} \right)
\sin \! \left(\frac{\vartheta}{2} \! + \! \vartheta_{0} \right)} \! - \!
\frac{2}{\sin \! \left(\frac{\vartheta}{2} \right)} \right) \\
&= \dfrac{\sqrt{3} \, (\epsilon b)^{1/3} \cos \! \left(\frac{\vartheta}{2}
\right) \! \left(2 \sin^{2} \! \left(\frac{\vartheta}{2} \right) \! + \!
\cos (\vartheta) \! - \! \cos (2 \vartheta_{0}) \right)}{\tau^{1/3} \sin \!
\left(\frac{\vartheta}{2} \right) \sin \! \left(\frac{\vartheta}{2} \! - \!
\vartheta_{0} \right) \sin \! \left(\frac{\vartheta}{2} \! + \! \vartheta_{0}
\right)}.
\end{align*}
One completes the proof upon noting that $2 \sin^{2}(\vartheta /2) \!
+ \! \cos (\vartheta) \! - \! \cos (2 \vartheta_{0}) \! = \! 2 \sin^{2}
\vartheta_{0} \! = \! 3$. \hfill $\qed$

Consider, now, more carefully, the cheese-like domain
$\widetilde{\mathcal{D}}_{u}$ (cf. Equation~\eqref{eqnotat4}). As follows
{}from the leading term of asymptotics~\eqref{eq:u-asympt-mult} for the case
$\Re (\widetilde{\nu} \! + \! 1) \! \neq \! 1/2$, where $\widetilde{\nu} \!
:= \! \widetilde{\nu}(0,0)$, the countable sets $\mathcal{P}_{u}^{\delta}$ and
$\mathcal{Z}_{u}^{\delta}$ are empty, that is, $\widetilde{\mathcal{D}}_{u}$
coincides with the strip $\mathcal{D}$ (cf. Equation~\eqref{eqdofd}). In this
case, $\delta$ is vacuous: strictly speaking, one can improve some of the
error estimates in Section~\ref{sec3} and this section by assuming, in lieu
of conditions~\eqref{eq3.11}, the following conditions for the coefficients:
\begin{equation*}
\hat{r}_{0}(\tau) \underset{\tau \to +\infty}{\sim} \mathcal{O}
(\tau^{-\frac{1}{3}+ \delta}), \qquad \hat{u}_{0}(\tau) \underset{\tau \to
+\infty}{\sim} \mathcal{O}(\tau^{-\frac{1}{3}+ \delta}), \qquad \hat{h}_{0}
(\tau) \underset{\tau \to +\infty}{\sim} \mathcal{O}(1).
\end{equation*}
To eschew this, we can formally substitute $\delta \! = \! 0$ in our
calculations; surely, though, we will not achieve the best possible error
estimates, but they are sufficient for our purposes. Consider the definition
of $\delta_{G}$ given in Equation~\eqref{newconds5}, where $\delta^{G}_{12}
\! = \! \delta^{G}_{21} \! = \! \delta^{0}$ (cf. Equations~\eqref{newconds4}
and~\eqref{eqtawonesixth}). Put $\delta \! = \! 0$ in the formulae for
$\overset{\infty}{\delta_{12}}$ and $\overset{\infty}{\delta_{21}}$,
and $\delta^{0}$ defined, respectively, in Lemma~\ref{lem3.3.1}, and
Theorem~\ref{theor3.3.1} and Lemma~\ref{lem3.3.2}. Analyzing the resulting
system for $\delta_{G}$, one has to choose $\varepsilon \! \in \! (0,1/9)$
so that $\delta_{G}$ attains its maximum; the result reads:
\begin{enumerate}
\item[(i)] for $\Re (\widetilde{\nu} \! + \! 1) \! \in \! (0,1/2)$ and
$\varepsilon \! = \! \delta_{1}/(7 \! + \! 6 \Re (\widetilde{\nu} \! + \!
1))$,
\begin{equation*}
0 \! < \! \delta_{G} \! < \! \delta_{1} \! \left(\dfrac{1 \! + \! 2 \Re
(\widetilde{\nu} \! + \! 1)}{7 \! + \! 6 \Re (\widetilde{\nu} \! + \! 1)}
\right);
\end{equation*}
\item[(ii)] for $\Re (\widetilde{\nu} \! + \! 1) \! \in \! (1/2,1)$ and
$\varepsilon \! = \! \delta_{1}/(9 \! + \! 2 \Re (\widetilde{\nu} \! + \!
1))$,
\begin{equation*}
0 \! < \! \delta_{G} \! < \! \delta_{1} \! \left(\dfrac{3 \! - \! 2 \Re
(\widetilde{\nu} \! + \! 1)}{9 \! + \! 2 \Re (\widetilde{\nu} \! + \! 1)}
\right).
\end{equation*}
\end{enumerate}

Together with the justification scheme of \cite{a20}, this completes the
proof of Theorem~\ref{theor2.1}.

The structure of the cheese-like domain $\widetilde{\mathcal{D}}_{u}$ for $\Re
(\widetilde{\nu} \! + \! 1) \! = \! 1/2$ is more complicated: to study this
case, it is convenient to introduce the following definition.
\begin{aaaaa}
Define $\hat{\vartheta}$ as the leading term of asymptotics of $\vartheta$
$($cf. Equation~\eqref{eqth1.3}$)$, that is,
\begin{equation*}
\vartheta \! = \! \hat{\vartheta} \! + \! \mathcal{O}(\tau^{-\delta_G}
\ln \tau).
\end{equation*}
\end{aaaaa}
\begin{bbbbb} \label{propastaw}
Suppose that
\begin{equation*}
g_{11}(\varepsilon_{1},\varepsilon_{2})g_{12}(\varepsilon_{1},\varepsilon_{2})
g_{21}(\varepsilon_{1},\varepsilon_{2})g_{22}(\varepsilon_{1},\varepsilon_{2})
\! \neq \! 0, \qquad \Re \! \left(\dfrac{\mi}{2 \pi} \ln (g_{11}
(\varepsilon_{1},\varepsilon_{2})g_{22}(\varepsilon_{1},\varepsilon_{2}))
\right) \! = \! \dfrac{1}{2},
\end{equation*}
and the branch of the function $\ln (\pmb{\cdot})$ is chosen so that
$\Im (\ln (-g_{11}(\varepsilon_{1},\varepsilon_{2})g_{22}(\varepsilon_{1},
\varepsilon_{2}))) \! = \! 0$. For $m \! \in \! \mathbb{N}$, consider the
following transcendental equations:
\begin{equation}
\dfrac{\hat{\vartheta}}{2} \! = \! \pi m, \qquad \dfrac{\hat{\vartheta}}{2}
\! + \! \vartheta_{0} \! = \! \pi m, \qquad \dfrac{\hat{\vartheta}}{2} \! -
\! \vartheta_{0} \! = \! \pi m. \label{eqtrans}
\end{equation}
There exist unique solutions of these equations, $\tau^{\infty}_{as,m}$,
$\tau^{+}_{as,m}$, and $\tau^{-}_{as,m}$, respectively\footnote{
$\tau^{\infty}_{as,m}$, $\tau^{+}_{as,m}$, and $\tau^{-}_{as,m}$ are the
poles and zeroes of the leading term of asymptotics of $u(\tau)$ $($cf.
Equation~\eqref{eq:u-asympt-mult}$)$.}, the signs of the real parts of
which correspond to $(-1)^{\varepsilon_{1}};$ furthermore,
\begin{align}
\tau_{as,m}^{\infty} \underset{m \to \infty}{=}& \, \me^{\mi \pi
\varepsilon_{1}} \left(\dfrac{2 \pi (-1)^{\varepsilon_{2}}m}{3 \sqrt{3} \,
(\epsilon b)^{1/3}} \right)^{3/2} \left(1 \! - \! \dfrac{3 \varrho_{1}
(\varepsilon_{1},\varepsilon_{2})}{4 \pi} \dfrac{\ln m}{m} \! - \! \dfrac{3
\varrho_{2}(\varepsilon_{1},\varepsilon_{2})}{4 \pi} \dfrac{1}{m} \! + \!
\mathcal{O} \! \left(\dfrac{\ln^{2}m}{m^{2}} \right) \right),
\label{eqtauinfas} \\
\tau_{as,m}^{\pm} \underset{m \to \infty}{=}& \, \me^{\mi \pi \varepsilon_{1}}
\left(\dfrac{2 \pi (-1)^{\varepsilon_{2}}m}{3 \sqrt{3} \, (\epsilon b)^{1/3}}
\right)^{3/2} \left(1 \! - \! \dfrac{3 \varrho_{1}(\varepsilon_{1},
\varepsilon_{2})}{4 \pi} \dfrac{\ln m}{m} \! - \! \dfrac{3}{4 \pi} \left(
\varrho_{2}(\varepsilon_{1},\varepsilon_{2}) \! \pm \! 2 \vartheta_{0} \right)
\dfrac{1}{m} \! + \! \mathcal{O} \! \left(\dfrac{\ln^{2}m}{m^{2}} \right)
\right), \label{eqtaupmas}
\end{align}
where
\begin{equation}
\varrho_{1}(\varepsilon_{1},\varepsilon_{2}) \! := \! \dfrac{1}{2 \pi}
\ln (-g_{11}(\varepsilon_{1},\varepsilon_{2})g_{22}(\varepsilon_{1},
\varepsilon_{2})) \quad (\in \! \mathbb{R}), \label{eqnrowone}
\end{equation}
\begin{align}
\varrho_{2}(\varepsilon_{1},\varepsilon_{2}) :=& \, \varrho_{1}
(\varepsilon_{1},\varepsilon_{2}) \ln (24 \pi) \! + \! (-1)^{\varepsilon_{2}}
\Re (a) \ln (2 \! + \! \sqrt{3}) \! + \! \dfrac{\pi}{2} \! - \! \dfrac{1}{2}
\arg \! \left(\dfrac{g_{11}(\varepsilon_{1},\varepsilon_{2})g_{12}
(\varepsilon_{1},\varepsilon_{2})}{g_{21}(\varepsilon_{1},\varepsilon_{2})
g_{22}(\varepsilon_{1},\varepsilon_{2})} \right) \nonumber \\
-& \, \arg \! \left(\Gamma \! \left(\tfrac{1}{2} \! + \! \mi \varrho_{1}
(\varepsilon_{1},\varepsilon_{2}) \right) \right) \! + \! \mi \! \left(
(-1)^{\varepsilon_{2}} \Im (a) \ln (2 \! + \! \sqrt{3}) \! + \! \dfrac{1}{2}
\ln \! \left(\left\vert \dfrac{g_{11}(\varepsilon_{1},\varepsilon_{2})g_{12}
(\varepsilon_{1},\varepsilon_{2})}{g_{21}(\varepsilon_{1},\varepsilon_{2})
g_{22}(\varepsilon_{1},\varepsilon_{2})} \right\vert \right) \right),
\label{eqnrowtwo}
\end{align}
and $\vartheta_{0}$ is given in Equation~\eqref{eq:theta0-def}.
\end{bbbbb}

\emph{Proof}. Follows {}from a successive approximations argument applied to
Equations~\eqref{eqtrans}. \hfill $\qed$
\begin{bbbbb} \label{propdeltagee}
Let the conditions on the monodromy data be the same as those in
Proposition~{\rm \ref{propastaw}}. Suppose that $0 \! < \! 13 \delta \! < \!
\delta_{1} \! < \! 1/3$. Then
\begin{equation*}
0 \! < \! \delta \! < \! \delta_{G} \! < \! \dfrac{\delta_{1} \! - \! 8
\delta}{5}.
\end{equation*}
\end{bbbbb}

\emph{Proof}. Consider the definition of $\delta_{G}$ given in
Equation~\eqref{newconds5}, where $\delta^{0}$ is defined in
Theorem~\ref{theor3.3.1} and Lemma~\ref{lem3.3.2},
$\overset{\infty}{\delta_{12}}$ and $\overset{\infty}{\delta_{21}}$ are
defined in Lemma~\ref{lem3.3.1}, and $\delta^{G}_{12} \! = \!
\delta^{G}_{21} \! = \! \delta^{0}$ (cf. Equations~\eqref{newconds4}
and~\eqref{eqtawonesixth}). Under the assumptions stated in the Proposition,
one can choose $\varepsilon$ which appears in these equations to satisfy
the inequalities
\begin{equation*}
\dfrac{3 \delta}{2} \! < \! \varepsilon \! < \! \dfrac{\delta_{1} \! + \! 2
\delta}{10};
\end{equation*}
then $\delta_{G} \! = \! 2(\varepsilon \! - \! \delta)$. \hfill $\qed$

This completes the proof of Theorem~\ref{theorem2.3} modulo the distribution
of the cheese-holes in $\widetilde{\mathcal{D}}_{u}$.
\begin{ddddd} \label{th:A1}
For any given solution $u(\tau)$ with the monodromy data satisfying the
conditions~\eqref{eqtheorem2.21}, consider the cheese-like domain
$\mathcal{D}_{u}$ $($cf. Remark~{\rm \ref{remark2.3}}$)$. There exists $T \!
> \! 0$ such that, $\forall$ $\tau^{\pm}_{as,m},\tau^{\infty}_{as,m}$ with
$\lvert \tau^{\pm}_{as,m} \rvert \! > \! T$ and $\lvert \tau^{\infty}_{as,m}
\rvert \! > \! T$ in each cheese-hole centered at $\tau^{\pm}_{as,m}$ $($or
$\tau^{\infty}_{as,m})$, there exists one, and only one, zero $\tau^{\pm}_m$
$($or pole $\tau^{\infty}_{m})$ of $u(\tau)$ located in this cheese-hole:
\begin{equation*}
\lvert \tau^{\pm}_{as,m} \! - \! \tau^{\pm}_{m} \rvert \! < \! C \lvert
\tau^{\pm}_{as,m} \rvert^{\frac{1}{3}-\delta}, \qquad \quad \lvert
\tau^{\infty}_{as,m} \! - \! \tau^{\infty}_{m} \rvert \! < \! C \lvert
\tau^{\infty}_{as,m} \rvert^{\frac{1}{3}-\delta},
\end{equation*}
where $C \! > \! 0$.
\end{ddddd}

\emph{Proof}. As in the main body of the paper, the case $\varepsilon_{1}
\! = \! \varepsilon_{2} \! = \! 0$ is considered in detail; in particular,
$\vartheta \! = \! \vartheta (0,0,\tau) \! = \! \widetilde{\vartheta}
(\tau) \! + \! \vartheta_{0} \! - \! \pi$: the remaining cases are proved
analogously upon applying the transformation $\mathscr{F}_{\varepsilon_{1},
\varepsilon_{2}}$ introduced at the beginning of Section~\ref{sec2}.
The strategy here is to integrate both sides of
Equations~\eqref{eq:Alog-u-derivative} and~\eqref{eqtheorem2.33}, and the
equation resulting upon multiplying Equation~\eqref{eqtheorem2.34} by the
factor $2 \tau^{-1}$, around the boundary of a cheese-hole centered at
$\tau_{*} \! = \! \tau^{\pm}_{as,m}$ (or $\tau^{\infty}_{as,m})$. Hereafter,
the symbol $\oint$ means that the above-mentioned integration is taken in
the counter-clockwise sense.

Denote by $n_{+}$, $n_{-}$, and $n_{p}$, respectively, the numbers of zeroes
of types $\tau_{+}$, $\tau_{-}$, and poles $\tau_{p}$ of a solution $u(\tau)$
inside the cheese-hole around whose boundary the integration is performed. The
integrals on the left-hand sides of Equations~\eqref{eq:Alog-u-derivative}
and~\eqref{eqtheorem2.33}, and the equation resulting upon multiplying
Equation~\eqref{eqtheorem2.34} by the factor $2 \tau^{-1}$, can be calculated
via the Cauchy Residue Theorem; the results, respectively, read:
\begin{equation} \label{eq:Aleftoint}
2 \pi \mi (n_{+}+n_{-}-2n_{p}), \qquad 2 \pi \mi (n_{-}+n_{p}), \qquad 2 \pi
\mi (n_{-}-n_{p}).
\end{equation}
To integrate the right-hand sides, note that the error estimates in the
above-mentioned equations are all of the form $\mathcal{O}(\tau^{-\frac{1}{3}
-\upsilon})$, for some $\upsilon \! > \! 0$; therefore,
\begin{equation*}
\oint \mathcal{O} \! \left(\dfrac{1}{\tau^{\frac{1}{3}+ \upsilon}}
\right) \md \tau \underset{\tau \to \infty}{\sim} \mathcal{O} \! \left(
\dfrac{C}{\tau_{*}^{\upsilon + \delta}} \right) \underset{\tau_{*} \to
\infty} \rightarrow 0.
\end{equation*}
Thus, for some $T \! > \! 0$, the integral is less than $2 \pi/3$. The leading
terms of the right-hand sides of Equations~\eqref{eq:Alog-u-derivative}
and~\eqref{eqtheorem2.33}, and the equation resulting upon multiplying
Equation~\eqref{eqtheorem2.34} by the factor $2 \tau^{-1}$, are meromorphic
functions in the strip $\mathcal{D}$, where $\mathcal{D}$ is defined by
Equation~\eqref{eqdofd}; therefore, one can also calculate $\oint$ by using
the Cauchy Residue Theorem\footnote{It follows from the derivation in
Section~\ref{sec3} and this section that the error estimate $\mathcal{O}
(\tau^{-\delta_{G}} \ln \tau)$ in the ``phase'', $\vartheta$, is a meromorphic
function in $\mathcal{D}$; this fact, however, is not used in the proof.}. To
do so, notice that the non-trivial part of the leading terms (non-vanishing
under integration) has the form $\tau^{-1/3}F(\vartheta)$, where the function
$F$ is expressed in terms of trigonometric functions:
\begin{align}
\oint 2 \sqrt{3} \, (\epsilon b)^{1/3} \tau^{-1/3}F(\vartheta) \, \md \tau
&= \oint \hat{\vartheta}^{\prime}F(\vartheta) \, \md \tau \! + \! \mathcal{O}
\! \left(\oint \dfrac{F(\vartheta)}{\tau} \, \md \tau \right) \nonumber \\
&= \oint \hat{\vartheta}^{\prime}F(\hat{\vartheta}) \, \md \tau \! + \!
\mathcal{O} \! \left(\oint \hat{\vartheta}^{\prime} \max (\lvert F^{\prime}
(\hat{\vartheta}) \rvert) \tau^{-\delta_{G}} \ln \tau \, \md \tau \right)
\! + \! \mathcal{O} \! \left(\oint \dfrac{F(\vartheta)}{\tau} \, \md \tau
\right) \nonumber \\
&= \oint F(\hat{\vartheta}) \, \md \hat{\vartheta} \! + \! \mathcal{O}
(\tau^{\delta -\delta_{G}}_{*} \ln \tau_{*}) \! + \! \mathcal{O}
(\tau_{*}^{-2/3}), \label{eq:Aoint-estimate}
\end{align}
where use has been made of the fact that $\delta \! < \! \delta_{G}$, $\max
(\lvert F^{\prime}(\hat{\vartheta}) \rvert) \! = \! \mathcal{O}(\tau_{*}^{2
\delta})$, and $F(\vartheta) \! = \! \mathcal{O}(\tau_{*}^{\delta})$ in the
annular region of width $\mathcal{O}(\tau^{-\delta_{G}} \ln \tau)$ around the
circle of integration. Again, for $\tau_{*} \! > \! T$, each estimate in
Equation~\eqref{eq:Aoint-estimate} is less than $2\pi/3$; thus, for large
enough $\tau_{*}$, the modulus of the contribution of the error estimates to
$\oint$ is less than $2 \pi$. (Their contribution is actually equal to zero,
since the results of the integration of the right-hand sides should coincide
with the results of the integration of the corresponding left-hand sides,
and, according to Equations~\eqref{eq:Aleftoint}, should be equal be
$2 \pi \mi$ multipled by an integer.) Now, calculating the integrals
$\oint F(\hat{\vartheta}) \, \md \hat{\vartheta}$ in
Equations~\eqref{eq:Alog-u-derivative} and \eqref{eqtheorem2.33}, and the
equation resulting upon multiplying Equation~\eqref{eqtheorem2.34} by the
factor $2 \tau^{-1}$, one arrives at the following systems: (i) for $\tau_{*}
\! = \! \tau^{\infty}_{as,m}$,
\begin{equation*}
2 \pi \mi (n_{+}+n_{-}-2n_{p}) \! = \! -4 \pi \mi, \qquad 2 \pi \mi (n_{-}+
n_{p}) \! = \! 2 \pi \mi, \qquad 2 \pi \mi (n_{-}-n_{p}) \! = \! -2 \pi \mi;
\end{equation*}
(ii) for $\tau_{*} \! = \! \tau^{+}_{as,m}$,
\begin{equation*}
2 \pi \mi (n_{+}+n_{-}-2n_{p}) \! = \! 2 \pi \mi, \qquad 2 \pi \mi (n_{-}
+n_{p}) \! = \! 0, \qquad 2 \pi \mi (n_{-}-n_{p}) \! = \! 0;
\end{equation*}
and (iii) for $\tau_{*} \! = \! \tau^{-}_{as,m}$,
\begin{equation*}
2 \pi \mi (n_{+}+n_{-}-2n_{p}) \! = \! 2 \pi \mi, \qquad 2 \pi \mi (n_{-}
+n_{p}) \! = \! 2 \pi \mi, \qquad 2 \pi \mi (n_{-}-n_{p}) \! = \! 2 \pi \mi.
\end{equation*}
One completes the proof upon solving these systems. \hfill $\qed$
\begin{fffff} \label{corayywone}
Define $\hat{\tau}_{as,m}^{\infty}$ and $\hat{\tau}_{as,m}^{\pm}$,
respectively, as the leading terms of asymptotics of $\tau^{\infty}_{as,m}$
and $\tau^{\pm}_{as,m}$ $($cf. Proposition~{\rm \ref{propastaw})}, that is,
\begin{equation*}
\tau^{\infty}_{as,m} \! =: \! \hat{\tau}_{as,m}^{\infty} \! + \! \mathcal{O}
\! \left(\dfrac{\ln^{2}m}{\sqrt{m}} \right), \quad \quad \tau^{\pm}_{as,m}
\! =: \! \hat{\tau}_{as,m}^{\pm} \! + \! \mathcal{O} \! \left(
\dfrac{\ln^{2}m}{\sqrt{m}} \right).
\end{equation*}
Under the conditions of Theorem~{\rm \ref{th:A1}},
\begin{equation*}
\tau^{\infty}_{m} \underset{m \to \infty}{=} \hat{\tau}_{as,m}^{\infty} \! +
\! \mathcal{O} \! \left(m^{\frac{1}{2}- \frac{3 \delta}{2}} \right), \quad
\quad \tau^{\pm}_{m} \underset{m \to \infty}{=} \hat{\tau}_{as,m}^{\pm} \! +
\! \mathcal{O} \! \left(m^{\frac{1}{2}- \frac{3 \delta}{2}} \right),
\end{equation*}
where $0 \! < \! \delta \! < \! 1/39$.
\end{fffff}

\emph{Proof}. As follows {}from Proposition~\ref{propdeltagee}, $0 \! < \!
\delta \! < \! 1/39$; thus, $\ln^{2}m/\sqrt{m} \! \ll \! m^{\frac{1}{2}-
\frac{3 \delta}{2}}$ as $m \! \to \! \infty$. \hfill $\qed$

This completes the proof of Theorem~\ref{theorem2.2}, which is equivalent to
the specification of the cheese-holes in $\widetilde{\mathcal{D}}_{u}$.

\vspace*{1.30cm}

\pmb{\large Acknowledgments}

The authors were supported, in part, by a College of Charleston (CofC)
Mathematics Department Summer Research Award.
\clearpage
\section*{Appendix A: Poles and Zeroes} \label{apndxca}
\setcounter{section}{1}
\setcounter{equation}{0}
\renewcommand{\thesection}{\Alph{section}}
\renewcommand{\theequation}{\Alph{section}.\arabic{equation}}
In this Appendix some basic information concerning the poles and zeroes of the
solution $u(\tau)$ of Equation~\eqref{eq1.1} and the associated functions
$\mathcal{H}(\tau)$ and $f(\tau)$ (cf. Equations~\eqref{eqh1} and~\eqref{eqf},
respectively) are presented.

Let $\tau_{\infty}$ be a pole of $u(\tau)$; then, its Laurent expansion reads
\begin{equation}
u(\tau) \! = \! -\dfrac{\tau_{\infty}}{4 \epsilon(\tau \! - \!
\tau_{\infty})^2} \! + \! a_{0} \! + \! \sum_{k=1}^{\infty}a_{k}(\tau \!
- \!\tau_{\infty})^{k}, \quad \epsilon \! = \! \pm 1, \label{equpole}
\end{equation}
where $a_{0}$ is a parameter, and the remaining coefficients, $a_{k}$, are
recursively and uniquely defined; for example,
\begin{gather*}
a_{1} \! = \! -\dfrac{a_{0}}{\tau_{\infty}}, \qquad \, a_{2} \! = \!
\dfrac{ab}{5 \tau_\infty} \! - \! \dfrac{12 \epsilon a_{0}^{2}}{5
\tau_{\infty}} \! + \! \dfrac{9a_{0}}{10 \tau_{\infty}^{2}}, \qquad \, a_{3}
\! = \! -\dfrac{8ab}{45 \tau_{\infty}^{2}} \! + \! \dfrac{24 \epsilon
a_{0}^{2}}{5 \tau_{\infty}^{2}} \! - \! \dfrac{4a_{0}}{5 \tau_{\infty}^{3}},
\\
a_{4} \! = \! -\dfrac{\epsilon b^{2}}{7 \tau_{\infty}} \! + \! \dfrac{32
a_{0}^{3}}{7 \tau_{\infty}^{2}} \! - \! \dfrac{4 \epsilon a_{0}ab}{7
\tau_{\infty}^{2}} \! + \! \dfrac{10ab}{63 \tau_{\infty}^{3}} \! - \!
\dfrac{47 \epsilon a_{0}^{2}}{7 \tau_{\infty}^{3}} \! + \! \dfrac{5a_{0}}{7
\tau_{\infty}^{4}}, \\
a_{5} \! = \! \dfrac{\epsilon b^{2}}{35 \tau_{\infty}^{2}} \! - \! \dfrac{96
a_{0}^{3}}{7 \tau_{\infty}^{3}} \! + \! \dfrac{124 \epsilon a_{0}ab}{105
\tau_{\infty}^{3}} \! - \! \dfrac{ab}{7 \tau_{\infty}^{4}} \! + \! \dfrac{57
\epsilon a_{0}^{2}}{7 \tau_{\infty}^{4}} \! - \! \dfrac{9a_{0}}{14
\tau_{\infty}^{5}}.
\end{gather*}
The Laurent expansion of the associated Hamiltonian, $\mathcal{H}(\tau)$,
reads
\begin{equation}
\mathcal{H}(\tau) \! - \! \dfrac{1}{2 \tau} \! \left(a \! - \! \dfrac{\mi}{2}
\right)^{2} \! = \! \dfrac{1}{\tau \!- \! \tau_{\infty}} \! + \! H_{0} \! +
\! \sum_{k=1}^{\infty}H_{k}(\tau \! - \! \tau_{\infty})^{k}, \label{eqhpole}
\end{equation}
where the coefficients, $H_{k}$, can be uniquely determined via Equation~(37)
of \cite{a1}; for example,
\begin{gather*}
H_{0} \!= \! 12 \epsilon a_{0}, \qquad \, H_{1} \! = \! -\dfrac{8 \epsilon
a_{0}}{\tau_{\infty}}, \qquad \, H_{2} \! = \! \dfrac{2 \mi \epsilon b}{
\tau_{\infty}} \! + \! \dfrac{6 \epsilon a_{0}}{\tau_{\infty}^{2}}, \\
H_{3} \! = \! -\dfrac{16a \epsilon b}{15 \tau_{\infty}^{2}} \! - \! \dfrac{16
a_{0}^{2}}{5 \tau_{\infty}^{2}} \! - \! \dfrac{24 \epsilon a_{0}}{5
\tau_{\infty}^{3}}, \qquad \, H_{4} \! = \! \dfrac{8 \mi ba_{0}}{\tau_{
\infty}^{2}} \! + \! \dfrac{8a \epsilon b}{9 \tau_{\infty}^{3}} \! + \!
\dfrac{8a_{0}^{2}}{\tau_{\infty}^{3}} \! + \! \dfrac{4 \epsilon a_{0}}{
\tau_{\infty}^{4}}, \\
H_{5} \! = \! \dfrac{52b^{2}}{35 \tau_{\infty}^{2}} \! - \! \dfrac{128ab
a_{0}}{35 \tau_{\infty}^{3}} \! + \! \dfrac{128 \epsilon a_{0}^{3}}{35
\tau_{\infty}^{3}} \! - \! \dfrac{8 \mi ba_{0}}{\tau_{\infty}^{3}} \! - \!
\dfrac{16a \epsilon b}{21 \tau_{\infty}^{4}} \! - \! \dfrac{468a_{0}^2}{35
\tau_{\infty}^{4}} \! - \! \dfrac{24 \epsilon a_{0}}{7 \tau_{\infty}^{5}}.
\end{gather*}
The function $f(\tau)$ also has a first-order pole at $\tau \! = \!
\tau_{\infty}$:
\begin{equation}
f(\tau) \! = \! -\dfrac{\tau_{\infty}}{2(\tau \! - \! \tau_{\infty})} \! + \!
f_{0} \! + \! \sum_{k=1}^{\infty}f_{k}(\tau \! - \! \tau_{\infty})^{k},
\label{eqhpole1}
\end{equation}
where the coefficients, $f_{k}$, can be uniquely determined via
Equation~\eqref{eqf}; for example,
\begin{gather*}
f_{0} \! = \! -\dfrac{\mi a}{2} \! - \! \dfrac{3}{4}, \qquad \, f_{1} \! = \!
-2a_{0} \epsilon, \qquad \, f_{2} \! = \! \mi \epsilon b \! + \! \dfrac{a_{0}
\epsilon}{\tau_{\infty}}, \\
f_{3} \! = \! \dfrac{8a_{0}^{2}}{5 \tau_{\infty}} \! - \! \dfrac{4a \epsilon
b}{5 \tau_{\infty}} \! + \! \dfrac{\mi \epsilon b}{\tau_{\infty}} \! - \!
\dfrac{3a_{0} \epsilon}{5 \tau_{\infty}^{2}}, \qquad \, f_{4} \! = \!
\dfrac{4 \mi a_{0}b}{\tau_{\infty}} \! - \! \dfrac{12a_{0}^{2}}{5
\tau_{\infty}^{2}} \! + \! \dfrac{4a \epsilon b}{45 \tau_{\infty}^{2}}
\! + \! \dfrac{2a_{0} \epsilon}{5 \tau_{\infty}^{3}}, \\
f_{5} \! = \! \dfrac{6b^{2}}{7 \tau_{\infty}} \! - \! \dfrac{64 \epsilon
a_{0}^{3}}{35 \tau_{\infty}^{2}} \! - \! \dfrac{48aba_{0}}{35 \tau_{
\infty}^{2}} \! - \! \dfrac{4a \epsilon b}{63 \tau_{\infty}^{3}} \! + \!
\dfrac{94a_{0}^{2}}{35 \tau_{\infty}^{3}} \! - \! \dfrac{2 \epsilon a_{0}}{7
\tau_{\infty}^{4}}.
\end{gather*}

There are two types of first-order zeroes of $u(\tau)$, denoted $\tau_{s}$,
$s \! = \! \pm 1$, which differ by their Taylor expansions:
\begin{equation}
u(\tau) \! = \! \sum_{k=1}^{\infty}b_{k}^{s}(\tau \! - \! \tau_{s})^{k},
\end{equation}
where, for example,
\begin{gather*}
b_{1}^{s} \! = \! \mi sb, \qquad b_{2}^{s} \! = \! -\dfrac{b}{\tau_{s}}
\! \left(a \! - \! \dfrac{\mi s}{2} \right), \qquad b_{3}^{s} \! \in \!
\mathbb{C}, \qquad b_{4}^{s} \! = \! \dfrac{1}{2 \tau_{s}} \! \left(4
\epsilon b^{2} \! + \! \mi sab_{3}^{s} \! - \! b_{3}^{s} \right), \\
b_{5}^{s} \! = \! \dfrac{1}{20 b \tau_{s}^{2}} \! \left(7bb_{3}^{s}
\! - \! 8b^{3} \epsilon \! - \! 6 \mi s(b_{3}^{s})^{2} \tau_{s}^{2}
\! + \! 48 \mi sab^{3} \epsilon \! - \! 10 \mi sabb_{3}^{s} \right).
\end{gather*}
These two types of first-order zeroes are ``analytically'' distinguished by
the Hamiltonian function; indeed, at $\tau \! = \! \tau_{+}$, the function
$\mathcal{H}(\tau)$ is holomorphic,
\begin{equation*}
\mathcal{H}(\tau) \! = \! H_{0}^{+} \! + \! \sum_{k=1}^{\infty}H_{k}^{+}
(\tau \! - \! \tau_{+})^{k},
\end{equation*}
where the first few coefficients are
\begin{gather*}
H_{0}^{+} \! = \! -\dfrac{3 \mi b_{3}^{+} \tau_{+}}{2b} \! + \! \dfrac{1}{2
\tau_{+}} \! \left(a \! - \! \dfrac{\mi}{2} \right)^{2} \! + \! \dfrac{\mi}{
\tau_{+}} \! \left(a \! - \! \dfrac{\mi}{2} \right), \qquad H_{1}^{+} \! =
\! \dfrac{1}{2 \tau_{+}^{2}} \left(a \! - \! \frac{\mi}{2} \right)^{2}, \\
H_{2}^{+} \! = \! -\dfrac{3b_{3}^{+}}{2b \tau_{+}} \! \left(a \! - \! \dfrac{
\mi}{2} \right) \! + \! \dfrac{1}{2 \tau_{+}^{3}} \! \left(a \! - \! \dfrac{
\mi}{2} \right)^{2} \! - \! \dfrac{\mi}{\tau_{+}^{3}} \! \left(a \! - \!
\dfrac{\mi}{2} \right)^{3},
\end{gather*}
while, at $\tau \! = \! \tau_{-}$, the function $\mathcal{H}(\tau)$ has a
first-order pole,
\begin{equation}
\mathcal{H}(\tau) \! - \! \dfrac{1}{2 \tau} \! \left(a \! - \! \dfrac{\mi}{2}
\right)^{2} \! = \! \dfrac{1}{\tau \!- \! \tau_{-}} \! + \! H_{0}^{-} \! + \!
\sum_{k=1}^{\infty}H_{k}^{-}(\tau \! - \! \tau_{-})^{k}, \label{eqhzero-}
\end{equation}
where, for example,
\begin{equation*}
H_{0}^{-} \! = \! \dfrac{3 \mi b_{3}^{-} \tau_{-}}{2b}, \qquad \, H_{1}^{-}
\! = \! -\dfrac{\mi b_{3}^{-}}{b}, \qquad \, H_{2}^{-} \! = \! -\dfrac{2
\mi \epsilon b}{\tau_{-}} \! + \! \dfrac{3 \mi b_{3}^{-}}{4b \tau_{-}}.
\end{equation*}
The zeroes of $u(\tau)$ are also distinguished by the function $f(\tau)$:
\begin{equation}
f(\tau) \! = \! \sum_{k=1}^{\infty}f_{k}^{+}(\tau \! - \! \tau_{+})^{k},
\label{eq:fzero+}
\end{equation}
where the first two coefficients are
\begin{equation*}
f_{1}^{+} \! = \! \dfrac{4a^{2} \! + \! 1}{8 \tau_{+}} \! - \! \dfrac{3 \mi
b_{3}^{+} \tau_{+}}{4b}, \qquad \, f_{2}^{+} \! = \! -2 \mi \epsilon b \!
- \! \dfrac{3b_{3}^{+}}{4b} \! \left(a \! - \! \dfrac{\mi}{2} \right) \!
- \! \dfrac{\mi (4a^{2} \! + \! 1)}{8 \tau_{+}^{2}} \! \left(a \! - \!
\frac{\mi}{2} \right),
\end{equation*}
and
\begin{equation}
f(\tau) \! = \! \dfrac{\tau_{-}}{2(\tau \! - \! \tau_{-})} \! + \! f_{0}^{-}
\! + \! \sum_{k=1}^{\infty}f_{k}^{-}(\tau \! - \! \tau_{-})^{k},
\label{eq:fzero-}
\end{equation}
where, for example,
\begin{equation*}
f_{0}^{-} \! = \! -\dfrac{\mi}{2} \! \left(a \! + \! \frac{\mi}{2} \right),
\qquad \, f_{1}^{-} \! = \! \dfrac{\mi b_{3}^{-} \tau_{-}}{4b}, \qquad \,
f_{2}^{-} \! = \! \mi \epsilon b \! - \! \dfrac{\mi b_{3}^{-}}{8b}.
\end{equation*}

As follows {}from the definitions of the associated functions $\mathcal{H}
(\tau)$ and $f(\tau)$, they have poles only at the poles and $\tau_{-}$-zeroes
of the function $u(\tau)$.
\section*{Appendix B: Comparison of Asymptotic Results} \label{apndxb}
\setcounter{section}{2}
\setcounter{equation}{0}
\setcounter{ddddd}{0}
\renewcommand{\thesection}{\Alph{section}}
\renewcommand{\theequation}{\Alph{section}.\arabic{equation}}
The ranges of the validity for asymptotics of $u(\tau)$ obtained in Part I
\cite{a1} of our studies and in this paper overlap: this enables us to resolve
the ambiguity in the choice of sign for the function $\hat{r}_{0}(\tau)$
discussed in the proof of Proposition~\ref{cor3.3.2}.

The asymptotics of $u(\tau)$ for $\lvert \Re (\widetilde{\nu}(\varepsilon_{1},
\varepsilon_{2}) \! + \! 1) \rvert \! < 1/6$, $\varepsilon_{1},\varepsilon_{2}
\! = \! 0,\pm 1$, are obtained in Theorem~3.1 of \cite{a1}, whilst the
asymptotics of $u(\tau)$ stated in Theorem~\ref{theor2.1} of this paper are
applicable for $\Re (\widetilde{\nu}(\varepsilon_{1},\varepsilon_{2}) \!
+ \! 1) \! \in \! (0,1) \setminus \lbrace 1/2 \rbrace$. In fact, if one is
concerned only with the leading exponent of $\cosh (\pmb{\cdot})$ in the
asymptotics of $u(\tau)$ presented in Theorem~3.1 of \cite{a1}, then
the range of its validity can be extended to $\lvert \Re (\widetilde{\nu}
(\varepsilon_{1},\varepsilon_{2}) \! + \! 1) \rvert \! < \! 1/2$. Analogously,
only the leading exponent in the expansion of $\sin^{-2}(\pmb{\cdot})$ in the
asymptotics of $u(\tau)$ obtained in Theorem~\ref{theor2.1} of this paper
is ``larger'' than the error correction term in case $\Re (\widetilde{\nu}
(\varepsilon_{1},\varepsilon_{2}) \! + \! 1) \! > \! 0$. Therefore, the
corresponding leading exponents of the asymptotics of $u(\tau)$ in
Theorem~3.1 of \cite{a1} and Theorem~\ref{theor2.1} of this paper should
coincide, provided
\begin{equation}
0 \! < \! \Re (\widetilde{\nu}(\varepsilon_{1},\varepsilon_{2}) \! + \! 1) \!
< \! \dfrac{1}{2}. \label{newrest}
\end{equation}
Recall the main asymptotic (as $\tau \! \to \! \infty)$ result for $u(\tau)$
obtained in Part I \cite{a1}\footnote{In order to facilitate the comparison,
the parameter $\varepsilon$ which appears in Theorem~3.1 of \cite{a1} has been
changed to $\epsilon$.}:
\begin{ddddd}[{\rm Theorem~3.1 \cite{a1}}] \label{theoa1}
Let $\varepsilon_{1},\varepsilon_{2} \! = \! 0,\pm 1$, $\epsilon b \! = \!
\lvert \epsilon b \rvert \me^{\mi \pi \varepsilon_{2}}$, and $u(\tau)$ be
a solution of the degenerate third Painlev\'{e} equation~\eqref{eq1.1}
corresponding to the monodromy data $(a,s^{0}_{0},s^{\infty}_{0},
s^{\infty}_{1},g_{11},g_{12},g_{21},\linebreak[4]
g_{22})$. Suppose that
\begin{equation*}
g_{11}(\varepsilon_{1},\varepsilon_{2})g_{12}(\varepsilon_{1},\varepsilon_{2})
g_{21}(\varepsilon_{1},\varepsilon_{2})g_{22}(\varepsilon_{1},\varepsilon_{2})
\! \neq \! 0, \qquad \left\vert \Re \! \left(\dfrac{\mi}{2 \pi}
\ln (g_{11}(\varepsilon_{1},\varepsilon_{2})g_{22}(\varepsilon_{1},
\varepsilon_{2})) \right) \right\vert \! < \! \dfrac{1}{6}.
\end{equation*}
Then, $\exists$ $\delta \! > \! 0$ such that $u(\tau)$ has the asymptotic
expansion
\begin{align}
u(\tau) \underset{\tau \to \infty \me^{\mi \pi \varepsilon_{1}}}{=}& \,
\dfrac{(-1)^{\varepsilon_{1}} \epsilon \sqrt{\lvert \epsilon b \rvert}}{
3^{1/4}} \! \left(\sqrt{\dfrac{\vartheta (\tau)}{12}}+ \! \sqrt{\widetilde{\nu}
(\varepsilon_{1},\varepsilon_{2}) \! + \! 1} \, \me^{\frac{3 \pi \mi}{4}}
\cosh \! \left(\vphantom{M^{M^{M}}} \mi \vartheta (\tau) \! + \!
(\widetilde{\nu}(\varepsilon_{1},\varepsilon_{2}) \! + \! 1) \ln \vartheta
(\tau) \right. \right. \nonumber \\
+&\left. \left. \, z(\varepsilon_{1},\varepsilon_{2}) \! + \!
o(\tau^{-\delta}) \vphantom{M^{M^{M}}} \right) \!
\vphantom{M^{M^{M^{M^{M^{M}}}}}} \right), \label{eqtheoa1.2}
\end{align}
where
\begin{equation*}
\vartheta (\tau) \! := \! 3 \sqrt{3} \, \lvert \epsilon b \rvert^{1/3}
\lvert \tau \rvert^{2/3}, \quad \quad \widetilde{\nu}(\varepsilon_{1},
\varepsilon_{2}) \! + \! 1 \! := \! \dfrac{\mi}{2 \pi} \ln (g_{11}
(\varepsilon_{1},\varepsilon_{2})g_{22}(\varepsilon_{1},\varepsilon_{2})),
\end{equation*}
\begin{align*}
z(\varepsilon_{1},\varepsilon_{2}) :=& \, \dfrac{1}{2 \pi} \ln 2 \pi \!
- \! \dfrac{\mi \pi}{2} \! - \! \dfrac{3 \pi \mi}{2}(\widetilde{\nu}
(\varepsilon_{1},\varepsilon_{2}) \! + \! 1) \! + \! (-1)^{\varepsilon_{2}}
\mi a \ln (2 \! + \! \sqrt{3}) \! + \! (\widetilde{\nu}(\varepsilon_{1},
\varepsilon_{2}) \! + \! 1) \ln 12 \\
-& \, \ln \! \left(\omega (\varepsilon_{1},\varepsilon_{2})
\sqrt{\widetilde{\nu}(\varepsilon_{1},\varepsilon_{2}) \! + \! 1} \, \Gamma
(\widetilde{\nu}(\varepsilon_{1},\varepsilon_{2}) \! + \! 1) \right),
\end{align*}
with
\begin{equation}
\omega (\varepsilon_{1},\varepsilon_{2}) \! := \! \dfrac{g_{12}
(\varepsilon_{1},\varepsilon_{2})}{g_{22}(\varepsilon_{1},\varepsilon_{2})},
\label{eqtheoa1.5}
\end{equation}
and $\Gamma (\pmb{\cdot})$ is the gamma function.

Let $\mathcal{H}(\tau)$ be the Hamiltonian function defined in
Equation~\eqref{eqh1} corresponding to the function $u(\tau)$ given above.
Then,
\begin{align*}
\mathcal{H}(\tau) \! - \! \dfrac{1}{2 \tau} \! \left(a \! - \!
\dfrac{(-1)^{\varepsilon_{2}} \mi}{2} \right)^{2} \underset{\tau \to \infty
\me^{\mi \pi \varepsilon_{1}}}{=}& \, 3(\epsilon b)^{2/3} \tau^{1/3}
\! + \! 2 \lvert \epsilon b \rvert^{1/3} \tau^{-1/3} \! \left(\left(
a \! - \! \dfrac{(-1)^{\varepsilon_{2}} \mi}{2} \right) \right. \\
-&\left. \, \mi 2 \sqrt{3} \, (\widetilde{\nu}(\varepsilon_{1},
\varepsilon_{2}) \! + \! 1) \! + \! o(\tau^{-\delta}) \right).
\end{align*}
\end{ddddd}
Now, taking into account condition~\eqref{newrest}, one shows that
asymptotics~\eqref{eqtheoa1.2} implies
\begin{align}
u(\tau) \underset{\tau \to \infty \me^{\mi \pi \varepsilon_{1}}}{=}& \,
\dfrac{(-1)^{\varepsilon_{2}} \epsilon \lvert \epsilon b \rvert^{2/3}}{2}
\lvert \tau \rvert^{1/3} \! + \! \dfrac{(-1)^{\varepsilon_{1}} \epsilon
\sqrt{\lvert \epsilon b \rvert}}{2 \cdot 3^{1/4}} \lvert \tau \rvert^{
\frac{2}{3} \Re (\widetilde{\nu}(\varepsilon_{1},\varepsilon_{2})+1)}
\me^{\Theta (\tau)} \! + \! \mathcal{O}(\tau^{-\frac{2}{3} \Re
(\widetilde{\nu}(\varepsilon_{1},\varepsilon_{2})+1)}) \nonumber \\
+& \, o(\tau^{\frac{2}{3} \Re (\widetilde{\nu}(\varepsilon_{1},
\varepsilon_{2})+1)- \delta}), \label{eqtheoa1.6}
\end{align}
where
\begin{align}
\Theta (\tau) :=& \, \mi 3 \sqrt{3} \, \lvert \epsilon b \rvert^{1/3}
\lvert \tau \rvert^{2/3} \! + \! \mi \Im (\widetilde{\nu}(\varepsilon_{1},
\varepsilon_{2}) \! + \! 1) \ln \lvert \tau \rvert^{2/3} \! + \!
(\widetilde{\nu}(\varepsilon_{1},\varepsilon_{2}) \! + \! 1) \ln
(36 \sqrt{3} \, \lvert \epsilon b \rvert^{1/3}) \nonumber \\
+& \, (-1)^{\varepsilon_{2}} \mi a \ln (2 \! + \! \sqrt{3}) \! + \!
\dfrac{\mi \pi}{4} \! - \! \dfrac{3 \pi \mi}{2}(\widetilde{\nu}
(\varepsilon_{1},\varepsilon_{2}) \! + \! 1) \! - \! \ln \! \left(
\dfrac{\omega (\varepsilon_{1},\varepsilon_{2}) \Gamma (\widetilde{\nu}
(\varepsilon_{1},\varepsilon_{2}) \! + \! 1)}{\sqrt{2 \pi}} \right).
\label{eqtheoa1.9}
\end{align}
Using the Euler formula for $\sin (\pmb{\cdot})$ to expand
asymptotics~\eqref{eqth1.2} (assuming condition~\eqref{newrest}) in
Theorem~\ref{theor2.1} of this work, one arrives, again, at
asymptotics~\eqref{eqtheoa1.6}, but with the error correction
\begin{equation*}
\mathcal{O}(\tau^{\frac{2}{3} \Re (\widetilde{\nu}(\varepsilon_{1},
\varepsilon_{2})+1)-\delta_{G}} \ln \tau) \! + \! \mathcal{O}
(\tau^{\frac{4}{3} \Re (\widetilde{\nu}(\varepsilon_{1},\varepsilon_{2})+1)
-\frac{1}{3}});
\end{equation*}
furthermore, there is a discrepancy in the definition of the parameter
$\omega (\varepsilon_{1},\varepsilon_{2})$ (cf. Equation~\eqref{eqtheoa1.9}):
instead of \eqref{eqtheoa1.5}, one finds that
\begin{equation}
\omega (\varepsilon_{1},\varepsilon_{2}) \! := \! g_{11}(\varepsilon_{1},
\varepsilon_{2})g_{12}(\varepsilon_{1},\varepsilon_{2}). \label{eqtheoa1.10}
\end{equation}
In fact, \eqref{eqtheoa1.10} is the correct definition of $\omega
(\varepsilon_{1},\varepsilon_{2})$ (see the discussion below); therefore,
\eqref{eqtheoa1.10} should be used in lieu of the definition of $\omega
(\varepsilon_{1},\varepsilon_{2})$ in Theorem~\ref{theoa1} and in Theorem~3.1
of \cite{a1}.

The root of this discrepancy is related with the fact that the right factor
of the second term in Equation~(79) of \cite{a1} should be changed, namely,
\begin{equation}
\ln \! \left(\dfrac{3(\varepsilon b)^{1/6}}{\sqrt{2}} \right) \! \to \! \ln
\! \left(\dfrac{3(\varepsilon b)^{1/6} \me^{-\mi \pi}}{\sqrt{2}} \right).
\label{eqtheoa1.11}
\end{equation}
As a consequence of \eqref{eqtheoa1.11}, corresponding changes should be made
to the following formulae in \cite{a1}: (i) to the right-hand sides of
Equations~(96), (98), and~(99) one must add the term $\mi \pi (\nu \! + \!
1)$; (ii) to the argument of the exponential function on the right-hand side
of Equation~(101) one must add the term $-2 \pi \mi (\nu \! + \! 1)$; (iii) to
the right-hand side of the formula for $z_{n}$ given in Corollary~4.3.1 one
must add the term $2 \pi \mi (\nu \! + \! 1)$; and (iv) the definition of
the parameter $\widehat{\omega}(\varepsilon_{1},\varepsilon_{2})$ appearing
in Theorem~A.1 should be changed\footnote{There are a few innocuous misprints
in \cite{a1}, which do not, however, affect the final results: (i) the factor
$3 \sqrt{3}- \! 2$ appearing on the right-hand sides of Equations~(76)
and~(92) should read $3(\sqrt{3}- \! 1)$; (ii) the factor $\mi
(\varepsilon b)^{1/3} \tau^{2/3}$ appearing on the right-hand side of
Equation~(98) and in the argument of the exponential function on the
right-hand side of Equation~(101) should read $\tfrac{3 \mi}{2}
(\varepsilon b)^{1/3} \tau^{2/3}$; and (iii) the factor $\mi (3 \sqrt{3}-
\! 1)(\varepsilon b)^{1/3} \tau^{2/3}$ appearing on the right-hand side of
Equation~(99) should read $\mi (3 \sqrt{3}- \! \tfrac{3}{2})(\varepsilon
b)^{1/3} \tau^{2/3}$.} to $\widehat{\omega}(\varepsilon_{1},\varepsilon_{2})
\! := \! \widehat{g}_{11}(\varepsilon_{1},\varepsilon_{2}) \widehat{g}_{12}
(\varepsilon_{1},\varepsilon_{2})$.

Making analogous expansions for the $\cot (\pmb{\cdot})$ functions in
Equation~\eqref{eqth1.6} of Theorem~\ref{theor2.1}, one obtains an
asymptotic formula for $\mathcal{H}(\tau)$ that coincides with the one
presented in Theorem~\ref{theoa1}.

For the other choice of $\hat{r}_{0}(\tau)$ (cf. the proof of
Proposition~\ref{cor3.3.2}) one would obtain the same
asymptotics~\eqref{eqtheoa1.6}--\eqref{eqtheoa1.9} but with the change
$\Theta (\tau) \! \to \! \Theta (\tau) \! + \! \mi \pi$.
\section*{Appendix C: Asymptotics for Imaginary $\tau$} \label{apndxc}
\setcounter{section}{3}
\setcounter{equation}{0}
\setcounter{ddddd}{0}
\setcounter{fffff}{0}
\setcounter{eeeee}{0}
\renewcommand{\thesection}{\Alph{section}}
\renewcommand{\theequation}{\Alph{section}.\arabic{equation}}
Here, asymptotics as $\tau \! \to \! \pm \mi \infty$ of the functions
$u(\tau)$, $\mathcal{H}(\tau)$, and $f(\tau)$ are presented. These results
are obtained by applying transformations~6.2.2 (changing\footnote{In
transformation~6.2.2, $\tau \! \to \! \tau$, that is, $\tau_{n} \!
= \! \tau_{o}$.} $a \! \to \! -a)$ and~6.2.3 (changing\footnote{In
transformation~6.2.3, $a \! \to \! a$, that is, $a_{n} \! = \! a_{o}$.}
$\tau \! \to \! \mi \tau)$ given in Section~6 of \cite{a1} to the asymptotic
results (for $\varepsilon_{1} \! = \! \varepsilon_{2} \! = \! 0)$ stated in
Theorems~\ref{theor2.1}--\ref{theorem2.3} of this paper. For this purpose,
it is convenient to introduce the auxiliary mapping\footnote{There
is a misprint on page~1202 (Appendix) of \cite{a1}: for items~(1)
and~(4) in the definition of the auxiliary mapping
$\widehat{\mathscr{F}}_{\varepsilon_{1},\varepsilon_{2}}$, the
change $a \! \to \! -a$ should be made everywhere.}
\begin{gather}
\widehat{\mathscr{F}}_{\varepsilon_{1},\varepsilon_{2}} \colon \mathscr{M}
\! \to \! \mathscr{M},\qquad
(a,s^{0}_{0},s^{\infty}_{0},s^{\infty}_{1},g_{11},
g_{12},g_{21},g_{22}) \mapsto \nonumber\\
((-1)^{1+\varepsilon_{2}}a,s^{0}_{0},
\widehat{s}^{\infty}_{0}(\varepsilon_{1},\varepsilon_{2}),
\widehat{s}^{\infty}_{1}(\varepsilon_{1},\varepsilon_{2}),\widehat{g}_{11}
(\varepsilon_{1},\varepsilon_{2}),\widehat{g}_{12}(\varepsilon_{1},
\varepsilon_{2}),\widehat{g}_{21}(\varepsilon_{1},\varepsilon_{2}),
\widehat{g}_{22}(\varepsilon_{1},\varepsilon_{2})), \varepsilon_{1} \! = \!
\pm 1, \varepsilon_{2} \! = \! 0,\pm 1, \nonumber
\end{gather}
which is equivalent to the formulae
for the monodromy data in transformations~6.2.2 and~6.2.3 in Section~6 of
\cite{a1} (with $l \! = \! \varepsilon_{1}$ and $p \! = \! \varepsilon_{2})$.
Define\footnote{$s^{0}_{0}(\varepsilon_{1},\varepsilon_{2}) \! = \!
s^{0}_{0}$.}:
\begin{enumerate}
\item[(1)] $\widehat{\mathscr{F}}_{-1,0}$: $\widehat{s}^{\infty}_{0}(-1,0) \!
= \! s^{\infty}_{1} \me^{\frac{3 \pi a}{2}}$, $\widehat{s}^{\infty}_{1}(-1,0)
\! = \! s^{\infty}_{0} \me^{\frac{\pi a}{2}}$, $\widehat{g}_{11}(-1,0) \! =
\! -g_{22} \me^{\frac{3 \pi a}{4}}$, $\widehat{g}_{12}(-1,0) \! = \! -(g_{21}
\! + \! s^{\infty}_{0}g_{22}) \me^{-\frac{3 \pi a}{4}}$, $\widehat{g}_{21}
(-1,0) \!= \! (s^{0}_{0}g_{22}-g_{12}) \me^{\frac{3 \pi a}{4}}$, and
$\widehat{g}_{22}(-1,0) \! = \! ( s^{0}_{0}(g_{21} \! + \! s^{\infty}_{0}g_{22})
-g_{11} \! -\! s^{\infty}_{0}g_{12}) \me^{-\frac{3 \pi a}{4}}$;
\item[(2)] $\widehat{\mathscr{F}}_{-1,-1}$: $\widehat{s}^{\infty}_{0}(-1,-1)
\! = \! s^{\infty}_{0} \me^{-\frac{\pi a}{2}}$, $\widehat{s}^{\infty}_{1}
(-1,-1) \! = \! s^{\infty}_{1} \me^{\frac{\pi a}{2}}$, $\widehat{g}_{11}
(-1,-1) \! = \! -\mi g_{21} \me^{-\frac{\pi a}{4}}$, $\widehat{g}_{12}(-1,-1)
\! = \! -\mi g_{22} \me^{\frac{\pi a}{4}}$, $\widehat{g}_{21}(-1,-1) \! =
\! -\mi (g_{11} \! - \! s^{0}_{0}g_{21}) \me^{-\frac{\pi a}{4}}$, and
$\widehat{g}_{22}(-1,-1) \! = \! -\mi (g_{12} \! - \! s^{0}_{0}g_{22})
\me^{\frac{\pi a}{4}}$;
\item[(3)] $\widehat{\mathscr{F}}_{-1,1}$: $\widehat{s}^{\infty}_{0}(-1,1)
\! = \! s^{\infty}_{0} \me^{-\frac{\pi a}{2}}$, $\widehat{s}^{\infty}_{1}
(-1,1) \! = \! s^{\infty}_{1} \me^{\frac{\pi a}{2}}$, $\widehat{g}_{11}(-1,1)
\! = \! g_{11} \me^{-\frac{\pi a}{4}}$, $\widehat{g}_{12}(-1,1) \! = \! g_{12}
\me^{\frac{\pi a}{4}}$, $\widehat{g}_{21}(-1,1) \! = \! g_{21}
\me^{-\frac{\pi a}{4}}$, and $\widehat{g}_{22}(-1,1) \! = \! g_{22}
\me^{\frac{\pi a}{4}}$;
\item[(4)] $\widehat{\mathscr{F}}_{1,0}$: $\widehat{s}^{\infty}_{0}(1,0) \!
= \! s^{\infty}_{1} \me^{\frac{\pi a}{2}}$, $\widehat{s}^{\infty}_{1}(1,0)
\! = \! s^{\infty}_{0} \me^{\frac{3 \pi a}{2}}$, $\widehat{g}_{11}(1,0) \!
= \! -\mi g_{12} \me^{\frac{\pi a}{4}}$, $\widehat{g}_{12}(1,0) \! = \!
-\mi (g_{11} \! + \! s^{\infty}_{0}g_{12}) \me^{-\frac{\pi a}{4}}$,
$\widehat{g}_{21}(1,0) \! = \! -\mi g_{22} \me^{\frac{\pi a}{4}}$, and
$\widehat{g}_{22}(1,0) \! = \! -\mi (g_{21} \! + \! s^{\infty}_{0}g_{22})
\me^{-\frac{\pi a}{4}}$;
\item[(5)] $\widehat{\mathscr{F}}_{1,-1}$: $\widehat{s}^{\infty}_{0}(1,-1)
\! = \! s^{\infty}_{0} \me^{\frac{\pi a}{2}}$, $\widehat{s}^{\infty}_{1}(1,-1)
\! = \! s^{\infty}_{1} \me^{-\frac{\pi a}{2}}$, $\widehat{g}_{11}(1,-1) \! =
\! g_{11} \me^{\frac{\pi a}{4}}$, $\widehat{g}_{12}(1,-1) \! = \! g_{12}
\me^{-\frac{\pi a}{4}}$, $\widehat{g}_{21}(1,-1) \! = \! g_{21}
\me^{\frac{\pi a}{4}}$, and $\widehat{g}_{22}(1,-1) \! = \! g_{22}
\me^{-\frac{\pi a}{4}}$;
\item[(6)] $\widehat{\mathscr{F}}_{1,1}$: $\widehat{s}^{\infty}_{0}(1,1) \!
= \! s^{\infty}_{0} \me^{\frac{\pi a}{2}}$, $\widehat{s}^{\infty}_{1}(1,1) \!
= \! s^{\infty}_{1} \me^{-\frac{\pi a}{2}}$, $\widehat{g}_{11}(1,1) \! = \!
\mi (g_{21} \! + \! s^{0}_{0}g_{11}) \me^{\frac{\pi a}{4}}$, $\widehat{g}_{12}
(1,1) \! = \! \mi (g_{22} \! + \! s^{0}_{0}g_{12}) \me^{-\frac{\pi a}{4}}$,
$\widehat{g}_{21}(1,1) \! = \! \mi g_{11} \me^{\frac{\pi a}{4}}$, and
$\widehat{g}_{22}(1,1) \! = \! \mi g_{12} \me^{-\frac{\pi a}{4}}$.
\end{enumerate}
\begin{ddddd} \label{theor2.3}
Let $\varepsilon_{1} \! = \! \pm 1$, $\varepsilon_{2} \! = \! 0,\pm 1$,
$\epsilon b \! = \! \vert \epsilon b \vert \me^{\mi \pi \varepsilon_{2}}$,
and $u(\tau)$ be a solution of Equation~\eqref{eq1.1} corresponding to the
monodromy data $(a,s^{0}_{0},s^{\infty}_{0},s^{\infty}_{1},g_{11},g_{12},
g_{21},g_{22})$. Suppose that
\begin{equation}
\widehat{g}_{11}(\varepsilon_{1},\varepsilon_{2}) \widehat{g}_{12}
(\varepsilon_{1},\varepsilon_{2}) \widehat{g}_{21}(\varepsilon_{1},
\varepsilon_{2}) \widehat{g}_{22}(\varepsilon_{1},\varepsilon_{2}) \! \neq
\! 0, \qquad \Re (\widehat{\nu}(\varepsilon_{1},\varepsilon_{2}) \! + \! 1)
\! \in \! (0,1) \setminus \left\lbrace \tfrac{1}{2} \right\rbrace,
\label{eqth4.1}
\end{equation}
where
\begin{equation}
\widehat{\nu}(\varepsilon_{1},\varepsilon_{2}) \! + \! 1 \! := \!
\dfrac{\mi}{2 \pi} \ln (\widehat{g}_{11}(\varepsilon_{1},\varepsilon_{2})
\widehat{g}_{22}(\varepsilon_{1},\varepsilon_{2})). \label{eqth4.4}
\end{equation}
Then there exist $\delta_{G}$ satisfying, for $0 \! < \! \Re (\widehat{\nu}
(\varepsilon_{1},\varepsilon_{2}) \! + \! 1) \! < \! \tfrac{1}{2}$, the
inequality
\begin{equation*}
0 \! < \! \delta_{G} \! < \! \dfrac{1}{3} \! \left(\dfrac{1 \! + \! 2 \Re
(\widehat{\nu}(\varepsilon_{1},\varepsilon_{2}) \! + \! 1)}{7 \! + \! 6
\Re (\widehat{\nu}(\varepsilon_{1},\varepsilon_{2}) \! + \! 1)} \right),
\end{equation*}
and, for $\tfrac{1}{2} \! < \! \Re (\widehat{\nu}(\varepsilon_{1},
\varepsilon_{2}) \! + \! 1) \! < \! 1$, the inequality
\begin{equation*}
0 \! < \! \delta_{G} \! < \! \dfrac{1}{3} \! \left(\dfrac{3 \! - \! 2 \Re
(\widehat{\nu}(\varepsilon_{1},\varepsilon_{2}) \! + \! 1)}{9 \! + \! 2
\Re (\widehat{\nu}(\varepsilon_{1},\varepsilon_{2}) \! + \! 1)} \right),
\end{equation*}
such that $u(\tau)$ has the asymptotic expansion
\begin{align}
u(\tau) \underset{\tau \to \infty \me^{\frac{\mi \pi \varepsilon_{1}}{2}}}{=}&
\, \dfrac{\me^{-\frac{\mi \pi \varepsilon_{1}}{2}} \epsilon (\epsilon
b)^{2/3}}{2}(\me^{-\frac{\mi \pi \varepsilon_{1}}{2}} \tau)^{1/3} \! \left(
1 \! - \! \dfrac{3}{2 \sin^{2}(\tfrac{1}{2} \beta (\varepsilon_{1},
\varepsilon_{2},\tau))} \right) \label{eqth4.2} \\
\underset{\tau \to \infty \me^{\frac{\mi \pi \varepsilon_{1}}{2}}}{=}& \,
\dfrac{\me^{-\frac{\mi \pi \varepsilon_{1}}{2}} \epsilon (\epsilon
b)^{2/3}}{2}(\me^{-\frac{\mi \pi \varepsilon_{1}}{2}} \tau)^{1/3} \dfrac{
\sin (\tfrac{1}{2} \beta (\varepsilon_{1},\varepsilon_{2},\tau) \! - \!
\vartheta_{0}) \sin (\tfrac{1}{2} \beta (\varepsilon_{1},\varepsilon_{2},\tau)
\! + \! \vartheta_{0})}{\sin^{2}(\tfrac{1}{2} \beta (\varepsilon_{1},
\varepsilon_{2},\tau))} , \label{eqth4.7}
\end{align}
where
\begin{align}
\beta (\varepsilon_{1},\varepsilon_{2},\tau) :=& \, \widehat{\phi}(\tau) \! -
\! \mi \left((\widehat{\nu}(\varepsilon_{1},\varepsilon_{2}) \! + \! 1) \! -
\! \dfrac{1}{2} \right) \ln \widehat{\phi}(\tau) \! - \! \mi \left((\widehat{
\nu}(\varepsilon_{1},\varepsilon_{2}) \! + \! 1) \! - \! \dfrac{1}{2} \right)
\ln 12 \! + \! (-1)^{1+ \varepsilon_{2}}a \ln (2 \! + \! \sqrt{3}) \nonumber
\\
+& \, \dfrac{\pi}{4} \! - \! \dfrac{3 \pi}{2}(\widehat{\nu}(\varepsilon_{1},
\varepsilon_{2}) \! + \! 1) \! + \! \mi \ln \! \left(\dfrac{\widehat{g}_{11}
(\varepsilon_{1},\varepsilon_{2}) \widehat{g}_{12}(\varepsilon_{1},
\varepsilon_{2}) \Gamma (\widehat{\nu}(\varepsilon_{1},\varepsilon_{2}) \! +
\! 1)}{\sqrt{2 \pi}} \right) \! + \! \mathcal{O}(\tau^{-\delta_{G}} \ln \tau),
\label{eqth4.3}
\end{align}
with
\begin{equation}
\widehat{\phi}(\tau) \! = \! 3 \sqrt{3} \, (-1)^{\varepsilon_{2}}
(\epsilon b)^{1/3}(\me^{-\frac{\mi \pi \varepsilon_{1}}{2}} \tau)^{2/3},
\label{eqth4.8}
\end{equation}
$\vartheta_{0}$ given in Equation~\eqref{eq:theta0-def}, and
$\Gamma (\pmb{\cdot})$ the Euler gamma function {\rm \cite{a24}}.

Let $\mathcal{H}(\tau)$ be the Hamiltonian function defined by
Equation~\eqref{eqh1} corresponding to the function $u(\tau)$ given
above. Then $\mathcal{H}(\tau)$ has the asymptotic expansion
\begin{align}
\mathcal{H}(\tau) \underset{\tau \to \infty \me^{\frac{\mi \pi
\varepsilon_{1}}{2}}}{=}& \, \me^{-\frac{\mi \pi \varepsilon_{1}}{2}} \left(
\vphantom{M^{M^{M^{M^{M}}}}} 3(\epsilon b)^{2/3}(\me^{-\frac{\mi \pi
\varepsilon_{1}}{2}} \tau)^{1/3} \! - \! \mi (-1)^{\varepsilon_{2}}4 \sqrt{3}
\, (\epsilon b)^{1/3}(\me^{-\frac{\mi \pi \varepsilon_{1}}{2}} \tau)^{-1/3}
\left(\vphantom{M^{M^{M}}} (\widehat{\nu}(\varepsilon_{1},\varepsilon_{2})+1)
\right. \right. \nonumber \\
-&\left. \left. \, \dfrac{1}{2} \! + \! \dfrac{1}{2 \sqrt{3}} \left(\mi
(-1)^{1+ \varepsilon_{2}}a \! + \! \dfrac{1}{2} \right) \! + \! \dfrac{\mi}{4}
\cot \left(\tfrac{1}{2} \beta (\varepsilon_{1},\varepsilon_{2},\tau) \right)
\! + \! \dfrac{\mi}{4} \cot \left(\tfrac{1}{2} \beta (\varepsilon_{1},
\varepsilon_{2},\tau) \! - \! \vartheta_{0} \right) \right. \right.
\nonumber \\
+&\left. \left. \, \mathcal{O}(\tau^{-\delta_{G}}) \vphantom{M^{M^{M}}}
\right) \vphantom{M^{M^{M^{M^{M}}}}} \! \right).
\label{eqth4.5}
\end{align}
The function $f(\tau)$ defined by Equation~\eqref{eqf} has the following
asymptotics:
\begin{equation}
f(\tau) \underset{\tau \to \infty \me^{\frac{\mi \pi \varepsilon_{1}}{2}}}{=}
-\dfrac{(-1)^{\varepsilon_{2}}(\epsilon b)^{1/3}}{2}(\me^{-\frac{\mi \pi
\varepsilon_{1}}{2}} \tau)^{2/3} \! \left(\mi \! + \! \dfrac{3}{\sqrt{2} \,
\sin (\tfrac{1}{2} \beta (\varepsilon_{1},\varepsilon_{2},\tau)) \sin
(\tfrac{1}{2} \beta (\varepsilon_{1},\varepsilon_{2},\tau) \! - \!
\vartheta_{0})} \right). \label{eqthasymptfmain}
\end{equation}
\end{ddddd}
\begin{eeeee} \label{appcremc1}
\textsl{Define the strip $($in the $\widehat{\phi}$-plane$)$
\begin{equation}
\widehat{\mathcal{D}} \! := \! \left\{\mathstrut \tau \! \in \! \mathbb{C}
\, \colon \, \Re (\widehat{\phi}(\tau)) \! > \! c_{1}, \, \lvert \Im
(\widehat{\phi}(\tau)) \rvert \! < \! c_{2} \right\}, \label{eqndhatremc1}
\end{equation}
where $\widehat{\phi}(\tau)$ is given in Equation~\eqref{eqth4.8}, and $c_{1},
c_{2} \! > \! 0$ are parameters. The asymptotics of $u(\tau)$, $\mathcal{H}
(\tau)$, and $f(\tau)$ presented in Theorem~{\rm \ref{theor2.3}} are actually
valid in the strip domain $\widehat{\mathcal{D}}$.}
\end{eeeee}
\begin{ddddd} \label{theorr2.4}
Let $\varepsilon_{1} \! = \! \pm 1$, $\varepsilon_{2} \! = \! 0,\pm 1$,
$\epsilon b \! = \! \vert \epsilon b \vert \me^{\mi \pi \varepsilon_{2}}$,
and $u(\tau)$ be a solution of Equation~\eqref{eq1.1} corresponding to the
monodromy data $(a,s^{0}_{0},s^{\infty}_{0},s^{\infty}_{1},g_{11},g_{12},
g_{21},g_{22})$. Suppose that
\begin{equation}
\widehat{g}_{11}(\varepsilon_{1},\varepsilon_{2}) \widehat{g}_{12}
(\varepsilon_{1},\varepsilon_{2}) \widehat{g}_{21}(\varepsilon_{1},
\varepsilon_{2}) \widehat{g}_{22}(\varepsilon_{1},\varepsilon_{2}) \!
\neq \! 0, \qquad \Re \! \left(\dfrac{\mi}{2 \pi} \ln (\widehat{g}_{11}
(\varepsilon_{1},\varepsilon_{2}) \widehat{g}_{22}(\varepsilon_{1},
\varepsilon_{2})) \right) \! = \! \dfrac{1}{2}. \label{eqth6.1}
\end{equation}
Let the branch of the function $\ln (\pmb{\cdot})$ be chosen\footnote{The
second condition of Equations~\eqref{eqth6.1} suggests that this branch
of $\ln (\pmb{\cdot})$ exists.} such that $\Im (\ln (-\widehat{g}_{11}
(\varepsilon_{1},\varepsilon_{2}) \widehat{g}_{22}(\varepsilon_{1},
\varepsilon_{2}))) \! = \! 0$. Define
\begin{equation}
\widehat{\varrho}_{1}(\varepsilon_{1},\varepsilon_{2}) \! := \!
\dfrac{1}{2 \pi} \ln (-\widehat{g}_{11}(\varepsilon_{1},\varepsilon_{2})
\widehat{g}_{22}(\varepsilon_{1},\varepsilon_{2})) \quad (\in \! \mathbb{R}).
\label{eqth6.2}
\end{equation}
Then $\exists$ $\delta \! \in \! (0,1/39)$ such that the function $u(\tau)$
has, for all large enough $m \! \in \! \mathbb{N}$, second-order poles,
$\widehat{\tau}_{m}^{\infty}$, accumulating at the point at infinity,
\begin{equation}
\widehat{\tau}_{m}^{\infty} \underset{m \to \infty}{=} \me^{\frac{\mi
\pi \varepsilon_{1}}{2}} \left(\dfrac{2 \pi (-1)^{\varepsilon_{2}}m}{3
\sqrt{3} \, (\epsilon b)^{1/3}} \right)^{3/2} \left(1 \! - \! \dfrac{3
\widehat{\varrho}_{1}(\varepsilon_{1},\varepsilon_{2})}{4 \pi}
\dfrac{\ln m}{m} \! - \! \dfrac{3 \widehat{\varrho}_{2}(\varepsilon_{1},
\varepsilon_{2})}{4 \pi} \dfrac{1}{m} \right) \! + \! \mathcal{O} \!
\left(m^{\frac{1}{2}- \frac{3 \delta}{2}} \right), \label{eqth6.3}
\end{equation}
where
\begin{align}
\widehat{\varrho}_{2}(\varepsilon_{1},\varepsilon_{2}) :=& \,
\widehat{\varrho}_{1}(\varepsilon_{1},\varepsilon_{2}) \ln (24 \pi) \! +
\! (-1)^{1+ \varepsilon_{2}} \Re (a) \ln (2 \! + \! \sqrt{3}) \! + \!
\dfrac{\pi}{2} \! - \! \dfrac{1}{2} \arg \! \left(\dfrac{\widehat{g}_{11}
(\varepsilon_{1},\varepsilon_{2}) \widehat{g}_{12}(\varepsilon_{1},
\varepsilon_{2})}{\widehat{g}_{21}(\varepsilon_{1},\varepsilon_{2})
\widehat{g}_{22}(\varepsilon_{1},\varepsilon_{2})} \right) \nonumber \\
-& \, \arg \! \left(\Gamma \! \left(\tfrac{1}{2} \! + \! \mi
\widehat{\varrho}_{1}(\varepsilon_{1},\varepsilon_{2}) \right) \right) \! +
\! \mi \! \left((-1)^{1+ \varepsilon_{2}} \Im (a) \ln (2 \! + \! \sqrt{3})
\! + \! \dfrac{1}{2} \ln \! \left(\left\vert \dfrac{\widehat{g}_{11}
(\varepsilon_{1},\varepsilon_{2}) \widehat{g}_{12}(\varepsilon_{1},
\varepsilon_{2})}{\widehat{g}_{21}(\varepsilon_{1},\varepsilon_{2})
\widehat{g}_{22}(\varepsilon_{1},\varepsilon_{2})} \right\vert \right)
\right);
\label{eqth6.5}
\end{align}
furthermore, the function $u(\tau)$ has, for all large enough $m \! \in
\! \mathbb{N}$, a pair of first-order zeroes, $\widehat{\tau}_{m}^{\pm}$,
accumulating at the point at infinity,
\begin{equation}
\widehat{\tau}_{m}^{\pm} \underset{m \to \infty}{=} \me^{\frac{\mi \pi
\varepsilon_{1}}{2}} \left(\dfrac{2 \pi (-1)^{\varepsilon_{2}}m}{3
\sqrt{3} \, (\epsilon b)^{1/3}} \right)^{3/2} \left(1 \! - \!
\dfrac{3 \widehat{\varrho}_{1}(\varepsilon_{1},\varepsilon_{2})}{4 \pi}
\dfrac{\ln m}{m} \! - \! \dfrac{3}{4 \pi}(\widehat{\varrho}_{2}
(\varepsilon_{1},\varepsilon_{2}) \! \pm \! 2 \vartheta_{0}) \dfrac{1}{m}
\right) \! + \! \mathcal{O} \! \left(m^{\frac{1}{2}- \frac{3 \delta}{2}}
\right), \label{eqth6.6}
\end{equation}
where $\vartheta_{0}$ given in Equation~\eqref{eq:theta0-def}.
\end{ddddd}
\begin{eeeee} \label{appcremc2}
\textsl{To present asymptotics of $u(\tau)$, $\mathcal{H}(\tau)$, and
$f(\tau)$ outside of neighborhoods of poles and zeroes, introduce the
cheese-like domain, $\widehat{\mathcal{D}}_{u}$, for a solution $u(\tau):$
\begin{equation*}
\widehat{\mathcal{D}}_{u} \! := \! \left\{\mathstrut \tau \! \in \!
\widehat{\mathcal{D}} \colon \, \lvert \widehat{\phi}(\tau) \! - \!
\widehat{\phi}(\widehat{\tau}_{m}^{\kappa}) \rvert \! \geqslant \! C
\lvert \widehat{\tau}_{m}^{\kappa} \rvert^{-\delta} \, \right\},
\end{equation*}
where the strip domain $\widehat{\mathcal{D}}$ is defined by
Equation~\eqref{eqndhatremc1}, $\widehat{\phi}(\tau)$ is given in
Equation~\eqref{eqth4.8}, $C \! > \! 0$ is a parameter, $\kappa \! = \!
\infty,\pm$ $(\widehat{\tau}_{m}^{\kappa}$ are the poles and zeroes introduced
in Theorem~{\rm \ref{theorr2.4}}$)$, and $0 \! < \! \delta \! < \! 1/39$.}
\end{eeeee}
\begin{ddddd} \label{theorrr2.4}
Let $\varepsilon_{1} \! = \! \pm 1$, $\varepsilon_{2} \! = \! 0,\pm 1$,
$\epsilon b \! = \! \vert \epsilon b \vert \me^{\mi \pi \varepsilon_{2}}$,
and $u(\tau)$ be a solution of Equation~\eqref{eq1.1} corresponding to the
monodromy data $(a,s^{0}_{0},s^{\infty}_{0},s^{\infty}_{1},g_{11},g_{12},
g_{21},g_{22})$. Suppose that conditions~\eqref{eqth6.1} are valid, the branch
of $\ln (\pmb{\cdot})$ is chosen as in Theorem~{\rm \ref{theorr2.4}}, and
$\widehat{\varrho}_{1}(\varepsilon_{1},\varepsilon_{2})$ is defined by
Equation~\eqref{eqth6.2}. Then there exist $\delta,\delta_{G} \! \in \!
\mathbb{R}_{+}$ satisfying the inequalities
\begin{equation*}
0 \! < \! \delta \! < \! \dfrac{1}{39}, \quad \quad 0 \! < \! \delta
\! < \! \delta_{G} \! < \! \dfrac{1}{15} \! - \! \dfrac{8 \delta}{5},
\end{equation*}
such that $u(\tau)$ has the asymptotic expansion
\begin{equation}
u(\tau) \underset{\underset{\scriptstyle \tau \in \widehat{\mathcal{D}}_{
u}}{\tau \to \infty \me^{\frac{\mi \pi \varepsilon_{1}}{2}}}}{=} \dfrac{
\me^{-\frac{\mi \pi \varepsilon_{1}}{2}} \epsilon (\epsilon b)^{2/3}}{2}
(\me^{-\frac{\mi \pi \varepsilon_{1}}{2}} \tau)^{1/3} \dfrac{\sin
(\tfrac{1}{2} \widehat{\beta}(\varepsilon_{1},\varepsilon_{2},\tau) \!
- \! \vartheta_{0}) \sin (\tfrac{1}{2} \widehat{\beta}(\varepsilon_{1},
\varepsilon_{2},\tau) \! + \! \vartheta_{0})}{\sin^{2}(\tfrac{1}{2}
\widehat{\beta}(\varepsilon_{1},\varepsilon_{2},\tau))}, \label{eqth5.1}
\end{equation}
where
\begin{align}
\widehat{\beta}(\varepsilon_{1},\varepsilon_{2},\tau) :=& \, \widehat{\phi}
(\tau) \! + \! \widehat{\varrho}_{1}(\varepsilon_{1},\varepsilon_{2})
\ln \widehat{\phi}(\tau) \! + \! \widehat{\varrho}_{1}(\varepsilon_{1},
\varepsilon_{2}) \ln 12 \! + \! (-1)^{1+ \varepsilon_{2}} \Re (a) \ln
(2 \! + \! \sqrt{3}) \! + \! \dfrac{\pi}{2} \nonumber \\
-& \, \dfrac{1}{2} \arg \! \left(\dfrac{\widehat{g}_{11}(\varepsilon_{1},
\varepsilon_{2}) \widehat{g}_{12}(\varepsilon_{1},\varepsilon_{2})}{\widehat{
g}_{21}(\varepsilon_{1},\varepsilon_{2}) \widehat{g}_{22}(\varepsilon_{1},
\varepsilon_{2})} \right) \! - \! \arg \! \left(\Gamma \! \left(\tfrac{1}{2}
\! + \! \mi \widehat{\varrho}_{1}(\varepsilon_{1},\varepsilon_{2}) \right)
\right) \! + \! \mi \! \left((-1)^{1+ \varepsilon_{2}} \Im (a) \ln (2 \!
+ \! \sqrt{3}) \right. \nonumber \\
+&\left. \, \dfrac{1}{2} \ln \! \left(\left\vert \dfrac{\widehat{g}_{11}
(\varepsilon_{1},\varepsilon_{2}) \widehat{g}_{12}(\varepsilon_{1},
\varepsilon_{2})}{\widehat{g}_{21}(\varepsilon_{1},\varepsilon_{2})
\widehat{g}_{22}(\varepsilon_{1},\varepsilon_{2})} \right\vert \right)
\right) \! + \! \mathcal{O}(\tau^{-\delta_{G}} \ln \tau), \label{eqth5.2}
\end{align}
with $\widehat{\phi}(\tau)$ and $\vartheta_{0}$ given, respectively, in
Equations~\eqref{eqth4.8} and~\eqref{eq:theta0-def}.

Let $\mathcal{H}(\tau)$ be the Hamiltonian function defined by
Equation~\eqref{eqh1} corresponding to the function $u(\tau)$ given above.
Then $\mathcal{H}(\tau)$ has the asymptotic expansion
\begin{align}
\mathcal{H}(\tau) \underset{\underset{\scriptstyle \tau \in
\widehat{\mathcal{D}}_{u}}{\tau \to \infty \me^{\frac{\mi \pi
\varepsilon_{1}}{2}}}}{=}& \, \me^{-\frac{\mi \pi \varepsilon_{1}}{2}}
\left(\vphantom{M^{M^{M^{M^{M}}}}} 3(\epsilon b)^{2/3}(\me^{-\frac{\mi \pi
\varepsilon_{1}}{2}} \tau)^{1/3} \! + \! (-1)^{\varepsilon_{2}}4 \sqrt{3} \,
(\epsilon b)^{1/3}(\me^{-\frac{\mi \pi \varepsilon_{1}}{2}} \tau)^{-1/3}
\left(\vphantom{M^{M^{M}}} \widehat{\varrho}_{1}(\varepsilon_{1},
\varepsilon_{2}) \right. \right. \nonumber \\
+&\left. \left. \, \dfrac{1}{2 \sqrt{3}} \left((-1)^{1+ \varepsilon_{2}}
a \! - \! \dfrac{\mi}{2} \right) \! + \! \dfrac{1}{4} \cot (\tfrac{1}{2}
\widehat{\beta}(\varepsilon_{1},\varepsilon_{2},\tau)) \! + \! \dfrac{1}{4}
\cot (\tfrac{1}{2} \widehat{\beta}(\varepsilon_{1},\varepsilon_{2},\tau)
\! - \! \vartheta_{0}) \right. \right. \nonumber \\
+&\left. \left. \, \mathcal{O}(\tau^{-\delta_{G}}) \vphantom{M^{M^{M}}}
\right) \vphantom{M^{M^{M^{M^{M}}}}} \! \right).
\label{eqth5.3}
\end{align}
The function $f(\tau)$ defined by Equation~\eqref{eqf} has the following
asymptotics:
\begin{equation}
f(\tau) \underset{\underset{\scriptstyle \tau \in \widehat{\mathcal{D}}_{u}}{
\tau \to \infty \me^{\frac{\mi \pi \varepsilon_{1}}{2}}}}{=} -\dfrac{
(-1)^{\varepsilon_{2}}(\epsilon b)^{1/3}}{2}(\me^{-\frac{\mi \pi
\varepsilon_{1}}{2}} \tau)^{2/3} \! \left(\mi \! + \! \dfrac{3}{\sqrt{2} \,
\sin (\tfrac{1}{2} \widehat{\beta}(\varepsilon_{1},\varepsilon_{2},\tau))
\sin (\tfrac{1}{2} \widehat{\beta}(\varepsilon_{1},\varepsilon_{2},\tau)
\! - \! \vartheta_{0})} \right). \label{eqth5.4}
\end{equation}
\end{ddddd}
\begin{eeeee} \label{appcremc3}
\textsl{For real, non-zero values of $b$, singular imaginary solutions
$u(\tau)$ $($for imaginary $\tau)$ of Equation~\eqref{eq1.1} are specified
by the following ``singular imaginary reduction'' for the monodromy
data\footnote{There exist regular imaginary solutions (cf. Part~I~\cite{a1},
Appendix) which are specified by another imaginary reduction.}$:$
\begin{equation}
\begin{gathered}
s^{0}_{0} \! = \! -\overline{s^{0}_{0}}, \qquad \widehat{s}^{\infty}_{0}
(\varepsilon_{1},\varepsilon_{2}) \! = \! -\overline{\widehat{s}^{\infty}_{1}
(\varepsilon_{1},\varepsilon_{2})} \, \me^{2 \pi a}, \qquad \widehat{g}_{11}
(\varepsilon_{1},\varepsilon_{2}) \! = \! -\overline{\widehat{g}_{22}
(\varepsilon_{1},\varepsilon_{2})}, \\
\widehat{g}_{12}(\varepsilon_{1},\varepsilon_{2}) \! = \! -\overline{
\widehat{g}_{21}(\varepsilon_{1},\varepsilon_{2})}, \qquad \Im (a) \! = \! 0.
\label{eqappc3.1}
\end{gathered}
\end{equation}
In this case, asymptotics of $\widehat{\tau}_{m}^{\infty}$,
$\widehat{\tau}_{m}^{\pm}$, $u(\tau)$, $\mathcal{H}(\tau)$, and $f(\tau)$
are as given in Equations~\eqref{eqth6.3}, \eqref{eqth6.6}, \eqref{eqth5.1},
\eqref{eqth5.3}, and~\eqref{eqth5.4}, respectively, but with the changes
$\widehat{\varrho}_{1}(\varepsilon_{1},\varepsilon_{2}) \! \to \!
\widehat{\varrho}_{0}(\varepsilon_{1},\varepsilon_{2})$,
$\widehat{\varrho}_{2}(\varepsilon_{1},\varepsilon_{2}) \! \to \!
\widehat{\varrho}_{0}^{\sharp}(\varepsilon_{1},\varepsilon_{2})$, and
$\widehat{\beta}(\varepsilon_{1},\varepsilon_{2},\tau) \! \to \!
\widehat{\beta}_{0}(\varepsilon_{1},\varepsilon_{2},\tau)$, where
\begin{align}
\widehat{\varrho}_{0}(\varepsilon_{1},\varepsilon_{2}) :=& \, \dfrac{1}{\pi}
\ln \lvert \widehat{g}_{11}(\varepsilon_{1},\varepsilon_{2}) \rvert,
\label{eqappc3.2} \\
\widehat{\varrho}_{0}^{\sharp}(\varepsilon_{1},\varepsilon_{2}) :=& \,
\widehat{\varrho}_{0}(\varepsilon_{1},\varepsilon_{2}) \ln (24 \pi) \! +
\! (-1)^{1+ \varepsilon_{2}} \Re (a) \ln (2 \! + \! \sqrt{3}) \! - \!
\dfrac{\pi}{2} \nonumber \\
-& \, \arg \! \left(\widehat{g}_{11}(\varepsilon_{1},\varepsilon_{2})
\widehat{g}_{12}(\varepsilon_{1},\varepsilon_{2}) \Gamma \! \left(
\tfrac{1}{2} \! + \! \mi \widehat{\varrho}_{0}(\varepsilon_{1},
\varepsilon_{2}) \right) \right), \label{eqappc3.3} \\
\widehat{\beta}_{0}(\varepsilon_{1},\varepsilon_{2},\tau) :=& \,
\widehat{\phi}(\tau) \! + \! \widehat{\varrho}_{0}(\varepsilon_{1},
\varepsilon_{2}) \ln \widehat{\phi}(\tau) \! + \! \widehat{\varrho}_{0}
(\varepsilon_{1},\varepsilon_{2}) \ln 12 \! + \! (-1)^{1+ \varepsilon_{2}}
\Re (a) \ln (2 \! + \! \sqrt{3}) \nonumber \\
-& \, \dfrac{\pi}{2} \! - \! \arg \! \left(\widehat{g}_{11}(\varepsilon_{1},
\varepsilon_{2}) \widehat{g}_{12}(\varepsilon_{1},\varepsilon_{2}) \Gamma
\! \left(\tfrac{1}{2} \! + \! \mi \widehat{\varrho}_{0}(\varepsilon_{1},
\varepsilon_{2}) \right) \right) \! + \! \mathcal{O}(\tau^{-\delta_{G}}
\ln \tau). \label{eqappc3.4}
\end{align}}
\end{eeeee}

\end{document}